\documentclass[a4paper,
11pt,  reqno]{amsart}
\usepackage[margin=1.2in,marginparwidth=1.5cm, marginparsep=0.5cm]{geometry}

\usepackage{amsmath,amssymb,amscd,amsthm,amsxtra, esint}
\usepackage{mathrsfs}
\usepackage[most]{tcolorbox}
\usepackage{color}
\usepackage{hyperref}
\usepackage{subcaption}
\usepackage{float}

\usepackage{bbm}
\usepackage{forest}
\newcommand{\red}{\color{red}}

\newcommand{\fia}{\mathbbm{1}_{|\Phi-\alpha|<M}}

\newcommand{\jap}[1]{\langle #1 \rangle}
\newcommand{\I}{\hspace{0.2mm}\textup{I}\hspace{0.2mm}}
\newcommand{\II}{\textup{I \hspace{-2.9mm} I} }
\newcommand{\III}{\textup{I \hspace{-2.9mm} I \hspace{-2.9mm} I}}
\newcommand{\IV}{\textup{I \hspace{-2.9mm} V}}

\newcommand{\noi}{\noindent}
\newcommand{\m}[1]{
	\ifdefequal{#1}{1}
	{\mathbbm{#1}}
	{\mathbb{#1}}
}
\newcommand{\sech}{\text{sech}}
\newcommand{\FL}{\mathcal{F} L}
\newcommand{\ft}{\widehat}
\newcommand{\ind}{\mathbbm 1}
\newcommand{\les}{\lesssim}
\newcommand{\dl}{\delta}
\newcommand{\F}{\mathcal{F}}
\newcommand{\eps}{\varepsilon}

\newcommand{\wt}{\widetilde}
\newcommand{\embeds}{\hookrightarrow}
\newcommand{\dx}{\partial_x}
\newcommand{\dt}{\partial_t}

\allowdisplaybreaks[2]

\newtheorem{theorem}{Theorem} [section]
\newtheorem*{theorem*}{Theorem} 

\newtheorem{lemma}[theorem]{Lemma}
\newtheorem{proposition}[theorem]{Proposition}

\newtheorem{definition}[theorem]{Definition}

\theoremstyle{remark}
\newtheorem{remark}[theorem]{Remark}

\usepackage{mathtools}
\mathtoolsset{showonlyrefs=true}

\newcommand{\R}{\mathbb{R}}
\newcommand{\bbC}{\mathbb{C}}

\newcommand{\bul}{\bullet}

\let\P= \undefined
\newcommand{\P}{\mathbf{P}}

\newcommand{\TT}{\mathcal{T}}
\newcommand{\Tr}{\mathrm{T}}

\newcommand{\al}{\alpha}
\newcommand{\be}{\beta}

\newcommand{\g}{\gamma}
\newcommand{\G}{\Gamma}

\newcommand{\ld}{\lambda}

\newcommand{\bfrak}{\mathfrak{b}}
\newcommand{\gfrak}{\mathfrak{g}}

\newcommand{\pb}{\mathbf{p}}

\newcommand{\sbf}{\mathbf{s}}

\newcommand{\Ft}{{\mathcal{F}}}

\newcommand{\ta}{\theta}

\renewcommand{\l}{\ell}
\renewcommand{\o}{\omega}

\newcommand{\ges}{\gtrsim}

\newcommand{\jb}[1]
{\langle #1 \rangle}

\newcommand{\ff}{\mathfrak{f}}

\newcommand{\M}{\mathcal{M}}
\newcommand{\MMb}{\mathbf{M}}

\newcommand{\N}{\mathbb{N}}

\DeclareMathOperator{\mul}{mul}
\DeclareMathOperator{\mulB}{mul_{\BB}}

\tikzset{
	dotA/.style={
		transform shape, fill,circle,inner sep=1.5pt, label distance=-1pt,
		font={\normalsize }
	},
	>=stealth,
}

\tikzset{
	dotAB/.style={
		transform shape, fill,circle,inner sep=1.5pt, label distance=-10pt,
		font={\small }
	},
	>=stealth,
}

\makeatletter
\def\DeclareSymbolNew#1#2#3#4#5{\expandafter\gdef\csname MH@symb@#1\endcsname{\tikz[baseline=#2,scale=0.5, every node/.style={dotA},
		level distance=15pt,
		level 1/.style={sibling distance=#3pt},
		level 2/.style={sibling distance=#4pt},
		level 3/.style={sibling distance=15pt}
		]{#5}}}
\def\<#1>{\csname MH@symb@#1\endcsname}
\makeatother

\def\distC{15}

\DeclareSymbolNew{1}{-7}{\distC}{\distC}
{
	\node   (root) {}
	child
	child
	;
	\node  at (root-1) {};
	\node at (root-2) {};
}

\DeclareSymbolNew{11}{-7}{\distC}{\distC}
{
	\node   (root) {}
	child
	child
	;
	\node  at (root-1) {};
	\node at (root-2) {};
}

\newcommand{\cherryS}{
	\resizebox{!}{1.3ex}{%
		\<1>%
	}
}

\newcommand{\cherry}[3]{
	\begin{tikzpicture}[
		baseline=-8, scale=0.6, every node/.style={dotA},
		level distance=15pt,
		level 1/.style={sibling distance=\distC pt}
		]
		\node  [label=above:{$#1$}]  (root) {}
		child
		child
		;
		\node [label= left:{$#2$}]  at (root-1) {};
		\node [label= right:{$#3$}] at (root-2) {};
	\end{tikzpicture}
}

\DeclareSymbolNew{01}{-7}{\distC}{\distC}
{
	\node  [circle,draw=black, fill=white, inner sep=0pt,minimum size=7pt]  (root) {}
	child  child;
	\node  at (root-1) {};
	\node at (root-2) {};
}

\DeclareSymbolNew{21}{-12}{\distC}{\distC}
{
	\node   (root) {}
	child {child child} child;
	\node  at (root-1) {};
	\node at (root-2) {};
	\node  at (root-1-1) {};
	\node at (root-1-2) {};
}

\DeclareSymbolNew{22}{-12}{\distC}{\distC}
{
	\node   (root) {}
	child child {child child};
	\node  at (root-1) {};
	\node at (root-2) {};
	\node  at (root-2-1) {};
	\node at (root-2-2) {};
}

\DeclareSymbolNew{31}{-12}{\distC}{\distC}
{
	\node   (root) {}
	child {child {child child} child} child;
	\node  at (root-1) {}; \node at (root-2) {};
	\node at (root-1-1) {}; \node at (root-1-2) {};
	\node  at (root-1-1-1) {}; \node at (root-1-1-2) {};
}

\DeclareSymbolNew{32}{-12}{\distC}{\distC}
{
	\node   (root) {}
	child {child child {child child} } child;
	\node  at (root-1) {}; \node at (root-2) {};
	\node at (root-1-1) {}; \node at (root-1-2) {};
	\node  at (root-1-2-1) {}; \node at (root-1-2-2) {};
}

\DeclareSymbolNew{33}{-12}{21}{\distC}
{
	\node   (root) {}
	child {child child} child {child child};
	\node  at (root-1) {}; \node at (root-2) {};
	\node at (root-1-1) {}; \node at (root-1-2) {};
	\node  at (root-2-1) {}; \node at (root-2-2) {};
}

\DeclareSymbolNew{34}{-12}{\distC}{\distC}
{
	\node   (root) {}
	child child {child {child child} child  } ;
	\node  at (root-1) {}; \node at (root-2) {};
	\node at (root-2-1) {}; \node at (root-2-2) {};
	\node  at (root-2-1-1) {}; \node at (root-2-1-2) {};
}

\DeclareSymbolNew{35}{-12}{\distC}{\distC}
{
	\node   (root) {}
	child child {child child {child child} } ;
	\node  at (root-1) {}; \node at (root-2) {};
	\node at (root-2-1) {}; \node at (root-2-2) {};
	\node  at (root-2-2-1) {}; \node at (root-2-2-2) {};
}

\newcommand{\Ical}{\mathcal{I}}

\newcommand{\Gc}{\mathcal{G}}

\numberwithin{equation}{section}
\numberwithin{theorem}{section}

\newcommand{\MM}{\mathrm{M}}
\newcommand{\NN}{\mathrm{N}}
\newcommand{\Ncal}{\mathcal{N}}
\newcommand{\BB}{\mathrm{B}}
\newcommand{\D}{\mathrm{D}}
\newcommand{\DD}{\mathcal{D}}

\title[Gauge transform for KdV and low regularity well-posedness]{Gauge transform for the Korteweg-de Vries equation and well-posedness below the $H^{-1}$-scale}
\author[A.~Chapouto \and S.~Correia \and J.~P.~Ramos]
{Andreia Chapouto \and Sim\~ao Correia \and Jo\~ao Pedro Ramos}

\date{\today}

\address{Andreia Chapouto,
	CNRS, Laboratoire de math\'ematiques de Versailles, UVSQ, Universit\'e Paris-Saclay, CNRS, 45 avenue des \'Etats-Unis, 78035 Versailles Cedex, France, 
    and School of Mathematics, Monash University, Clayton, VIC 3800, Australia}

\email{andreia.chapouto@monash.edu}

\address{Sim\~ao Correia,
	Center for Mathematical Analysis, Geometry and Dynamical Systems,
	Department of Mathematics,
	Instituto Superior T\'ecnico, Universidade de Lisboa
	Av. Rovisco Pais, 1049-001 Lisboa, Portugal
}

\email{simao.f.correia@tecnico.ulisboa.pt}

\address{Jo\~ao Pedro Ramos,
	IMPA, Rio de Janeiro, Brazil}

\email{joao.ramos@impa.br}

\subjclass[2020]{35A01, 35A22, 35Q53, 70K45}

\keywords{Korteweg-de Vries equation, gauge transform,
	well-posedness.
}

\begin{document}

\begin{abstract}

We propose a new formulation of the Korteweg-de Vries equation (KdV) on the real line, via a gauge transform.
While KdV and the gauged equation are equivalent for smooth solutions, the latter is better behaved at low regularity in Fourier-Lebesgue spaces. In particular, the admissible regularities go beyond the $H^{-1}$-scale, which is a well-known threshold for KdV. As a byproduct, by reversing the gauge transform, we are able to improve on the known theory for KdV and derive sharp local well-posedness in Fourier-Lebesgue spaces with large integrability exponent.
Our strategy is based on an infinite normal form reduction and Fourier restriction estimates, together with a thorough exploitation of algebraic cancellations. Additionally, our method is totally independent of the KdV completely integrable structure, and extends to other non-integrable models with quadratic nonlinearities.
\end{abstract}

	\maketitle

    \section{Introduction}\label{SEC:intro}

    \subsection{Korteweg-de Vries equation}

We consider the Korteweg-de Vries equation (KdV) on the real line $\R$:
	\begin{align}
		\label{kdv}\tag{KdV}
		\dt u + \dx^3 u = \dx(u^2),
	\end{align}
	which was first introduced by Boussinesq in 1887 \cite{Boussinesq} and rediscovered by Korteweg and de Vries in 1895 \cite{KortewegDeVries}. It was proposed as a model for the unidirectional propagation of solitary waves in a channel, which had been observed by John Scott Russell, as mentioned in his 1882 book \cite{Russell}. Since then, \eqref{kdv} has appeared in the context of fluid mechanics, oceanography, acoustics, geophysics, plasma physics and interacting particle fields \cite{Fordy, CCR_particle,  Kuznetsov_Nakoryakov_Pokusaev_Shreiber_1978,    Nezlin, OlsonChristensen, StewartCorones,WashimiTaniuti}, making it one of the most ubiquitous models for the nonlinear propagation of waves. For more details, we refer to the excellent review by Crighton \cite{Crighton}.


From the analytical viewpoint, \eqref{kdv} is a canonical dispersive equation, with a rich structure, a variety of special solutions (such as solitons and multi-solitons), and it is completely integrable \cite{ZakharovFadeev}. As such, it admits a Lax pair formulation and infinitely many conservation laws, one in each $L^2$-based Sobolev space $H^k$ for $k\in\N$. Some of the most relevant conserved quantities include:
\begin{align*}
\text{Mass:}& \int_\R u dx , \\
\text{Momentum:}& \int_\R u^2 dx , \\
\text{Energy:}& \int_\R \frac{(\dx u)^2}{2} + \frac{u^3}{3}  dx.
\end{align*}
Another important feature of \eqref{kdv} is its scaling invariance: given a solution $u$ to \eqref{kdv}, then
\begin{equation}\label{scaling}
u_\lambda(t,x)=\lambda^{2}u(\lambda^3 t, \lambda x)
\end{equation}
is also a solution, for any $\lambda>0$. We note that this scaling leaves the $\dot H^{-\frac32}(\R)$ norm invariant. As such, we refer to $s=-\frac32$ as the scaling-critical regularity (in the context of $H^s(\R)$ spaces).


Our focus in this work is the low-regularity local well-posedness problem for \eqref{kdv}, which requires a brief overview of the known literature{\footnote{Which is by no means exhaustive, as it spans more than half a century.}}. In the 1970s, local well-posedness for \eqref{kdv} in $H^s(\R)$, $s>\frac32$, was shown by viewing \eqref{kdv} as a quasilinear hyperbolic equation \cite{BonaSmith, Kato_kdv, SautTemam, TsutsumiMukasa}. These works did not rely on the dispersive character of \eqref{kdv}. The first {author} to incorporate dispersion was Kato \cite{Kato_kdv_disp}, later expanded upon in the works of Kenig, Ponce, and Vega \cite{KPV91, KPV_gkdv}, which reached the threshold $s>\frac34$.

Soon after, Bourgain \cite{bourg2} reached $s=0$ by introducing the Fourier restriction norm method,
based on spaces tailored to the dispersion of the equation. This seminal paper was later improved in \cite{kpv_bilin} to the threshold $s>-\frac34$. The endpoint case $s=-\frac34$ was first attained by Christ, Colliander, and Tao \cite{CCT_ill} and revisited by Guo \cite{Guo_gwp} and Kishimoto \cite{Kishimoto_lwp_crit}. 
The aforementioned results are based on contraction arguments, which can be extended globally-in-time by means
 of (modified) conservation laws associated with \eqref{kdv}  \cite{bourg2, CKSTT_sharp_kdv, Guo_gwp}. For $s<-\frac34$, the initial-value problem is \emph{ill-posed} in the sense of failure of uniform continuity of the data-to-solution map \cite{KPV_ill}. As such, $s=-\frac34$ is a threshold for contraction arguments, which lies substantially above the
scaling-critical regularity $s=-\frac32$.

In order to surpass the $s=-\frac34$ threshold, exploiting more deeply the structure of \eqref{kdv}, namely its complete integrability, has
led to remarkable improvements to the well-posedness theory. Kappeler \cite{Kappeler_rough} and Tsutsumi \cite{Tsutsumi_kdv_measures} constructed solutions with Radon measures as initial data. In \cite{BuckmasterKoch}, Buckmaster and Koch proved the existence of $H^{-1}(\R)$ weak solutions. Later, Killip and Vi\c{s}an \cite{KillipVisan_kdv} showed the (global) well-posedness of \eqref{kdv} in $H^{-1}(\R)$, where the solution is understood as the unique limit of Schwartz solutions. This result is \emph{optimal}, as it is not possible to define a continuous flow in $H^s(\R)$, for $s<-1$ \cite{Molinet_ill_line}. In the periodic case, complete integrability has led to even more substantial improvements to the local result of Kenig, Ponce, and Vega \cite{kpv_bilin}, {as the reader may find in }\cite{KMMT16, KappelerMolnar18, KT06, KillipVisan_kdv}.

As illustrated above, for particular dispersive PDEs such as \eqref{kdv}, complete integrability techniques
have emerged as one of the most powerful tools in the field.
However, these methods are rigid and dependent on the delicate integrability structure, which is often lost under small perturbations of the dispersion and nonlinearity.
The breakdown of these techniques forces one to go back to perturbative methods.
Consequently,
as a primary postulate of this work, we
avoid 
the complete integrability of \eqref{kdv} and pursue a Fourier analytic approach instead,
building upon the contraction methods which allow us to encompass perturbations of \eqref{kdv} and offer new insights into its rich dynamics.

To bridge the gap between the sharp well-posedness threshold in $H^s(\R)$ ($s=-1$) and the critical regularity $s=-\frac32$, one may consider a different functional framework.
One such choice, amenable to Fourier analytic techniques, is the Fourier-Lebesgue spaces $\FL^{s,p}(\R)$, defined through the following norm\footnote{
In this text, we use the convention
$
 \Ft(f)(\xi) := (2\pi)^{-\frac12} \int_\R f(x) e^{-i x \xi} dx.
$
For space-time functions $f: \R^2 \to \R$, we may write $\Ft_t (f)$ and $\Ft_x(f)$ to denote the Fourier transform with respect to $t$ and $x$, respectively. When clear from context, we also use $\ft{f}$ to denote either the spatial Fourier transform or the space-time Fourier transform. }
	\begin{equation}
		\|u \|_{\FL^{s,p}(\R)}=\|\jb{\xi}^{s}\ft{u  }(\xi)\|_{L^p(\R)},
	\end{equation}
    where $\jb{\,\cdot\,} = (1+|\cdot|^2)^\frac12$.
These spaces agree with the classical $H^s(\R)$ spaces when $p=2$ 
and were used in many works \cite{adams, BO94, Chapouto_FL_mkdv, Christ_power,CC24,Grunrock_mkdv, Grunrock_ldv_hier, HerrGrunrock_FL,  GrunrockVega}, where larger values of $p$ tend to close the gap between the critical and the known well-posedness regularity.

A useful heuristic comes from the behaviour of the homogeneous ${\FL}^{s,p}$-norm under the natural scaling \eqref{scaling} of \eqref{kdv}.
Since
$$
\big\| |\xi|^s  \ft{u_\ld}(\xi) \big\|_{L^p(\R)} \sim \lambda^{1+s+\frac1p}
\big\| |\xi|^s \ft{u} (\xi) \big\|_{L^p(\R)},
$$
we say that the Fourier-Lebesgue space $\FL^{s,p}(\R)$ with
$s=s_c(p)=-1-\frac1p$ is \emph{critical}, as the corresponding norm is left invariant under scaling. For comparison, we define $s_{H^{-1}}(p) = - \frac12 - \frac1p$ as the Fourier-Lebesgue regularity at the $H^{-1}$-scale.
 With such a perspective in mind, we note that \eqref{kdv} in Fourier-Lebesgue spaces was studied by Gr\"unrock \cite{Grunrock_ldv_hier},  where he proves local well-posedness of \eqref{kdv} in $\mathcal{F}L^{s,p}(\R)$ for
	\begin{equation}
		2<p<\infty,\qquad s>\max\left\{ -\tfrac12 - \tfrac{1}{2p}, -\tfrac14 - \tfrac{11}{8p} \right\}.
\label{eq:sp_grunrock}
	\end{equation}
    In Appendix \ref{sec:-14}, we extend this result to $p=\infty$. 
The result in \cite{Grunrock_ldv_hier} remains the current state-of-the-art in this functional framework, leaving a substantial gap from the scaling critical regularity $s_c(p)$ and even from the $H^{-1}$-scale $s_{H^{-1}}(p)$. The regularity constraint \eqref{eq:sp_grunrock} comes from the main bilinear estimate, in particular, the high$\times$high$\to$low frequency interaction $|\xi_1| \sim |\xi_2| \gg \big||\xi_1| - |\xi_2|\big|$ (see Appendix~\ref{sec:-14}).
To go beyond this result, our first main contribution is the introduction of a gauge transform that
successfully removes these interactions and allows for a local well-posedness theory via a contraction argument at lower regularity, surpassing even the $H^{-1}$-scale, for the gauged equation. Afterwards, by reversing the gauge transform, we derive the sharp range for well-posedness of \eqref{kdv} in $\FL^{s,p}(\R)$.

We summarize our main contributions for the \eqref{kdv} well-posedness theory in the next (informal) statement. 
For precise statements, see Theorem \ref{thm:lwp_rkdv}, Theorem \ref{thm:wp_kdv_intro}, and Proposition~\ref{prop:ill_intro}.

\begin{theorem}\label{thm:main}
Fix $2\le p \le \infty$.

\smallskip

\noi
{\rm(i)} \eqref{kdv} is (analytically) locally well-posed in $\FL^{s,p}(\R)$ for
$$
s>\max\left\{-\tfrac12-\tfrac1p, -1+\tfrac1{2p}\right\}.
$$

\noi
{\rm(ii)} \eqref{kdv} is ill-posed in $\FL^{s,p}(\R)$ (in the sense of non-existence of a continuous flow)~for
$$
s<-\tfrac12-\tfrac1p.
$$

\noi
{\rm(iii)} 
If $2\le p <\infty$, then
\eqref{kdv} is ill-posed in $\FL^{s,p}(\R)$ (in the sense  of failure of local uniform continuity of the flow)  for
$$
s<-1+\tfrac1{2p}.
$$

\noi
{\rm(iv)} There exists a gauged KdV equation \eqref{rkdv} (see also \eqref{z2}) satisfying the following:
\begin{enumerate}
\item[(a)] For $s' \gg 1$, \eqref{kdv} and \eqref{rkdv} are equivalent in $\FL^{s',p}(\R)$. In particular, there exists an invertible map $\Gc: \FL^{s',p}(\R) \to \FL^{s',p}(\R)$ (see \eqref{gauge}) such that for small data $u_0 \in \FL^{s',p}(\R)$, $u$ solves \eqref{kdv} with $u(0) = u_0$ if and only if $z:= \Gc[u]$ solves \eqref{rkdv} with $z(0) = \Gc[u_0]$. 

\item[(b)] \eqref{rkdv} is (analytically) well-posed for small data in $\FL^{s,p}(\R)$ for
\begin{align*}
s> -1 + \tfrac{1}{2p}
.
\end{align*}
\end{enumerate}

\end{theorem}

\begin{figure}[h]
		\centering
		\includegraphics[height=4.4cm]{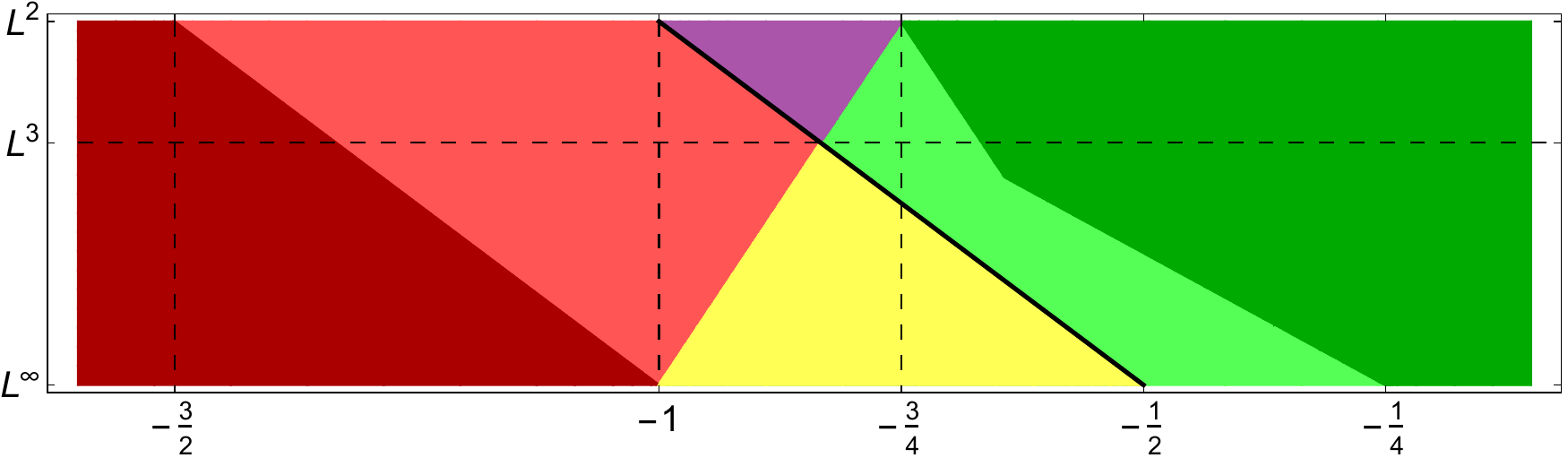}
		\caption{$\FL^{s,p}(\R)$ well-posedness  theory in the $(s,\frac1p)$-plane.
        Dark-red denotes the supercritical regime. In the light-red and yellow regions, \eqref{kdv} is ill-posed (for continuous flows). In the purple range, \eqref{kdv} is ill-posed (for uniformly continuous flows). The solid black line is the $H^{-1}$-scale.
        The dark-green region is the result by Grünrock \cite{Grunrock_ldv_hier}. The light-green region is the new (analytic) well-posedness range for \eqref{kdv}. In the yellow region, \eqref{kdv} is ill-posed while \eqref{rkdv} is well-posed.}\label{fig:kdv_intro}
	\end{figure}

Theorem~\ref{thm:main} provides a definite picture of the well-posedness of \eqref{kdv} in Fourier-Lebesgue spaces. 
Firstly, for any $2 \le p \le \infty$, Part (i) improves the result by Gr\"unrock in \cite{Grunrock_mkdv}. 
Secondly, for $3 \le p \le \infty$, we establish the $H^{-1}$-scale as the \textit{optimal} range for \eqref{kdv}: well-posedness holds for $s>s_{H^{-1}}(p) = - \frac12 - \frac1p$ (Part (i)) and there is no continuous extension of the data-to-solution map if $s<s_{H^{-1}}(p)$ (Part (ii)).
Lastly, we propose the gauged KdV equation \eqref{rkdv} as the correct model to study at low regularity, since \eqref{rkdv} and \eqref{kdv} are \emph{equivalent} for smooth solutions while \eqref{rkdv} is better behaved at low regularity. 
In particular, for $2 \le p < 3$, we have well-posedness of \eqref{rkdv} for $-1 + \frac{1}{2p} \ge s>s_{H^{-1}}(p)$, while \eqref{kdv} is analytically ill-posed, and for $3 \le p \le \infty$, \eqref{rkdv} is well-posed \emph{below} the $H^{-1}$-scale; see Part (iv).
See Figure~\ref{fig:kdv_intro} for details on the parameter ranges of Theorem~\ref{thm:main}.

We arrive at the gauged equation \eqref{rkdv} via an infinite iteration of normal form reductions and exploiting crucial cancellations introduced by the gauge transform $\Gc$; see \eqref{gauge}. 
As a result, the nonlinearity is much more involved than that for \eqref{kdv}.
Although our gauge transform 
$\Gc$ resembles the classical Miura map \eqref{eq:defiMiura} (see Subsection~\ref{SUBSEC:miura}), 
Theorem~\ref{thm:main} is of \emph{semilinear} nature and extends to non-integrable perturbations of \eqref{kdv} (see Remark \ref{rem:dgbo}).

\subsection{Gauged KdV equation}\label{SEC:intro-g}

To handle the obstructions found in \cite{Grunrock_ldv_hier} (which lead to the constraints \eqref{eq:sp_grunrock}), we introduce the gauge transform $\Gc:\FL^{s,p}(\R) \to \FL^{s,p}(\R)$, whose inverse is formally given~by
    \begin{equation}
\label{gauge}
\Gc^{-1}[f]
:= f + \Ft_x^{-1}\bigg( - \frac{1}{3} \int_{\xi=\xi_1+\xi_2} \ind_{A^c}(\xi_1,\xi_2) \frac{1}{\xi_1\xi_2} \ft{f}(\xi_1) \ft{f}(\xi_2) d\xi_1 \bigg),
    \end{equation}
    where $\Ft^{-1}_x$ denotes the inverse Fourier transform and
$A^c=\{(\xi_1,\xi_2)\in \R^2: |\xi_1|, |\xi_2|\ge 1\}$ includes the problematic high$\times$high$\to$low frequency interactions.
The map $\Gc$ and its inverse are well-defined for small data in $\FL^{s,p}(\R)$ if and only if $s>s_{H^{-1}}(p)=-\frac12-\frac1p$ (Lemma \ref{LEM:D}), i.e., above the $H^{-1}$-scale.
    Moreover,
the gauge transform connects \eqref{kdv} to a novel KdV-type equation: given a smooth solution $u$ to \eqref{kdv}, then $z$, defined as
    \begin{equation}\label{defi_z_intro}
z=\mathcal{G}[u],
    \end{equation}
solves the \emph{gauged KdV equation}:
\begin{equation}
\label{rkdv}\tag{$\Gc$-KdV}
\partial_t z+ \partial_x^3 z =
\sum_{j=2}^\infty
N_j[z].
	\end{equation}
Here,
each $N_j$, $j\ge 2$, is a nonlocal nonlinear operator with homogeneity $j$ in $z$. In particular, the quadratic nonlinear terms $N_2[z]$ in \eqref{rkdv} do not contain the bad high$\times$high$\to$low frequency interactions, at the cost of more involved cubic and higher-order nonlinearities.
  Despite the apparent increased complexity of the nonlinearity of \eqref{rkdv}, it is in fact \emph{better-behaved} than \eqref{kdv}, and we can push its well-posedness theory to lower regularities, namely beyond the $H^{-1}$-scale.
The improved regularity and the particular form of the nonlinear terms in \eqref{rkdv} come from an infinite normal form reduction
accounting for careful algebraic cancellations together with $\Gc^{-1}$ in \eqref{gauge}.

    Crucially, \eqref{kdv} and \eqref{rkdv} are \textit{equivalent} for small\footnote{The smallness assumption is harmless in the context of \eqref{kdv} since we are working in the scaling-subcritical regime $s>s_c(p)$.} smooth solutions, due to the mapping properties of the gauge transform $\Gc$ in \eqref{gauge} (see Proposition \ref{PRO:equiv}):
    $$
\eqref{kdv} \overset{\Gc}{\longleftrightarrow} \eqref{rkdv},
    $$
while \eqref{rkdv} is better behaved at lower regularity.
The idea of using a gauge transform to write a given dispersive equation in an equivalent way (for smooth data) while improving the well-posedness theory has been exploited in other contexts; see 
\cite{ 
takaoka_dnls, 
CKSTT_dnls,
tao_bo,  
Herr_dnls, 
MolinetPilod_bo, 
IfrimTataru,
Chapouto_FL_mkdv, 
KishimotoTsutsumi_kdnls}, for example.
To the best of our knowledge, the gauge transform $\Gc$ is a new transform for \eqref{kdv} on the real line, with the advantage of being independent of
any integrable structure, and thus being easily generalized to other dispersive equations with quadratic nonlinearities, such as 
\eqref{gBO} and Novikov-Veselov equations; see Remark~\ref{rem:dgbo}.

 Our first result is the local well-posedness for the gauged KdV equation \eqref{rkdv} in the Fourier-Lebesgue class.

	\begin{theorem}\label{thm:lwp_rkdv}
Let $2 \le p \le \infty$ and $s> - 1 + \frac1{2p}$. Given $T>0$, there exists $\dl = \dl(T)>0$, such that for $z_0 \in \FL^{s, p}(\R)$ with
\begin{align}
    \|z_0\|_{\FL^{s,p}(\R)} < \dl,
\end{align}
 the following holds:

\vspace{1mm}

\noi{\rm(i)}
When $ 2 \le p <\infty$,
there exist $b>\frac{1}{p'}$ and  a unique strong integral solution of \eqref{rkdv} on $[0,T]$
	$$z\in X^{s,b}_{p, p}(0,T)\hookrightarrow C([0,T]; \mathcal{F}L^{s,p}(\R)),$$
    with absolute summability of the nonlinear terms (see \eqref{Xsb} for the definition of $X^{s,b}_{p,q}(0,T)$);

\vspace{1mm}

\smallskip
\noi{\rm(ii)}
For $p=\infty$, there exist $b=1-$, $\tilde{b}=0+$ and a unique strong integral solution of \eqref{rkdv} on $[0,T]$
	$$z\in X^{s,b}_{\infty,\infty}(0,T)\cap X^{s,\tilde{b}}_{\infty,1}(0,T)
    \quad
    \text{with}
    \quad
    e^{t\partial_x^3}z \in C([0,T]; \mathcal{F}L^{s,\infty}(\R))
    ,
    $$
    with absolute summability of the nonlinear terms, where $e^{-t\dx^3}$ denotes the linear Airy group at time $t$.

\smallskip
\noi
    Moreover, in \textup{(i)-(ii)}, the solution depends analytically on the initial data.
\end{theorem}

\begin{figure}[h]
		\centering
		\includegraphics[height=4.4cm]{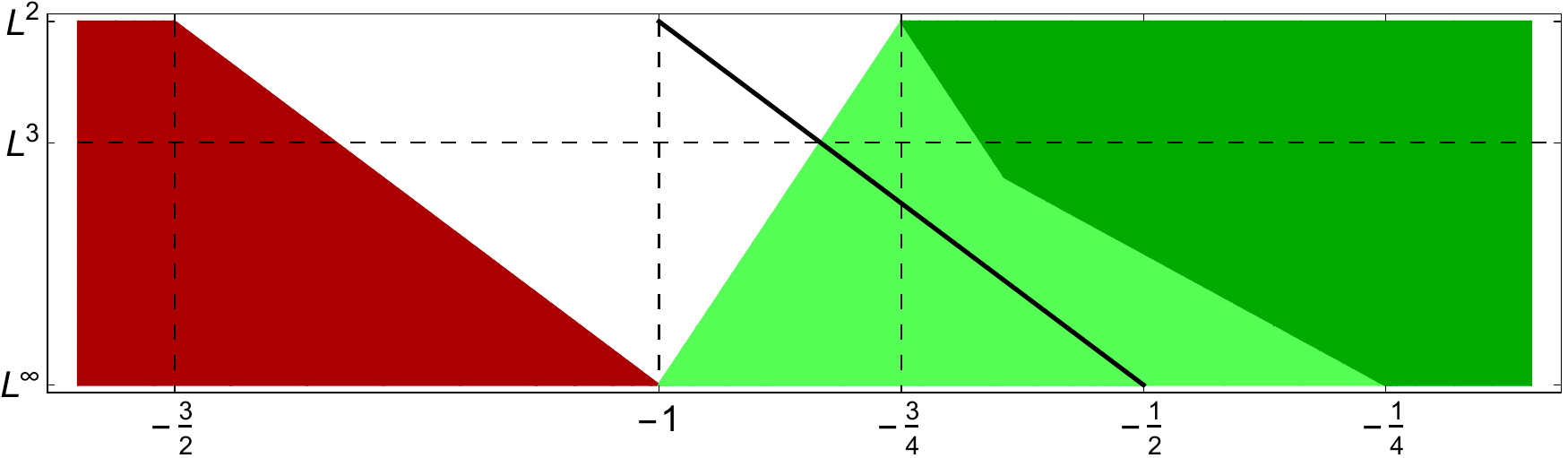}
		\caption{Well-posedness for \eqref{rkdv} in the $(s,\frac1p)$-plane. The dark-red region corresponds to the supercritical region. The black line represents the $H^{-1}$-scale, which is the threshold for the boundedness of the gauge transform. The dark-green region corresponds to the well-posedness result of \cite{Grunrock_ldv_hier}, while the light green corresponds to Theorem \ref{thm:lwp_rkdv}.}\label{fig:rkdv}
	\end{figure}

Theorem~\ref{thm:lwp_rkdv} establishes well-posedness of \eqref{rkdv} beyond the $H^{-1}$-scale when $p$ is sufficiently large, namely $p>{3}$. See Figure~\ref{fig:rkdv} for a visual depiction of the range covered in the result, compared with the scaling indices of \eqref{kdv}.
Also, 
the threshold in Theorem~\ref{thm:lwp_rkdv} is always below that in \eqref{eq:sp_grunrock}, and thus extends the well-posedness result for \eqref{rkdv} obtained from combining \cite{Grunrock_ldv_hier} with the gauge in \eqref{gauge}.

\begin{remark}
\noi(i)
Unlike \eqref{kdv}, \eqref{rkdv} does not enjoy the scaling invariance \eqref{scaling} (due to the restriction to the set $A^c$ appearing in the nonlinear terms $N_j$), and thus we cannot trivially extend the result in Theorem~\ref{thm:lwp_rkdv} from small to large data.

\smallskip

\label{remark:Linfty_continuity}
\noi (ii) When $p=\infty$, the time-continuity of the solution holds only at the level of the interaction representation $e^{t\dx^3}z(t)$. In fact, due to the high oscillations of the linear group at high frequency, $z$ is actually \emph{discontinuous} in $\FL^{s,\infty}(\R)$. This subtlety is not present when $p$ is finite (because of dominated convergence) and is a generic feature of well-posedness in $\FL^{s,\infty}(\R)$, independent of the underlying dispersive equation (see Section~\ref{sec:fourier_restr}, \cite[Section~7]{CC24}, and \cite[Remark~1.4.6]{Forlano}).

\end{remark}

\bigskip

The proof of Theorem \ref{thm:lwp_rkdv} follows a standard fixed-point argument in Fourier restriction spaces adapted to the Fourier-Lebesgue setting,  once we derive the necessary multilinear estimates for each of the nonlinear terms in \eqref{rkdv} (with some nontrivial adaptations in the case $p=\infty$, see Section \ref{sec:fourier_restr}), whose complexity increases with $j$.
We do so systematically by reducing these to Fourier restricted estimates, as introduced by the second author in~\cite{COS}.
The main novelty in this work comes from \emph{deriving} the gauged equation \eqref{rkdv}, whose structure is essential to break the $H^{-1}$-regularity barrier.
This process relies on an infinite normal form reduction
and on a thorough exploitation of algebraic cancellations coming from the phase function and the  gauge transform $\Gc$.
These cancellations crucially remove the bad frequency interactions that create the regularity bottleneck in \cite{Grunrock_ldv_hier} for \eqref{kdv}.

\subsection{From gauged KdV back to KdV}
With the local well-posedness for \eqref{rkdv} in Theorem~\ref{thm:lwp_rkdv}, we return our attention to \eqref{kdv} and improve the regularity constraint \eqref{eq:sp_grunrock} from \cite{Grunrock_ldv_hier}. 
The gauge transform $\Gc$ is a local bijection over $B_\delta\subset \FL^{s,p}(\R)$, for $s>s_{H^{-1}}(p) = -\frac12-\frac1p$, where
$$
B_\delta=\left\{ f \in \FL^{s,p}(\R): \|f\|_{\FL^{s,p}(\R)}<\delta\right\},\qquad 0<\delta\ll 1.
$$
As such, given initial data $u_0\in \FL^{s,p}(\R)$ for \eqref{kdv}, we must first use the scaling invariance \eqref{scaling} to reduce to sufficiently small data, move to \eqref{rkdv} using $\Gc$, construct the solution via Theorem \ref{thm:lwp_rkdv}, and then undo both the gauge and the rescaling.
To rigorously formalize this procedure,  given $\lambda>0$, we define the scaling map
   \begin{equation}
    \label{SR}
    S_\lambda[f](t,x) := \lambda^2 f(\lambda^3 t, \lambda x),
   \end{equation}
which we extend in a natural way to time-independent functions, keeping (with a slight abuse) the same notation.

\begin{theorem}\label{thm:wp_kdv_intro} Let $2\le p\le \infty$,
\begin{equation}
\label{eq:reg-kdv}
s>\max\left\{-\tfrac12 - \tfrac1p, -1+\tfrac1{2p}\right\}, 
\end{equation}
and $\Theta$ be the data-to-solution map for \eqref{rkdv} over the time interval $[0,1]$ defined in Theorem~\ref{thm:lwp_rkdv}. Then, the following hold:

\vspace{1mm}
{
\noi{\rm(i)} If $2 \le p<\infty$, given $R>0$, there exist $T, \ld>0$, depending on $R$, such that the map
\begin{align}\label{wtTa}
u_0\in B_R\subset\FL^{s,p}(\R)\mapsto u&
=\wt\Theta[u_0]
\\
&:= S_{\lambda^{-1}}\circ\mathcal{G}^{-1} \circ\Theta  \circ \Gc \circ S_\lambda[u_0]\in C([0,T];\FL^{s,p}(\R)),
\end{align}
is the unique continuous extension of the \eqref{kdv} data-to-solution map, defined for smooth initial data.
}

\smallskip

\noi{\rm(ii)} If $p=\infty$, given $R>0$, there exist $T, \ld>0$, depending on $R$, such that the map
\begin{equation}
u_0\in B_R\subset\FL^{s,\infty}(\R)\mapsto \wt\Theta[u_0]  \in L^\infty([0,T];\FL^{s,\infty}(\R))
\end{equation}
is a continuous extension of
the data-to-solution map for \eqref{kdv}, defined for smooth initial data, and it satisfies
\begin{equation}
e^{t\dx^3} u = e^{t\dx^3} \wt\Theta[u_0] \in C([0,T] ; \FL^{s,\infty}(\R))
.
\end{equation}

\noi
Moreover, the data-to-solution map $\wt{\Theta}$ is analytic.
	\end{theorem}

\begin{remark}
 In the case $p=\infty$, we cannot show that $\wt{\Theta}$ is the \emph{unique} extension of the flow defined for smooth data, simply because smooth decaying functions are not dense in $\FL^{s,\infty}(\R)$.
\end{remark}

Firstly, 
Theorem~\ref{thm:wp_kdv_intro} extends the result in \cite{Grunrock_ldv_hier} to lower regularities, as the threshold in \eqref{eq:reg-kdv} always lies below that in \eqref{eq:sp_grunrock}.
Secondly, since $\Gc^{-1}$ is unbounded on $\FL^{s,p}(\R)$ with $s< s_{H^{-1}}(p) $, we cannot go beyond the $H^{-1}$-scale for \eqref{kdv}, as we do for \eqref{rkdv}. This stark contrast in the well-posedness results is a property of the equations \eqref{kdv}-\eqref{rkdv} and not a deficiency of our method, as shown by the ill-posedness claims in Theorem \ref{thm:main}~(ii)-(iii).
These claims  exploit two known ill-posedness mechanisms. 
For (iii), we use the breather  solutions to \eqref{kdv} , namely special solutions which are time-periodic and spatially decaying, already used in \cite{KPV_ill} to show ill-posedness of \eqref{kdv} in $H^s(\R)$, $s<-\frac34$. Our argument is a simple extension of  this argument to the Fourier-Lebesgue class. For (ii), which concerns ill-posedness below the $H^{-1}$-scale, we adapt the arguments of Molinet \cite{Molinet_ill_line}, which takes advantage of the unboundedness of the Miura map (see Subsection~\ref{SUBSEC:miura}).

\begin{remark}[Extensions to other dispersive models]
\label{rem:dgbo}

Given $a \in [1,2]$, consider the dispersion-generalized Benjamin-Ono equation
\begin{align}
\label{gBO}\tag{dgBO}
\dt u + |\dx|^a \dx u = \dx (u^2)
\end{align}
where $|\dx|^a$ denotes the fractional Laplacian given by $|\dx|^a f = \Ft_{x}^{-1} (|\xi|^a \ft f (\xi))$.
In particular, \eqref{gBO} with $a=1$ corresponds to the Benjamin-Ono equation (BO), while $a=2$ gives \eqref{kdv}. 
Unlike the special cases of $a=1,2$, \eqref{gBO} is not completely integrable for $1<a<2$.
In the $H^s(\R)$-class, Herr \cite{herr_bo} proved the local well-posedness of \eqref{gBO} in $H^s (\R) \cap \dot H^{\frac12-\frac 1a}(\R)$ for $s>\frac{3}{4}(1-a)$, using a contraction argument in Fourier restriction spaces. The restriction to $\dot H^{\frac12-\frac 1a}(\R)$, caused by high$\times$low$\to$high interactions, has been recently lifted in \cite{AiLiu} through a gauge transform analogous to that of Tao \cite{tao_bo} for BO. For other works concerning \eqref{gBO} and its variants, we refer to \cite{CamposLinaresSantos, FarahLinaresPastor_bo, Jockel, KenigMartelRobbiano}.
As our results are perturbative by nature, the well-posedness theory developed throughout this work can also be adapted to the \eqref{gBO} equation. 
 We point out the necessary adaptations in Subsection \ref{subsec:gbo}.

This example shows the generality of our methodology, which extends to dispersive equations with amenable quadratic nonlinearities (namely, satisfying the key divisibility condition in \eqref{eq:divisible}), outside of the scope of integrable systems. 
Moreover, the  framework presented in this work is not restricted to 1$d$ models or to the underlying spatial geometry.
In the companion paper \cite{CCR-NV}, we explore its applicability to the Novikov-Veselov equation, a 2$d$-analogue of \eqref{kdv}, reaching sharp local well-posedness in $H^s(\R^2)$.

\end{remark}

\subsection{Derivation of the gauged equation}
\label{SUBSEC:method}

    In the following, we restrict our attention to the analysis of \eqref{kdv} in $\FL^{s,\infty}(\R)$, for simplicity. 
	To pinpoint the precise obstacles in \cite{Grunrock_ldv_hier},  it is convenient to introduce the interaction representation~\cite{Ginibre_interaction_rep}
	\begin{equation}\label{inter}
		v(t) : = e^{t\partial_x^3}u(t).
	\end{equation}
This classical change of variables emphasizes the multilinear dispersion effects in the nonlinearity.
Indeed, we may recast \eqref{kdv} as\footnote{Here and throughout this work, we use the notation $d\vec{\xi}$ to represent the volume form over the frequency integration hyperplane.}
	\begin{align}
		\label{dtv}
		\dt{\widehat{v}}(t,\xi) = \int_{\xi = \xi_1+\xi_2} e^{-it (\xi^3 -\xi_1^3 -
			\xi^3_2)} i \xi \widehat{v}(t,\xi_1) \widehat{v}(t,\xi_2)  d\vec{\xi}
		,
	\end{align}
where the oscillatory integral depends on the resonance function
	\begin{equation}\label{phase-intro}
		\Phi = \Phi(\xi,\xi_1,\xi_2)=-\xi^3+\xi_1^3+\xi_2^3 = -3\xi\xi_1\xi_2, \qquad \xi=\xi_1+\xi_2.
	\end{equation}
One can see that $\Phi(\xi, \xi_1, \xi-\xi_1)$ has two stationary interactions
(corresponding to $\partial_{\xi_1}\Phi=0$), namely $\xi_1=\xi_2$ and $\xi_1=-\xi_2$. Away from these points, we can reach scaling criticality in $\FL^{s,\infty}(\R)$, via a change of variables $\xi_1\mapsto \Phi$, in analogy with the computations of \cite{kpv_bilin}.
In fact, it is the stationary interactions that cause the regularity obstruction of $s>-\frac14$ (compare the proof of Proposition \ref{prop:bourgLinfty} for the Type $\I$ term with Lemma \ref{LEM:bil0}).

	In order to see if these conditions are merely technical or truly a qualitative threshold for a fixed-point argument, one may look at the {first Picard iterate}, on the Fourier side:\footnote{Observe that the unboundedness of the first Picard iterate can be used to show that the KdV flow cannot be $C^2$, which implies the failure of the contraction arguments. See  \cite{bou_ill_kdv, MolTzeSaut_ill, Tzvetkov_ill}.}
	\begin{equation}\label{picard1}
\Ft_x\bigg(
	\int_0^t e^{-(t-t')\dx^3}\partial_x((e^{-t'\dx^3}u_0)^2)dt'
\bigg) (\xi)
=
e^{it \xi^3}
\int_{\xi=\xi_1+\xi_2}\frac{e^{-3it\xi\xi_1\xi_2}-1}{-3\xi_1\xi_2} \widehat{u}_0( \xi_1) \widehat{u}_0(\xi_2) d\vec{\xi}.
\end{equation}

\noi
	For $\xi_1, \xi_2$ large,
	a simple computation shows that this term is bounded in $\mathcal{F}L^{s,\infty}(\R)$ (for $u_0\in \mathcal{F}L^{s,\infty}(\R)$) if and only if $s>-\frac12$, thus going beyond the restriction imposed by the stationary interactions for \eqref{dtv}.

	The above argument shows the existence of a gap between the direct well-posedness theory in Fourier restriction spaces (which requires $s>-\frac14$ for initial data in $\mathcal{F}L^{s,\infty}(\R)$) and the range of boundedness of the first Picard iterate. There are several examples of dispersive equations where this gap occurs, such as the KP-II equation \cite{hadac, HHK_kpii}, the one-dimensional quadratic nonlinear Schrödinger equation \cite{taobejenaru, OS12}, and the NLS-KdV system~\cite{CLS}.
To bridge this gap, one can construct a space adapted to the solution and/or exploit further structure of the equation.
In the first direction,
in \cite{taobejenaru, hadac, HHK_kpii}, the authors adapt the Fourier restriction space with specific weights that already take into account the nonlinear effects present in the equation.
In the second direction, a now classical approach to improve the observable properties of a given dispersive equation goes back to the work of Babin, Ilyin, and Titi \cite{bit}, where the authors use the idea of differentiation-by-parts to show global well-posedness and Lipschitz continuity for periodic KdV.
This technique, often referred to as a (Poincaré-Dulac) \emph{normal form reduction},
has been further developed in the context of unconditional uniqueness \cite{GuoKwonOh_nls, kishimoto, KwonOh_uu_mkdv, KOY20}, local well-posedness \cite{CLS}, global well-posedness \cite{CGKO_normal_form, IfrimTataru}, and nonlinear smoothing \cite{CKV, CS20,  ErdoganTzirakis_kdv, ErdoganTzirakis_nls, McConnel_nls}.
The differentiation-by-parts technique is also a key aspect of the space-time resonances method of Germain, Masmoudi, and Shatah \cite{GMS_nls_3d}, which can be used to prove small data global well-posedness and scattering \cite{katopusateri, GMS_gravity, GermainMasmoudi_euler_maxwell,  anjolras}. 
 Other methods to expand the nonlinear terms and avoid certain problematic interactions include the power series method of Christ \cite{Christ_power}, and the Birkhoff normal form \cite{Bambusi_birkhoff, bourg_birkhoff_normal, CollianderKwonOh_birkhoff}, which relies on the Hamiltonian structure of the underlying PDE.

To bridge the regularity gap between the bilinear term in \eqref{dtv} and \eqref{picard1},
we apply this idea of differentiation-by-parts to \eqref{dtv}.
Our goal is to rewrite the equation around the problematic  high$\times$high$\to$low frequency interactions contained in set $A^c$, which include the stationary region where the obstruction of $s>-\frac14$ appears. By considering \eqref{dtv}, splitting the region of integration between non-stationary and stationary sets $A$ and $A^c$, respectively, and differentiating-by-parts in the latter only, we get
	\begin{align}
		\dt\widehat{v}(t,\xi)&= \int_{A\cup A^c}e^{-it (\xi^3 -
			\xi_1^3 - \xi_2^3) } i\xi \widehat{v}(t, \xi_1) \widehat{v}(t,\xi_2) d\vec{\xi}
		\\&
=\int_0^t \int_{A}e^{-it (\xi^3 -
			\xi_1^3 - \xi_2^3) } i\xi \widehat{v}(t, \xi_1) \widehat{v}(t,\xi_2) d\vec{\xi}  + \dt\left[ \int_{A^c}\frac{e^{-it (\xi^3 -
				\xi_1^3 - \xi_2^3) }}{-3\xi_1\xi_2}  \widehat{v}(t, \xi_1) \widehat{v}(t,\xi_2) d\vec{\xi}\right] \\&\quad- \int_{A^c}\frac{e^{-it (\xi^3 -
				\xi_1^3 - \xi_2^3) }}{-3\xi_1\xi_2}  \partial_t\left(\widehat{v}(t, \xi_1) \widehat{v}(t,\xi_2)\right) d\vec{\xi}dt
		.
\label{NF00}
	\end{align}
We can then replace \eqref{dtv} in $\dt \ft{v}$ terms in the last contribution and integrate in time, thus finding
	\begin{align}
			\widehat{v}(t,\xi)&=  \int_0^t \int_{A}e^{-it' (\xi^3 -
				\xi_1^3 - \xi_2^3) } i\xi \widehat{v}(t', \xi_1) \widehat{v}(t',\xi_2) d\vec{\xi} dt'
\\&\quad+ \bigg[ \int_{A^c}\frac{e^{-it' (\xi^3 -
					\xi_1^3 - \xi_2^3) }}{-3\xi_1\xi_2}  \widehat{v}(t', \xi_1) \widehat{v}(t',\xi_2) d\vec{\xi}\bigg]^{t'=t}_{t'=0}
                    \label{NF0}
                    \\
			&\quad
			+\frac23 i
\int_0^t
			\int_{A^c, \xi_1=\xi_{11}+\xi_{12}}\frac{e^{-it' (\xi^3 -
					\xi_{11}^3-
					\xi_{12}^3 - \xi_2^3) }}{\xi_2}  \widehat{v}(t', \xi_{11})\widehat{v}(t', \xi_{12}) \widehat{v}(t',\xi_2) d\vec{\xi} dt'
			.
	\end{align}
We identify two obstacles in the new formulation of \eqref{dtv} in \eqref{NF0}:

\smallskip
\noi{\rm(i)} The second term in \eqref{NF0}, which we refer to as the boundary term, relates to the first Picard iterate in \eqref{picard1}, and lives in $\FL^{s,\infty}(\R)$ for $s>-\frac12$ and each fixed $t$. However, due to the lack of a time integral, it is not amenable to estimates in the Fourier restriction framework.

\smallskip
\noi{\rm(ii)} The cubic integral term can be shown to have regularity $s>-1$ in the non-stationary region $(\xi_{11}, \xi_{12})\in A$, but once again we find the restriction $s>-\frac12$ when $(\xi_{11}, \xi_{12})\in A^c$.
The improved regularity of the cubic term comes from the two-derivative gain generated by the differentiation-by-parts, which introduces a factor of $|\Phi(\xi,\xi_1, \xi_2)|^{-1} \sim |\xi\xi_1\xi_2|^{-1}$ in the stationary region $(\xi_1,\xi_2)\in A^c$.

\medskip

Setting aside the boundary term and focusing on (ii), we repeat the process of isolating the stationary region, where we perform a further differentiation-by-parts to mimic the derivative gain in the cubic term, observed in the first iteration. However, the same obstruction as in (ii) appears at each step of this iteration, and thus we must perform this reduction \textit{infinitely} many times.
This requires us to proceed carefully, as the differentiation-by-parts in time hinges on the identity
\begin{align}
\label{exp}
\partial_t\frac{e^{it\Psi}}{i\Psi} = {e^{it\Psi}},
\end{align}
where $\Psi$ is the total resonance function appearing in the oscillatory integrals, which may introduce singularities when $\Psi=0$.
Namely, we argue iteratively as follows:

    \smallskip
\noi{\bf Step 1.}
    Combine all the nonlinear integral terms with the ``same" structure, that is,
with the same homogeneity in $v$ and with spatial frequency integral over the exact same convolution hyperplanes,  and lying over the stationary set.
As we will see, the resulting multiplier can always be written as
\begin{align}
\text{controllable terms}
\ + \mathcal{K}\times \Psi
.
\end{align}
The first term can be estimated for $s>-1$, while  $\mathcal{K}$ is a smoothing operator analogous to the gain of two derivatives observed in the first iteration. Being well-behaved, the first contribution is left aside.

    \smallskip
\noi{\bf Step 2.}
    Perform a differentiation-by-parts on the second term, cancelling the $\Psi$ factor and therefore avoiding the introduction of singularities. The resulting boundary term is also set aside.

    \smallskip
\noi{\bf Step 3.}
Replace any $\dt \ft{v}$ terms by \eqref{dtv} and split any new frequency integrations into stationary and non-stationary regions $A^c \cup A$. The latter can be shown to be bounded for $s>-1$ and we iterate the process for the stationary terms.

    \medskip

    In other words, we must exploit in an optimal way the cancellations between comparable nonlinear terms before exploring further time oscillations. This procedure is possible due to the simple algebraic fact that, over the stationary set $A^c$, the quotient between the multiplier for the \eqref{kdv} nonlinearity and the resonance $\Phi$ in \eqref{phase-intro} is bounded:
\begin{equation}\label{eq:divisible}
		\left|\frac{\xi}{\Phi}\right|=\frac{1}{3|\xi_1\xi_2|}\lesssim 1\qquad \mbox{over }A^c.
\end{equation}
In fact, $\mathcal{K}$ will be a product of such quotients.
		Importantly, the property in \eqref{eq:divisible}, a key part of our method, obviously holds in non-completely integrable contexts.

At each finite iteration of Steps 1-3, we see remainder terms whose estimate morally reduces to that of the bilinear term in \eqref{dtv} restricted to the stationary set $A^c$,  imposing $s>-\frac12$. Thus, we must iterate this procedure infinitely many times,  relegating the regularity obstructions to a nonlinear term of increasingly higher-order which, in the limit, formally converges  to zero.
    In the end, we will have rewritten \eqref{dtv} as an infinite normal form equation
\begin{equation}\label{eq:nfe_intro}
		\dt\widehat{v}(t,\xi)=  \sum_{j=1}^\infty \left(\dt\left(\mbox{boundary terms at step } j\right)+ \mbox{integral terms at step }j\right),
	\end{equation}
	whose integral terms are now bounded in the Fourier restriction space associated with $\mathcal{F}L^{s,\infty}(\R)$ for $s>-1$.

\begin{remark}
\rm
(i) The algebraic procedure in Steps 1-3 is \textit{independent} of the set $A$, and therefore we can always formally write a version of \eqref{eq:nfe_intro} for any choice of $A$.
Our particular choice of $A$ is motivated a posteriori by the analytical properties of nonlinear terms appearing in \eqref{eq:nfe_intro}.

\smallskip

\noi(ii) The first work introducing an infinite iteration of normal form reductions is due to Guo, Kwon, and Oh \cite{GuoKwonOh_nls} 
and has since been further developed in various contexts \cite{CS20, kishimoto, Kishimoto_nls,KwonOh_uu_mkdv, KOY20}.
In the arguments therein, in order to avoid singularities in using \eqref{exp}, the differentiation-by-parts is performed over the \emph{non-resonant set} $|\Psi|\gg 1$, which varies along the iteration process. This requires one to keep track of what occurred at every single step, and thus each nonlinear term (in the words of \cite{GuoKwonOh_nls}) ``has to know how it grew''. In contrast, by exploiting in a refined way the algebraic structure of the equation, we can perform each iteration over the \emph{same} stationary region, which allows us to combine all terms with the same structure, regardless of ``how they grew''. The fact that the problematic region does not change along the procedure proves essential in the referred algebraic cancellations and will also impact the treatment of the boundary terms (Remark \ref{remark:factorial}).

\end{remark}

We must now turn to face the boundary terms.
Due to the presence of a time derivative in front of each boundary term in \eqref{eq:nfe_intro} and of the problematic frequency interactions in each and every single term (as observed in \eqref{NF0}),  we are unable to handle them with the Fourier restriction framework needed for the integral terms. As such, we must introduce a gauge transform
    capable of absorbing \emph{all} boundary terms at once.
    First, let $\D[v,v]$ represent the  boundary term (at time) $t$ in~\eqref{NF0},
\begin{equation}
	(\D[v,v])^\wedge(\xi) =  \int_{A^c}\frac{e^{-it (\xi^3 -
			\xi_1^3 - \xi_2^3) }}{-3\xi_1\xi_2}  \widehat{v}(t, \xi_1) \widehat{v}(t,\xi_2) d\vec{\xi}.\label{eq:D_intro}
\end{equation}

	We claim that, as a consequence of the careful normal form reduction described in Steps 1-3, all boundary terms can be written as successive compositions of $\D$. Indeed, the boundary terms at the second iteration are
	$$
	-\D\big[\D[v,v],v\big]
\quad
\text{and}
\quad
-\D\big[v,\D[v,v]\big]
	$$
	the ones of the third iteration are
\begin{gather}
\D\Big[\D\big[\D[v,v],v\big],v\Big],
\quad \D\Big[\D\big[v,\D[v,v]\big],v\Big],
\quad \D\big[\D[v,v],\D[v,v]\big],
\\
\D\Big[\D\big[v, \D[v,v],v\big]\Big],\quad
\text{and}
\quad  \D\Big[v,\D\big[v,\D[v,v]\big]\Big],
\end{gather}
	and so on. The infinite sum of these compositions corresponds to the bilinear Neumann series of the formal inverse of the nonlinear map
\begin{equation}\label{eq:renorm}
		w\mapsto v=w+\D[w,w].
	\end{equation}
    This can be easily seen by succesively rewriting
	$$
	w=v-\D[w,w]
	$$
by replacing $w$ on the right-hand side by its expression above:
	\begin{align}
		w
&=v-\D[w,w]=v-\D\big[v-\D[w,w],v-\D[w,w]\big]
\\
&=v-\D[v,v]+\D\big[\D[w,w],v\big]+\D\big[v,\D[w,w]\big]-\D\big[\D[w,w],\D[w,w]\big]=\dots.
	\end{align}
	If we replace $v$ in \eqref{eq:nfe_intro} using \eqref{eq:renorm}, we find
	\begin{equation}\label{eq:nfe_intro2}
		\dt \widehat{w}(t,\xi)= \sum_{j\ge 2} (\mbox{integral terms in }w\mbox{ of order }j),
	\end{equation}
where the resulting integral terms differ from those in \eqref{eq:nfe_intro}. In fact, their algebraic complexity is significantly reduced, as \eqref{eq:renorm} also introduces cancellations at this level.

    Now, we obtain the gauged KdV equation \eqref{rkdv} from \eqref{eq:nfe_intro2} for $z(t)=e^{-t\partial_x^3}w(t)$:
    \begin{equation}\label{eq:z_intro_2}
\dt z + \dx^3 z = \sum_{j= 2}^\infty N_j[z],\quad N_j\mbox{ with homogeneity }j \mbox{ in }z. 
    \end{equation}
    For the exact form, see \eqref{z1}.
The map connecting  $u$ and $z$  is precisely the gauge transform $\Gc$, defined in \eqref{gauge}, since
    \begin{equation}
u=e^{-t\dx^3}v= e^{-t\dx^3}\left(w+\D[w,w]\right)= z + e^{-t\dx^3}\D[e^{t\dx^3}z, e^{t\dx^3}z] = \Gc^{-1}[z].
    \end{equation}

	\medskip

	Having derived \eqref{rkdv}, we move on to the proof of Theorem \ref{thm:lwp_rkdv}, which hinges on suitable multilinear bounds for each $N_j$, $j\ge 2$. This requires two more ingredients:

    \smallskip
    \noi\textit{\underline{Algebraic cancellations $\I$.}} After the change of variables $v\mapsto w$, new cancellations appear among terms with the same ``structure'', which may have risen from the replacement of $v$ either by $w$ or $\D[w,w]$ (see \eqref{eq:renorm}). We describe this new algebraic miracle in Section \ref{SUBSEC:gauge}:
for each $j\ge 2$, the number of terms of order $j$ grows exponentially, with each term having a very particular structure (see also Figure \ref{FIG:types}).

    \smallskip
    
    \noi\textit{\underline{Algebraic cancellations $\II$.}} For \eqref{kdv}, there exists a further level of cancellations in the gauged equation, 
    related to the symmetry properties of the multiplier $\frac\xi\Phi$. In the limit case $A=\emptyset$, these cancellations would kill all non-cubic terms, reducing \eqref{rkdv} to 
    the modified Korteweg-de Vries equation (see Subsection \ref{SUBSEC:miura} and Remark \ref{remark:miura} below).

    \smallskip
    
    \noi\textit{\underline{Frequency-restricted estimates.}} 
    We systematically prove multilinear bounds via frequency-restricted estimates, in the spirit of \cite{CC24, CLS, COS}. These estimates are an effective way to prove sharp bounds for nonlinear terms of arbitrary order, in both $H^s(\R)$ (where we find the classical regularity restriction $s>-\frac34$) and in $\FL^{s,\infty}(\R)$ (which imposes $s>-1$). The general case $\FL^{s,p}(\R)$, $2<p<\infty$ follows from a standard interpolation argument.

\begin{remark}
In other models, such as \eqref{gBO} discussed in Subsection~\ref{subsec:gbo}, the step of \textit{Algebraic cancellations $\II$} above does not necessarily apply.
Nevertheless, it is possible to prove frequency-restricted estimates without using this extra cancellation, which still yield improved regularity assumptions for the gauged equation, in the high-integrability setting. For instance, for \eqref{rkdv}, this already yields the regularity range
$$
s>\max\left\{-1+\tfrac1{2p}, -\tfrac34-\tfrac1{4p}\right\}, 
$$
which is sufficient to cover the full well-posedness region for \eqref{kdv} in Theorem \ref{thm:wp_kdv_intro} (See Figure \ref{fig:cancel}).
\begin{figure}[h]
		\centering
		\includegraphics[height=4.4cm]{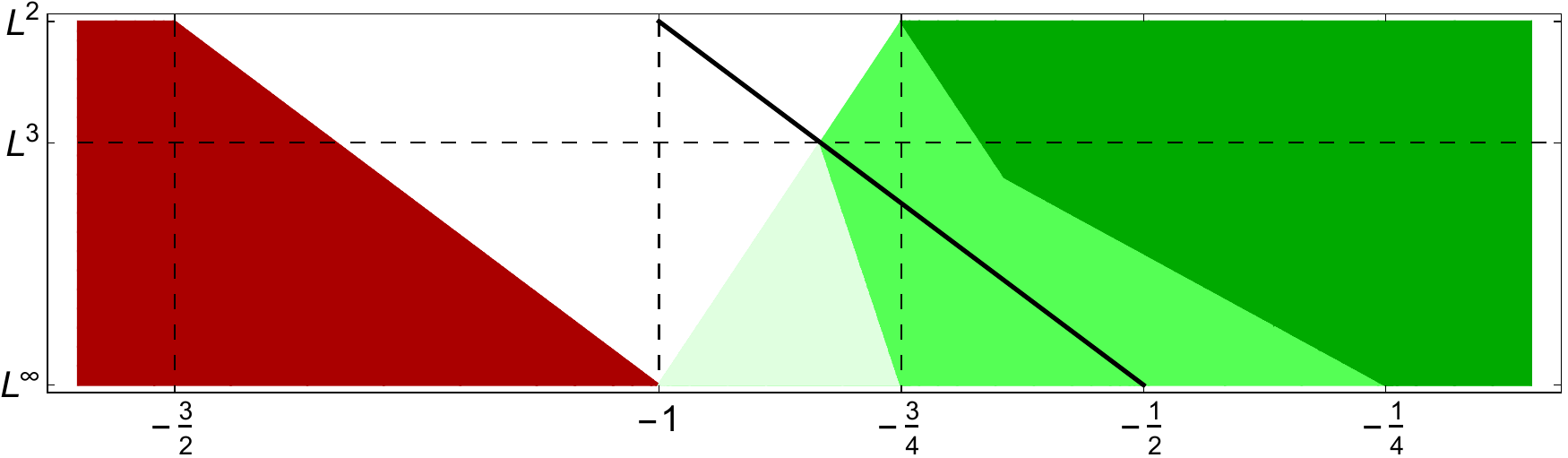}
		\caption{Well-posedness range for \eqref{rkdv} vs. algebraic cancellations. The dark-green region corresponds to the direct result by Grünrock \cite{Grunrock_ldv_hier}. The green region corresponds to the well-posedness for \eqref{rkdv}using just the \textit{Algebraic cancellations $\I$} (which do not use any particular further structure of \eqref{kdv}). The light-green region is the improved regularity range after exploiting the \textit{Algebraic cancellations $\II$}  specific to \eqref{kdv} }\label{fig:cancel}
\end{figure}
In particular, the extra algebraic cancellations are only needed to find the full well-posedness range for \eqref{rkdv} in Theorem~\ref{thm:lwp_rkdv}. 
\end{remark}

	Once we have the multilinear bounds, the small-data local well-posedness of \eqref{rkdv} (Theorem \ref{thm:lwp_rkdv}) follows from a standard fixed-point argument in Fourier restriction spaces. Afterwards, using the properties of the gauge transform (Section \ref{sec:gauge_prop}), we go back to \eqref{kdv}. This requires us to prove that, at high regularity, \eqref{kdv} and \eqref{rkdv} are equivalent, with local time of existence depending only on the smallness of the weaker $\FL^{s,p}$-norm. This is achieved by first truncating the nonlinearity, then showing the equivalence between the truncated equations, and finally proving the convergence in the truncation parameter. Once this equivalence is settled, Theorem \ref{thm:wp_kdv_intro} follows via a standard density argument.

\begin{remark}
\rm

(i) Due to the  algebraic cancellations discussed above, \eqref{rkdv} has an exponential growth of the number of nonlinear terms of order $j\ge 2$. Consequently, we impose smallness of the initial data $z_0$ in Theorem~\ref{thm:lwp_rkdv} to counter-act this growth, not due to any criticality of the functional space.
\\
Alternatively, we could add a sufficiently small constant to the definition of the stationary sets $A^c$, depending on the size of the initial data, and obtain a version of Theorem~\ref{thm:lwp_rkdv} for large data $z_0$. 
We choose to present the derivation with the small data assumption for clarity of exposition.

\smallskip
\noi
(ii) In the infinite normal form reduction in \cite{GuoKwonOh_nls, kishimoto}, as the nonlinear terms in the equation need to encode ``how they grew'', the number of terms of order $j$ grows factorially with $j$. Consequently, to obtain summability in $j$, a smallness condition on the initial data is insufficient (as this gives only exponential decay). Instead, in setting the condition on the phase $|\Psi| \gg 1$ at the $j$-th iteration, they require the underlying constants to grow sufficiently fast with $j$.
This prevents an analogue of our algebraic miracle for the boundary terms, as the latter relies crucially on the restriction of all frequencies to the same set $A^c$.

\end{remark}

\begin{remark}\label{remark:dgbo} In the case of \eqref{gBO}, 
since the quadratic nonlinearity still satisfies the divisibility property \eqref{eq:divisible}, we can define the gauge transform $\Gc$ analogously to  \eqref{gauge}, and obtain an analogue of \eqref{rkdv} 
and associated well-posedness results. 
 The \textit{Algebraic cancellations $\I$} hold with no further modification, reinforcing the level of applicability of our approach. Also, the frequency-restricted estimates (Section \ref{sec:multi}) can be easily adapted to \eqref{gBO}, as they rely on the homogeneity of the phase function, instead of its specific formula. As such, one must simply adapt the proofs for $a=2$  to the general $a\in(1,2)$ case (see Subsection \ref{subsec:gbo}).

\end{remark}

    \subsection{Connection between the gauge transform and the Miura map}\label{SUBSEC:miura}
		As observed by Miura \cite{Miura}, there is a strong link between \eqref{kdv} and the modified KdV~equation:
		\begin{equation}
			\tag{mKdV}\label{mkdv}
			\partial_t z + \partial_x^3 z = \tfrac29 \partial_x(z^3).
		\end{equation}
		More precisely, if $z$ is a solution to \eqref{mkdv}, then
		\begin{equation}\label{eq:defiMiura}
			\MMb[z]=\partial_x z +\tfrac{1}{3}z^2
		\end{equation}
		is a solution to \eqref{kdv}. 
        The transform $\MMb$, often called the Miura map, can be used to transfer results from \eqref{kdv} to \eqref{mkdv}: an example is the global well-posedness for \eqref{mkdv} in $H^s(\R)$, $s>\frac14$, which is a consequence of the corresponding result for \eqref{kdv} in $H^s$, $s>-\frac34$, proved by Colliander, Keel, Staffilani, Takaoka, and Tao \cite{CKSTT_sharp_kdv}. Unfortunately, the opposite direction is considerably harder due to the lack of invertibility of $\MMb$. In the $H^s$-setting, this was thoroughly studied in  \cite{BuckmasterKoch, KPST_miura}, relying on spectral analytic tools (in connection to the complete integrability of~\eqref{kdv}).

		Since the Miura map formally satisfies
		\begin{equation}\label{eq:miura2}
			{\MMb[\partial_x^{-1}z]}^\wedge(\xi) = \widehat{z}(\xi)-\frac{1}{3}\int_{\xi=\xi_1+\xi_2} \frac{\widehat{z}(\xi_1)\widehat{z}(\xi_2)}{\xi_1\xi_2}d\xi_1
            , 
		\end{equation}
         we have the following agreement with the inverse gauge transform $\Gc^{-1}$:
        \begin{equation}
\MMb[\dx^{-1} z]=\Gc^{-1}[z]\quad \mbox{ when }A=\emptyset\mbox{ in }\eqref{gauge}.
        \end{equation}

\noi
Both the Miura map $\MMb$ and $\Gc^{-1}$ relate \eqref{kdv} with another KdV-type equation, and each have their own advantages.
On the one hand, $\MMb$ is connected with the integrable structure of \eqref{kdv} and links it to \eqref{mkdv}, which has a single cubic nonlinearity. However, it fails to be invertible due to the low frequencies $|\xi_1|, |\xi_2|\ll 1$.
 This defect was corrected in \cite{CCT_ill}, via a generalized Miura map, relating \eqref{kdv} to a coupled KdV-mKdV system. 
In contrast, $\Gc^{-1}$ is invertible and independent of the underlying integrable structure, but connects \eqref{kdv} to a far more complex equation (not a system), with nonlinear terms of any order $j\ge2$, which are nevertheless amenable to improved estimates.

        \begin{remark}\label{remark:miura}
For $A=\emptyset$ in \eqref{gauge}, the gauge transform would eliminate all non-cubic terms in \eqref{rkdv} and reduce it to (an equivalent form of) \eqref{mkdv}. This is the miracle responsible for the \textit{Algebraic cancellations $\II$} for \eqref{kdv}, which is potentially lost for perturbed equations.
Hence, we can see the Miura map $\MMb$ in \eqref{eq:miura2} as a particular case of the (inverse) gauge transform $\Gc^{-1}$ coming from the normal form reduction described above, which remains valid in non-integrable cases.

        \end{remark}

\subsection{Organization and notation} This manuscript is organized as follows: in Section \ref{sec:infr}, we introduce the tree formalism (Subsection \ref{SUB:tree}) that will guide our infinite normal form reduction and make the \textit{Algebraic cancellations $\I$} clear, derive the normal form equation for \eqref{kdv} (Subsection \ref{sec:nfe}), and apply the gauge transform $\Gc$ to obtain the precise form of the gauged equation \ref{rkdv} (Subsection \ref{SUBSEC:gauge}). 
The \textit{Algebraic cancellations $\II$} are shown in Section~\ref{sec:cancel_kdv}.
In Section \ref{sec:functional}, we introduce the Fourier restriction space framework (Subsection \ref{sec:fourier_restr}) and reduce the multilinear estimates to frequency-restricted estimates (Subsection \ref{sec:fre}), showing the latter in Section \ref{sec:multi} for \eqref{rkdv}. 
Section \ref{sec:lwp} is devoted to the properties of the gauge transform, the equivalence between \eqref{kdv} and \eqref{rkdv} (for smooth solutions) and to the proofs of Theorems \ref{thm:lwp_rkdv} and \ref{thm:wp_kdv_intro}.

\bigskip

We conclude this introduction with some notation that will be used throughout this work. Given real numbers $A,B$, we write $A \les B$ if there exists $C>0$ such that $A \le CB$, and $A \ll B$ if $0<C< 10^{-10}$. If $A\les B$ and $B\les A$, we say that $A\sim B$. Additionally, $|A| \simeq |B|$ if $\max(|A|,|B|) \gg \big| |A| - |B| \big|$.
We often use $A \land B = \min(A,B)$ and $A \lor B = \max(A,B)$.
Lastly,
the notations $\al +$ and $\al-$ stand for $\al+\eps$ and $\al-\eps$ for any $\eps>0$, respectively.

\subsection{Acknowledgements} A.C. was partially supported by  CNRS-INSMI through a grant ``PEPS Jeunes chercheurs et jeunes chercheuses 2025". S.C. was partially supported by Funda\c{c}\~ao para a Ci\^encia e Tecnologia (FCT), through CAMGSD (project UID/04459/2025). A.C. and S.C. were also supported by FCT, through the grant FCT/Mobility/1330917020 /2024-25). A.C. is thankful for the kind hospitality of Instituto Superior Técnico, where part of this work was developed. S.C. and J.P.R. were also supported by FCT through project SHADE (project 2023.17881.ICDT, 
DOI 10.54499/ 2023.17881.ICDT).

\section{Infinite normal form equation and gauge transform}\label{sec:infr}

In this section, we describe the procedure to derive the \eqref{rkdv} equation:
first, one writes an infinite normal form equation for \eqref{kdv} (that is, after infinitely
many differentiations-by-parts in time) and then, applying the gauge transform $\Gc$ and exploiting further cancellations, one arrives at \eqref{rkdv}.

\subsection{Tree formalism}
\label{SUB:tree}

To properly describe the normal form procedure, we start by introducing suitable combinatoric structures, namely  binary trees, which we use to track the reduction process and index the terms in the equation.
We also refer to \cite{Bruned, Christ_power, GuoKwonOh_nls,kishimoto, KOY20} and references therein for further details on the use of (decorated) trees to describe infinite expansions in the context of dispersive PDEs.

In this subsection, we recall relevant notations related to binary trees, which we use to describe our infinite normal form equation.

\begin{definition}
\label{DEF:tree}
\rm
A \emph{binary tree} $\TT$ is a partially ordered set (with partial order $\le$) satisfying the following properties:

\begin{enumerate}
\item[(i)] $\TT$ has a unique minimal element called the \textit{root node} labelled by $\emptyset$ (i.e., $\emptyset \le a$ for all $a \in \TT$).

\smallskip
\item[(ii)] Given $a\in \TT\setminus \{ \emptyset \}$,
there exists a unique $b \in \TT$ such that: $b\ne a$, $b \le a$, and $b \le c \le a$ implies $c\in\{a,b\}$. In that case, we say that $b$ is the \textit{parent} of $a$ and $a$ is a \textit{child} of $b$.

\smallskip
\item[(iii)] A node $a\in \TT$ is called \textit{terminal} or a \textit{leaf} if it has no children. A \textit{non-terminal node} or a \textit{parent} $a\in\TT$ is a node with exactly two children, one called \textit{left child} and the other \textit{right child}.\footnote{Introducing the notation of left and right children of a node  reflects the lack of ambiguity in the usual planar representation of binary trees.}
We write $\TT^0, \TT^\infty$ for the sets of parent and terminal nodes of $\TT$, respectively.  Note that $\TT$ is the disjoint union of $\TT^0$ and $\TT^\infty$.

\smallskip
\item[(iv)]
Each node in $\TT$ is labelled with a word in the following way: the root node has word~$\emptyset$; and for a given parent node $a\in \TT^0$, its left and right children have labels $a1$ and $a2$,  obtained by concatenating $1$ or $2$ to the right of $a$, respectively. Note that children of the root node have labels 1 and 2.

\end{enumerate}

\noi
We interchangeably use $\TT$ to denote the representation of the tree as a finite connected graph and to denote the set of words used to label its nodes. Note that the latter is also a partially ordered set, with words $a,b$ satisfying $a\le b$ if there exists a word $c$ such that $b = ac$.

\end{definition}

Given $j\in\N$, we use the notation $\Tr_j$ to denote the family of all binary trees
with $j$ non-terminal nodes, which we call the $j$-th generation of trees. From Definition~\ref{DEF:tree}, it follows that
\begin{align}
| \TT^0 |  = j \quad \text{and} \quad |\TT^\infty| = j+1.
\end{align}

\noi
It is not difficult to see that the number of distinct trees in $\Tr_j$ is given by the $j$-th Catalan number
\begin{align}\label{eq:catalan}
| \Tr_j | = \frac{1}{j+1} \binom{2j}{j}  = \frac{(2j)!}{(j+1)! j !}.
\end{align}
In particular, the first 3 generations of trees are succintly written below, identified by their graphs:
\begin{align*}
\Tr_1 & = \big\{ \<1> \big\}, \\
\Tr_2 & = \left\{ \<21>, \<22> \right\}, \\
\Tr_3 &  = \left\{ \<31>, \<32>, \<33>, \<34>, \<35> \right\}.
\end{align*}
Similarly, we can represent the trees above by the following partially ordered sets of words used to label their nodes (as per Definition \ref{DEF:tree}(iv)):
\begin{align}
\<1>
&
= \{\emptyset\}
,
\\
\<21>
&
=
\{ \emptyset, 1 ,2, 11,12  \}
,
\\
\<22>
&
=
\{\emptyset, 1,2, 21,22\}
,
\\
\<31>
&
=
\{
\emptyset, 1,2, 11,12, 111, 112
\}
, \ldots
\end{align}

\smallskip

For readibility, we will reserve the symbols $\bul$, $\star$, and $\ast$ to denote nodes in a tree and their respective associated word.
Next, we introduce some additional notation to identify different nodes in a tree.

\begin{definition}\rm
\label{DEF:tree2}
Let $\TT$ be a binary tree in the sense of Definition~\ref{DEF:tree} and
$\bul, \star$ be two nodes in $\TT$.

\begin{enumerate}

\smallskip
\item[(i)]
If $\bul \in \TT\setminus \{ \emptyset\}$, we write $\pb(\bul)$ to denote the \emph{parent node} of $\bul$.

\smallskip
\item[(ii)]
If $\bul, \star\in \TT \setminus\{ \emptyset \}$ are distinct and share the same parent, they are called \textit{siblings}. We write $\sbf(\bul)$ for the sibling of $\bul$.

\smallskip
\item[(iii)]
Given $\bul\in \TT^0$, we say $\bul$ is a \textit{final parent} if its children $\bul1, \bul2 \in \TT^\infty$. The set of final parents is denoted by $\TT^{0,\ff}$.
.

\end{enumerate}

\end{definition}

\begin{remark}\rm
Note that from Definition~\ref{DEF:tree}, the maps $\pb, \sbf$, for the parent and sibling of a given node, are well-defined in their respective domains, since if either the parent or sibling node exists, it is unique.

\end{remark}

We next introduce concatenation and subtraction of trees.

\begin{definition}
\label{DEF:tree3}
\rm

Let $\TT$ be a binary tree in the sense of Definition~\ref{DEF:tree}.

\begin{enumerate}

\smallskip
\item[(i)] (Concatenation)
Given $ \star \in \TT^\infty$, we define
$$
\TT \overset{\star}{\oplus} \<1>
: = \TT \cup \{\star1, \star2\}
$$
as the result of adding the children $\star1,\star2$ of $\star$ to the tree $\TT$.
Note that this operation now makes $\star$ a final parent in the tree $\TT \overset{\star}{\oplus} \<1>$.

\smallskip
\item[(ii)] (Subtraction)
Given a final parent $\star  \in \TT^{0,\ff}$, we define the tree
$$
\TT \overset{\star }{\ominus } \<1>
: = \TT \setminus \{ \star1, \star2  \}
$$
as the result of eliminating the children $\star1, \star2$ of $\star$
from $\TT$. Note that $\star$ becomes a terminal node in the new tree $\TT \overset{\star}{\ominus } \<1>$.
More generally, given $C\subseteq \TT^{0,f}$, we define
\begin{align}
\label{treeSubtr}
\TT \overset{C}{\ominus} \<1> : =
 \bigcap_{\star \in C} ( \TT \overset{\star}{\ominus}  \<1>)
= \bigcap_{\star\in C} (\TT \setminus \{\star1,\star2\})
= \TT \setminus \Big( \bigcup_{\star\in C} \{\star1,\star2\} \Big),
\end{align}
which corresponds to the set where all the children of the final parents $\star\in C$ have been removed from the tree.

\end{enumerate}

\end{definition}

The tree structure defined above and its word-labelling will help in keeping track of the nonlinear terms appearing in our normal form reduction procedure. Namely, after an differentiation-by-parts step, a time derivative falls on a $v$-factor as in \eqref{NF00}, after which we use equation \eqref{dtv} to replace it by a quadratic term, generating a higher degree nonlinear contribution. We can then use binary trees to label the frequencies of the corresponding factors, in the following way.

\begin{definition}\rm

Let $j\in\N$ and $\TT \in \Tr_j$.

\begin{enumerate}
\item[(i)] (Index function)
We define the index function $\vec\xi : \TT\to \R^\TT$ such that $\vec\xi = (\xi_\bul)_{\bul\in\TT}$, where
\begin{align}
\label{parent}
\xi_\bul = \xi_{\bul1} + \xi_{\bul2}, \quad \bul\in \TT^0,
\end{align}
and $\bul1$ and $\bul2$ denote the children of $\bul$. Moreover, we define the set $\G(\TT)$ as the set of convolution relations on the tree $\TT$:
\begin{align}
\label{Gamma}
\G(\TT) : = \big\{ \vec\xi \in \R^\TT: \ \xi_\bul = \xi_{\bul1} + \xi_{\bul2}, \ \bul \in \TT^0 \big\}.
\end{align}

\smallskip

\item[(ii)] (Phase function)
We use $\Psi(\TT)$ to denote the \textit{phase function} associated with the tree $\TT$, given by
\begin{align}
\label{psi-tree}
\Psi(\TT)
:=
- \xi^3 + \sum_{\bul\in \TT^\infty} \xi_\bul^3,
\end{align}
where $\TT^\infty$ denotes the set of leaves of $\TT$. Given a parent node $\star\in\TT^0$, we also introduce the partial phase function for the tree $\cherry{\star}{}{}=\{\star, \star1, \star2\}$:
\begin{align}
\label{psi-circ}
\Psi(\star) = \Psi(\cherry{\star}{}{}) : = - \xi^3_\star + \xi^3_{\star1} + \xi^3_{\star2}.
\end{align}
Note that the following additivity property holds for $\Psi$:
\begin{equation}\label{eq:additivity}
\Psi(\TT) = \sum_{\star \in \TT^0} \Psi(\star).
\end{equation}

\noi
Lastly, for $\star\in \TT^0$, we define the multiplier $K_\star$ given by
\begin{equation}
\label{Kmul}
K_{\star} := \frac{\xi_{\star}}{\Psi(\star)} = - \frac{1}{3 \xi_{\star1} \xi_{\star2}},
\end{equation}
where the last equality holds under the condition in \eqref{parent}.

\end{enumerate}

\end{definition}

To describe the infinite normal form reduction, we introduce a correspondence between binary trees and oscillatory integrals.
Given $j\in\N$, a binary tree $\TT \in \Tr_j$, a multiplier $m:\R^{\TT}\to \bbC$ on $\TT$, and a
function $v: \R\times \R \to\R$, we define the associated integral term
\begin{equation}\label{eq:defiI}
\Ical(\TT, m; v)(t) :=   i \int_{\G(\TT)}
e^{it\Psi(\TT)}m(\vec{\xi})
\prod_{\bul\in \TT^\infty}  \ft v_\bul (t) d\vec{\xi},
\end{equation}
where
we use the short-hand notation  $\ft v_\bul(t) := \ft v(t,\xi_\bul)$ for a label $\bul$.

\subsection{Deriving the normal form equation}\label{sec:nfe}

Using the tree formalism in Subsection~\ref{SUB:tree}, we now derive the normal form equation for \eqref{kdv}.
We apply a normal form reduction in the form of ``differentiation-by-parts" as in \cite{bit}, only in a problematic region where we cannot show the nonlinear estimate at low regularity
(namely, in the region that imposes $s>-\frac14$ in the bilinear estimate in Lemma~\ref{LEM:bil0}).
In particular, we consider the following sets $A$ and $A^c = \G({\<1>}) \setminus A$:
\begin{align}
\label{A}
\begin{split}
A & := \big\{ (\xi_{1}, \xi_{2}) \in \R^2: \ \xi_{1}+ \xi_{2} = \xi, \  |\xi_{1}| \land |\xi_{2}| \le 1  \big\}, \\
A^c & := \big\{ (\xi_{1}, \xi_2) \in \R^2: \ \xi_{1}+ \xi_{2} = \xi,  \ |\xi_{1}| ,  |\xi_{2}| > 1  \big\},
\end{split}
\end{align}
where $A^c$ denotes the region where we perform a reduction.

Using the tree-indexing in Subsection~\ref{SUB:tree} and the notations in \eqref{psi-tree}, \eqref{psi-circ}, and \eqref{Kmul}, we rewrite \eqref{dtv} as follows
\begin{align}
\label{duhamel2}
\dt\ft v(t) =  i \int_{\xi = \xi_1 + \xi_2} e^{it \Psi(\<1>)}  K \Psi(\<1>)   (\widehat{v}_{1} \widehat{v}_{2}) (t)
d\xi_{1}.
\end{align}
Then, splitting the region of integration into $A\cup A^c$ as in \eqref{A}, we may write \eqref{duhamel2} as
\begin{align}
\dt\ft v
&
=
\Ical\big( \<1>, \ind_A K \Psi(\<1>); v
\big)
+ \Ical\big( \<1>, \ind_{A^c} K \Psi(\<1>) ;v
\big)
\\
&
=: \Ft_x[ \NN ( \<1>; v )] + \Ft_x[ \MM(\<1>; v) ]
\label{NFE000}
\end{align}
where $\Ical$ is as in \eqref{eq:defiI}. Performing a differentiation-by-parts in time on the $\MM$-term restricted to the bad region $A^c$ using the relation
$$
\dt\frac{e^{it\Psi(\<1>)}}{i\Psi(\<1>)}=e^{it\Psi(\<1>)},
$$
we obtain
\begin{align}
\label{NFE01}
\begin{split}
\dt\ft v(t)
&=
\Ft_x[ \NN ( \<1>; v )] (t)
+ \dt\bigg[ \int_{\xi = \xi_1 + \xi_2} e^{it \Psi(\<1>)}
[\ind_{A^c}
 K ]
(\ft v_{1} \ft v_{2} )(t)
d\xi_{1}
\bigg]
\\
&
\quad
-
   \int_{\xi = \xi_1 + \xi_2} e^{it \Psi(\<1>)}
[ \ind_{A^c} K ]
\partial_{t} (\ft v_{1} \ft v_{2} )(t)
d\xi_{1}
\\
& =
\Ft_x[ \NN ( \<1>; v )] (t)
+  \dt\bigg[ \int_{\xi = \xi_1 + \xi_2} e^{it \Psi(\<1>)}
[\ind_{A^c}
 K ]
(\ft v_{1} \ft v_{2} )(t)
d\xi_{1}
\bigg]
\\
&
\quad
- i
 \int_{\substack{\xi = \xi_1 + \xi_2 \\ \xi_1 = \xi_{11}+ \xi_{12}  }} \exp\Big(it \Psi(\<21>) \Big)
[ \ind_{A^c} K K_1 \Psi(1) ]
 (\ft v_{11} \ft v_{12} \ft v_2 )(t)
d\xi_{1} d\xi_{11}
\\
&
\quad
-i  \int_{\substack{\xi = \xi_1 + \xi_2 \\ \xi_2 = \xi_{21}+ \xi_{22}  }} \exp\Big(it \Psi(\<22>) \Big)
[ \ind_{A^c} K K_2 \Psi(2) ]
 (\ft v_{1} \ft v_{21} \ft v_{22} )(t)
d\xi_{1} d\xi_{21}
\end{split}
\end{align}
where we replaced $\dt \ft v_j$ by \eqref{duhamel2} to obtain the last equality. This can be rewritten as
\begin{align}
\label{NFE00}
\begin{split}
\dt\ft v
&
=
 \Ft_x[ \NN ( \<1>; v )]
+
\dt\Ical \big( \<1>,  - i \ind_{A^c} K; v \big)
\\
&
\quad
+
\Ical \Big( \<21>, - \ind_{A^c} KK_1\Psi(1) ; v  \Big)
+
\Ical \Big( \<22>, - \ind_{A^c} KK_2\Psi(2) ; v  \Big)
.
\end{split}
\end{align}

In order to iterate the procedure taking \eqref{duhamel2} into \eqref{NFE00}, we use the following lemma, which provides a systematic way to perform a normal form reduction on integral terms of the form $\Ical$ in \eqref{eq:defiI}. This translates to an operation on trees (based on tree concatenation in Definition~\ref{DEF:tree3}~(ii)) and suitable modifications of the multipliers.

\begin{lemma}[Differentiation-by-parts on trees]
\label{LEM:IBP}
Let $j\in\N$, $\TT\in\Tr_j$, $m: \R^{\TT}\to \bbC$ a multiplier on $\TT$,
and $v$ be a smooth function such that $\ft v$ satisfies  \eqref{dtv}. Then, the following identity holds
\begin{align}\label{IBP}
\begin{split}
\Ical (\TT,m; v)
= \dt \Ical \Big(\TT,\frac{m}{ i \Psi(\TT)} ; v \Big)
+ \sum_{\star \in \TT^\infty} \Ical
\Big(\TT \overset{\star}{\oplus} \<1>
, - K_\star \Psi(\star) \frac{m}{\Psi(\TT)} ; v \Big)
,
\end{split}
\end{align}
where $\Ical$ is the integral  functional in \eqref{eq:defiI}, $\TT \overset{\star}{\oplus} \<1> $ denotes tree concatenation at the final node $\star\in \TT^\infty$ as in Definition~\ref{DEF:tree3}~(ii), and $\Psi(\TT), \Psi(\star),  K_\star$ are as in \eqref{psi-tree}, \eqref{psi-circ}, and~\eqref{Kmul}.

\end{lemma}
\begin{proof}
Fix $j\in\N$, $\TT\in\Tr_j$, the multiplier $m$, and $v$ solving \eqref{dtv}.
From \eqref{eq:defiI}, we have
\begin{align}
\Ical(\TT, m; v)(t)
&
= i  \int_{\G(\TT)} e^{it\Psi(\TT)} m(\vec \xi) \prod_{\bul\in \TT^\infty} \ft v_\bul (t) d\vec\xi
\\
& =  i  \int_{\G(\TT)} \dt \bigg( \frac{e^{it\Psi(\TT)}}{ i \Psi(\TT)} \bigg) m(\vec \xi) \prod_{\bul\in \TT^\infty} \ft v_\bul (t) d\vec\xi \\
& = \dt\bigg[ i \int_{\G(\TT)} e^{it \Psi(\TT)} \frac{m(\vec\xi)}{i \Psi(\TT)}  \prod_{\bul\in \TT^\infty} \ft v_\bul (t) d\vec\xi  \bigg]
\\
& \qquad - \int_{\G(\TT)} e^{it \Psi(\TT)} \frac{m(\vec\xi)}{\Psi(\TT)} \partial_{t}  \bigg(\prod_{\bul\in \TT^\infty} \ft v_\bul \bigg)(t) d\vec\xi.
\end{align}
From \eqref{dtv}, we have
\begin{align}
\Ical(\TT, m; v)(t)
& = \dt \Ical\Big( \TT, \frac{m}{i \Psi(\TT)} ;v \Big)(t)\\
& \quad - \sum_{\star \in \TT^\infty} i  \int_{\G ({\TT \overset{\star}{\oplus} \cherryS})} e^{it[\Psi(\TT) + \Psi(\star) ]} K_\star \Psi(\star) \frac{m(\vec\xi) }{\Psi(\TT)} \bigg( \ft v_{\star1} \ft v_{\star2} \prod_{\bul\in\TT^\infty\setminus\{\star\}} \ft v_\bul \bigg)(t) d\vec\xi
\\
& =
\dt \Ical\Big( \TT, \frac{m}{i \Psi(\TT)} ;v \Big) (t)
+
\sum_{\star\in \TT^\infty}
 \Ical
\Big( \TT \overset{\star}{\oplus} \<1>,  - K_\star \Psi(\star) \frac{m }{\Psi(\TT)} ; v \Big)(t),
\end{align}
completing the proof.
\qedhere

\end{proof}

\noi
As in \eqref{NFE01}, we split each region of frequency integration between a good and a bad set, and only perform a differentiation-by-parts when restricted to the latter. Generalizing the definition of $A, A^c$ in \eqref{A}, given a tree $\TT$ and a parent node $\bul \in \TT^0$, we define the analogous sets $A_\bul$ and $A_\bul^c$ as
\begin{align}
A_\bul&
:=
 \big\{ (\xi_{\bul1}, \xi_{\bul2}) \in \R^2: \ \xi_\bul =
\xi_{\bul1}
+ \xi_{\bul2}, \ |\xi_{\bul1}| \land  |\xi_{\bul2}|  \le 1
 \big\} , \\
A_\bul^c & := \big\{ (\xi_{\bul1}, \xi_{\bul2}) \in \R^2: \ \xi_\bul =
\xi_{\bul1} +
\xi_{\bul2},\
|\xi_{\bul1}| ,  |\xi_{\bul2}| > 1
\big\}
.
\label{defA}
\end{align}
Note that, when $\bul = \emptyset$ is the root node, $A_\bul =A$ as in \eqref{A}.
We will often omit the dependence on the children frequencies $\xi_{\bul 1}, \xi_{\bul2}$ of $\xi_\bul$, and use the short-hand notation $\ind_{A_\bul}=\ind_{A_\bul} (\xi_{\bul1},\xi_{\bul2})$.

Returning to \eqref{NFE00}, we continue to derive the second iteration of the normal form equation.
Unlike the first normal form reduction from equation \eqref{duhamel2} to \eqref{NFE00}, where we split the integral term into $A$ and $A^c$ (see \eqref{A}) and differentiate-by-parts in the latter region, we now have an additional step, which we describe in the following. Considering the two last terms in \eqref{NFE00}, we first split their regions of integration $(\xi_{\bul1}, \xi_{\bul2}) \in A_\bul \cup A_\bul^c$, where $\bul = 1,2$, respectively:
\begin{align}
\label{NFE001}
\dt\ft v
&
=
 \Ft_x[ \NN ( \<1>; v )]
+
\dt \Ical \big( \<1>, -i  \ind_{A^c} K; v \big)
\\
&
\quad
+
\Ical \Big( \<21>, - \ind_{A^c \cap A_1} KK_1\Psi(1) ; v  \Big)
+
\Ical \Big( \<21>, - \ind_{A^c \cap A_1^c} KK_1\Psi(1) ; v  \Big)
\\
&
\quad
+
\Ical \Big( \<22>, - \ind_{A^c \cap A_2} KK_2\Psi(2) ; v  \Big)
+
\Ical \Big( \<22>, - \ind_{A^c \cap A_2^c} KK_2\Psi(2) ; v  \Big) .
\end{align}
As before, we keep the integral contributions localized to a ``good set" $A_\bul$. Regarding the contributions localized to $A_\bul^c$ for all $\bul\in \TT^0$, despite the restriction to the ``bad regions'', we will keep a suitable part of these terms before performing another differentiation-by-parts. In particular, by using the additivity of the phase $\Psi$ in \eqref{eq:additivity}, we write
\begin{align}
\Psi(\star)
= \big[ \Psi(\star) + \Psi(\<1>) \big]- \Psi(\<1>)
=
\Psi(\TT) - \Psi(\<1>),
\quad (\star, \TT)\in \left\{ (1, \<21>) , (2, \<22>) \right\},
\end{align}
and keep the contributions with partial phase $-\Psi(\<1>)$, only differentiating-by-parts the terms with total phase $\Psi(\TT)$. Then, from \eqref{NFE001}, we obtain the second iteration of the normal form equation
\begin{align}
\dt\ft v
&
=
 \Ft_x[ \NN ( \<1>; v )]
+
\dt \Ical \big( \<1>, - i  \ind_{A^c} K; v \big)
\\
&
+
\Ical \Big( \<21>, -[ \ind_{A^c \cap A_1} \Psi(1) - \ind_{A^c \cap A_1^c} \Psi(\<1>)  ]KK_1 ; v  \Big)
+
\Ical \Big( \<21>, - \ind_{A^c \cap A_1^c} KK_1\Psi(\<21>) ; v  \Big)
\\
&
+
\Ical \Big( \<22>, - [\ind_{A^c \cap A_2} \Psi(2) - \ind_{A^c \cap A_2^c} \Psi(\<1>)  ] KK_2; v  \Big)
+
\Ical \Big( \<22>, - \ind_{A^c \cap A_2^c} KK_2\Psi(\<22>) ; v  \Big)
\\
&
= :
 \Ft_x[ \NN ( \<1>; v )]
+
 \Ft_x[
\BB(\<1>;v)
]
\\
&
\quad
+
 \Ft_x[
\NN(\<21>; v )
]
+
 \Ft_x[
\MM(\<21>; v )
]
+
 \Ft_x[
\NN(\<22>; v)
]
+
 \Ft_x[
\MM(\<22>;v)
]
.
\label{NFE002}
\end{align}
When writing the Duhamel formulation for $v$, based on \eqref{NFE002}, $\NN$ will correspond to time-integral terms, $\BB$ to boundary terms, and $\MM$ are remainders, to which we apply \eqref{IBP} at the next step. Thus, in the following, we refer to them as integral terms, boundary terms, and remainders, respectively.

In the following iterations of the normal form equation, as in \eqref{NFE00} and \eqref{NFE002}, for each tree $\TT$, we want to have a unique integral term $\NN(\TT)$ and boundary term $\BB(\TT)$,
where this must encode all the possible ways of arriving at such a tree structure.
From Lemma~\ref{LEM:IBP}, we note that given a tree $\TT \in \Tr_j$ for $j\ge2$, after differentiation-by-parts, the resulting contributions associated with the tree $\TT$ can come from any term of the form $\Ical(\TT \overset{\star}{\ominus} \<1>)$ where $\star \in \TT^{0, \ff}$, i.e., that the children of $\star$ are leaves in $\TT$. Consequently, in the following, we need to perform differentiation-by-parts on the remainder terms $\MM(\TT)$ simultaneously and combine any contributions associated with the same tree before preparing the integral terms to be kept and those which will be differentiated in the following step. The simplest example appears in the third iteration of the normal form reduction, with the contributions associated with the tree $\<33>$, for which we need to combine the terms coming from $\<21>$ and $\<22>$ in \eqref{NFE002}.
In particular, from \eqref{NFE002} and \eqref{IBP}, we have
\begin{align}
&
 \Ft_x[
\MM(\<21>;v)
]
+
 \Ft_x[
\MM(\<22>;v)
]
\\
&
=
\dt \Ical\left(\<21>, i \ind_{A^c \cap A_1^c} KK_1 ; v\right)
+
\dt \Ical\left(\<22>, i \ind_{A^c \cap A_2^c} KK_2 ; v\right)
\\
&
\quad
+ \Ical\left( \<31> ,  \ind_{A^c\cap A_1^c} KK_1 K_{11} \Psi(11)   ; v \right)
+ \Ical\left( \<32> ,  \ind_{A^c\cap A_1^c} KK_1 K_{12}  \Psi(12) ; v \right)
\\
&
\quad
+ \Ical\left( \<33> ,  [ \ind_{A^c\cap A_1^c} \Psi(2)+ \ind_{A^c \cap A_2^c} \Psi(1) ] KK_1 K_2    ; v \right)
\\
&
\quad
+ \Ical\left( \<34> ,  \ind_{A^c \cap A_2^c } K K_2 K_{21} \Psi(21) ; v \right)
+ \Ical\left( \<35> , \ind_{A^c \cap A_2^c } K K_2 K_{22}  \Psi(22)   ; v \right).
\end{align}
After decomposing the integral terms as in \eqref{NFE002} and using the additivity property \eqref{eq:additivity},
\begin{align}
&
 \Ft_x\left[
\MM\left(\<21>;v\right)
\right]
+
 \Ft_x\left[
\MM\left(\<22>;v\right)
\right]
\\
&
=
\dt \Ical\left(\<21>, i \ind_{A^c \cap A_1^c} KK_1 ; v\right)
+
\dt\Ical\left(\<22>, i \ind_{A^c \cap A_2^c} KK_2 ; v\right)
\\
&
\quad
+ \Ical\left( \<31> ,  [\ind_{A^c\cap A_1^c \cap A_{11}} \Psi(11) - \ind_{A^c\cap A_1^c \cap A_{11}^c} \Psi(\<21>) ] KK_1 K_{11}  ; v \right)
\\
&
\quad
+
 \Ical\left( \<32> , [ \ind_{A^c\cap A_1^c \cap A_{12}} \Psi(12) - \ind_{A^c \cap A_1^c \cap A_{12}^c} \Psi(\<21>)  ]KK_1 K_{12}  ; v \right)
\\
&
\quad
+ \Ical\left( \<33> ,  [ \ind_{A^c\cap A_1^c \cap A_{2}} \Psi(2)+ \ind_{A^c \cap A_1 \cap A_2^c} \Psi(1)  - \ind_{A^c \cap A_1^c \cap A^2_c} \Psi(\<1>)] KK_1 K_2    ; v \right)
\\
&
\quad
+ \Ical\left( \<34> , [ \ind_{A^c \cap A_2^c \cap A_{21}} \Psi(21) - \ind_{A^c \cap A_2^c \cap A_{21}^c} \Psi(\<22>) ]K K_2 K_{21}  ; v \right)
\\
&
\quad
+ \Ical\left( \<35> ,[ \ind_{A^c \cap A_2^c \cap A_{22} } \Psi(22) - \ind_{A^c \cap A_2^c \cap A_{22}^c} \Psi(\<22>) ] K K_2 K_{22}    ; v \right)
\\
&
\quad
+ \sum_{\TT\in\Tr_3}
 \Ft_x[
\MM(\TT;v)
]
\\
&
=:
\sum_{\TT\in \Tr_2}
\Ft_x[
\BB(\TT;v)
]
 + \sum_{\TT \in \Tr_3}
 \Ft_x[
\NN(\TT;v)
]
 + \sum_{\TT\in\Tr_3}
 \Ft_x[
\MM(\TT;v)
]
.
\label{M3}
\end{align}
Finally, replacing \eqref{M3} in \eqref{NFE002}, we obtain the third iteration normal form equation, succintly written as follows
\begin{align}
\dt v
&
= \sum_{j=1}^2 \sum_{\TT \in \Tr_j} \BB(\TT ; v)
+
\sum_{j=1}^3 \sum_{\TT\in \Tr_j} \NN(\TT; v)
+
\sum_{j=1}^3 \sum_{\TT\in \Tr_j} \MM(\TT; v),
\end{align}
where the boundary terms $\BB$, the integral terms $\NN$, and the remainders $\MM$ are as in \eqref{NFE000}, \eqref{NFE002}, and \eqref{M3}.

The following result provides a systematic way of repeating the process described above, with an explicit description of the contributions and multipliers involved at each finite normal form equation.
We note that we see one boundary term and an integral term for each tree of a given generation $\Tr_j$, with the complexity of the contributions being encoded in the corresponding multipliers.

\begin{proposition}
\label{PRO:NFE}
Let $J\in\N$, $u$ be a smooth solution to \eqref{kdv}, and $v$ as in \eqref{inter} be a solution to \eqref{dtv}.
Then,
the $J$-th normal form equation for $v$ is given by
\begin{align}
\label{NFE1}
\dt v=
\sum_{j=1}^{J-1} \sum_{\TT\in \Tr_j} \BB(\TT; v)
+ \sum_{j=1}^{J} \sum_{\TT\in \Tr_j} \NN(\TT; v)
+ \sum_{\TT\in \Tr_{J}}
\MM(\TT; v),
\end{align}
where the operators are defined via their spatial Fourier transforms
\begin{align}
\label{NFE1a}
\begin{split}
\Ft_x [\BB(\TT;v)]  & = \dt\Ical(\TT, - i \mulB(\TT); v),\\
\Ft_x [\NN(\TT;v)]  & = \Ical(\TT, \mul(\TT) ; v  ), \\
\Ft_x [\MM(\TT;v)] & = \Ical(\TT, \Psi(\TT)\mulB(\TT); v),
\end{split}
\end{align}

\noi
with $\Ical(\TT, m; v)$ as in \eqref{eq:defiI}.
Moreover,
the multipliers $\mul(\TT)$ and $\mulB(\TT)$ for $\TT\in \Tr_j$ are given by
\begin{align}
\label{NFE0b}
\begin{split}
\mulB(\TT)
& =  (-1)^{j+1}  \prod_{\bul \in \TT^0}  \ind_{A^c_{\bul}} K_\bul
,
\\
\mul(\TT) & =  (-1)^{j+1} \bigg(\prod_{\bul \in \TT^0} K_\bul \bigg) \bigg\{  \sum_{\star \in \TT^{0,f}}
\ind_{A_{\star}} \Psi(\star) \bigg(\prod_{\substack{\bul \in \TT^0 \setminus\{\star\}
}} \ind_{A^c_{\bul}} \bigg)\\
&
\hspace{4cm}
-
\bigg( \prod_{\bul \in \TT^0} \ind_{A_\bul^c} \bigg)
\sum_{\star\in \TT^0 \setminus \TT^{0,\ff}} \Psi(\star)
 \bigg\},
\end{split}
\end{align}
where $K_\bul$ and $\Psi(\star)$ are as in \eqref{Kmul} and \eqref{psi-circ}.

\end{proposition}

\begin{proof}

We proceed by induction on $J$.
The claim \eqref{NFE1} for $J=1,2$ holds due to \eqref{NFE000} and \eqref{NFE002}.
To show \eqref{NFE1} for $J\ge3$, it suffices to show the following identity for $J\ge 2$:
\begin{equation}\label{eq:teseinducao}
\sum_{\TT\in\Tr_{J}} \MM(\TT; v)
=
\sum_{\TT\in \Tr_{J}} \BB(\TT; v)
+
\sum_{\TT\in \Tr_{J+1}} \NN(\TT; v)
+ \sum_{\TT\in \Tr_{J+1}} \MM(\TT; v),
\end{equation}
along with the definitions in \eqref{NFE1a}-\eqref{NFE0b}.

Let $J\in \N$ with $J\ge2$, and assume that \eqref{eq:teseinducao} holds for $1\le M \le J$,
with $\BB(\TT; v)$ and $\MM(\TT; v)$ given in \eqref{NFE1a}, and multipliers as in \eqref{NFE0b}, for any $\TT\in \Tr_M$, $1\le M \le J$, and similarly for $\NN(\TT; v)$ for $\TT\in \Tr_M$ for $1\le M \le J+1$.
For a fixed $\TT\in\Tr_{J}$, using the differentiation-by-parts formula \eqref{IBP}, we obtain
\begin{align*}
\Ft_x [\MM(\TT; v) ]
&
=\Ical(\TT, \Psi(\TT)\mulB(\TT); v)
\\
& =
\dt\Ical(\TT, -i \mulB(\TT); v) +
\sum_{\star\in\TT^\infty}
\Ical (\TT \overset{\star}{\oplus} \<1>, -\mulB(\TT) K_\star \Psi(\star) ; v )
\\
&
= \Ft_x[\BB(\TT; v) ]
+
\sum_{\star\in\TT^\infty} \Ical (\TT \overset{\star}{\oplus} \<1>, -\mulB(\TT) K_\star \Psi(\star) ; v).
\end{align*}

\noi
Thus, summing over the choices of $\TT \in \Tr_{J}$, we have
\begin{align}
\label{mul2a}
\begin{split}
&
\sum_{\TT\in \Tr_{J}} \Ft_x [\MM(\TT; v) ]
\\
&
= \sum_{\TT\in \Tr_{J}} \Ft_x[\BB(\TT; v) ]
+
\sum_{\TT\in \Tr_{J}}
\sum_{\star\in\TT^\infty}
\Ical
\big(\TT \overset{\star}{\oplus} \<1>,
-\mulB(\TT) K_\star \Psi(\star);  v \big)
\\
&
=
 \sum_{\TT\in \Tr_{J}} \Ft_x[\BB(\TT; v) ]
+
\sum_{\TT\in \Tr_{J+1}}
\Ical
\bigg(\TT,
-\sum_{\star\in\TT^{0,f}}
\mulB(\TT \overset{\star}{\ominus}  \<1> ) K_\star \Psi(\star) ;
v
\bigg)
\end{split}
\end{align}
 since the map $(\TT,\star)\mapsto (\TT \overset{\star}{\oplus} \<1>, \star)$ from  $\{(\TT,\star): \TT\in \Tr_J, \ \star \in \TT^\infty\}$ to $\{(\TT,\star): \TT\in \Tr_{J+1}, \ \star \in \TT^{0,\ff}\}$  is a bijection, and
where the last equality in \eqref{mul2a} follows from the linearity of $\Ical$ with respect to the multiplier, as seen in \eqref{eq:defiI}. Since for $\TT\in \Tr_{J+1}$, the tree $\TT \overset{\star}{\ominus} \<1>\in \Tr_{J}$, by the induction hypothesis, using \eqref{NFE0b}, we can write the multiplier of the last term as follows

\noi
\begin{align}
\label{mul2b}
\begin{split}
&
-\sum_{\star\in\TT^{0,f}}
\mulB(\TT \overset{\star}{\ominus}  \<1> ) K_\star \Psi(\star)
\\
&
=
(-1)^{ N+2} \sum_{\star\in \TT^{0,\ff}} K_\star \Psi(\star)\bigg( \prod_{\bul \in \TT^0 \setminus \{\star\}} \ind_{A_\bul^c} K_\bul  \bigg)
\\
& =
(-1)^{ N+2} \bigg( \prod_{\bul \in \TT^0} K_\bul \bigg)
\bigg\{ \sum_{\star \in \TT^{0,\ff}} \Psi(\star) \bigg(\prod_{\bul \in \TT^0 \setminus\{\star\}} \ind_{A^c_\bul}  \bigg) \bigg\}
\\
& = (-1)^{ N+2} \bigg( \prod_{\bul \in \TT^0} K_\bul \bigg)
\bigg\{ \sum_{\star \in \TT^{0,\ff}} \ind_{A_\star} \Psi(\star) \bigg(\prod_{\bul \in \TT^0 \setminus\{\star\}} \ind_{A^c_\bul}\bigg)\\
& \hspace{4cm}
+ \bigg(\prod_{\bul\in \TT^0} \ind_{A_\bul^c} \bigg) \bigg(\sum_{\star\in \TT^{0,\ff}}\Psi(\star) - \Psi(\TT) + \Psi(\TT)
\bigg)  \bigg\}
\\
& =
(-1)^{ N+2} \bigg( \prod_{\bul \in \TT^0} K_\bul \bigg) \bigg\{ \sum_{\star \in \TT^{0,\ff}} \ind_{A_\star} \Psi(\star) \bigg( \prod_{\bul \in \TT^{0} \setminus \{\star\}} \ind_{A_\bul^c} \bigg) - \bigg(\prod_{\bul\in \TT^0} \ind_{A_\bul^c} \bigg) \sum_{\star \in \TT^{0} \setminus \TT^{0,\ff}} \Psi(\star) \bigg\} \\
&
\quad + (-1)^{ N+2} \bigg(\prod_{\bul \in \TT^0} \ind_{A^c_\bul} K_\bul \bigg) \Psi(\TT) \\
&  = \mul(\TT) + \Psi(\TT)\mulB(\TT),
\end{split}
\end{align}
where we used that the set of parent nodes $(\TT \overset{\star}{\ominus} \<1>)^0$ is given by $\TT^0 \setminus\{\star\}$, since per Definition~\ref{DEF:tree3} we have removed the children nodes of $\star$.

Combining \eqref{mul2a} and \eqref{mul2b} gives
\begin{align}
& \sum_{\TT\in \Tr_{J+1}}
\Ical
\bigg(\TT,
-\sum_{\star\in\TT^{0,\ff}}
\mulB(\TT \overset{\star}{\ominus}  \<1> ) K_\star \Psi(\star) ;
v
\bigg)
\\
&
=
\sum_{\TT\in \Tr_{J+1}}  \Ical \big(\TT, \mul(\TT) + \Psi(\TT)\mulB(\TT) ; v \big)
 \\
&
=
 \sum_{\TT\in \Tr_{J+1}}  \Ft_x [\NN(\TT ; v ) ]
+  \sum_{\TT\in \Tr_{J+1}} \Ft_x [\MM(\TT ; v ) ]
,
\end{align}
from linearity of $\Ical$ in the multiplier (see \eqref{eq:defiI}), and \eqref{NFE1a}. Thus, \eqref{eq:teseinducao} holds and the result follows. \qedhere

\end{proof}

Motivated by Proposition \ref{PRO:NFE}, we can formally take the limit $J\to \infty$ to derive the infinite normal form equation for $v$:
\begin{align}
\label{NFE}
\dt v= \sum_{j=1}^{\infty} \sum_{\TT\in \Tr_j} \left(\NN(\TT;v) +
\BB(\TT;v) \right)
,
\end{align}
with $\NN(\TT;v)$ and $\BB(\TT;v)$ as in \eqref{NFE1a}.
In Subsection \ref{sec:gauge_prop}, we will prove that \eqref{dtv} and \eqref{NFE} are equivalent for sufficiently smooth solutions.

\subsection{Change of variables and modified normal form equation}
\label{SUBSEC:gauge}

In this subsection, we derive the modified normal form equation for the new variable $w$ related to $v$ via a nonlinear change of variables. In particular, the equation for $w$ (see \eqref{w2}) includes no boundary terms and all the integral terms can be estimated via frequency-restricted estimates.

In the normal form equation \eqref{NFE}, the quadratic boundary term $\BB(\<1>;v)$ is unbounded in $\FL^{s,\infty}(\R)$ for $s<-\frac12$; see the proof of Lemma~\ref{LEM:D} together with \eqref{w} and \eqref{gauge1}. This motivates the introduction of a nonlinear change of variables which removes this contribution and leads to a better-behaved normal form equation at low regularity.

Before defining our change of variables, we introduce a multilinear version of the functional $\Ical$ in \eqref{eq:defiI}.
Let $j\in \N$ and $\TT\in \Tr_j$. Consider the $j+1$ leaves of $\TT$ ordered lexicographically with respect to their words:
\begin{align}
\label{lexic}
\bul_1\le \bul_2\le \dots \le \bul_{j+1}.
\end{align}
Graphically, in a
planar representation of the tree, this corresponds to ordering the leaves of
$\TT$ from left to right. Given space-time functions $v_1,\dots, v_{j+1}: \R \times \R \to \R$, we
define the multilinear operator
\begin{equation}\label{eq:defiImulti}
\Ical(\TT, m; v_1,\dots,v_{j+1}) (t)=i
\int_{\Gamma_\TT}
e^{it\Psi(\TT)}m(\vec{\xi}) \prod_{\ell
=1}^{j+1} \ft v_\l (t,\xi_{\bul_\ell}) d\vec{\xi}
.
\end{equation}
Note that
when $v_1 = \cdots = v_{j+1} = v$, the definition above agrees with that of $\Ical(\TT, m; v )$ in \eqref{eq:defiI}.  Moreover, we define the multilinear analogues of $\BB(\TT)$ and $\NN(\TT)$ in \eqref{NFE1a} in Proposition~\ref{PRO:NFE}, as follows:

\noi
\begin{align*}\Ft_x [\BB(\TT; v_1, \ldots, v_{j+1}) ]  &=
\dt\Ical (\TT,- i \mulB(\TT) ; v_1,\dots,v_{j+1}),\\
\Ft_x [\NN(\TT; v_1, \ldots, v_{j+1} ) ]
&=
\Ical(\TT,\mul(\TT) ; v_1, \ldots,
v_{j+1} )
.
\end{align*}

We now introduce our change of variables $v\mapsto w$, defined via
\begin{align}
\label{w}
v
=
w + \Ical(\<1>, - i \mulB(\<1>); w)
=:
w + \D[w, w]
,
\end{align}
where $\Ical$ is as in \eqref{eq:defiI}. Then, we have
\begin{align}
\label{wa}
\dt v
=
\dt w  + \BB(\<1>; w),
\end{align}
for $\BB$ as in
\eqref{NFE1a} in
Proposition~\ref{PRO:NFE}; see also \eqref{NFE002}.
Combining \eqref{wa} with the normal form equation \eqref{NFE}, and replacing $v$ with the change of variables \eqref{w}, we obtain a closed equation for~$w$:
\begin{align}
\dt w & = -  \BB(\<1>;w) + \BB(\<1>;v) + \NN(\<1>;v)
+ \sum_{j=2}^\infty \sum_{\TT\in\Tr_j} \Big( \NN(\TT; v)  +
\BB(\TT; v)
 \Big) \nonumber\\
& =
 \BB\big(\<1>; w, \D[w, w] \big)
+ \BB \big(\<1>; \D[w, w],w\big)
+
\BB\big(\<1>; \D[w,w], \D[w, w]\big)
  \nonumber\\
& \qquad
+ \NN(\<1>; w ,w )  + \NN\big(\<1>; w, \D[w, w]\big)
+
\NN\big(\<1>; \D[w, w], w\big)
+ \NN\big(\<1>; \D[w, w],\D[w, w]\big)
 \nonumber\\
& \qquad
+ \sum_{j=2}^\infty \sum_{\TT\in \Tr_j}
\Big( \NN\big(\TT; w+
\D[w, w]
\big)
+  \BB\big(\TT; w + \D[w, w]\big) \Big).
\label{w1}
\end{align}

Note that unlike \eqref{NFE} for $v$, the equation in \eqref{w1} has multiple contributions with the same homogeneity in $w$, coming from operators $\BB$ and $\NN$ evaluated at different trees,
since the operator $\D[w, w]$ has
two copies of $w$.
For example, consider the time-integral terms in \eqref{w1} with homogeneity 4 associated with the tree $\<33>$:
\begin{align}
\label{Bcherry}
\NN(\<33>;w) , \ \NN(\<21>; w, \D[w,w] ), \ \NN(\<22>;
\D[w,w], w ) , \ \text{and} \ \NN( \<1>; \D[w, w], \D[w,w] ).
\end{align}
By using a red tree ${\red \<1>}$ to represent $\D[w, w]$, the contributions in \eqref{Bcherry} are associated with folowing coloured trees, respectively,
\begin{align}
\label{Bcherry2}
\<33>
\
,
\quad
\begin{tikzpicture}[
baseline=-12, scale=0.5, every node/.style={dotA},
level distance=15pt,
level 1/.style={sibling distance=21 pt},
level 2/.style={sibling distance=\distC pt},
level 3/.style={sibling distance=15pt}
]
\node   (root) {}
child {child child} child {child child};
\node  at (root-1) {}; \node [color=red] at (root-2) {};
\node at (root-1-1) {}; \node at (root-1-2) {};
\node [color=red]  at (root-2-1) {}; \node [color=red] at (root-2-2) {};
\path  (root-2) edge [thick, color=red] (root-2-1) ;
\path  (root-2) edge [thick, color=red] (root-2-2) ;
\end{tikzpicture}
\ ,
\quad
\begin{tikzpicture}[
baseline=-12, scale=0.5, every node/.style={dotA},
level distance=15pt,
level 1/.style={sibling distance=21 pt},
level 2/.style={sibling distance=\distC pt},
level 3/.style={sibling distance=15pt}
]
\node   (root) {}
child {child child} child {child child};
\node [color=red] at (root-1) {}; \node at (root-2) {};
\node [color=red]  at (root-1-1) {}; \node [color=red]  at (root-1-2)
{};
\node  at (root-2-1) {}; \node at (root-2-2) {};
\path  (root-1) edge [thick, color=red] (root-1-1) ;
\path  (root-1) edge [thick, color=red] (root-1-2) ;
\end{tikzpicture}
\ , \quad
\text{and}
\quad
\quad
\begin{tikzpicture}[
baseline=-12, scale=0.5, every node/.style={dotA},
level distance=15pt,
level 1/.style={sibling distance=21 pt},
level 2/.style={sibling distance=\distC pt},
level 3/.style={sibling distance=15pt}
]
\node   (root) {}
child {child child} child {child child};
\node [color=red] at (root-1) {}; \node [color=red] at (root-2) {};
\node [color=red] at (root-1-1) {}; \node [color=red] at (root-1-2) {};
\node [color=red] at (root-2-1) {}; \node [color=red] at (root-2-2) {};
\path  (root-2) edge [thick, color=red] (root-2-1) ;
\path  (root-2) edge [thick, color=red] (root-2-2) ;
\path  (root-1) edge [thick, color=red] (root-1-1) ;
\path  (root-1) edge [thick, color=red] (root-1-2) ;
\end{tikzpicture}
\ .
\end{align}

\noi Analogously to the $v$ equation, we rewrite \eqref{w1} by combining terms with the same homogeneity in $w$ and associated with the same tree, thus moving the complexity to the corresponding multipliers. In this way, we do not need to consider coloured trees as in \eqref{Bcherry2}. More importantly, this process leads to a fundamental cancellation, since \textit{all the boundary terms cancel}, as shown in the following lemma.

\begin{lemma}
\label{LEM:wbd}

Let $j\in\N$, $\TT\in \Tr_j$,
$\star \in \TT^{0,\ff}$, and $w_1, \ldots, w_{j+1}: \R\times \R \to \R$ be any smooth functions. Moreover, consider the lexicographical ordering of the leaves of $\TT$ as in \eqref{lexic},
and let $1 \le k \le j$ be such that $(\bul_k, \bul_{k+1})$ are the children of $\star$ in $\TT$.
Then, we have
\begin{align}
\label{wbd0}
\BB\big(\TT ; w_1, \ldots, w_{j+1} \big)
+ \BB\big(\TT\overset{\star}{\ominus} \<1>; w_1,  \ldots, w_{k-1} , \D[w_k,w_{k+1}] , w_{k+2}, \ldots,  w_{j+1} \big)
\equiv 0.
\end{align}
Consequently, all the boundary terms in \eqref{w1} vanish:
\begin{multline}
\label{wbd1}
  \BB \big(\<1>; w, \D[w, w] \big)
+ \BB\big(\<1>; \D[w, w],w\big)
\\
+
\BB\big(\<1>;\D[w, w], \D[w, w]\big)
+ \sum_{j=2}^\infty \sum_{\TT\in \Tr_j}  \BB\big(\TT; w + \D[w, w]\big)  \equiv 0  .
\end{multline}

\end{lemma}

\begin{proof}

Let $j\in\N$, $\TT\in\Tr_j$, and $\star \in \TT^{0,\ff}$, and  consider the terminal nodes of $\TT$ ordered as in \eqref{lexic}.
We start by showing \eqref{wbd0}. From \eqref{NFE1a} in Proposition~\ref{PRO:NFE} and \eqref{w}, we have
\begin{align}
&
\F_x\big[
\BB\big(\TT ; w_1, \ldots, w_{j+1} \big)
+ \BB\big(\TT\overset{\star}{\ominus} \<1>; w_1,  \ldots, w_{k-1} , \D[w_k,w_{k+1}] , w_{k+2}, \ldots,  w_{j+1} \big)
\big](t)
\\
&
=
-i\dt\bigg[
\int_{\G(\TT)} e^{it \Psi(\TT)} \mulB(\TT) \prod_{\l=1}^{j+1} \ft{w}_{\l} (t, \xi_{\bul_\l}) d\vec\xi
\\
&
\quad
+
\int_{\G(\TT\overset{\star}{\ominus}\<1>)} \int_{\xi_{\star} = \xi_k + \xi_{k+1}} e^{it \Psi(\TT\overset{\star}{\ominus}\<1>) + it\Psi(\star)} \mulB(\TT\overset{\star}{\ominus}\<1>) \mulB(\cherry{\star}{}{}) \prod_{\l=1}^{j+1} \ft{w}_{\l} (t, \xi_{\bul_\l})  \, d\xi_k d\vec\xi
\bigg]
\\
&
=
-i\dt\int_{\G(\TT)}
e^{it\Psi(\TT)}
\big\{
\mulB(\TT)
+
 \mulB(\TT\overset{\star}{\ominus}\<1>) \mulB(\cherry{\star}{}{})
\big\}
\prod_{\l=1}^{j+1} \ft{w}_{\l} (t, \xi_{\bul_\l}) d\vec\xi
\end{align}
where we used the additivity of $\Psi$ in \eqref{eq:additivity}. Replacing \eqref{NFE0b} in the multiplier above gives
\begin{align}
&
\mulB(\TT)
+
 \mulB(\TT\overset{\star}{\ominus}\<1>) \mulB(\cherry{\star}{}{})
\\
&
=
(-1)^{j+1} \prod_{\bul \in \TT^0} \ind_{A^c_\bul} K_\bul
+
\bigg( (-1)^{j} \prod_{\bul\in \TT^0 \setminus \{\star\}} \ind_{A^c_\bul} K_\bul \bigg)
\big(
(-1)^2 \ind_{A^c_\star} K_\star
\big)
=0,
\end{align}
since $(\TT \overset{\star}{\ominus} \<1>)^{0} = \TT^0 \setminus \{\star\}$, from which \eqref{wbd0} follows.

To show \eqref{wbd1}, let $j\in \N$ and $\TT \in \Tr_j$. Then, gathering all the boundary terms in \eqref{w1}
associated with the tree $\TT$ (when viewed as coloured trees as in \eqref{Bcherry2}), we have
\begin{align*}
\wt{\BB}(\TT; w)  & =
\sum_{ P\subseteq \TT^{0,f}  } \BB\Big(\TT\overset{P}{\ominus} \<1>;  f_{\bul} = w, \ \bul
\in (\TT\overset{P}{\ominus} \<1>)^\infty \setminus P , f_{\star} = \D[w,w] , \star \in P
\Big).
\end{align*}

Using \eqref{NFE1a}, \eqref{eq:defiI}, and \eqref{NFE0b}, we see that
\begin{align}
\wt{\BB}(\TT;w)   = \dt \Ical(\TT, \wt{\mul}_\BB (\TT);w )
\end{align}
where
\begin{align}
\label{wbd1a}
\begin{split}
\wt{\mul}_\BB(\TT)
 &
= \sum_{P\subseteq \TT^{0,\ff}} \mul_{\BB}\Big(\TT\overset{P}{\ominus} \<1>\Big)
\prod_{\star \in P} \mulB(\cherry{\star}{}{})
\\
&
= \sum_{P\subseteq \TT^{0,\ff}} (-1)^{|\TT^0| - |P| + 1} \bigg( \prod_{\bul \in \TT^0 \setminus P} \ind_{A_\bul^c } K_\bul \bigg) \bigg( \prod_{\star \in P} \ind_{A_\star^c} K_\star \bigg) \\
&
=(-1)^{|\TT^0|+1} \bigg( \prod_{\bul \in \TT^0} \ind_{A_\bul^c} K_\bul \bigg) \sum_{P \subseteq \TT^{0,\ff}} (-1)^{|P|} \\
&
=(-1)^{|\TT^0|+1} \bigg( \prod_{\bul \in \TT^0} \ind_{A_\bul^c} K_\bul \bigg) \sum_{k=0}^{|\TT^{0,\ff}|} \binom{|\TT^{0,\ff}|}{k} (-1)^{k}  =0 ,
\end{split}
\end{align}
and the result follows. \qedhere

\end{proof}

\begin{remark} \label{remark:factorial}
To reach the algebraic cancellations in Lemma \ref{LEM:wbd}, it is clear that the restriction imposed on $(\xi_{\bul1}, \xi_{\bul2})$ in the set $A^c_\bul$ must be independent on the parent node~$\bul$.
This is in stark contrast with the normal form approach in \cite{kishimoto, GuoKwonOh_nls}, where the non-resonant set (where the differentiation-by-parts is performed), imposing a restriction on $(\xi_{\bul1}, \xi_{\bul2})$, depends on all the nodes added in the previous iterations.
Moreover, this requires the usage of ordered trees, which retain the information on how they ``grew''. As the number of such trees grows factorially with $j$, to ensure summability, one must introduce a sequence of small parameters (depending on $j$) in the definition of the resonant set. This is yet another obstacle to the algebraic cancellations described in Lemma \ref{LEM:wbd}.
\end{remark}

We can then further rewrite the equation \eqref{w1} by grouping integral terms in a similar fashion. Similarly to Lemma~\ref{LEM:wbd} which establishes the cancellation of \textit{all} boundary terms in the $w$ equation \eqref{w1}, the next lemma also shows some cancellation for integral terms.

\begin{lemma}
\label{LEM:wint}
Let $w:\R^2\to \R$ be a smooth function and $\NN$ be as in \eqref{NFE1a}.
Then, we~have
\begin{align}
\label{wint0a}
\sum_{j=1}^\infty \sum_{\TT\in \Tr_j} \NN\big(\TT; w+
\D[w, w]
\big)
=
\NN_1(\<1>;w) + \sum_{j=2}^\infty \sum_{\TT\in\Tr_j} \sum_{\ell=2}^4  \NN_\ell(\TT; w),
\end{align}
where the nonlinear contributions are defined via their spatial Fourier transform
\begin{align}
\label{wint0b}
\begin{split}
\Ft_x \big[ \NN_\ell (\TT; w) \big] & =
\Ical \big(\TT, {\mul}_\l(\TT); w \big),\quad \l=1,\dots, 4,
\end{split}
\end{align}
for $\Ical(\TT, m; w)$ as in \eqref{eq:defiI},
where the multipliers are given by
\begin{align}
    \mul_1(\<1>)&= \ind_A K \Psi(\<1>),\\
    \mul_2(\TT)&=(-1)^{j+1} \ind_{A_\star} K_\star\Psi(\star) \bigg( \prod_{\bul \in \TT^0 \setminus\{\star\}}  K_\bul  \ind_{A_\bul^c} \bigg),
    && \text{if} \quad  \TT^{0,\ff} = \{\star\}, \ \star\neq \emptyset,
    \\
    \mul_3(\TT)&=(-1)^{j}  K_{\pb(\star)}\Psi(\pb(\star)) \bigg( \prod_{\bul \in \TT^0 \setminus\{\pb(\star)\}}  K_\bul  \ind_{A_\bul^c} \bigg),
    &&
    \text{if} \quad  \TT^{0,\ff} = \{\star\},  \ \star\neq \emptyset
    \label{wint0c}
    \\
     \mul_4(\TT)&=(-1)^{j+1}  K_\star\Psi(\star) \bigg( \prod_{\bul \in \TT^0 \setminus\{\star\}}  K_\bul  \ind_{A_\bul^c} \bigg),
    && \text{if} \quad  \TT^{0,\ff} = \{\star1, \star2\},
    \\
    \mul_{\l}(\TT) & = 0, \ \l=2,3,4 && \text{otherwise},
\end{align}
with $\TT\in \Tr_j$,  $j \ge 2$, and $K$ and $\Psi$ as in \eqref{Kmul} and \eqref{psi-tree}, respectively.
Moreover, for $\l=1,\ldots, 4$, we will refer to the contributions $\NN_\l(\TT) $ as terms of Type $\I,\II,\III,\IV$, respectively.

\end{lemma}

\begin{proof}

We note that there is a unique element on l.h.s. of \eqref{wint0a} with homogeneity 2 in $w$, namely $\NN(\<1>; w,w)$.
Therefore, $\NN_1(\<1>;w) = \NN(\<1>;w)$ and $\mul_1(\<1>) = \mul(\<1>)$ as in \eqref{NFE1a}-\eqref{NFE0b}.

Next, we consider the contributions on the l.h.s. of \eqref{wint0a} with underlying tree structure of $\<21>$, $\<22>$, and $\<33>$.
Namely, from \eqref{wint0a}, \eqref{NFE1a}, and \eqref{NFE0b}, we have
\begin{align}
&
		\NN(\<1>; \D[w, w], w) + \NN(\<21>; w) \\
		&
		=  \Ical \Big(\<21>, - K K_{1} \big(   \ind_{A^c\cap A_{1}} \Psi(1) - \ind_{A^c\cap A_{1}^c} \Psi(\<1>)  \big); w \Big)
		+
		\Ical \Big(\<21>, K K_{1}  \ind_{A\cap A_{1}^c} \Psi(\<1>); w \Big)
\\
&
=  \Ical \Big(\<21>, - K K_{1}   \ind_{A^c\cap A_{1}} \Psi(1) ; w \Big)
		+
		\Ical \Big(\<21>, K K_{1}  \ind_{ A_{1}^c} \Psi(\<1>); w \Big)
\\
&
= : \NN_2(\<21>; w) + \NN_3(\<21>; w )
,
\\
&
		\NN(\<1>; w, \D[w, w] ) + \NN(\<22>; w) \\
		&
		=  \Ical \Big(\<22>, - K K_{2} \big(   \ind_{A^c\cap A_{2}} \Psi(2) - \ind_{A^c\cap A_{2}^c} \Psi(\<1>)  \big); w \Big)
		+
		\Ical \Big(\<22>, K K_{2}  \ind_{A\cap A_{2}^c} \Psi(\<1>); w \Big)
\\
&
=
\Ical \Big(\<22>, - K K_{2}   \ind_{A^c\cap A_{2}} \Psi(2) ; w \Big)
		+
		\Ical \Big(\<22>, K K_{2}  \ind_{ A_{2}^c} \Psi(\<1>); w \Big)
\\
&
=:
\NN_2 (\<22>; w) + \NN_3 (\<22>; w)
,
\\
&
		\NN(\<1>;  \D[w, w], \D[w, w] ) + \NN(\<21>; w,w,  \D[w, w])  + \NN(\<22>; \D[w, w],w, w) + \NN(\<33>;w)  \\
		&
		=
		\Ical \Big( \<33>, KK_1K_2 \Psi(\<1>) \ind_{A\cap A_{1}^c\cap A_{2}^c}  \Big)
		+\Ical \Big( \<33>, KK_1K_2 \Psi(\<1>) \ind_{A^c\cap A_{1}^c\cap A_{2}^c}  \Big)
        \\
        &
		=
		\Ical \Big( \<33>, KK_1K_2 \Psi(\<1>) \ind_{ A_{1}^c\cap A_{2}^c}  \Big)
\\
&
=: \NN_4 (\<33>; w)
,
	\end{align}

\noi
where we note the agreement with \eqref{wint0b}-\eqref{wint0c}.

	Now let $j\in\N$ with $j\ge3$ and $\TT\in \Tr_j$, such that $\TT \neq \<33>$.
	We write the resulting term corresponding to tree $\TT$ by  gathering the integral terms in \eqref{wint0a}
	associated with the tree $\TT$ without colourings:
	\begin{align*}
		\sum_{ C \subseteq \TT^{0,f}  } \NN\big(\TT \overset{C}{\ominus} \<1> ;
		f_{\bul} = w \ \text{ if }\ \bul \in (\TT\overset{C}{\ominus}\<1>)^\infty \setminus C , \
		 f_{\bul} = \D[w, w] \ \text{ if } \ \bul \in C
		\big)
		 = \Ical ( \TT, \wt{\mul}(\TT); w),
	\end{align*}
	where $\Ical$ is as in \eqref{eq:defiI} and the multiplier $\wt{\mul}(\TT)$ is given by
	\begin{align}
			\wt{\mul}(\TT)
		& = \sum_{C \subseteq \TT^{0,\ff}} \mul(\TT\overset{C}{\ominus} \<1>)
			\bigg( \prod_{\star \in C } \mulB(\cherry{\star}{}{}) \bigg)
			\\
			& = \sum_{C \subseteq \TT^{0,\ff}} (-1)^{ j + 1 + |C|}
			\bigg( \prod_{\bul\in \TT^0 \setminus C    } K_\bul \bigg)
			\bigg(\prod_{\star \in C} \ind_{A_\star^c} K_\star \bigg) \\
			&  \times \bigg\{ \sum_{\star \in (\TT \overset{C}{\ominus} \cherryS)^{0,\ff}} \ind_{A_\star} \Psi(\star)
			\bigg(\prod_{\bul \in \TT^0 \setminus(C \cup \{\star\})} \ind_{A_\bul^c} \bigg)
			-
			\bigg(\prod_{\bul \in \TT^0 \setminus C} \ind_{A_\bul^c}\bigg)
			\sum_{\star \in (\TT \overset{C}{\ominus} \cherryS)^0  \setminus (\TT \overset{C}{\ominus} \cherryS)^{0,\ff}} \Psi(\star)
			\bigg\} \\
			& =: {\mul}_{\gfrak}(\TT) + {\mul}_{\bfrak}(\TT),
		\label{wint1a}
	\end{align}
	where we used \eqref{NFE0b}, and the fact that $(\TT\overset{C}{\ominus}\<1>)^0 = \TT^0 \setminus C$. It remains to simplify the multipliers above.

	We first consider ${\mul}_{\bfrak}(\TT)$, which we write as
	\begin{align}
		{\mul}_{\bfrak} (\TT)
&
= (-1)^j \bigg( \prod_{\bul \in \TT} \ind_{A_\bul^c} K_\bul \bigg)
\sum_{\star \in \TT^0 } \Psi(\star) \sum_{ C \subseteq \TT^{0,\ff}} (-1)^{|C|} \ind_{\star\in (\TT \overset{C}{\ominus} \cherryS)^0  \setminus (\TT \overset{C}{\ominus} \cherryS)^{0,\ff}},
	\label{wint2a}
	\end{align}
where it remains to simplify the latter sum, by considering different cases for $\star$.
Note that if $\star \in \TT^{0,\ff}$ or $\star\in \TT^0\setminus \TT^{0,\ff}$ with at least one child in $\TT^0 \setminus \TT^{0,\ff}$, then sum over all $C\subseteq \TT^{0,\ff}$ in \eqref{wint2a} gives a zero contribution. Thus, it remains to consider $\star \in \TT^0 \setminus \TT^{0,\ff}$ and $\{\star1, \star2\} \subseteq \TT^{0,\ff} \cup \TT^\infty$, noting that we cannot have both children being leaves.

\smallskip

\noi
\underline{\textbf{Case 1}:  $(\star1, \star2) \in \TT^{0,\ff} \times \TT^\infty$ or $(\star1, \star2) \in \TT^\infty \times \TT^{0,\ff}$}

If
$(\star1, \star2) \in \TT^{0,\ff} \times \TT^\infty$, the only non-zero contributions in the sum over $C$ in \eqref{wint2a}~are
\begin{align}
\sum_{C\subseteq \TT^{0,\ff} \setminus \{\star1\}} (-1)^{|C|}
=  \sum_{j=0}^{|\TT^{0,\ff}|-1} \binom{|\TT^{0,\ff}|-1}{j} (-1)^j
=  \ind_{|\TT^{0,\ff}| = 1},
\end{align}
and we must have $\TT^{0,\ff} = \{\star 1\}$ to have a non-zero term. The remaining case is analogous by symmetry between $\star1$ and $\star2$.

\smallskip
\noi
\underline{\textbf{Case 2}:  $\star1,\star2\in \TT^{0,\ff}$}

In this case, $C$ can contain at most one element in $\{\star1,\star2\}$, thus we get
\begin{align}
 \sum_{C\subseteq \TT^{0,\ff} \setminus\{\star1,\star2\}}[ (-1)^{|C|}  + (-1)^{|C|+1} + (-1)^{|C|+1}]
=
-  \ind_{|\TT^{0,\ff}| = 2}
\end{align}
which is nonzero when $\TT^{0,\ff} = \{\star1, \star2\}$.

Combining the different cases, we get
\begin{align}\label{mulb}
\mul_{\bfrak}(\TT) = (-1)^{j+1} \bigg( \prod_{\bul \in \TT^0} \ind_{A_\bul^c} K_\bul \bigg) \sum_{\star\in\TT^0\setminus\TT^{0,\ff}}\Psi(\star)
\big(
- \ind_{\TT^{0,\ff} = \{\star1\}} - \ind_{\TT^{0,\ff} = \{\star2\}} + \ind_{\TT^{0,\ff} = \{\star1, \star2\}}
\big).
\end{align}

	We now consider the first multiplier ${\mul}_\gfrak$ in \eqref{wint1a}, which we rewrite as follows
	\begin{align}
    \begin{split}
		 {\mul}_\gfrak (\TT)
		&
		=
(-1)^{j+1}
		\bigg( \prod_{\bul \in \TT^0} K_\bul \bigg)
		\sum_{ C \subseteq \TT^{0,\ff}} (-1)^{|C|}  \sum_{\star \in (\TT \overset{C}{\ominus} \cherryS)^{0,\ff}}  \ind_{A_\star} \Psi(\star) \bigg(\prod_{\bul \in \TT^0\setminus\{\star\}} \ind_{A_\bul^c} \bigg)\label{wint3a}
\\
&
=
(-1)^{j+1}
		\bigg( \prod_{\bul \in \TT^0} K_\bul \bigg)
 \sum_{\star \in \TT^0} \ind_{A_\star} \Psi(\star)
		\bigg( \prod_{\bul\in\TT^0 \setminus\{\star\}} \ind_{A_\bul^c}\bigg)
		\sum_{C\subseteq \TT^{0,\ff}} (-1)^{|C|} \ind_{\star \in (\TT\overset{C}{\ominus}\cherryS)^{0,\ff} } .
        \end{split}
	\end{align}
As before, we consider different cases depending on $\star$ and focus on the inner sum in $C$.
Note that any $\star \in \TT^{0}\setminus \TT^{0,\ff}$ with at least one child in $\TT^0\setminus\TT^{0,\ff}$ gives a zero contribution in \eqref{wint3a}. Thus, we consider $\star \in \TT^{0,\ff}$ or $\{\star1,\star2\} \subseteq \TT^{0,\ff} \cup \TT^\infty$.

\smallskip

\noi
\underline{\textbf{Case 1}: $\star\in \TT^{0,\ff}$  }

The inner summand in \eqref{wint3a} is non-zero when $C\subseteq \TT^{0,\ff} \setminus \{\star\}$:
\begin{align}
\sum_{C\subseteq \TT^{0,\ff} \setminus\{\star\}} (-1)^{|C|} = \ind_{|\TT^{0,\ff}|=1},
\end{align}
thus we need $\TT^{0,\ff} = \{\star\}$ to get a non-zero contribution.

\smallskip
\noi
\underline{\textbf{Case 2}:  $(\star1, \star2) \in \TT^{0,\ff} \times \TT^\infty$ or $(\star1, \star2) \in \TT^\infty \times \TT^{0,\ff}$}

If $(\star1, \star2) \in \TT^{0,\ff} \times \TT^\infty$,
the inner summand in \eqref{wint3a} is non-zero when $ \star1 \in C$, giving
\begin{align}
\sum_{C\subseteq \TT^{0,\ff} \setminus\{\star1\}} (-1)^{|C|+1} = - \ind_{|\TT^{0,\ff}|=1},
\end{align}
thus we need $\TT^{0,\ff} = \{\star1\}$ to get a non-zero contribution. This is analogous for the other case by symmetry.

\smallskip

\noi
\underline{\textbf{Case 3}:  $\star1,\star2\in \TT^{0,\ff}$}

In this case, we need $C$ in \eqref{wint3a} containing both $\star1$ and $\star2$, which gives
\begin{align}
\sum_{C\subseteq \TT^{0,\ff} \setminus\{\star1, \star2\}} (-1)^{|C|+2} = \ind_{|\TT^{0,\ff}|=2},
\end{align}
thus we need $\TT^{0,\ff} = \{\star1, \star2\}$ to get a non-zero contribution.

Combining the different terms, we get
\begin{align}\label{mulg}
{\mul}_\gfrak (\TT)
& =
(-1)^{j+1}
		\bigg( \prod_{\bul \in \TT^0} K_\bul \bigg)
 \sum_{\star \in \TT^0} \ind_{A_\star} \Psi(\star)
		\bigg( \prod_{\bul\in\TT^0 \setminus\{\star\}} \ind_{A_\bul^c}\bigg)
\\
&
\quad \times
 \big(
\ind_{ \TT^{0,\ff}  = \{\star\}}
- \ind_{ \TT^{0,\ff}  = \{\star1\}}
- \ind_{ \TT^{0,\ff}  = \{\star2\}}
+\ind_{ \TT^{0,\ff}  = \{\star1, \star2\}}
\big).
\end{align}

Then, from \eqref{mulb} and \eqref{mulg}, we obtain
\begin{align}
    &
     \mul_\bfrak (\TT) + \mul_\gfrak(\TT)
     \\
&
=
(-1)^{j+1}
\sum_{\star\in \TT^0} K_\star \Psi(\star)
\bigg( \prod_{\bul \in \TT^0 \setminus \{\star\}} K_\bul \ind_{A^c_\bul} \bigg)
\bigg(
\ind_{A_\star^c}
\big(
- \ind_{\TT^{0,\ff} = \{\star1\}} - \ind_{\TT^{0,\ff} = \{\star2\}} + \ind_{\TT^{0,\ff} = \{\star1, \star2\}}
\big)
\\
&
\hspace{3cm}
+
\ind_{A_\star}
\big(
\ind_{ \TT^{0,\ff}  = \{\star\}}
- \ind_{ \TT^{0,\ff}  = \{\star1\}}
- \ind_{ \TT^{0,\ff}  = \{\star2\}}
+\ind_{ \TT^{0,\ff}  = \{\star1, \star2\}}
\big)
\bigg)
\\
&
=
(-1)^{j+1}
\sum_{\star\in \TT^0} K_\star \Psi(\star)
\bigg( \prod_{\bul \in \TT^0 \setminus \{\star\}} K_\bul \ind_{A^c_\bul} \bigg)
\bigg(
\ind_{A_\star}
\ind_{\TT^{0,\ff} = \{\star\}}
-
\ind_{\TT^{0,\ff}=\{\star\l\}, \l=1,2 }
+
\ind_{\TT^{0,\ff} = \{\star1,\star2\}}
\bigg)
\\
&
= : \mul_2(\TT) + \mul_3(\TT) + \mul_4(\TT)
\label{mul234}
\end{align}

\noi
which can be rewritten as in \eqref{wint0c}.
\qedhere

\end{proof}

\begin{remark}\label{rm:tree_type} \rm
(i)
     From \eqref{wint0c}, we see that the only trees with nonzero contributions are, up to permutations, as depicted in Figure \ref{FIG:types}.

\begin{figure}[!h]
\begin{subfigure}[b]{0.2\textwidth}
\centering
\begin{tikzpicture}[
baseline=-15, scale=1, every node/.style={dotA},
level distance=17pt,
sibling distance=20pt,
level 4/.style={sibling distance=15 pt}
]
\node  [label=above:{}]  (root) {}
child
child ;
\node  at (root-1) {};
\node  at (root-2) {};
\path  (root-1) edge [ dotted, color=black] (root-1-1) ;
\node
at (root) {};
\end{tikzpicture}
\caption{Type $\I$.}
\end{subfigure}
\begin{subfigure}[b]{0.2\textwidth}
\centering
\begin{tikzpicture}[
baseline=-15, scale=1, every node/.style={dotA},
level distance=17pt,
sibling distance=20pt,
]
\node  [label=above:{}]  (root) {}
child {
	child
	{
		edge from parent[draw=none]
		child
			{
			child
			child
			}
		child
	}
	child
	{edge from parent[draw=none]}
}
child {};
\node at (root) {};
\node  at (root-1) {};
\node  at (root-2) {};
\node at (root-1-1) {};
\node at (root-1-1-2) {};
\node at (root-1-1-1-2) {};
\node  at (root-1-1-1-1) {};
\node
at (root-1-1-1) {};
\path  (root-1) edge [ dotted, color=black] (root-1-1) ;
\end{tikzpicture}
\caption{Types $\II$-$\III$.}
\end{subfigure}
\begin{subfigure}[b]{0.2\textwidth}
\centering
\begin{tikzpicture}[
baseline=-15, scale=1, every node/.style={dotA},
level distance=17pt,
sibling distance=20pt,
level 4/.style={sibling distance=15 pt}
]
\node  [label=above:{}]  (root) {}
child {
	child
	{
		edge from parent[draw=none]
		child
			{
			child
			child
			}
		child
			{
			child
			child
			}
	}
	child
	{edge from parent[draw=none]}
}
child {};
\node at (root) {};
\node  at (root-1) {};
\node  at (root-2) {};
\path  (root-1) edge [ dotted, color=black] (root-1-1) ;
\node
at (root-1-1) {};
\node at (root-1-1-2) {};
\node at (root-1-1-1-2) {};
\node  at (root-1-1-1-1) {};
\node at (root-1-1-2-2) {};
\node  at (root-1-1-2-1) {};
\node
at (root-1-1-1) {};
\node at (root-1-1-2) {};
\end{tikzpicture}
\caption{Type $\IV$.}
\end{subfigure}
\caption{Graphical representation of tree of Type $\I$, and examples of trees of Types $\II$-$\IV$. }
\label{FIG:types}
\end{figure}
\noi Observe that, given $j\ge 2$ the number of trees of Types $\II$-$\III$ (with $j$ parent nodes) and of Type $\IV$ (with $j+2$ parent nodes) is exactly $2^{j-1}$. This is significantly smaller than the total number of trees with $j$ parent nodes given in \eqref{eq:catalan}, which grows as $\sim 4^j/j^{\frac32}$. However, the main gain in the cancellations in Lemma \ref{LEM:wint} is not the number of trees but rather their simple structure. This translates into a strong restriction on the frequency interactions in the nonlinear terms, which heavily simplifies the subsequent multilinear estimates.

\medskip

\noi (ii) From the definition of the multipliers in \eqref{wint0c}, note that the terms of Type $\III$ and $\IV$ do not have a restriction on the frequencies in the generation with parent $\pb(\star)$ and $\star$, respectively. This is a consequence of the calculation in \eqref{mul234}, where we see that they come from combining terms from $\mul_\bfrak(\TT)$ and $\mul_{\gfrak}(\TT)$, which come with $\ind_{A_{\star}^c} $ and $\ind_{A_\star}$, respectively.

\end{remark}

Combining \eqref{w1}, Lemma~\ref{LEM:wbd}, and Lemma~\ref{LEM:wint}, we derive the modified normal form equation for the new variable $w$:
\begin{align}
\label{w2}
\dt w & =  \NN_1(\<1>;w) + \sum_{j=2}^\infty \sum_{\TT \in \Tr_j} \sum_{\l=2}^4  \NN_\l(\TT;w),
\end{align}
where the operators $\NN_\l$ are as in \eqref{wint0b}.
We can regard $w$ as the interaction representation of another variable $z$, namely
\begin{align}
\label{z0}
z(t) : = e^{-t\dx^3} w(t),
\end{align}
analogous to the relation between $u$ and $v$ in \eqref{inter}. In particular, $z$ satisfies the gauged KdV equation \eqref{rkdv}, whose exact form is given by
\begin{align}
\label{z1}
\dt z
+ \dx^3 z
&
=
e^{-t\dx^3}
\NN_1(\<1>;e^{t\partial_x^3}z)
+
\sum_{j=2}^\infty \sum_{\TT \in \Tr_j} \sum_{\l=2}^4
e^{-t\dx^3}
\NN_\l(\TT;e^{t\partial_x^3}z)
\\
&
=:
\Ncal_1(\<1>;z) + \sum_{j=2}^\infty \sum_{\TT \in \Tr_j} \sum_{\l=2}^4  \Ncal_\l(\TT;z)
.
\end{align}

{
The computations leading \eqref{kdv} to \eqref{z1} are, for now, completely formal. Nevertheless, \eqref{z1} is well-defined \emph{per se} and we may move to the proof of Theorem \ref{thm:lwp_rkdv} (Subsection \ref{SUBSEC:lwp}). This will require us to introduce the functional framework (Section \ref{sec:functional}) and prove the necessary multilinear estimates (Section \ref{sec:multi}). Afterwards, using Theorem \ref{thm:lwp_rkdv} and the gauge properties (see Subsection \ref{sec:gauge_prop}), we will prove the equivalence between both equations for sufficiently smooth solutions (Subsection \ref{SUBSEC:lwp}). This equivalence will then allow us to prove Theorem \ref{thm:wp_kdv_intro}.
}

	\section{Functional framework}\label{sec:functional}

	In this section, we introduce the function spaces over which the solution to \eqref{rkdv} will be built and prove linear estimates associated with the Duhamel integral formulation (Section \ref{sec:fourier_restr}). The construction of solutions reduces to proving specific multilinear estimates, via frequency-restricted estimates (Section \ref{sec:fre}).

	\subsection{Fourier restriction spaces}\label{sec:fourier_restr}

	We first recall the Fourier restriction spaces,
     first introduced by Bourgain \cite{bourg1, bourg2} and Klainerman-Machedon in \cite{KlainermanMachedon}, and adapted to the Fourier-Lebesgue setting in \cite{Grunrock_mkdv, HerrGrunrock_FL}.

		Let $s,b\in \R$ and $1\le p,q \le \infty$. We consider the Fourier restriction spaces $X^{s,b}_{p,q}$ defined through the norm
		\begin{equation}\label{Xsb}
			\| u \|_{X^{s,b}_{p,q}} := \| \jap{\xi}^s \jap{\tau-\xi^3}^b \F_{t,x}(u) (\tau, \xi) \|_{L^p_\xi L^q_\tau}
            =
            \| \jb{\xi}^s \jb{\tau}^b \F_{t,x}(e^{t \dx^3} u ) (\tau, \xi) \|_{L^p_\xi L^q_\tau}.
		\end{equation}
When $p=q$, we may also abbreviate $X^{s,b}_{p,q}$ as $X^{s,b}_p$.
		The case $p=q=2$ corresponds to the classic Fourier restriction spaces \cite{bourg1, bourg2}.
        Given $T>0$, we also recall the time-localized space $X^{s, b}_{p,q} (0,T)$ defined via the following norm
        \begin{align}
            \| u\|_{X^{s, b}_{p,q}(0,T)} : = \inf
            \big\{
            \| v\|_{X^{s, b}_{p,q}} : \ v\in X^{s,b}_{p,q}, \ v\vert_{[0,T]} = u
            \big\},
        \end{align}
        where the infimum is taken over all extensions of $u$ on $[0,T]$. See \cite[Chapter 2.6]{tao_book}, \cite[Chapter 7.4]{linares_ponce_book}), and  \cite{HerrGrunrock_FL} for further details on these spaces.
        Let  $\psi\in C_c^\infty(\R)$ be a smooth cut-off function with $\psi\equiv 1$ over $[-1,1]$ and $\psi\equiv 0$ on $\R\setminus [-2,2]$, and given $T>0$, we set $\psi_T(t)=\psi(t/T)$.

\begin{lemma}
\label{lem:linear_estimates_qfinite}
			Let $s\in\R$, $2\le p <\infty$, and $b,b'\in\R$ with $b > \frac{1}{p'}$ and $-\frac1p<b' \le 0 \le b \le 1+b'$.

			\noi{\rm(i)} The following embedding holds: $X^{s,b}_{p} \hookrightarrow C(\R;\FL^{s,p}(\R))$.

			\smallskip
			\noi{\rm(ii)} For all $f \in \FL^{s,p}(\R)$ and $F \in X^{s, b'}_p$, we have
			\begin{align}
				\label{lin1_gh}
				\| \psi e^{-t\dx^3} f \|_{X^{s,b}_{p}} &\les \| f \|_{\FL^{s,p}},
                \\
				\bigg\| \psi_T (t) \int_0^t e^{-(t-t')\dx^3 } F(t') \, dt' \bigg\|_{X^{s,b}_{p} }
                &
                \les T^{1+b'-b} \| F \|_{X^{s,b'}_{p}}.
                \label{lin2a_gh}
			\end{align}

		\end{lemma}

		It is important to stress that, in the above result, the case $p=\infty$ is out-of-reach {for (i) and \eqref{lin2a_gh} in (ii)}. Nevertheless, this delicate case can still be addressed by exploiting a different integrability exponent for the time-frequency variable, as the next lemma shows. For a similar result, see \cite[Section 7]{CC24}.

		\begin{lemma}\label{lem:linear_estimates}
			Let $s\in\R$, $b,\wt b, b', \wt b'\in\R$ satisfying
            \begin{align}
                \label{bs}
                -1 < b' \le 0 \le b < 1 + b'
                \quad
                \text{and}
                \quad
                -1 < \wt b' \le 0 \le \wt b < 1 + \wt b',
            \end{align}
        and $\wt u(t) := e^{t\dx^3 } u(t)$.

            \noi
			{\rm(i)} For $b>0$,  given $u\in X^{s,b}_{\infty,1}$, we have $\wt{u}\in C(\R;\FL^{s,\infty}(\R))$ with
			\begin{equation}\label{lin}
				\|\wt{u}\|_{L^\infty_t\FL^{s, \infty}_x}\lesssim \|u\|_{X^{s,b}_{\infty,1}},
			\end{equation}
			and the following embedding holds
			\noi
			\begin{align}
				X^{s,b}_{\infty, 1} &
                \embeds C(\R; \FL_{\mathrm{loc}}^{s,\infty}(\R))\cap L^\infty(\R; \FL^{s,\infty}(\R)),
			\end{align}
            where $f \in \FL^{s, \infty}_{\mathrm{loc}}(\R)$ if for any compact set $K\subset \R$, $\P_K f \in \FL^{s, \infty}(\R)$, with $\P_K$  denoting the Dirichlet projection onto spatial frequencies $ \xi \in K$.

			\smallskip

			\noi {\rm(ii)}
            For $ 1 \le q \le \infty$, we have
			\begin{align}
				\label{lin1}
				\| \psi e^{-t\dx^3} f \|_{X^{s,b}_{\infty,q}} \les \| f \|_{\FL^{s,\infty}}.
			\end{align}

			\noi{\rm(iii)}
            The following inhomogeneous estimates hold

			\noi
			\begin{align}
				\label{lin2a}
				\bigg\| \psi_T (t) \int_0^t e^{-(t-t')\dx^3 } F(t') \, dt' \bigg\|_{X^{s,\wt{b}}_{\infty,1} } & \les T^{1+\wt{b}'-\wt{b}} \| F \|_{X^{s,\wt{b}'}_{\infty,1}}
                ,
            \\
				\label{lin2b}
				\bigg\| \psi_T (t) \int_0^t e^{-(t-t')\dx^3} F(t') \, dt' \bigg\|_{X^{s,b}_{\infty, \infty} }
                &
                \les T^{1+b'-b} \| F\|_{X^{s, b'}_{\infty, \infty}} + T^{2+\wt{b}' -b} \|F\|_{X^{s, \wt{b}'}_{\infty, 1}}.
			\end{align}

		\end{lemma}

		\begin{proof}

			We begin with (i). Given $u\in X^{s,b}_{\infty,1}$,
			$$
            \Ft_x(\wt u)(t,\xi) = \int_\R e^{it\tau}\ft{u}(\tau+\xi^3,\xi)d\tau,
			$$
            from which we have
$$\| u(t)\|_{\FL^{s,\infty}_x} = \|\wt{u}(t)\|_{\FL^{s,\infty}_x} \lesssim \|\jap{\xi}^s\ft {u}\|_{L^\infty_\xi L^1_\tau} \les \|u\|_{X^{s,b}_{\infty,1}},
			$$
            for $b\ge 0$.
			Moreover, given $t,t'\in\R$,
			\begin{align*}
				\|\wt{u}(t)-\wt{u}(t')\|_{\FL_x^{s,\infty}} &\lesssim \bigg\|\int_\R \frac{e^{it\tau}-e^{it'\tau}}{\jap{\tau}^b}\jap{\xi}^s\jap{\tau}^b\ft{u}(\tau+\xi^3,\xi)d\tau\bigg\|_{L^\infty_\xi}\\&\lesssim \sup_{\tau} \bigg| \frac{e^{it\tau}-e^{it'\tau}}{\jap{\tau}^b}\bigg| \|u\|_{X^{s,b}_{\infty,1}} \to 0 \quad \mbox{as }t\to t',
			\end{align*}
            for any $b>0$,
			and thus $\wt{u}\in C_b(\R;\FL^{s,\infty}(\R))$. Finally, considering a compact set $K\subset \R$, we have
			\begin{align*}
				\| \P_K [u(t)-u(t')]\|_{\mathcal{F}L^{s,\infty}_x}
                &
                \lesssim \|\wt{u}(t)-\wt{u}(t')\|_{\FL^{s,\infty}_x}
                + \|\jb{\xi}^s (1-e^{-i(t-t')\xi^3})\ft{u}(t')\|_{L^\infty_\xi (K) } \\
                &\lesssim \|\wt{u}(t)-\wt{u}(t')\|_{\FL^{s,\infty}_x} + \sup_{\xi\in K}\big|1 -e^{-i (t-t')\xi^3}\big| \| u \|_{X^{s,b}_{\infty,1}} \to 0\quad \mbox{as }t\to t',
			\end{align*}
            for $b\ge0$ together
            with the earlier convergence,
			which implies that $u\in C(\R;\mathcal{F}L^{s,\infty}_{\mathrm{loc}}(\R))$.

            The proof of (ii) is a direct consequence of
			$$
			\|\psi e^{-t\partial_x^3}f\|_{X^{s,b}_{\infty,q}} = \|\jap{\xi}^s\jap{\tau}^b\ft{\psi}(\tau)\ft{f}(\xi)\|_{L^\infty_\xi L^q_\tau} \lesssim \|\psi\|_{\FL^{b,q}} \|f\|_{\FL^{s,\infty}}.
			$$

			\smallskip

			We now move to (iii), where we proceed as in the proof of \cite[Lemma~2]{Grunrock_mkdv}.
            Note that \eqref{lin2a}-\eqref{lin2b} are equivalent to
            \begin{align}
                \bigg\| \psi_T (t) \int_0^t G(t') \, dt' \bigg\|_{\FL^{s, p}_x \FL^{b_1, q_1 }_t }
                \les T^{1+ q_2 - q_1} \| G \|_{\FL^{s, p}_x \FL^{b_2 , q_2}_t },
            \end{align}
            for $G = e^{t\partial_x^3} F$, for suitable choices of the parameters. Thus, we start by expanding the quantity on the l.h.s.:
			\begin{equation}
				\psi_T(t)\int_0^t G(t') dt' =  c \psi_T (t)\int_\R \frac{e^{it\eta}-1}{i\eta } (\mathcal{F}_t G )(\eta,x) d\eta = \I + \II + \III,
				\label{lin3a}
			\end{equation}
			where
			\begin{equation}
				\label{lin3}
				\begin{aligned}
					\I & :=c\psi_T(t)\sum_{k= 1}^\infty \frac{t^k}{k!}\int_{|\eta| T \le 1} (i\eta)^{k-1} (\mathcal{F}_t G )(\eta,x) d\eta, \\
					\II& :=-c\psi_T (t)\int_{|\eta|T \ge 1} (i\eta)^{-1} (\mathcal{F}_t G )(\eta,x) d\eta , \\
					\III & :=c\psi_T (t)\int_{|\eta| T  \ge 1} e^{it\eta}(i\eta)^{-1} (\mathcal{F}_t G )(\eta,x) d\eta,
				\end{aligned}
			\end{equation}
			for some constant $c\in\R$.

            The space-time Fourier transform of $\I$ in \eqref{lin3} satisfies
			\begin{equation}
				|(\mathcal{F}_{t,x}\I)(\tau, \xi)| \les \sum_{k= 1}^\infty \frac{|\widehat{t^k\psi_T}(\tau)|}{k!} \int_{|\tau|T  \le 1} |\eta|^{k-1} |\ft{G }(\eta,\xi)| d\eta,
				\label{lin4a}
			\end{equation}
            which combined with H\"older's inequality gives
			\begin{align}
					\| \jap{\xi}^s \jap{\tau}^b \F_{t,x} \I \|_{L^\infty_\xi L^\infty_\tau} & \les \sum_{k=1}^\infty \frac{1}{k!} \| \jap{\tau}^b\ft{t^k \psi_T} (\tau) \|_{L^\infty_\tau} \bigg( \int_{|\eta|T \le 1} |\eta|^{k-1} \jap{\eta}^{-b'} d\eta \bigg) \| \jap{\xi}^s \jap{\eta}^{b'} \ft{G}(\eta, \xi) \|_{L^\infty_\xi L^\infty_\eta} \\
					& \les_\psi T^{1+ b' -b} \| G \|_{\FL^{s, \infty}_x \FL^{b', \infty}_t },
				\label{lin4b}
			\end{align}
			\noi
			since $k-1-b' \ge0$ for $b' \le 0$.
			Similarly, since $\tilde{b}'\le 0$,
			\begin{equation}
				\begin{aligned}
					\| \jap{\xi}^s \jap{\tau}^{\tilde{b}} \F_{t,x} \I \|_{L^\infty_\xi L^1_\tau} & \les \sum_{k=1}^\infty \frac{1}{k!} \| \jap{\tau}^{\wt{b}}\ft{t^k \psi_T} (\tau) \|_{L^1_\tau} \bigg( \sup_{|\eta|T \le 1} |\eta|^{k-1} \jap{\eta}^{-\wt{b}'}  \bigg) \| \jap{\xi}^s \jap{\eta}^{\wt{b}'} \ft{G}(\eta, \xi) \|_{L^\infty_\xi L^1_\eta} \\
					& \les_\psi T^{1+ \wt{b}' -\wt{b}} \| G \|_{\FL^{s, \infty}_x \FL^{\wt b ', 1}_t},
				\end{aligned}
				\label{lin4c}
			\end{equation}

			For $\II$ in \eqref{lin3}, by H\"older's inequality, we have
			\begin{equation}
				\begin{aligned}
					\| \jap{\xi}^s \jap{\tau}^b \F_{t,x} \II \|_{L^\infty_\xi L^\infty_\tau} & \les \bigg\| \jap{\tau}^b \ft{\psi}_T(\tau) \int_{|\eta|T >1} \jap{\xi}^s |\eta|^{-1} | \ft{G}(\eta, \xi)| d\eta \bigg\|_{L^\infty_\xi L^\infty_\tau} \\
					& \les \| \jap{\tau}^b \ft{\psi_T}(\tau) \|_{L^\infty_\tau} \Big( \sup_{|\eta| T >1} |\eta|^{-1} \jap{\eta}^{-\wt{b}'} \Big) \| \jap{\xi}^s \jap{\eta}^{\wt{b}'} \ft{G}(\eta, \xi) \|_{L^\infty_\xi L^1_\eta} \\
					& \les_\psi T^{2 + \wt{b}' - b} \| G  \|_{\FL^{s, \infty}_x \FL^{\wt b ', 1}_t},
				\end{aligned}
				\label{lin5a}
			\end{equation}
			given that $\wt b ' >-1$.
			Analogously,
			\begin{equation}
				\begin{aligned}
					\| \jap{\xi}^s \jap{\tau}^{\wt{b}} \F_{t,x} \II \|_{L^\infty_\xi L^1_\tau} & \les \bigg\| \jap{\tau}^{\wt{b}} \ft{\psi}_T(\tau) \int_{|\eta|T >1} \jap{\xi}^s |\eta|^{-1} | \ft{G}(\eta, \xi)| d\eta \bigg\|_{L^\infty_\xi L^\infty_\tau} \\
					& \les \| \jap{\tau}^{\wt{b}} \ft{\psi_T}(\tau) \|_{L^1_\tau} \Big( \sup_{|\eta| T >1} |\eta|^{-1} \jap{\eta}^{-\wt{b}'} \Big) \| G  \|_{\FL^{s, \infty}_x \FL^{\wt b', 1}_t}
                    \\
						& \les_\psi T^{1 + \wt{b}' - \wt{b}} \| G \|_{\FL^{s, \infty}_x \FL^{\wt b ', 1}_t},
				\end{aligned}
				\label{lin5c}
			\end{equation}
            for $\wt b'>-1$.

			Lastly, for $\III$, from \eqref{lin3} and H\"older's inequality, we have
			\begin{equation}
				\begin{aligned}
					\| \jap{\xi}^s \jap{\tau}^b \F_{t,x} \III \|_{L^\infty_\xi L^\infty_\tau} & \les \bigg\| \jap{\xi}^s \int_{|\eta|T >1} \jap{\eta}^b |\eta|^{-1} |\ft{G}(\eta, \xi) \ft{\psi}_T(\tau-\eta) | d\eta \bigg\|_{L^\infty_\xi L^\infty_\tau} \\
					& \quad + \bigg\| \jap{\xi}^s \int_{|\eta|T >1}  |\eta|^{-1} |\ft{G}(\eta, \xi)| \jap{\tau - \eta}^b | \ft{\psi}_T(\tau-\eta) | d\eta \bigg\|_{L^\infty_\xi L^\infty_\tau} \\
					& \les \Big( \sup_{|\eta|T>1} \jap{\eta}^{b-b'} |\eta|^{-1} \Big) \| \jap{\xi}^s \jap{\eta}^{b'} \ft{G}(\eta, \xi) \|_{L^\infty_\xi L^\infty_\eta} \| \ft{\psi}_T \|_{L^1_\tau} \\
					& \quad + \Big( \sup_{|\eta|T>1} |\eta|^{-1} \jap{\eta }^{-b'} \Big)  \| \jap{\xi}^s \jap{\eta}^{b'} \ft{G}(\eta, \xi) \|_{L^\infty_\xi L^\infty_\eta} \| \jap{\tau}^{b} \ft{\psi}_T \|_{L^1_\tau } \\
					& \les T^{1+b'-b} \| G\|_{\FL^{s, \infty}_x \FL^{b', \infty}_t}
				\end{aligned}
				\label{lin6a}
			\end{equation}
			given that $b-b' < 1 $ and $b'>-1$.
			Similarly,
			\begin{equation}
				\begin{aligned}
					\| \jap{\xi}^s \jap{\tau}^{\wt{b}} \F_{t,x} \III \|_{L^\infty_\xi L^1_\tau}
					& \les \Big( \sup_{|\eta|T>1} \jap{\eta}^{\wt{b}-\wt{b}'} |\eta|^{-1} \Big) \| \jap{\xi}^s \jap{\eta}^{\wt{b}'} \ft{G}(\eta, \xi) \|_{L^\infty_\xi L^1_\eta} \| \ft{\psi}_T \|_{L^1_\tau} \\
					& \quad + \Big( \sup_{|\eta|T>1} |\eta|^{-1} \jap{\eta }^{-\wt{b}'} \Big)  \| \jap{\xi}^s \jap{\eta}^{\wt{b}'} \ft{G}(\eta, \xi) \|_{L^\infty_\xi L^1_\eta} \| \jap{\tau}^{\wt{b}} \ft{\psi}_T \|_{L^1_\tau } \\
					& \les T^{1+\wt{b}'-\wt{b}} \| \jap{\xi}^s \jap{\eta}^{\wt{b}'} \ft{G}(\eta, \xi) \|_{L^\infty_\xi L^1_\eta},
				\end{aligned}
				\label{lin6c}
			\end{equation}
			for that $\wt b- \wt b' < 1 $ and $\wt b'>-1$.
            Combining the estimates above, we obtain \eqref{lin2a}-\eqref{lin2b}, completing the proof.
            \qedhere

		\end{proof}

	\subsection{Reduction to frequency-restricted estimates}\label{sec:fre}

	The proof of well-posedness of the gauged KdV equation \eqref{z1} (Theorem~\ref{thm:lwp_rkdv}) reduces to establishing multilinear estimates on terms of the form
		\begin{equation}\label{eq:multilinear_geral}
			\Ft_x (\mathcal{N}[z_1,\dots, z_k])(\xi)=\int_{\xi=\xi_1+\dots+\xi_k} m(\xi_1,\dots,\xi_k)\prod_{j=1}^k\ft{z}_j(\xi_j) d\xi_1\cdots d\xi_{k-1},
		\end{equation}
where $m$ denotes a Fourier multiplier.
In particular, in this section, we describe how to reduce such multilinear estimates in Fourier restriction spaces to frequency-restricted estimates as in \cite{CLS}. Let  $\Phi = \Phi(\xi, \xi_1, \ldots, \xi_k)$ denote the phase function of the term in \eqref{eq:multilinear_geral}, given by
\begin{equation}\label{Phi}
 \Phi=-\xi^3+\sum_{j=1}^k\xi_j^3.
\end{equation}

\noi
We start by showing a bilinear estimate in $X^{s,b}_{\infty, \infty} \cap X^{s, \wt b}_{\infty, 1}$.

		\begin{lemma}\label{lem:reduction_fre_k2}
			Let $m= m(\xi_1,\xi_2)$ be a multiplier and $\Ncal$ a bilinear operator as in \eqref{eq:multilinear_geral} with $k=2$.
Suppose that there exist $s\in\R$ and $\beta,\lambda\ge 0$ with $\beta+\lambda<1$ such that
			\begin{equation}\label{fre}
				\sup_{\xi} \int_{\xi=\xi_1+\xi_2}\frac{|m(\xi_1,\xi_2)|\jap{\xi}^s}{\jap{\xi_1}^s\jap{\xi_2}^s} \mathbbm{1}_{|\Phi-\alpha|<M} d\xi_1 \lesssim \jap{\alpha}^\lambda M^\beta,
			\end{equation}
for any $\alpha\in \R$ and $M>1$,
where $\Phi$ is as in \eqref{Phi}.
Then, choosing $b,\wt b\in\R$ satisfying
\begin{align}\label{freb}
 \beta<b<1-\lambda
\quad
\text{and}
\quad
 0<\wt{b}<b-\beta,
\end{align}
and $ 0 < T \le 1$, the following estimate holds
\begin{equation}\label{eq:bilin_infty}
    \left\|\psi_T (t)\int_0^t e^{-(t-t')\partial_x^3}\Ncal[z_1,z_2](t')dt' \right\|_{X^{s,b}_{\infty, \infty} \cap X^{s,\tilde{b}}_{\infty,1}}  \les T^{0+}\|z_1\|_{X^{s,b}_{\infty, \infty} \cap X^{s,\tilde{b}}_{\infty,1}}\|z_2\|_{X^{s,b}_{\infty, \infty} \cap X^{s,\tilde{b}}_{\infty,1}}.
\end{equation}

		\end{lemma}

		\begin{proof}
Assume that \eqref{fre} holds and take $b,\wt b$ satisfying \eqref{freb} and $b', \wt b'$ such that
			$$
			 b < 1+b'< 1-\lambda
             \quad
             \text{and}
             \quad
             -1+\wt{b}<\wt{b}'<b-\beta-1.
$$
From the linear estimates \eqref{lin2a} and \eqref{lin2b} in Lemma~\ref{lem:linear_estimates}, we have
        $$
        \left\|\psi_T(t)\int_0^t e^{-(t-t')\partial_x^3}\Ncal[z_1,z_2](t')dt' \right\|_{X^{s,b}_{\infty, \infty} \cap X^{s,\wt{b}}_{\infty,1}}  \lesssim (T^{1+b'-b} + T^{2+\wt{b}'-b})\left\|\Ncal[z_1,z_2] \right\|_{X^{s,b'}_{\infty, \infty} \cap X^{s,\wt{b}'}_{\infty,1}}.
        $$

We now estimate the norms of $\Ncal$ above. Let $\M = \M(\vec\xi)$ given by
\begin{equation}
\mathcal{M}(\vec\xi)=\frac{|m(\vec\xi)|\jap{\xi}^s}{\jap{\xi_1}^s\jap{\xi_2}^s},
\end{equation}
and note that
\begin{equation}
(\tau - \xi^3) - \Phi = (\tau_1 - \xi_1^3 ) + (\tau_2 - \xi_2^3),
\end{equation}
for $\tau=\tau_1+ \tau_2$, $\xi=\xi_1+\xi_2$, and $\Phi$ as in \eqref{Phi}. Then, by symmetry, we assume that $|\Phi - (\tau-\xi^3) | \les |\tau_1 - \xi_1^3|$, from which get
	\begin{align}
				\|\Ncal[z_1,z_2] \|_{X^{s,b'}_{\infty, \infty}} &\lesssim \sup_{\tau,\xi} \int_{\substack{\tau=\tau_1+\tau_2 \\ \xi = \xi_1 + \xi_2}}
|m(\vec\xi) | \jb{\xi}^s
\jap{\tau-\xi^3}^{b'}
|\ft{z}_1(\tau_1, \xi_1)
\ft{z}_2(\tau_2, \xi_2)
|
d\tau_1 d\xi_1
\\
&
\lesssim
 \|z_1\|_{X^{s,b}_{\infty,\infty}}
 \|z_2\|_{X^{s,0}_{\infty, 1}}
\sup_{\xi}
\int_{\xi = \xi_1 + \xi_2}
\frac{|m(\vec\xi)| \jb{\xi}^s}{\jb{\xi_1}^s \jb{\xi_2}^s}
\frac{\jb{\tau-\xi^3}^{b'}}{\jb{\Phi - (\tau-\xi^3)}^{b}}
d\xi_1
\\
&
\lesssim
 \|z_1\|_{X^{s,b}_{\infty,\infty}}
 \|z_2\|_{X^{s,0}_{\infty, 1}}
\sup_{\alpha,\xi}
\int_{\xi = \xi_1 + \xi_2}  \mathcal{M}(\vec\xi)\frac{\jap{\alpha}^{b'}}{\jap{\Phi-\alpha}^b  } d\xi_1 .
\label{FREinfty1}
			\end{align}

		\noi
Then, the estimate for the $X^{s, b'}_{\infty, \infty}$-norm follows once we control the last factor above, which holds due to \eqref{fre}:
			\begin{align}
			&
            \sup_{\alpha,\xi} \int \mathcal{M}(\vec\xi)\frac{\jap{\alpha}^{b'}}{\jap{\Phi-\alpha}^b  } d\xi_1
            \\
&
\les
\sup_{\alpha,\xi}
\ \jap{\alpha}^{b'}
\bigg(
\int \mathcal{M}(\vec\xi)\mathbbm{1}_{|\Phi-\alpha|\le 1} d\xi_1
+
\sum_{M\ge 1 \text{ dyadic}} \frac{1}{M^{b}}\int \mathcal{M}(\vec\xi)\mathbbm{1}_{|\Phi-\alpha|\sim M} d\xi_1
\bigg)
\\
&\les
\sup_{\alpha,\xi}\sum_{M \text{ dyadic}} \frac{\jap{\alpha}^{b'}}{M^{b}}\int \mathcal{M}(\vec\xi)\mathbbm{1}_{|\Phi-\alpha|<M} d\xi_1
\\
&
\lesssim \sup_{\alpha,\xi}\sum_{M \text{ dyadic}} \frac{\jap{\alpha}^{b'+\lambda}}{M^{b-\beta}}
<\infty,
			\end{align}
			since $b'\le -\lambda$ and $b>\beta$.

For the $X^{s,\wt{b}'}_{\infty, 1}$-norm of $\Ncal$, we have
\begin{equation}
				\|\Ncal[z_1,z_2]\|_{X^{s,\wt{b}'}_{\infty,1}} \lesssim \sup_\xi
\int  \int_{\substack{\tau=\tau_1+\tau_2 \\ \xi = \xi_1 + \xi_2}}
|m(\vec\xi)|
\jb{\xi}^s
 \jap{\tau-\xi^3}^{\wt{b}'}  | \ft{z}_1 (\tau_1, \xi_1) \ft{z}_2(\tau_2, \xi_2)|  d\tau_1 d\xi_1 d\tau.
			\end{equation}

        \noi
			If $|\tau-\xi^3|\gtrsim |\Phi|$, then
			\begin{align}
				\|\Ncal[z_1,z_2]\|_{X^{s,\wt{b}'}_{\infty,1}}
&
\lesssim
\bigg(
\prod_{j=1, 2}\|z_j\|_{X^{s,\wt{b}}_{\infty,1}}
\bigg)
\sup_\xi \int_{\xi=\xi_1+\xi_2} \mathcal{M}(\vec\xi) \jap{\Phi}^{\wt{b}'} d\xi_1  \\
&
\lesssim
\bigg(
\prod_{j=1, 2}\|z_j\|_{X^{s,\wt{b}}_{\infty,1}}
\bigg)
\sum_{M \text{ dyadic}} M^{\wt{b}'} \sup_\xi \int \mathcal{M} (\vec\xi)
\mathbbm{1}_{|\Phi|<M} d\xi_1
 \lesssim
\prod_{j=1, 2}\|z_j\|_{X^{s,\wt{b}}_{\infty,1}} ,
			\end{align}
			by \eqref{fre} and $\wt{b}'+\beta<b-1<0$.
Otherwise, $|\tau-\xi^3|\ll |\Phi|$ and we may assume that $|\tau_1-\xi_1^3|\gtrsim |\Phi|$. Then, \eqref{fre} gives

			\begin{align}
				\|\Ncal[z_1,z_2]\|_{X^{s,\wt{b}'}_{\infty,1}}
&
\lesssim \sup_\xi \int \int_{\substack{\tau=\tau_1+\tau_2\\ \xi=\xi_1+\xi_2}}
|m(\vec\xi)|
\jb{\xi}^s \frac{\jap{\tau-\xi^3}^{\wt{b}'}}{\jap{\tau_1-\xi_1^3}^b} \jap{\tau_1-\xi_1^3}^b
|\ft{z}_1(\tau_1, \xi_1)
\ft{z}_2(\tau_2, \xi_2)
|
 d\tau_2 d\xi_1 d\tau
 \\
&
\lesssim
\|z_1\|_{X^{s,b}_{\infty,\infty}} \|z_2\|_{X^{s,\wt{b}}_{\infty,1}}
\sup_\xi \int \frac{\mathcal{M}(\vec\xi)}{\jap{\Phi}^{b-\wt{b}'-1-}\jap{\tau-\xi^3}^{1+} } d\tau d\xi_1
\\
&
\lesssim
\|z_1\|_{X^{s,b}_{\infty,\infty}} \|z_2\|_{X^{s,\wt{b}}_{\infty,1}}
\sum_{M \text{ dyadic}} M^{1-b+\wt{b}'+} \sup_\xi \int \mathcal{M}(\vec\xi) \mathbbm{1}_{|\Phi|<M} d\xi_1
\\
&
 \lesssim\|z_1\|_{X^{s,b}_{\infty,\infty}}  \|z_2\|_{X^{s,\wt{b}}_{\infty,1}}
			\end{align}

            \noi
			since $b- \be - 1 >\wt b '$, completing the proof.

		\end{proof}

We recall the next lemma from \cite{CLS}, on bilinear estimates in $L^2$-based spaces.
		\begin{lemma}[{\cite[Lemma 3.1]{CLS}}]\label{lem:fre_bi}
            Let $s\in\R$.
			Suppose that, for all $\al\in \R$ and $M>1$, there exists $0<\beta<1$, such that one of the following holds:
			\begin{enumerate}
				\item[(i)] for all distinct $ k, \l\in\{\emptyset,1,2\}$\footnote{The empty set symbol $\emptyset$ here corresponds to the frequency component $\xi$, without index.}
                and $1<M\lesssim |\alpha|$,
				\begin{equation}\label{eq:frequad}
					\sup_{\xi_k}\int_{\xi=\xi_1+\xi_2}
                    \frac{|m(\xi_1,\xi_2)|^2\jap{\xi}^{2s}}{\jap{\xi_{1}}^{2s}\jap{\xi_{2}}^{2s}}\mathbbm{1}_{|\Phi-\alpha|<M} d\xi_{\l} \lesssim \jap{\alpha}^\beta M^\beta;
				\end{equation}

				\item[(ii)]
                there exist distinct $k, \l\in\{1,2\}$ such that
				\begin{equation}\label{eq:frequad2}
					\sup_{\xi_k}\int_{\xi=\xi_1+\xi_2}
                    \frac{|m(\xi_1,\xi_2)|^2\jap{\xi}^{2s}}{\jap{\xi_{1}}^{2s}\jap{\xi_{2}}^{2s}}\mathbbm{1}_{|\Phi-\alpha|<M} d\xi_{\l}  \lesssim \jap{\alpha} M^\beta;
				\end{equation}

				\item[(iii)] or we have that
				\begin{equation}\label{eq:frequad3}
					\sup_{\xi}\int_{\xi=\xi_1+\xi_2}\frac{|m(\xi_1,\xi_2)|^2\jap{\xi}^{2s}}{\jap{\xi_{1}}^{2s}\jap{\xi_{2}}^{2s}}\mathbbm{1}_{|\Phi-\alpha|<M} d\xi_{1}  \lesssim \jap{\alpha}^\beta M;
				\end{equation}
			\end{enumerate}

            \noi with $\Phi$ as in \eqref{Phi} with $k=2$.
			Then, for $b=\frac12+$ and $0<T \le 1$, we have
			\begin{equation}
				\left\|\psi_T (t)\int_0^te^{-(t-t')\partial_x^3} \Ncal[z_1,z_2](t')dt'\right\|_{X^{s,b}_{2}} \lesssim T^{0+}\|z_1\|_{X^{s,b}_{2}} \|z_2\|_{X^{s,b}_{2}},
			\end{equation}
            where $\Ncal$ is as in \eqref{eq:multilinear_geral}.
		\end{lemma}

		We now consider $\FL^{\infty}$-estimates for multilinear operators of the form \eqref{eq:multilinear_geral} with higher homogeneity.

		\begin{lemma}\label{lem:reduction_fre_kge3}
			Fix $2\le \ell\le k$ and $0<T \le 1$. Suppose that there exist $0\le \beta<1$ and multipliers $\mathcal{M}_1:\R^\l \to \R^+$ and $\mathcal{M}_2:\R^{k+1-\ell}\to \R^+$ such that
			\begin{equation}
				\frac{|m(\xi_1, \ldots, \xi_k)|\jap{\xi}^s}{\prod_{j=1}^k\jap{\xi_j}^s} \lesssim \mathcal{M}_1(\xi_1,\dots, \xi_{\l})\mathcal{M}_2(\xi, \xi_{\ell+1},\dots, \xi_{{k}}),
			\end{equation}

            \noi
            satisfying the following bounds
            \begin{align}
                \sup_{\xi, \xi_{\ell+1},\dots , \xi_{k},\alpha} \int_{\xi=\xi_1 + \cdots + \xi_k} \mathcal{M}_1(\xi_1,\dots,\xi_\l) \mathbbm{1}_{|\Phi-\alpha|<M} d\xi_1\cdots d\xi_{\ell-1} &
                \lesssim
                M^\beta
                 \label{fre_kge3},
                 \\
                 \sup_\xi\int \mathcal{M}_2(\xi, \xi_{\ell+1},\ldots, \xi_{{k}}) d\xi_{\ell+1}\cdots d\xi_k
                &
                \les 1,
                \label{fre_kge32}
            \end{align}
            for $\Phi$ as in \eqref{Phi}.
            Then, given $b,\wt b\in \R$ with
            $$
            1-\frac{1-\beta}{\ell}<b<1,\quad 0<\wt{b}<b-\beta,
            $$
            the following  bound holds:
			\begin{align}\label{frek}
				\left\|\psi_T (t)\int_0^te^{-(t-t')\partial_x^3}\Ncal[z_1,\dots, z_k](t')dt' \right\|_{X^{s,b}_{\infty, \infty}\cap X^{s,\wt{b}}_{\infty, 1}} & \les T^{0+}\prod_{j=1}^k\|z_j\|_{X^{s,b}_{\infty, \infty} \cap X^{s,\wt{b}}_{\infty,1}}
                ,
			\end{align}
			for $\Ncal$ as in \eqref{eq:multilinear_geral}.
		\end{lemma}

		\begin{proof}

        Let $b', \wt b' \in \R$ such that \eqref{bs} holds
        {and $\wt b' < b - 1 - \be$.}
        From Lemma~\ref{lem:linear_estimates}~(iii),
        \begin{align}
            \bigg\|\psi_T (t)\int_0^t e^{-(t-t')\partial_x^3}\Ncal[z_1,\ldots, z_k](t')dt' \bigg\|_{X^{s,b}_{\infty, \infty} \cap X^{s, \wt b}_{\infty, 1}}
            &
            \les T^{0+} \| \Ncal[z_1, \ldots, z_k] \|_{X^{s, b'}_{\infty, \infty} \cap X^{s, \wt b '}_{\infty, 1} } .
        \end{align}

        We first estimate the $X^{s, b'}_{\infty, \infty}$-norm. Since $b'<0<\wt b$ and using \eqref{fre_kge32}, we have
\begin{align*}
				&\|\Ncal[z_1,\dots, z_k] \|_{X^{s,b'}_{\infty, \infty}}
                \\
                & \les
                \bigg(\prod_{j=1}^\ell \|z_j\|_{X_{\infty,\infty}^{s,b}}
                \bigg)
                \sup_{\tau, \xi} \int_{\substack{\tau=\tau_1 + \cdots + \tau_k \\ \xi = \xi_1 + \cdots + \xi_k}} \frac{|m(\xi_1, \ldots, \xi_k)|\jap{\xi}^s}{\prod_{j=1}^k\jap{\xi_j}^s \times \prod_{j=1}^\ell \jap{\tau_j-\xi_j^3}^b}\prod_{j=\ell+1}^k \jap{\xi_j}^s|\ft{z}_j| d\tau_2\cdots d\tau_k d\xi_2 \cdots d\xi_k
				\\
                &\les
                \bigg(\prod_{j=1}^\ell \|z_j\|_{X_{\infty,\infty}^{s,b}}
                \bigg)
                \sup_{\substack{\tau, \tau_{\ell+1},\dots,\tau_k\\ \xi, \xi_{\ell+1}, \xi_k}} \left(\int_{\xi=\xi_1+\cdots+\xi_k} \mathcal{M}_1(\xi_1, \ldots, \xi_\l) \frac{1}{\prod_{j=1}^\ell \jap{\tau_j-\xi_j^3}^b}d\tau_2\dots d\tau_\ell d\xi_2\dots d\xi_\ell\right) \times
                \\ & \qquad \qquad \qquad \times \sup_\xi\int \mathcal{M}_2(\xi,\xi_{\l+1}, \ldots, \xi_k) \prod_{j=\ell+1}^k \jap{\xi_j}^s|\ft{z}_j| d\tau_{\ell+1}\cdots d\tau_k d\xi_{\ell+1} \cdots d\xi_k
				\\
                &
                \les
                \bigg(\prod_{j=1}^\ell \|z_j\|_{X_{\infty,\infty}^{s,b}}
                \bigg)
                \bigg( \prod_{j=\l+1}^k \|z_j\|_{X^{s, \wt b}_{\infty, 1}} \bigg)
                \sup_\xi\int  \mathcal{M}_2 (\xi, \xi_{\l+1}, \ldots, \xi_k)
                d\xi_{\l+1} \cdots d\xi_{k}
				\\
                & \qquad \qquad \qquad \times
                \sup_{\substack{\tau, \tau_{\ell+1},\dots,\tau_k\\ \xi, \xi_{\ell+1}, \xi_k}} \int_{\xi = \xi_1 + \cdots + \xi_k} \frac{ \mathcal{M}_1(\xi_1, \ldots, \xi_\l)
                }{ \jap{- \Phi + \tau-\xi^3 - \sum_{j=\ell+1}^k (\tau_j-\xi_j^3)}^{\ell b - \ell +1}} d\xi_2\cdots d\xi_\ell  \\
                &
                \les
                \bigg(\prod_{j=1}^k \|z_j\|_{X_{\infty,\infty}^{s,b} \cap X^{s, \wt b}_{\infty, 1}}
                \bigg)
                \sup_{{\alpha,  \xi, \xi_{\ell+1}, \xi_k}} \int_{\xi = \xi_1 + \cdots + \xi_k}  \frac{\mathcal{M}_1(\xi_1, \ldots, \xi_\l)}{ \jap{\Phi -\alpha}^{\ell b - \ell +1}} d\xi_2\cdots d\xi_\ell
                ,
			\end{align*}
for $1>b> 1- \frac1\l$.
			To control the last factor, we dyadically decompose  $\Phi-\alpha$ and use \eqref{fre_kge3}, to obtain
			\begin{align*}
				&\sup_{{\alpha,  \xi, \xi_{\ell+1}, \xi_k}} \int  \frac{\mathcal{M}_1(\xi_1, \ldots, \xi_\l)}{ \jap{\Phi -\alpha}^{\ell b - \ell +1}} d\xi_2\dots d\xi_\ell \\\lesssim\ & \sum_{M \text{ dyadic}} M^{-\ell b +\ell-1} \sup_{{\alpha,  \xi, \xi_{\ell+1}, \xi_k}} \int \mathcal{M}_1(\xi_1, \ldots, \xi_\l) \fia d\xi_2\dots d\xi_\ell \\\lesssim\ & \sum_{M \text{ dyadic}} M^{-\ell b +\ell-1+\beta} <\infty,
			\end{align*}
            for $b > 1 - \frac{1-\be}{\l}$.
			This concludes the proof of the first estimate. To estimate the $X^{s, \wt b'}_{\infty,1}$-norm, one proceeds as in the case $k=2$ in the proof of Lemma~\ref{lem:reduction_fre_k2}.
            {If $|\tau - \xi^3| \ges |\Phi|$, we place all factors in $X^{s, \wt b}_{\infty, 1}$, while when $|\tau-\xi^3| \ll |\Phi|$, assuming $|\Phi| \les |\tau_1-\xi_1^3|$, we place $z_1$ in $X^{s, b}_{\infty, \infty}$ and the remaining factors in $X^{s, \wt b}_{\infty, 1}$. We omit further details.}
\qedhere

		\end{proof}

        The multilinear estimate analogous to \eqref{frek} in $L^2$-based spaces follows from
        Lemma~\ref{lem:linear_estimates_qfinite}~(ii) and \cite[Lemma 2]{COS}.

		\begin{lemma}[{\cite[Lemma 2]{COS}}]\label{lem:fre_tri}
			Let $s\in\R$, $0<T\le 1$, $B\subset \{\emptyset, 1, \ldots, k\}$ be a proper non-empty subset, and  $\M_1 , \M_2: \R^{k} \to \R^+$ be multipliers satisfying
            \begin{align}
            \frac{|m(\xi_1, \ldots, \xi_k)|^2\jap{\xi}^{2s}}{\prod_{j=1}^k\jap{\xi_{j}}^{2s}}
            \les
                \M_1(\vec\xi) \M_2(\vec\xi)
                ,
            \end{align}
            where $\vec\xi=(\xi_1, \ldots, \xi_k)$.
            Then, if there exists $0<\be<1$ such that
			\begin{equation}\label{L2hypo}
			\sup_{\xi_{j\in B}} \int_{\xi=\xi_1+\ldots + \xi_k} \mathcal{M}_1 (\vec\xi) \mathbbm{1}_{|\Phi-\alpha|<M} d\xi_{j\notin B} + \sup_{\xi_{j\notin B}} \int_{\xi=\xi_1+\ldots+\xi_k} \mathcal{M}_2 (\vec\xi) \mathbbm{1}_{|\Phi-\alpha|<M} d\xi_{j\in B} \lesssim M^\beta,
			\end{equation}
            for all $M\ge 1$ and $\al\in\R$,
            where $\Phi$ is as in \eqref{Phi},
			then the following estimate holds for {$b=\frac12+$}, depending only on $\be$,
            and $\Ncal$ as in \eqref{eq:multilinear_geral}:
            \begin{align}
				\bigg\|\psi_T(t)\int_0^te^{-(t-t')\partial_x^3}\Ncal[z_1,\dots, z_k](t')dt' \bigg\|_{X^{s,b}_{2}} & \les T^{0+}\prod_{j=1}^k\|z_j\|_{X^{s,b}_{2}}.
			\end{align}
		\end{lemma}

		We conclude this subsection with some technical lemmas which are useful in showing frequency-restricted estimates.
        The first claims that one can extend the condition $\be<1$ in Lemmas \ref{lem:reduction_fre_k2} to \ref{lem:fre_tri} to include the endpoint $\be=1$, at the cost of an $\eps$-loss of the largest frequency.

		\begin{lemma}[{\cite[Lemma 3.3]{CLS}}]\label{lem:eta1}
			Let $K:\R^n\times \R^m\times \R\to \R^+$ be a positive measurable function such that, for some $N\in \mathbb{N}$,
			$$
			|K(x,y,M)|\lesssim \max\{1,|x|,|y|\}^N,\quad \forall x\in\R^n,\ y\in \R^m, \ M\in \R.
			$$
			Suppose that
			$$
			\sup_y \int K(x,y,M)\max\{1,|x|,|y|\}^{0+} dx \lesssim M,\quad \mbox{for all }M>1.
			$$
			Then, there exists $0<\beta<1$ such that
			$$
			\sup_y \int K(x,y,M) dx \lesssim M^{\beta},\quad \mbox{for all } M>1.
			$$
		\end{lemma}

Given an integral of the form
$$
\int \mathcal{M}(\vec\xi)\mathbbm{1}_{|\Phi-\alpha|<M}d\vec{\xi},
$$
with $\Phi$ as in \eqref{Phi},
our approach follows the guidelines:
\begin{itemize}
\item If $|\nabla \Phi|$ is large (non-stationary case), then we may perform a change of variables $\xi_j\mapsto \Phi$ for some $j$. The integration in $\Phi$ yields directly a factor of $M$, while the remaining integral have to be shown to be bounded.
\item If $|\nabla \Phi|$ is small (stationary case), then we are near a stationary point of the resonance function $\Phi$. Assuming that the stationary point is non-degenerate, we can replace the resonance function by a quadratic form by means of Morse's Lemma.

\end{itemize}

		\begin{lemma}[Morse's Lemma with parameters, {\cite[Lemma C.6.1]{Hormander}}]
			\label{LEM:morse}
		Let $\Phi: \R^2\times \R^d \to \R$, $\Phi=\Phi(x,y;z)$ be a smooth function such that
        $$D_{xy} \Phi(x_0,y_0;z_0) = 0,\qquad\det(D^2_{xy} \Phi)(x_0,y_0;z_0) \ne 0.$$
Then, for $|z-z_0|$ small, there exist a unique critical point $(x_0(z),y_0(z))$ near $(x_0,y_0)$ and  a smooth local change of variables $f^z: (p,q) \mapsto (x,y) $ such that $f^z(0,0) = (x(z),y(z))$ and
\begin{align}
\Phi(x,y;z) = \Phi(x_0(z), y_0(z)) + \ld_1 p^2 + \ld_2 q^2,
\end{align}
where $\ld_1,\ld_2 \in \{\pm1\}$ are the signs of the eigenvalues of $D^2\Phi(x_0,y_0)$.
		\end{lemma}

       After application of Morse's Lemma, the frequency-restricted estimates will often be reduced to sublevel estimates for quadratic forms.

		\begin{lemma}[{\cite[Lemma 5]{COS}}]
        \label{lem:quadraticas}
			Given $\alpha\in \R$ and $M\ge 1$, the following estimates hold:
			\begin{align}
            \label{eq:posdef}
				\iint_{|q_1|, |q_2|\le 1} \m 1_{|q_1^2\pm q_2^2-\alpha|<M}dq_1dq_2 &\lesssim M^{1-}, \\
                \label{eq:quadratica}
				\int_{-1}^{1} \m 1_{|q^2  -\alpha | < M} dq
                &\lesssim \sqrt M
                .
			\end{align}

		\end{lemma}

\section{Multilinear estimates}\label{sec:multi}

This section is dedicated to
estimating the nonlinear terms $\Ncal_\l(\TT;z)$, $\l=1, \ldots, 4$,  appearing in the gauged equation \eqref{z1} for $z$, in suitable Fourier restriction spaces $X^{s, b}_{p}$.
To that end, we consider the associated multilinear operators defined as
\begin{align}
\label{Ncal}
    \Ft_x\big[ \Ncal_\l (\TT; z_1, \ldots, z_{j+1} ) \big](\xi)
    &
    =
    i \int_{\G(\TT)} \mul_\l(\TT)(\vec\xi) \prod_{k=1}^{j+1} \ft{z}_k(t, \xi_{\bul_k}) d\vec\xi,
\end{align}
with $\G(\TT)$ in \eqref{Gamma} and considering the lexicographical ordering of the final nodes in $\TT$ in \eqref{lexic}.
We start with the case $p=\infty$.

\begin{proposition}[$\FL^\infty$-estimates]\label{prop:bourgLinfty}
Let $0< T \le 1$, $j \in\N$, $\TT \in \Tr_j$, $\Ncal_\l(\TT)$, $\l=1,\dots, 4$, as in \eqref{Ncal}, and $s$ satisfying
\begin{align}\label{eq:defisinfty}
s>s_\infty(\TT):=
\begin{cases}
- 1 , & \text{if }\l= 1 , 2 , \text{ or } \l=3, |\TT^0| = 2, 
\\
-\frac34, & \text{if }  \l=4 \text{ or }  \l=3, |\TT^0|\ge 3.
\end{cases}
\end{align}
Then, there exist
$b=1-$ and $\wt{b}=0+$, such that
	\begin{multline}
		\left\| \int_0^t e^{-(t-t')\dx^3}\Ncal_\l(\TT; z_1,\dots, z_{j+1})(t')dt'\right\|_{X^{s,b}_{\infty, \infty} (0,T)\cap X^{s,\wt b}_{\infty, 1}(0,T)}
        \\
        \le T^{0+} C^j\prod_{k=1}^{j+1}\|z_k\|_{X^{s,b}_{\infty, \infty} (0,T) \cap X^{s,\wt{b}}_{\infty,1}(0,T)}, \label{eq:Linftybourg1}
	\end{multline}
for some constant $C>0$ independent of $j$ and $\TT$.

\end{proposition}

We also prove analogous multilinear estimates in $L^2$-based spaces.
\begin{proposition}[$L^2$-estimates]\label{prop:bourgL2}
Let $0< T \le1$, $j \in\N$, $\TT \in \Tr_j$, $\Ncal_\l(\TT)$, $\l=1,\dots,4$, as in \eqref{Ncal}, and $s$ satisfying 
\begin{align}\label{eq:defis2}
s>s_2(\TT):=
\begin{cases}
- 1 , & \text{if }\l= 1 , 2 ,
\\
-\frac34, & \text{if }   \l=3, |\TT^0| =2, 
\\
-\frac78, & \text{if }  \l=4 \text{ or }  \l=3, |\TT^0|\ge 3.
\end{cases}
\end{align}
	Then, there exists $b=\frac12+$ such that
	\begin{align}\label{eq:bourg_bi_l2}
		\left\|\int_0^t e^{-(t-t')\dx^3} \Ncal_\l(\TT; z_1,\dots, z_{j+1})(t') dt'  \right\|_{X^{s,b}_{2}(0,T)}  \le T^{0+} C^j\prod_{k=1}^{j+1}\|z_k\|_{X^{s,b}_{2} (0,T)}
	\end{align}
    for some constant $C>0$ independent of $j$ and $\TT$.
\end{proposition}

Lastly, by interpolation between the results in \eqref{eq:Linftybourg1} and \eqref{eq:bourg_bi_l2}, we obtain the following estimates in $\FL^{p}$-based Fourier restriction norm spaces, for $2 < p <\infty$.

\begin{proposition}[$\FL^p$-estimates]\label{prop:bourgLp}
Given $2 < p <\infty$, let $0<T \le 1$, $j\in\N$,  $\TT\in \Tr_{j}$, and $\Ncal_\l(\TT)$, $\l=1,\dots, 4$, as in \eqref{z1}, and $s$ satisfying
$$
s>s_p(\TT):=\begin{cases}
-1, &  \text{if }\l= 1 , 2 , 
\\
-1+\frac1{2p}, & \text{if }   \l=3, |\TT^0| =2, 
\\
-\frac34-\frac1{4p}, & \text{if }  \l=4 \text{ or }  \l=3, |\TT^0|\ge 3.
\end{cases}
$$ 
Then
there exists $b=\frac{1}{p'}+$ such that
	\begin{align}\label{eq:bourg_bi_lp}
		\left\|\int_0^t e^{-(t-t')\dx^3} \Ncal_\l(\TT; z_1,\dots, z_{j+1})(t') dt'  \right\|_{X^{s,b}_{p}(0,T)}
        \le T^{0+} C^j\prod_{k=1}^{j+1}\|z_k\|_{X^{s,b}_{p}(0,T)}
	\end{align}
    for some constant $C>0$ independent of $j$ and $\TT$.
\end{proposition}

We postpone the proofs of Propositions~\ref{prop:bourgLinfty}, \ref{prop:bourgL2}, and \ref{prop:bourgLp} to Subsections~\ref{SUB:Linfty}, \ref{SUB:L2}, and \ref{SUBSEC:bourgainLp}, respectively.

\subsection{Proof of Proposition~\ref{prop:bourgLinfty} ($\FL^\infty$-estimates)}
\label{SUB:Linfty}

In this subsection, we show \eqref{eq:Linftybourg1} in Proposition~\ref{prop:bourgLinfty} via a reduction to frequency-restricted estimates.
Observe that the multilinear operator
$\Ncal_\l(\TT; z_1, \ldots, z_{j+1})$ (defined in \eqref{Ncal}) is of the form $\Ncal(z_1, \ldots, z_{j+1})$ in \eqref{eq:multilinear_geral} with multiplier $m=i \mul_\l(\TT)$, as in \eqref{wint0c}.
 Thus, estimate  \eqref{eq:Linftybourg1}  follows from the inhomogeneous estimates \eqref{lin2a}-\eqref{lin2b}, together with Lemma~\ref{lem:reduction_fre_k2} and Lemma~\ref{lem:reduction_fre_kge3}, once we verify the respective hypotheses.

 For the sake of simplicity, we only present the proof for $s<0$, and we estimate each $\Ncal_\l(\TT)$ contribution based on the tree classification in Lemma~\ref{LEM:wint}.

\smallskip
\noi
$\bul$ \textbf{Type I.}
From \eqref{wint0c}, the only tree for which $\mul_1(\TT)\neq 0$ is $\TT=\<1>$. Then, we have
\begin{align}\label{type1}
\jap{\xi}^s
| {\mul}_1(\TT) |
\prod_{\bul\in \TT^\infty} \jap{\xi_\bul}^{-s}
=
\frac{\jap{\xi}^s |\xi|}{\jap{\xi_1}^s \jap{\xi_2}^s} \ind_{A}.
\end{align}
By symmetry, we assume that $|\xi_1| \ge |\xi_2|$ and $|\xi_2| \le 1$ due to the restriction to the set $A$ in \eqref{A}.
We consider two cases: (i) $|\xi_1|\les 1$, or (ii) $|\xi| \sim |\xi_1| \sim |\xi_1-\xi_2| \gg 1 \ge |\xi_2|$.
If (i) holds, the multiplier in \eqref{type1} is bounded since all frequencies all small, and
$$
\sup_{\al, M}
\int_{\xi=\xi_1+\xi_2}\frac{ \jb{\xi}^{s+ } }{\jb{\xi_1}^s \jb{\xi_2}^s } \mul_{1} (\TT) \ind_{|\Psi(\TT) - \al| < M }d\xi_1 \lesssim 1.
$$
Thus, we consider (ii), where we note that $\phi = \Psi(\<1>) = -3 \xi \xi_1(\xi-\xi_1)$ satisfies
\begin{align}\label{partial-phi}
| \partial_{\xi_1} \phi | = 3 | \xi (\xi_1 - \xi_2)| \sim |\xi|^2.
\end{align}
Combining \eqref{type1} with the change of variables $\xi_1 \mapsto \phi $, we obtain
	\begin{align*}
    \sup_\xi
\int_{\xi=\xi_1+\xi_2}
\frac{ \jb{\xi}^{s+ } }{\jb{\xi_1}^s \jb{\xi_2}^s } \mul_{1} (\TT) \ind_{|\Psi(\TT) - \al| < M }d\xi_1
 &
\les
\sup_\xi
\int_{A} \frac{\jap{\xi}^{1 + s + }}{\jap{\xi_{1}}^{s}\jap{\xi_{2}}^{s}} \ind_{|\Psi(\TT) - \al| < M } d\xi_{1}
\\&\les
\sup_\xi
 \int \frac{\ind_{|\phi-\al| < M} }{ \jap{\xi}^{1 - }} d\phi \lesssim M,
	\end{align*}
	for $s>-1$. Therefore, \eqref{fre} (with $\lambda=0$ and $\beta=1-$) follows from Lemma~\ref{lem:eta1}.

\medskip
\noi
$\bul$ \textbf{Type $\II$.} From \eqref{wint0c}, we only need to consider $\TT\in \Tr_j$, $j\ge 2$, with $\TT^{0,\ff}=\{\star\}$. We control the multiplier by
\begin{align}\label{L22}
\jb{\xi}^s |\mul_2(\TT)| \prod_{\bul \in \TT^\infty} \frac{1}{\jb{\xi_\bul}^s}
&
\le
{C^j}
\frac{ \ind_{A_\star} \ind_{A^c_{\pb(\star)}} }{ \jb{\xi_{\star1}}^s \jb{\xi_{\star2}}^s  \jb{\xi_{\sbf(\star)}}^{s+1} } \prod_{\bul \in \TT^{\infty}_\star } \frac{1}{\jb{\xi_{\bul}}^{s+1} \jb{\xi_{\sbf(\bul)}}}
\\
&
=:
{C^j}
\mathcal{M}_1(\xi_{\star1}, \xi_{\star2}, \xi_{\sbf(\star)})\mathcal{M}_2(\xi, \xi_{\bul\in \TT^\infty_\star})
,
\end{align}
 where
 $\TT^\infty_\star : = \TT^\infty\setminus \{\sbf(\star), \star1, \star2\}$,
 recalling that $\sbf(\star)$ is the sibling node of $\star$.
Since $s>-1$, we can handle the $\M_2$ part of the multiplier
\begin{equation}\label{eq:estMlinha}
    \sup_\xi\int \mathcal{M}_2(\xi, \xi_{\bul\in \TT^\infty_\star})  d\xi_{\bul\in\TT^\infty_\star} \lesssim \bigg(
\prod_{\bul \in \TT^\infty_\star }
\int \frac{1}{ \jb{\wt\xi_\bul}^{s+2} } d\wt\xi_{\bul} \bigg)  \le C^j, 
\end{equation}
where $|\wt\xi_\bul| = |\xi_\bul| \land |\xi_{\sbf(\bul)}|$.
Hence, {for all $\al\in\R$ and $M>1$, }it remains to bound
\begin{align}
\label{IIbb}
\sup_{\xi, \xi_{\bul\in\TT^\infty_\star}}
\int_{\G(\TT)} \mathcal{M}_1(\xi_{\star1}, \xi_{\star2}, \xi_{\sbf(\star)})\ind_{|\Psi(\TT) - \al| < M} d\xi_{\star1} d\xi_{\star2}. 
\end{align}
Assume that
$|\xi_{\star1} | \ge |\xi_{\star2}|$, and thus $|\xi_{\star2}| \le 1$ from the restriction to $A_\star$.
If $|\xi_{\star1}| \les 1$, then  \eqref{IIbb} is trivially bounded since all frequencies are $\les 1$. Thus, we assume that $|\xi_{\star1}|\gg1$ in the following. 

If $|\xi_{\star1}|\not\simeq|\xi_{\sbf(\star)}|$, for fixed $\xi_{\star2}$, we perform the change of variables $\xi_{\star1} \mapsto \phi = \Psi(\TT)$, noting that $|\partial_{\xi} \phi|\sim (| \xi_{\star1}| \lor |\xi_{\sbf(\star)}|)^2$:
\begin{align}
\eqref{IIbb} 
&
\les 
\int \frac{\ind_{|\phi - \al| < M } \ind_{|\xi_{\star2}| \le 1}}{ \jb{\xi_{\star1}}^{s+2-}   } d\phi d\xi_{\star2}
\les M , 
\end{align}
for $s>-1$. 
Otherwise, we have $|\xi_{\star1}| \simeq | \xi_{\sbf(\star)}|$. Since $|\xi_{\star1}|\not\simeq |\xi_{\star2}|$, we instead fix $\xi_{\sbf(\star)}$ when doing the change of variables:
\begin{align}
\eqref{IIbb} 
&
\les 
\int \frac{\ind_{|\phi - \al| < M }}{ \jb{\xi_{\star1}}^{s+2-} \jb{\xi_{\sbf(\star)}}^{s+1}   } d\phi d\xi_{\sbf(\star)}
\les 
M \int \jb{\xi_{\star2}}^{-2s-3+} d\xi_{\sbf(\star)}
\les M , 
\end{align}
for $s>-1$.
With Lemma~\ref{lem:eta1}, we conclude that the assumptions \eqref{fre_kge3}-\eqref{fre_kge32} of Lemma~\ref{lem:reduction_fre_kge3} hold (for~$\be=1-$).

\smallskip

\noi
$\bul$ \textbf{Type $\III$, $j=2$.} 
From \eqref{wint0c}, we can restrict to $\TT= \<21>$, where the corresponding multiplier is bounded by
\begin{align}\label{mul3}
    \jb{\xi}^s | \mul_3 (\TT)|
    \prod_{\bul \in \TT^\infty} \frac{1}{\jb{\xi_\bul}^s}
    &
    \les
    \frac{\jb{\xi}^{s} |\xi|  \ind_{A_1^c} }{ \jb{\xi_{11}}^{s+1} \jb{\xi_{12}}^{s+1} \jb{\xi_{2}}^s}
    .
\end{align}
By symmetry, we assume $|\xi_{11} | \ge |\xi_{12}|$.

\smallskip
\noi
\underline{\textbf{Case 1:} $|\xi_{2}| \ges |\xi_{11}|$}

\smallskip
\noi
\underline{\textbf{Subcase 1.1:} $|\xi_{2}| \not\simeq |\xi_{12}|$}

Here, it holds that $|\xi| \les |\xi_2|$ and fixing $\xi, \xi_{11}$, we have
\begin{align}\label{IIIb-new}
\begin{split}
|\partial_{\xi_{2}} \Phi|
= 3 | (\xi-\xi_{11}-\xi_2)^2 - \xi_2^2|
=3|\xi_{12}^2 - \xi_2^2|
\sim |\xi_2|^2.
\end{split}
\end{align}
Thus, we perform the change of variables $\xi_2 \mapsto \phi = \Phi$ as in \eqref{IIIb-new}, to obtain
\begin{align}
\int_{\xi=\xi_{11} + \xi_{12} + \xi_2} \frac{ \ind_{A_1^c} \jb{\xi}^{s+1 }  }{ \jb{\xi_{11}}^{s+1  } \jb{\xi_{12}}^{s+1} \jb{\xi_{2}}^{s-} }
\ind_{|\Phi - \al| < M}
 d\xi_2 d\xi_{11}
&
\sim \int \frac{ \ind_{|\phi - \al| < M } }{ \jb{\xi_{11}}^{s+1} \jb{\xi_2}^{1-}   } d\phi d\xi_{11}
\\
&
\les  \int \frac{\ind_{|\phi - \al|< M }}{ \jb{\xi_{11}}^{s+2 - } } d\phi d\xi_{11} \les M,
\end{align}
given that $s>-1$. Applying Lemma~\ref{lem:eta1}, we find \eqref{fre_kge3} with $\be=1-$.

\smallskip
\noi
\underline{\textbf{Subcase 1.2:} $|\xi_{2}| \simeq |\xi_{12}|$}

Here, $|\xi_2| \sim |\xi_{11}| \sim |\xi_{12}|$. 
If $|\xi_\bul| \not\simeq |\xi_\circ|$ for a pair of distinct $\bul, \circ \in \{11,12,2\}$, then we can proceed as in Subcase 1.1, integrating in $\{\xi_{11},\xi_{12}, \xi_2\} \setminus \{\xi_{\circ}\}$, with the change of variables $\xi_\bul \mapsto \phi$. 
Otherwise, we have $|\xi_2| \simeq |\xi_{11}| \simeq |\xi_{12}|$ and $|\xi_2| \sim |\xi|$. The last claim is obvious if $\xi_{11}, \xi_{12}, \xi_2$ have the same sign, while if only two have the same sign, say $\xi_{11}, \xi_{12}>0>\xi_2$, then
\begin{align}
|\xi| = |\xi_{11} + \xi_{12} + \xi_2| 
\ge 
|\xi_{11}| - |\xi_{12} + \xi_2| 
= |\xi_{11}| - \big| |\xi_{12}| - |\xi_2| \big|
\sim |\xi_{11}|, 
\end{align}
since $\big| |\xi_{12}| - |\xi_2| \big| \ll |\xi_{12}| \sim |\xi_{11}|$.
We consider the normalization
\begin{align}
\Phi(\xi, \xi_{11}, \xi_2) = - \xi^3 + \xi_{11}^3 +  \xi_2^3 + (\xi - \xi_{11} - \xi_2)^3  &= \xi^3 \Phi(1, \tfrac{\xi_{11}}{\xi}, \tfrac{\xi_2}{\xi}) 
\\
&
=: \xi^3 \Phi(1, p_{11}, p_2),
\label{Phi1}
\end{align}
and use Morse's Lemma to find a suitable change of variables.
From the frequency assumptions, with $p_{12}:= 1 - p_{11} - p_2$, we have
\begin{align}
\big| |p_{11}| - |p_{12}| \big|,
\big| |p_{11}| - |p_2| \big|,
\big| |p_{12}| - |p_2| \big|
 \ll |p_{11}| \sim |p_{12} | \sim |p_2|  \sim 1,
\end{align}
from which we conclude that
$(p_{11}, p_2)$ must be around one of the four stationary points of $\Phi(1, p_{11}, p_2)$ in \eqref{Phi1}:
\begin{align}
(1, 1),
\quad
(1, -1),
\quad
(-1,1),
\quad
\text{and}
\quad
(\tfrac13, \tfrac13).
\end{align}

\noi
Looking at the Hessian $D^2 \Phi(1, p_{11}, p_2)$, we note that
{
$\det (D^2 \Phi) = 36 (p_{11} p_{12} + p_{11} p_2 + p_{12} p_2) \ne0 $
}
on the critical points above. Thus, applying Morse's Lemma (Lemma~\ref{LEM:morse}) around each point implies the existence of a local change of coordinates $(p_{11},p_2) \mapsto (q_{11}, q_2)$ such that $|q_{11}|, |q_2| \les 1$ and
\begin{align}
\Phi(1, p_{11},p_2) = \Phi(1, c) + \lambda_1q_{11}^2 + \lambda_2q_2^2, \quad \lambda_1,\lambda_2\in\{\pm1\},
\end{align}

\noi
{and $c$ being one of the critical points.}
Consequently, applying the change of variables $(\xi_{11},\xi_2) \mapsto (q_{11}, q_2)$, together with \eqref{eq:posdef}, we obtain
\begin{align}
&
\sup_{\xi, \al} \int_{\xi = \xi_{11} + \xi_{12} + \xi_2} \frac{\ind_{A_1^c}  \jb{\xi}^{s+1} }{ \jb{\xi_{11}}^{s+1} \jb{\xi_{12}}^{s+1} \jb{\xi_2}^s}
\ind_{|\Phi - \al| < M}
d\xi_{11} d\xi_2
\\
&
\les
\sup_{\xi, \al}
\jb{\xi}^{ - 2 s - 1   } \int \ind_{|\xi^3 \Phi(1, p_{11}, p_2) - \al| < M} \xi^2 dp_{11} dp_2
\\
&
\les
\sup_{\xi, \al}
\jb{\xi}^{- 2s + 1  } \int_{|q_{11}|, |q_2| \les 1} \ind_{|  \Phi(1,c) + \lambda_1q_{11}^2 +  \lambda_2q_2^2 - \al | < \frac{M}{|\xi|^3}  } dq_{11} d q_2
\\
&
\les
M^{1-}
\sup_{\xi}
\jb{\xi}^{- 2s -2 +  }  \les M^{1-}
\end{align}
for $s>-1$.

\smallskip
\noi
\underline{\textbf{Case 2:} $|\xi_{2}| \ll |\xi_{11}|$}

In this case, we have $|\xi_{11}|\sim|\xi|$ and $|\xi_{11}| \not\simeq |\xi_2|$ and can proceed as in Subcase 1.1 for $s>-1$, with the change of variables $\xi_{11}\mapsto \phi$:
\begin{align}
\int_{\xi=\xi_{11} + \xi_{12} + \xi_2} \frac{ \ind_{A_1^c} \jb{\xi}^{s+1 }  \ind_{|\Phi - \al| < M}  }{ \jb{\xi_{11}}^{s+1 - } \jb{\xi_{12}}^{s+1} \jb{\xi_{2}}^s }
 d\xi_{11} d\xi_{12}
&
\sim \int \frac{ \ind_{|\phi - \al| < M } }{ \jb{\xi_{11}}^{s+2-} \jb{\xi_{12}}^{s+1}   } d\phi d\xi_{12}
\\
&
\les  \int \frac{\ind_{|\phi - \al|< M }}{ \jb{\xi_{12}}^{2s+3 - } } d\phi d\xi_{12} \les M,
\end{align}
completing the estimate for Type $\III$ terms with $|\TT^0| = 2$.

\smallskip

\noi
$\bul$ \textbf{Type $\III$, $j\ge3$.} 
We can restrict to $\TT\in \Tr_j$, $j\ge 3$,  with $\TT^{0,\ff}=\{\star\}$. In this case,
\begin{align}
\jb{\xi}^s | \mul_3 (\TT)|
    \prod_{\bul \in \TT^\infty} \frac{1}{\jb{\xi_\bul}^s}
    &
    \le
    {C^j}
   \frac{ \ind_{A^c_{\star} \cap A^c_{\pb(\pb(\star))} } }{ \jb{\xi_{\star1}}^{s+1}  \jb{\xi_{\star2}}^{s+1}  \jb{\xi_{\sbf(\star)}}^{s}  \jb{\xi_{\sbf(\pb(\star))}}^{s+1}  } \times 
\prod_{\bul \in \TT^\infty_\star} \frac{1}{\jb{\xi_\bul}^{s+1} \jb{\xi_{\sbf(\bul)}}}
    \\
    &=: {C^j}
    \mathcal{M}_1(\xi_{\star1}, \xi_{\star2}, \xi_{\sbf(\star)}, \xi_{\sbf(\pb(\star))} ) \times \mathcal{M}_2(\xi,\xi_{\bul\in \TT^\infty_\star})
    ,
\end{align}
where $\TT_\star^\infty=\TT^\infty\setminus\{\star1,\star2, \sbf(\star), \sbf(\pb(\star))\}$ and $\sbf(\bul)$ (resp. $\pb(\bul)$) denotes the sibling (resp. parent) node of $\bul$. Since \eqref{eq:estMlinha} holds,  it remains to bound
\begin{align}
\label{IIIbb}
\sup_{ \xi, \alpha}
\int_{\xi=\xi_2 + \xi_{12} + \xi_{111} + \xi_{112}} 
\mathcal{M}_1(\xi_{111}, \xi_{112}, \xi_{12}, \xi_{2} )
\ind_{|\Psi - \al| < M} d\xi_{2} d\xi_{12} d\xi_{112}, 
\end{align}
with $\Psi = \xi^3 - \xi_{2}^3 - \xi_{12}^3 - \xi_{111}^3 - \xi_{112}^3$. 
 From the localization to $A^c_{11} \cap A^c$ in \eqref{wint0c}, we have $|\xi_{112}|, |\xi_{2}| >1$ and,  by symmetry, we assume $|\xi_{111}| \ge |\xi_{112}|$.

\smallskip
\noi
\underline{\textbf{Case 1:} $|\xi_{12}| \ges |\xi_{2}| \ges |\xi_{111}| $.}

\smallskip
\noi
\underline{\textbf{Subcase 1.1:} $|\xi_{12}| \not\simeq |\xi_{112}|$.}

For fixed $\xi_{2}, \xi_{112}$, we change  variables $\xi_{12}\mapsto \phi = \Psi$, noting that $|\partial_{\xi_{12}} \phi| \sim |\xi_{12}|^2$:
\begin{align}
\eqref{IIIbb}
&
\les 
\sup_{ \xi, \alpha}
\int
\frac{  \ind_{|\phi - \al| < M}}{ \jb{\xi_{111}}^{s+1}  \jb{\xi_{112}}^{s+1}  \jb{\xi_{12}}^{s+2- } \jb{\xi_{2}}^{s+1}  }
d \phi
 d\xi_{2} d\xi_{112}
\\
&
\les 
M 
\int \frac{1}{   \jb{\xi_{112}}^{2(s+1) + \frac12-}  \jb{\xi_{12}}^{2 (s+1) + \frac12- }   }
 d\xi_{2} d\xi_{112}
\les M
\label{typeIIIa}
\end{align}
for $s>-\frac34$.

\smallskip
\noi
\underline{\textbf{Subcase 1.2:} $|\xi_{12}| \simeq |\xi_{112}|$.}

Here, we have $|\xi_{2}| \sim |\xi_{12}| \sim |\xi_{111}| \sim | \xi_{112}|$. If $|\xi_\bul|\not\simeq|\xi_\circ|$ for some distinct $\bul, \circ \in \{2, 12, 111, 112\}$, then we proceed as in Subcase 1.1, integrating in $\{ \xi_{2}, \xi_{12}, \xi_{111},  \xi_{112}  \} \setminus\{\xi_\circ\}$ and with the change of variables $\xi_\bul \mapsto \phi$. 
It remains to consider the case $|\xi_{2}| \simeq |\xi_{12}| \simeq |\xi_{111}| \simeq | \xi_{112}|$, where we proceed as in Subcase 1.2 for the Type $\III$ terms with $j=2$, namely defining variables $p_\bul = \frac{\xi_\bul}{\xi_{112}}$, $\bul \in \{2, 12, 112\}$, rewriting the partial phase as in \eqref{Phi1}, with $\xi_{112}$ replacing $\xi$, which with Morse's Lemma (Lemma~\ref{LEM:morse}) and Lemma~\ref{lem:quadraticas} gives
\begin{align}
\eqref{IIIbb}
&
\les 
\sup_{\xi, \al}
\int \ind_{|q_2|, |q_{12}| \les 1, \ | \Phi(1, q_2, q_{12}) - \al| < \frac{M}{|\xi_{112}|^3}} \jb{\xi_{112}}^{-4s-1+}
dq_{2} dq_{12} d\xi_{112}
\\
&
\les M^{1-} \int \jb{\xi_{112}}^{-4s-4+} d\xi_{112} \les M^{1-}
\label{typeIIIb}
\end{align}
for $s>-\frac34$.

\smallskip
\noi
\underline{\textbf{Case 2:} $|\xi_{12}| \ll |\xi_{2}|$ or $|\xi_{2}|\ll |\xi_{111}| $.}

Note that in this case, we always have $|\xi_{\max}| \not\simeq |\xi_{12}|$ where $|\xi_{\max}| = \max(|\xi_2|, |\xi_{12}|, |\xi_{111}|, |\xi_{112}|)$. Moreover, noting that the multiplier in \eqref{IIIbb} is better in this case, as we see an extra negative power of the largest frequency, as compared to Case 1, we can proceed as in Subcase~1.1, changing variables $\xi_{\max} \mapsto \phi$ and integrating in $\xi_{\max}, \xi_{12}$ and another frequency.

\smallskip
\noi
$\bul$ \textbf{Type IV, $j=3$.} In this case, the only relevant tree is $\TT=\<33>$ and the corresponding multiplier can be controlled as
$$
\jb{\xi}^s  |\mul_4(\TT)|
\prod_{\bul \in \TT^\infty} \frac{1}{\jb{\xi_\bul}^s}
\le
C^{j}
\frac{ \jb{\xi}^{1+s} \ind_{A^c_{1} \cap A^c_{2}} }{
\jb{\xi_{11}}^{s+1} \jb{\xi_{12}}^{s+1} \jb{\xi_{21}}^{s+1} \jb{\xi_{22}}^{s+1}}=:\mathcal{M}_1(\xi_{11}, \xi_{12}, \xi_{21}, \xi_{22}),
$$
Since we will not make use of any further relations between $\xi_{11}, \xi_{12}, \xi_{21}$ and $\xi_{22}$, we assume that $|\xi_{11}|\ge |\xi_{12}|\ge |\xi_{21}|\ge |\xi_{22}|>1$, which implies $|\xi_{11}|\gtrsim |\xi|$.

\smallskip
\noi
\underline{\textbf{Case 1:} $|\xi_{11}| \not\simeq |\xi_{22}|$.}

We perform the change of variables $\xi_{1}\mapsto \phi= \Psi(\TT)$, taking into account that $|\partial_{\xi_11}\phi|\sim |\xi_{11}|^2$:
\begin{align}
\sup_{\xi,\alpha} \int \jap{\xi_{11}}^{0+}\mathcal{M}_1\ind_{|\Psi(\TT)-\alpha|<M}d\xi_{11} d\xi_{12}d\xi_{21} &\lesssim \sup_{\xi,\alpha} \int \frac{\ind_{|\phi-\alpha|<M}}{\jap{\xi_{12}}^{s+1}\jap{\xi_{21}}^{s+1}|\xi_{11}|^{2-}} d\phi d\xi_{12} d\xi_{21}\\&\lesssim M,\label{eq:Tipo43}
\end{align}
for $s>-1$.

\smallskip
\noi
\underline{\textbf{Case 2:} $|\xi_{11}| \simeq |\xi_{22}|$.}

Fixing $\xi, \xi_{12}$, we consider
the normalization $p_\bul = \frac{\xi_\bul}{\xi_{12}}$ for $\bul\in\{ \emptyset, 11,  21, 22\}$ and write
\begin{align}
    \Psi(\TT)
    &
    =  - \xi^3 + (\xi-\xi_{12} - \xi_{21} - \xi_{22})^3 + \xi_{12}^3 +  \xi_{21}^3 +  \xi_{22}^3
    \\
    &
    =
     - \xi_{12}^3 \big( p^3 - (p-1 - p_{21}-p_{22})^3 - 1 - p_{21}^3 - p_{22}^3 \big)
    \\
    & = : \xi_{12}^3 P(p_{21}, p_{22} ).
\end{align}
From the frequency assumptions, we have
$$1\sim |p-1-p_{21} - p_{22}| \simeq  |p_{21} |\simeq |p_{22}|,$$
and the variables $(p_{21}, p_{22})$ are around one of the stationary points of $P$:
\begin{align}\label{crit}
    (p-1, p-1), \quad (\pm(p-1), \mp(p-1), \quad \tfrac13(p-1, p-1).
\end{align}
Calculating the Hessian $D^2 P(p_{21}, p_{22})$, we note that
\begin{align}
 \det(D^2 P)(p_{21}, p_{22})
 &
 = 36 (p_{11} p_{21} + p_{11} p_{22} + p_{21} p_{22} ) \sim 1
\end{align}
near each of the stationary points.
Consequently, by Morse's Lemma with parameters (Lemma~\ref{LEM:morse}), there exists a local change of coordinates $( p_{21}, p_{22}) \mapsto (q_{21}, q_{22})$ such that $|q_{21}|, |q_{22}| \les 1 $ and
\begin{align}
P(p_{21}, p_{22})
 = P(c) + \ld_{21} q_{21}^2 + \ld_{22} q_{22}^2,
\end{align}
for $c$ being one of the critical points in \eqref{crit} and $\ld_{21}, \ld_{22} \in \{\pm1\}$ depending on the sign of the eigenvalues of $D^2P$ at $c$.
With this diffeomorphism in hand and using Lemma~\ref{lem:quadraticas}, we obtain
\begin{align}
&\sup_{\xi,\alpha} \int \jap{\xi_{11}}^{0+}\mathcal{M}_1\ind_{|\Psi(\TT)-\alpha|<M}d\xi_{12} d\xi_{21}d\xi_{22} 
\\&
\les
\sup_{\xi, \al}
\int \frac{ 1 }{ \jb{\xi_{12}}^{3s+3} }
\bigg(
\int \ind_{|\xi_{12}^3 P( p_{21}, p_{22}) - \al| < M } \xi_{12}^2 dp_{21} dp_{22}
\bigg)
d\xi_{12}
\\
&
\les
\sup_{\xi, \al}
\int \frac{1}{\jb{\xi_{12}}^{3s+1}}
\bigg(
\int_{|q_{21}| ,  |q_{22}| \les 1}
\ind_{|\xi_{12}^3 (P(c) + \ld_{21} q_{12} + \ld_{22} q_{22}) - \al| < M  }
d q_{21} dq_{22}
\bigg)
d\xi_{12}
\\
&
\les M^{1-} \int \frac{1}{\jb{\xi_{12}}^{3s+4-}} d \xi_{12}
\les M^{1-},\label{eq:Tipo43Morse}
\end{align}
given that $s> -1$,

\smallskip
\noi
$\bul$ \textbf{Type IV, $j\ge 4$.}
From \eqref{wint0c}, it suffices to take
$\TT\in\Tr_j$ such that $\TT^{0,\ff} = \{\star1,\star2\}$ for some $\star\in \TT^0$, where $\star1,\star2$ denote the children of~$\star$.
Then, we have
\begin{align}\label{L24}
&\jb{\xi}^s  |\mul_4(\TT)|
\prod_{\bul \in \TT^\infty} \frac{1}{\jb{\xi_\bul}^s}
\\&\le
C^{j}
\frac{ \jb{\xi_\star}^{-s+} \ind_{A^c_{\star1} \cap A^c_{\star2}}}{
\jb{\xi_{\star11}}^{s+1} \jb{\xi_{\star12}}^{s+1} \jb{\xi_{\star21}}^{s+1} \jb{\xi_{\star22}}^{s+1}   }\cdot \frac{1}{\jb{\xi_\star}^{-s+}\jap{\xi_{\sbf(\star)}}^{1+s}} \prod_{\substack{\bul\in \TT^\infty_\star\\\bul\neq \sbf(\star)}} \frac{1}{\jb{\xi_\bul}^{s+1} \jb{\xi_{\sbf(\bul)}}}, \\ &=:
{C^j}
\mathcal{M}_1(\xi_{\star11}, \xi_{\star12}, \xi_{\star21}, \xi_{\star22})
\M_2(\xi,\xi_{\bul\in\TT^\infty_\star})
\end{align}
where $\TT^\infty_\star = \TT^\infty \setminus\{ \star11, \star12, \star21, \star22\}$. Similarly to \eqref{eq:estMlinha}, since $s>-1$, 
$$
\sup \xi \int \mathcal{M}_2(\xi,\xi_{\bul\in\TT^\infty_\star}) d \xi_{\bul\in\TT^\infty_\star} \lesssim C^j.
$$
As such, it suffices to derive a bound on the integral of $\M_1$, that is,
\begin{align}\label{eq:frequadd}
    \sup_{\xi, \al}
\int_{\xi=\xi_{11} + \xi_{12} + \xi_{21} + \xi_{22}} \mathcal{M}_1(\xi_{11}, \xi_{12}, \xi_{21}, \xi_{22}) \mathbbm{1}_{|\Psi(\TT')-\alpha|<M}d\xi_{11}d\xi_{12}d\xi_{21},
\end{align}
where $\TT'=\<33>$. We argue as for $j=3$ and assume that $|\xi_{11}|\ge |\xi_{12}|\ge |\xi_{21}|\ge |\xi_{22}|>1$.

\smallskip
\noi
\underline{\textbf{Case 1:} $|\xi_{11}| \not\simeq |\xi_{22}|$.} As in \eqref{eq:Tipo43}, the change of variables yields
\begin{align}
\eqref{eq:frequadd} 
&\les \sup_{ \xi, \al }
\int \frac{ \ind_{|\phi-\al|<M} }{ \jb{\xi_{11}}^{3+2s-}  \jb{\xi_{12}}^{s+1} \jb{\xi_{21}}^{s+1} } d\phi d\xi_{21} d\xi_{22} \\& \sup_{ \xi, \al }
\int \frac{ \ind_{|\phi-\al|<M} }{ \ \jb{\xi_{12}}^{2s+\frac52-} \jb{\xi_{21}}^{2s+\frac52} } d\phi d\xi_{21} d\xi_{22}\lesssim M,
\end{align}
for $s>-\frac34$.

\smallskip
\noi
\underline{\textbf{Case 2:} $|\xi_{11}|\simeq |\xi_{22}|$.}

Proceeding as in Case 2 for the Type $\IV$ term with $j=3$, the analogue of \eqref{eq:Tipo43Morse}~is
\begin{align}
&\sup_{\xi,\alpha} \int \jap{\xi_{11}}^{0+}\mathcal{M}_1\ind_{|\Psi(\TT)-\alpha|<M}d\xi_{12} d\xi_{21}d\xi_{22} 
\\&
\les
\sup_{\xi, \al}
\int \frac{ \jap{\xi}^{-s} }{ \jb{\xi_{12}}^{4s+4-} }
\bigg(
\int \ind_{|\xi_{12}^3 P( p_{21}, p_{22}) - \al| < M } \xi_{12}^2 dp_{21} dp_{22}
\bigg)
d\xi_{12}
\\
&
\les
\sup_{\xi, \al}
\int \frac{1}{\jb{\xi_{12}}^{5s+2-}}
\bigg(
\int_{|q_{21}| ,  |q_{22}| \les 1}
\ind_{|\xi_{12}^3 (P(c) + \ld_{21} q_{12} + \ld_{22} q_{22}) - \al| < M  }
d q_{21} dq_{22}
\bigg)
d\xi_{12}
\\
&
\les M^{1-} \int \frac{1}{\jb{\xi_{12}}^{5s+5-}} d \xi_{12}
\les M^{1-},
\end{align}
for $s>-\frac45$,
 proving \eqref{eq:frequadd}.

\smallskip

Combining the estimates for the terms of Types $\I$-$\IV$, we find \eqref{eq:Linftybourg1} with a choice of $b$ and $\wt{b}$ independently of the tree\footnote{Observe that, in the application of Lemma \ref{lem:reduction_fre_kge3}, $\ell\le 4$ and $\beta$ depends only on $\ell$ and not on the size of the tree. See Remark~\ref{REM:b} for more details.} $\TT$,  completing the proof of Proposition~\ref{prop:bourgLinfty}. ~\qed

\bigskip

\begin{remark}\rm
\label{REM:b}
    When applying Lemma \ref{lem:reduction_fre_kge3} to an integral term parametrized by a tree $\TT\in \Tr_j$, $j\in\N$, if one had chosen $\ell=k = j+1$, then the parameter $b$ would depend on $j$ and converge to $1$  as $j\to \infty$. As such, no uniform choice could be made for all terms in the infinite normal form equation \eqref{z1} in order to ensure that $b<1$ (as required by Lemma~\ref{lem:linear_estimates}).
\end{remark}

\begin{remark} The estimates in Proposition \ref{prop:bourgLp} rely solely on the algebraic cancellations induced by the gauge transformation and thus can be extended to other dispersive models. As we will see in Section \ref{sec:cancel_kdv}, in the specific case of \eqref{kdv}, one may first recombine type $\III$ and $\IV$ terms before proceeding with the multilinear estimate (Lemma \ref{LEM:cancel}).

\end{remark}

\subsection{Proof of Proposition~\ref{prop:bourgL2} ($L^2$-estimates)}
\label{SUB:L2}

In this subsection, we present the proof of Proposition~\ref{prop:bourgL2}
on the $X^{s, b}_{2}$-estimates.
Since the frequency-restricted estimates needed to show \eqref{eq:bourg_bi_l2} differ from those in $L^\infty$-spaces, we present full details below.

As in Subsection~\ref{SUB:Linfty},
the estimate \eqref{eq:bourg_bi_l2} is a consequence of the inhomogeneous linear estimate \eqref{lin2a_gh} in Lemma~\ref{lem:linear_estimates_qfinite} together with Lemma~\ref{lem:fre_bi} and Lemma~\ref{lem:fre_tri}, reducing it to showing suitable frequency-restricted  estimates.
Given $j\in \N$ and $\TT\in \Tr_j$,
based on Lemma~\ref{LEM:wint}, it suffices to consider $\TT$ and contributions $\Ncal_
\l(\TT)$, $\l=1, \ldots, 4$, as in \eqref{z1} for $\TT$ of Type $\I$-$\IV$ (see also \eqref{Ncal}). 
For simplicity, we present the proof only in the case~$s<0$.

\smallskip

\noi
$\bul$ \textbf{Type $\I$.}
Here, we have $\TT = \<1>$, and
from \eqref{wint0c} and Lemma~\ref{lem:fre_bi}, we consider the multiplier
\begin{align}
\label{L21}
 \frac{\jb{\xi}^{2s} |\xi|^2 }{\jb{\xi_1}^{2s}\jb{\xi_2}^{2s}}
    \ind_A \ind_{|\Psi(\TT) - \al| < M } ,
\end{align}
 for $A$ as in \eqref{A}. We assume $|\xi_1| \ge |\xi_2|$ by symmetry, and thus $|\xi_2|\le 1$.
Fixing $\xi$ and changing variables $\xi_1 \mapsto \phi = 3 \xi\xi_1(\xi-\xi_1)$, since $|\partial_{\xi_1}\phi| \sim |\xi|^2$, we get
\begin{align}\label{eq:estL21}
    \sup_\xi \int_{\xi=\xi_1+\xi_2} \eqref{L21} d\xi_1
    \les \sup_\xi \int |\xi|^2 \ind_{|\phi - \al| < M } \frac{d\phi}{|\xi|^2 } \les M ,
\end{align}
thus condition (iii) in Lemma~\ref{lem:fre_bi} holds and the estimate follows, ragardless of $s$.

\medskip

For the remaining terms, the multilinear estimate \eqref{eq:bourg_bi_l2} will follow from Lemmas \ref{lem:fre_tri}-\ref{lem:eta1}.

\medskip

\noi
$\bul$ \textbf{Type $\II$.}
From \eqref{wint0c}, we only consider $\TT\in \Tr_j$ with $j\ge 2$ and $\TT^{0,\ff} = \{\star\}$. 
From \eqref{L22}, we can control the multiplier as follows
\begin{align}\label{L22_1}
\jb{\xi}^s |\mul_2(\TT)|  \prod_{\bul \in \TT^\infty} \frac{1}{\jb{\xi_\bul}^s}
&
\sim 
\frac{ \jb{\xi}^{s} \ind_{A_\star \cap A_{\pb(\star)}^c} }{  \jb{\xi_{\star1}}^s  \jb{\xi_{\star2}}^s  \jb{\xi_{\sbf(\star)}}^{s+1}}
\times
\prod_{\bul\in \TT^\infty_\star} \frac{1}{\jb{\xi_\bul}^{s+1} \jb{\xi_{\sbf(\bul)}} }
\\
&
=:
\M_1^\frac12(\vec\xi) \times \M_2^\frac12(\vec\xi)
,
\end{align}
where $\TT^\infty_\star = \TT^\infty \setminus \{\sbf(\star), \star1, \star2\}$, with $\sbf(\star)$ denoting the sibling node of $\star$. Without loss of generality, we assume that $|\xi_{\star1}| \ge |\xi_{\star2}| $ and thus $|\xi_{\star2}| \le 1$ from the restriction to $A_{\star}$. We consider only the case $|\xi_{\star1}|\ge |\xi_{\sbf(\star)}|$, as this is the worst-case scenario.

First, for fixed $\xi_{\star1}, \xi_{\star2}, \xi_{\sbf(\star)}$,  
we have
\begin{align}
\sup_{\xi_{\star1}, \xi_{\star2}, \xi_{\sbf(\star)} , \al} \int_{\G(\TT)} \M_2(\vec\xi) 
\ind_{|\Psi(\TT)-\al| < M}
d\xi_{\bul\in\TT^\infty_\star} 
&
\\\les 
\int
\prod_{\bul\in \TT^\infty_\star} \frac{1}{\jb{\xi_\bul}^{2(s+1)} \jb{\xi_{\sbf(\bul)}}^2 } d\xi_{\bul\in\TT^\infty_\star} 
&\les C^j\label{eq:int_M2}
\end{align}
for $s>-1$. For the remaining integral, we fix $\xi, \xi_{\bul \in \TT^\infty_\star}$, and consider some cases. If $|\xi_{\star1}| \les 1$, then the estimate for the $\M_1$ term in \eqref{L2hypo} holds trivially for $s\ge -1$, thus assume that $|\xi_{\star1}| \gg 1$. Then, we have $|\xi_{\star1}| \not\simeq |\xi_{\star2}|$ and for $\xi, \xi_{\bul\in\TT^\infty_\star}, \xi_{\sbf(\star)}$ fixed, we consider the change of variables $\xi_{\star1} \mapsto \phi = \Psi(\TT)$, which satisfies $|\partial_{\xi_{\star1}} \phi| \sim | \xi_{\star1}^2 - \xi_{\star2}^2| \sim |\xi_{\star1}|^2$:
\begin{align*}
    &\sup_{\xi, \xi_{\bul\in \TT^\infty_\star}, \alpha} \int_{\Gamma(\TT)} \jap{\xi_{\star1}}^{0^+}\M_1(\vec\xi) \mathbbm{1}_{|\Psi(\TT)
    -\alpha|<M}d\xi_{\star1} d\xi_{\sbf(\star)}
\\\lesssim \ &
\sup_{\xi, \xi_{\bul\in \TT^\infty_\star}, \alpha} \int \jap{\xi_{\sbf(\star)}}^{-2s-3} \mathbbm{1}_{|\phi-\alpha|<M}d\phi d\xi_{\sbf(\star)}
\lesssim M, 
\end{align*}
again for $s>-1$. 
This completes the proof for Type $\II$ terms.

\smallskip

\noi
$\bul$ \textbf{Type $\III$ for $j=2$.}
Here, up to permutations, $\TT=\<21>$. 
From \eqref{wint0c}, the relevant multiplier is controlled as in \eqref{mul3}. As we will no longer use any particular relation between $\xi, \xi_{11}, \xi_{12}, \xi_{2}$ (apart from the global convolution relation), we consider only the worst-case scenario $|\xi_2|\ge |\xi|\ge |\xi_{11}|\ge |\xi_{12}|$.

\smallskip
\noi
\underline{\textbf{Case 1.} $|\xi_2| \not\simeq |\xi_{11}|$, $|\xi| \not\simeq |\xi_{12}|$.}
Here, we consider the change of variables $\xi_{12} \mapsto \phi = \Psi(\TT)$ and $\xi_{11}\mapsto \phi = \Psi(\TT)$, respectively, noting that $|\partial_{\xi_{12}} \phi| \sim |\partial_{\xi_{11}} \phi| \sim |\xi_2|^2$: 
\begin{align}
\sup_{\xi_{12}, \xi, \al} \int_{\G(\TT)}
|\xi_2|^2 \ind_{|\Psi(\TT)-\al|<M}
d\xi_{11}
&
\les  
\sup_\al \int \ind_{|\phi - \al|<M} d\phi 
\les 
M
, 
\\
\sup_{\xi_{11}, \xi_{2}, \al} \int_{\G(\TT)}
\frac{\jb{\xi}^{2(s+1)}
\ind_{|\Psi(\TT)-\al|<M}
}{ \jb{\xi_{11}}^{2(s+1)} \jb{\xi_{12}}^{2(s+1)} \jb{\xi_{2}}^{2(s+1) -} } 
d \xi_{12}
&
\les 
\sup_{\xi_2, \al} 
\int \jb{\xi_{2}}^{-2+} \ind_{|\phi-\al|<M} d\phi \les M 
, 
\end{align}
given that $s>-1$.

\smallskip

\noi
 \underline{\textbf{Case 2:} $|\xi_2| \not\simeq |\xi_{11}|$, $|\xi| \simeq|\xi_{12}|$.}
Then all frequencies are comparable and we proceed as in Case 1, performing the change of variables only in the first integral:
\begin{align}
\sup_{\xi_{12}, \xi, \al} \int_{\G(\TT)}
|\xi_2|^{2} \ind_{|\Psi(\TT)-\al|<M}
d\xi_{11}
&
\les  
M
, 
\\
\sup_{\xi_{11}, \xi_{2}, \al} \int_{\G(\TT)}
\frac{\jb{\xi}^{2(s+1)}
\ind_{|\Psi(\TT)-\al|<M}
}{ \jb{\xi_{11}}^{2(s+1)} \jb{\xi_{12}}^{2(s+1)} \jb{\xi_{2}}^{2(s+1)-} } 
d \xi_{12}
&
\les 
\sup_{\xi_2, \al} 
\int \jap{\xi_{12}}^{-4(s+1)+}d\xi_{12}\lesssim 1
, 
\end{align}
for $s> -\frac34$. 

\smallskip

\noi
\underline{\textbf{Case 3:} $|\xi_2|\simeq |\xi_{11}|$, $|\xi_2| \not\simeq |\xi_{12}|$.}
Then, 
$|\xi_2| \sim |\xi|\sim|\xi_{11}|\sim|\xi_{12}|$ and we can proceed as in Case 2, exchanging the roles of $\xi_{11}$ and $\xi_{12}$.

\smallskip

\noi
\underline{\textbf{Case 4:} $|\xi_2| \simeq |\xi| \simeq |\xi_{11}| \simeq |\xi_{12}|$.}
All frequencies are comparable,  at least two among $\xi_{11}, \xi_{12}, \xi_2$ have the same sign, and all frequencies are away from zero.
Without loss of generality, suppose that $\xi_{12} \xi_2 >0$. Then, set
\begin{align}
    p_\bul = \frac{\xi_\bul}{\xi}, \quad P(p_{12}) = - 1 + p_{11}^3 + p_{12}^3 + p_2^3.
\end{align}
From the frequency assumptions, we have
\begin{align}
    \big| 1 - |p_{12}| \big| ,
    \big| 1 - |p_{11}| \big| ,
    \big| 1 - |p_{2}| \big| \ll 1,
\end{align}
and for fixed $\xi_{11}, \xi$, it holds that $\partial_{p_{12}} P(p_{12}) = 3 (p_{12} - p_2) (p_{12} + p_2) $ and $\partial^2_{p_{12}}P(p_{12}) = 6 (p_{12} + p_2) \ne 0$ since $p_{12}p_2 >0$, thus we are considering $p_{12}$ near a non-degenerate critical point of $P$. Consequently, there exists a change of coordinates $p_{12} \mapsto q_{12}$ such that
\begin{align}
    P(p_{12}) = c \pm q_{12}^2,
\end{align}
for some $c\in\R$
and with $|q_{12}| \ll 1$. Therefore, we have
\begin{align}
\sup_{\xi_{11}, \xi} \int
\frac{\jb{\xi}^{s+} |\xi|  }{ \jb{\xi_{11}}^{s+1} \jb{\xi_{12}}^{s+1} \jb{\xi_{2}}^s}
    \ind_{|\Psi(\TT)- \al | < M } d\xi_{12}
    &
    \les \sup_{\xi_{11}, \xi}
    \jb{\xi}^{-1-2s+}
    \int \ind_{|\xi^2P(p_{12}) - \al | < M } |\xi| d p_{12}
    \\
    &
    \les \sup_{\xi_{11}, \xi}
    \jb{\xi}^{-2s+}
    \int \ind_{|\xi^3(c \pm p_{12}^2) - \al| < M} dq_{12}
    \\
    &
    \les \sup_\xi \jb{\xi}^{-2s - \frac32 + } M^\frac12 \les M^\frac12,
\end{align}
for $s>-\frac34$, using \eqref{eq:quadratica} in  Lemma~\ref{lem:quadraticas}. For fixed $\xi_{12}, \xi_2$, we can perform a similar change of variables and evaluate the other integral analogously.

\smallskip
\noi
$\bul$ \textbf{Type $\III$ for $j\ge 3$.}
From \eqref{wint0c}, we assume that $\TT \in \Tr_j$ with $j\ge3$ and $\TT^{0,\ff} = \{\star\}$. Then,  we have the following bound for the multiplier:
\begin{align}\label{mul3j}
    \jb{\xi}^s | \mul_3 (\TT)|
    \prod_{\bul \in \TT^\infty} \frac{1}{\jb{\xi_\bul}^s}
    &
    \le 
    \frac{\jb{\xi}^{s} \ind_{A_\star^c} }{ \jb{\xi_{\star1}}^{s+1} \jb{\xi_{\star2}}^{s+1} \jb{\xi_{\sbf(\star)}}^s \jb{\xi_{\sbf(\pb(\star))}}^{s+1}} 
	\prod_{\bul \in \TT^\infty_\star} \frac{1}{\jb{\xi_\bul}^{s+1} \jb{\xi_{\sbf(\bul)}} }
\end{align}
where $\TT^\infty_\star : = \TT^\infty \setminus \{\star1, \star2, \sbf(\star), \sbf(\pb(\star))\}$, recalling that $\pb(\bul)$ and $\sbf(\bul)$ denote the parent and sibling nodes of $\bul$, respectively.
Without loss of generality, we take $-\frac78<s<-\frac12$ and, since we will rely solely on the global convolution relation, we can focus on the worst-case scenario $|\xi_{\sbf(\star)}|\ge |\xi_{\sbf(\pb(\star))}|\ge |\xi_{\star1}|\ge |\xi_{\star2}|\ge |\xi|$.

\smallskip
\noi
\underline{\textbf{Case 1:} $|\xi_{\star2}|\ge |\xi_{\star1}|$}

If $|\xi_{\sbf(\star)}|\not\simeq |\xi_{\star2}|$, we change variables $\xi_{\sbf(\star)} \mapsto \phi = \Psi(\TT)$:
\begin{align*}
& \sup_{\xi_{\sbf(\pb(\star))}, \xi_{\star1}, \al} \int_{\G(\TT)} |\xi_{\sbf(\star)}|^2 \jb{\xi}^{2s}  \ind_{|\Psi(\TT) - \al|<M} \prod_{\bul \in \TT^\infty_{\star}} \frac{1}{\jb{\xi_\bul}^{2(s+1)} 
\jb{\xi_{\sbf(\bul)}}^2 } d\xi_{\bul\in\TT^\infty_\star} d\xi_{\sbf(\star)} d\xi 
\\
& 
\les \sup_\al \int \ind_{|\phi-\al|<M} \jb{\xi}^{2s} d\phi d\xi \les M 
, 
\\
&
\sup_{\xi, \xi_{\bul\in\TT^\infty_\star}, \xi_{\sbf(\star)}, \xi_{\star2}} \int_{\G({\TT})}  
\frac{1 }{ \jb{\xi_{\star1}}^{2(s+1)} \jb{\xi_{\star2}}^{2(s+1)} \jb{\xi_{\sbf(\star)}}^{2(s+1+)} \jb{\xi_{\sbf(\pb(\star))}}^{2(s+1)}} 
d\xi_{\star1}
\\
&
\les \int \frac{1}{\jb{\xi_{\star1}}^{8(s+1)-}} d\xi_{\star1} \les 1
,
\end{align*}
for $s> - \frac78$.
If instead $|\xi_{\sbf(\star)}| \simeq |\xi_{\star2}|$, then $|\xi_{\sbf(\star)}| \sim |\xi_{\sbf(\pb(\star))} |\sim |\xi_{\star1}| \sim |\xi_{\star2}|$, and we consider two cases. If we have $|\xi_{\bul}| \not\simeq |\xi_{\circ}|$ for any two distinct $\bul, \circ\in\{\sbf(\star), \sbf(\pb(\star)), \star1 , \star2\}$, then we can proceed as above, with the change of variables $\xi_\bul \mapsto \phi$.

It remains to consider $|\xi_{\star1}| \simeq |\xi_{\star2}| \simeq |\xi_{\sbf(\star)}| \simeq |\xi_{\sbf(\pb(\star))}| $.
Fixing $\xi, \xi_{\sbf(\pb(\star))}, \xi_{\bul\in\TT^\infty_\star}$, let $p_\bul = \frac{\xi_\bul}{\xi_{\sbf(\pb(\star))}}$ and $P(p_{\star1}, p_{\star2}) = p^3_{\star1} + p^3_{\star2} + p_{\sbf(\star)}^3$. The restrictions on the frequencies imply that $(p_{\star1},p_{\star2})$ are near one of the critical points of $P$,
$$(c_0,c_0), \quad (\pm c_0, \mp c_0), \quad (\tfrac{ c_0}{3}, \tfrac{ c_0}{3}), \quad \text{with} \quad c_0 = p - \sum_{\bul \in \TT^\infty_\star } p_\bul - p_{\sbf(\pb(\star))}.  $$
Since
\begin{align}
    | D^2 P(p_{\star1}, p_{\star2}) |
    &
    = 36 |p_{\star1} p_{\star2} + p_{\star1} p_{\sbf(\star)} + p_{\star2} p_{\sbf(\star)} |
    \sim 1,
\end{align}
by Morse's lemma (Lemma~\ref{LEM:morse}), there exists a change of variables $(p_{\star1}, p_{\star2}) \mapsto (q_{\star1}, q_{\star2})$ with $|q_{\star1}|, |q_{\star2}| \ll 1$ and
\begin{align}
    P(p_{\star1}, p_{\star2}) = c + \ld_{\star1} q_{\star1}^2 + \ld_{\star2} q_{\star2}^2,
\end{align}
for some $\ld_{\star1}, \ld_{\star2} \in \{\pm1 \}$,
{where $c$ is the value of $P$ at each of the critical points above.}
Consequently, setting $\wt\al = \al + \xi^3 - \sum_{\bul \in \TT^\infty_\star} \xi_\bul^3$, we have
\begin{align}
    &
    \sup_{\xi,  \xi_{\sbf(\pb(\star))}, \xi_{\bul\in\TT^\infty_\star}}
     \jb{\xi_{\sbf(\pb(\star))}}^{1-}
     \int
    \ind_{|\Psi(\TT) - \al| < M } d \xi_{\star1} d\xi_{\star2}
    \\
    &
    \quad
    \les
    \sup_{\xi, \xi_{\sbf(\pb(\star))}, \xi_{\bul\in\TT^\infty_\star}}
     \jb{\xi_{\sbf(\pb(\star))}}^{3-}
     \int
    \ind_{|\xi_{\sbf(\pb(\star))}^3 P(p_{\star1}, p_{ \star2}) - \wt\al| < M } d p_{\star1} dp_{\star2}
    \\
    &
    \quad
    \les
    \sup_{\xi, \xi_{\sbf(\pb(\star))}, \xi_{\bul\in\TT^\infty_\star}}
     \jb{\xi_{\sbf(\pb(\star))}}^{3-}
     \int
    \ind_{|\xi_{\sbf(\pb(\star))}^3 ( c + \ld_{\star1} q_{\star1}^2 + \ld_{\star2} q_{\star2}^2) - \wt\al| < M } d q_{\star1} dq_{\star2} \les
    M^{1-} ,
    \\
    &
    \sup_{\xi_{\star1}, \xi_{\star2}, {\xi_{\sbf(\star)} }}
    \int
    \jb{\xi}^{2s} \jb{\xi_{\star1}}^{-8s - 7 +} \ind_{|\Psi(\TT) - \al| < M } d \xi \prod_{\bul\in \TT^\infty_\star} \frac{1}{\jb{\xi_\bul}^{2s+4}} d\xi_\bul
    \\
    &
    \quad
    \les
    \sup_{\xi_{\star1}, \xi_{\star2}, {\xi_{\sbf(\star)} }}
    \int
    \jb{\xi}^{-6s - 7 + }  d \xi \prod_{\bul\in \TT^\infty_\star} \frac{1}{\jb{\xi_\bul}^{2s+4}} d\xi_\bul \les 1,
\end{align}
where we used \eqref{eq:posdef} in Lemma~\ref{lem:quadraticas} on the first estimate and $s>-\frac78$ in the second estimate.

\smallskip
\noi
\underline{\textbf{Case 2:} $   |\xi_{\star1}| \gg |\xi_{\star2}| $}

As $|\xi_{\sbf(\star)}|\gg |\xi|$ and $|\xi_{\pb(\star)}|\gg |\xi_{\star2}|$, we may perform a change of variables on both sides of the interpolation:
\begin{align*}
& \sup_{\xi_{\sbf(\pb(\star))}, \xi_{\star1}, \xi_{\star2}, \al} \int_{\G(\TT)} |\xi_{\sbf(\star)}|^2 \jb{\xi}^{2s}  \ind_{|\Psi(\TT) - \al|<M} \prod_{\bul \in \TT^\infty_{\star}} \frac{1}{\jb{\xi_\bul}^{2(s+1)} 
\jb{\xi_{\sbf(\bul)}}^2 } d\xi_{\bul\in\TT^\infty_\star} d\xi_{\sbf(\star)} 
\\
& 
\les \sup_\al \int \ind_{|\phi-\al|<M}  d\phi  \les M 
, 
\\
&
\sup_{\xi, \xi_{\bul\in\TT^\infty_\star}, \xi_{\sbf(\star)},\alpha } \int_{\G({\TT})}  
\frac{\ind_{|\Psi(\TT) - \al|<M}  }{ \jb{\xi_{\star1}}^{2(s+1)} \jb{\xi_{\star2}}^{2(s+1)} \jb{\xi_{\sbf(\star)}}^{2(s+1+)} \jb{\xi_{\sbf(\pb(\star))}}^{2(s+1)}} 
d\xi_{\star1} d\xi_{\star2}
\\
&
\les \sup_\alpha \int \frac{1}{\jb{\xi_{\star1}}^{6s+8-}} \ind_{|\phi-\al|<M}  d\phi d\xi_{\star1} \les M
,
\end{align*}

\medskip

\noi
$\bul$ \textbf{Type $\IV$, $j=3$.}

Considering the multiplier in \eqref{wint0c}, for $\TT = \<33>$, we have 
\begin{align}\label{mul43}
\jb{\xi}^s \mul_4(\TT)  \prod_{\bul\in \TT^\infty} \frac{1}{\jb{\xi_\bul}^s} \les 
\frac{\jb{\xi}^{s+1}}{\jb{\xi_{11}}^{s+1} \jb{\xi_{12}}^{s+1} \jb{\xi_{21}}^{s+1} \jb{\xi_{22}}^{s+1} }
\end{align}
Since $s>-1$, we can assume the worst-case scenario $|\xi|\ge|\xi_{11}| \ge |\xi_{12}|\ge |\xi_{21}| \ge |\xi_{22}|$, which implies immediately that $|\xi|\not\simeq |\xi_{22}|$. As such, we may perform the change of variables $\xi \mapsto \phi = \Psi(\TT)$, since $|\partial_\xi \phi| \sim |\xi|^2$:
\begin{align}\label{L2-43-11}
\sup_{\xi_{11}, \xi_{21},\xi_{12} ,\al } \int_{\G(\TT) } |\xi|^2 \ind_{|\Psi(\TT) - \al| < M } d\xi 
&
\les M 
, 
\\
\sup_{\xi,   \xi_{22}} \int_{\G(\TT)} 
\frac{\jb{\xi}^{2(s+1)-2}}{\jb{\xi_{11}}^{2(s+1)} \jb{\xi_{21}}^{2(s+1)} \jb{\xi_{12}}^{2(s+1)} } d\xi_{11} d\xi_{12}
&
\les 
\int 
\frac{1}{\jb{\xi_{11}}^{1+}  \jb{\xi_{12}}^{4(s+1) + 1 -} } d\xi_{11} d\xi_{12}
\les 1
, 
\end{align}
for $s>-1$.

\medskip

\noi
$\bul$ \textbf{Type $\IV$, $j\ge 4$.}

From \eqref{wint0c}, with $\TT^\infty_\star = \TT^\infty \setminus\{ \star11 , \star12, \star21, \star22, \sbf(\star)\}$, we control the multiplier as follows
\begin{align}
\jb{\xi}^s |\mul_4(\TT)| \prod_{\bul\in\TT^\infty} \frac{1}{\jb{\xi_\bul}^s} \les \frac{\jb{\xi}^s}{  \jb{\xi_{\star11}}^{s+1} \jb{\xi_{\star12}}^{s+1} \jb{\xi_{\star21}}^{s+1} \jb{\xi_{\star22}}^{s+1}     \jb{\xi_{  \sbf(\star)}}^{s+1} }
\prod_{\bul\in\TT^\infty_\star } \frac{1}{\jb{\xi_\bul}^{s+1} \jb{\xi_{\sbf(\bul)}}}
.
\end{align}
We suppose, without loss in generality, that $s<-\frac12$. By splitting the multiplier, we first have 
\begin{align}
&
\sup_{\xi_{\star11}, \xi_{\star12}, \xi_{\star21},\xi_{\star22},\xi_{\sbf(\star)}}  \int_{\G(\TT) } \jb{\xi}^{2s} \prod_{\bul\in \TT^\infty_\star } \frac{1}{\jb{\xi_\bul}^{2(s+1)} \jb{\xi_{\sbf(\bul)}}^2} d\xi_{\bul\in\TT^\infty_\star} d\xi  
\les \int \jb{\xi}^{2s} d\xi \les 1, 
\end{align}
where we used $s>-1$. 

To estimate the other integral, by symmetry, 
we assume that $|\xi_{\star11}| \ge |\xi_{\star12}|\ge |\xi_{\star21}| \ge |\xi_{\star22}|\ge |\xi_{\sbf(\star)}|$  and consider different cases.

\smallskip
\noi
\underline{\textbf{Case 1:} $|\xi_{\star11}| \not\simeq |\xi_{\sbf(\star)}|$}

Changing variables $\xi_{\star11} \mapsto \phi = \Psi(\TT)$, with $|\partial_{\xi_{\star11}} \phi| \sim |\xi_{\star11}|^2$, we get
\begin{align}
&
\sup_{\xi_{\bul\in \TT^\infty_\star}, \xi, \al} \int_{\G(\TT)} 
\frac{ \ind_{|\Psi(\TT)-\al|<M} }{  \jb{\xi_{\star11}}^{2(s+1)-} \jb{\xi_{\star12}}^{2(s+1)} \jb{\xi_{\star21}}^{2(s+1)} \jb{\xi_{\star22}}^{2(s+1)}     \jb{\xi_{  \sbf(\star)}}^{2(s+1)} }
d\xi_{\star11}
d\xi_{\star12}d\xi_{\star 21} d\xi_{\star22}
\\
&
\les \sup_\al \int 
\frac{ \ind_{|\phi-\al|<M} }{  \jb{\xi_{\star12}}^{1+}  \jb{\xi_{\star21}}^{1+}  \jb{\xi_{\star 22}}^{8(s+1)}   }
d\phi
d\xi_{\star12}d\xi_{\star 21} d\xi_{\star22} \lesssim M
\end{align}
for $s>-\frac78$.

\smallskip
\noi
\underline{\textbf{Case 2:} $|\xi_{\star11}| \simeq  |\xi_{\sbf(\star)}|$}

Here, using Morse's Lemma (Lemma~\ref{LEM:morse}) and Lemma~\ref{lem:quadraticas} as before, we obtain
\begin{align}
&
\sup_{\xi_{\bul\in \TT^\infty_\star}, \xi, \al} \int_{\G(\TT)} 
\frac{ \ind_{|\Psi(\TT)-\al|<M} }{  \jb{\xi_{\star11}}^{2(s+1)-} \jb{\xi_{\star12}}^{2(s+1)} \jb{\xi_{\star21}}^{2(s+1)} \jb{\xi_{\star22}}^{2(s+1)}     \jb{\xi_{  \sbf(\star)}}^{2(s+1)} }
d\xi_{\star11}
d\xi_{\star12}d\xi_{\star 21} d\xi_{\star22}
\\
&
\les 
\sup_{ \al } \int \jb{\xi_{\star22}}^{-10(s+1) + 2 + } \ind_{| \ld_{11} q_{11}^2 + \ld_{12} q_{12}^2 - \al | < \frac{M}{|\xi_{\star22}|^3}} d q_{11} dq_{12} d\xi_{\star21}d\xi_{\star22}
\\
&
\les  M^{1-} \int \jb{\xi_{\star22}}^{-10(s+1) -1 + } d\xi_{\star21}d\xi_{\star22} \les M^{1-}, 
\end{align}
for $s>-\frac9{10}$. 

This completes the proof of Proposition~\ref{prop:bourgL2}. \qed

\subsection{Proof of Proposition \ref{prop:bourgLp} ($\FL^p$-estimates, $2<p<\infty$)}
\label{SUBSEC:bourgainLp}

The proof of \eqref{eq:bourg_bi_lp} follows from an interpolation argument between the $L^2$ and $\FL^\infty$ estimates \eqref{eq:bourg_bi_l2}-\eqref{eq:Linftybourg1}. However, since the latter involves the $X_{\infty,1}^{s,\wt{b}}$ norm, the argument is slightly more delicate, presented below.

Fix $j\in\N$, a tree $\TT\in \Tr_j$,  $\ell=1,\dots,4$, and $\Ncal_\l(\TT)$ as in \eqref{Ncal}. Define $s_2=s_2(\TT)$ and $s_\infty=s_\infty(\TT)$ as in \eqref{eq:defis2}-\eqref{eq:defisinfty}. First, we observe that, since the proof of \eqref{eq:bourg_bi_l2} results from the application of the inhomogeneous estimate in $X_2^{s,b}$-spaces together with either Lemma \ref{lem:fre_bi} or Lemma \ref{lem:fre_tri}, it actually implies that
\begin{equation}\label{eq:L2_interp}
		\left\| \Ncal_\l(\TT; z_1,\dots, z_{j+1})\right\|_{X^{s,b_2'}_{2}(0,T)}  \le  C^j\prod_{k=1}^{j+1}\|z_k\|_{X^{s,b_2}_{2}(0,T)}
\end{equation}
for $s>s_2$, $b_2=\frac12 +$ and $b_2'=(b_2-1)+$. Similarly, the proof of \eqref{eq:Linftybourg1} also gives
\begin{align*}
    \left\| \Ncal_\l(\TT; z_1,\dots, z_{j+1})\right\|_{X^{s,b'}_{\infty, \infty}(0,T)} &\le  C^j\prod_{k=1}^{j+1}\|z_k\|_{X^{s,b}_{\infty, \infty} (0,T) \cap X^{s,\wt{b}}_{\infty,1}(0,T)},\quad s>s_\infty.
\end{align*}
Here, $(b, b', \wt{b})$ must satisfy
\begin{equation}\label{eq:cond_lp}
    b-1<b'<0
    \quad
    \text{and}
    \quad
    0<\wt{b}<1-\beta + b
\end{equation}
(as per \eqref{bs} and \eqref{freb}),
where $\beta<1$ is related to the power of the parameter $M$ in the frequency-restricted estimates (which, as mentioned in Remark \ref{REM:b}, can be made uniform in $\TT$, $j$, and $\l$).

For $p>2$ fixed, we choose $(b, b'_\infty, \wt{b}_\infty)$ satisfying \eqref{eq:cond_lp} and
\begin{equation}
    1+\wt{b}_\infty <\tfrac{2}{p-2}(1-b_2+b_2') + 1 + b'_\infty.
\end{equation}
Then, given $b_\infty$ satisfying
\begin{equation}
    1+\wt{b}_\infty<b_\infty <\tfrac{2}{p-2}(1-b_2+b_2') + 1 + b'_\infty,
\end{equation}
since $b_\infty>\wt{b}_\infty+1>b$,
\begin{align}
    \left\| \Ncal_\l(\TT; z_1,\dots, z_{j+1})\right\|_{X^{s,b_\infty'}_{\infty, \infty}(0,T)} &\le  C^j\prod_{k=1}^{j+1}\|z_k\|_{X^{s,b}_{\infty, \infty}\cap X^{s,\wt{b}_\infty}_{\infty,1}(0,T)}\\&\le  C^j\prod_{k=1}^{j+1}\|z_k\|_{X^{s,b_\infty}_{\infty, \infty}(0,T)}.\label{eq:linfty_interp}
\end{align}
Interpolating \eqref{eq:linfty_interp} and \eqref{eq:L2_interp}, we find
\begin{equation}\label{eq:Lp_interp}
		\left\| \Ncal_\l(\TT; z_1,\dots, z_{j+1})\right\|_{X^{s,b_p'}_{p}}  \le  C^j\prod_{k=1}^{j+1}\|z_k\|_{X^{s,b_p}_{p}},
\end{equation}
for
$$
b_p=\tfrac{2}{p}b_2+\tfrac{p-2}{p}b_\infty,
\quad
b'_p=\tfrac{2}{p}b'_2+\tfrac{p-2}{p}b'_\infty,
\quad
\text{and}
\quad
s> \tfrac2p s_2 + \tfrac{p-2}{p} s_\infty
.
$$
It is now a simple computation to check that
$$
b_p>\tfrac1{p'},
\quad b_p-1<b_p',
\quad
-\tfrac{1}{p}<b_p'<0,
\quad
\text{and}
\quad
s> s_p(\TT).
$$
In particular, applying Lemma \ref{lem:linear_estimates_qfinite}, we obtain
	\begin{align}
		\bigg\|\int_0^t e^{-(t-t')\dx^3} \Ncal_\l(\TT; z_1,\dots, z_{j+1})(t') dt'  \bigg\|_{X^{s,b_p}_{p}(0,T)}
        &	\lesssim T^{1+b_p'-b_p}\big\|\Ncal_\l(\TT; z_1,\dots, z_{j+1})  \big\|_{X^{s,b_p'}_{p} (0,T)}
        \\
        &\le T^{1+b_p'-b_p} C^j\prod_{k=1}^{j+1}\|z_k\|_{X^{s,b_p}_{p}(0,T)},
	\end{align}
as intended.
\qed

\section{Additional cancellation for KdV}
\label{sec:cancel_kdv}

Due to the \eqref{kdv} structure, we observe an additional cancellation for the nonlinear terms $\Ncal_\l(\TT;z)$ appearing in \eqref{z1} when $\l=3,4$, namely for terms of Types $\III$-$\IV$ as classified in Lemma~\ref{LEM:wint}. 

\begin{lemma}\label{LEM:cancel}
Let $z\in H^\infty(\R)$. Then, the following hold:

\noi{\rm(i)} For any $j\ge 3$ and $\TT \in \Tr_j$, we can write
\begin{align}
\Ncal_3(\TT;z)
&
= 
\Ncal_3^{(1)}(\TT;z)
+
\Ncal_3^{(2)}(\TT;z)
+
\Ncal_3^{(3)}(\TT;z), 
\label{Ncal3}
\end{align}
where $\Ncal_3^{(k)}(\TT;z)$ are of the form \eqref{Ncal} with multiplier $i \mul_3^{(k)}(\TT)$, $k=1,2,3$, given by 
\begin{align}
\mul_3^{(1)}(\TT)
& = 
-  \frac{ \xi_{\pb(\pb(\star))}   }{36   \xi_{\sbf(\pb(\star))} \xi_{ \sbf(\star)} \xi_{\star1} \xi_{\star2 }} 
\ind_{|\xi_{\sbf(\pb(\star))} |, |\xi_{\sbf(\star)}| , |\xi_{\star1}|, |\xi_{\star2}| > 1 } \prod_{\bul\in \TT^0_\star} K_\bul \ind_{A_\bul^c}
, 
\\
\mul_3^{(2)}(\TT)
&
=
 \frac{ \xi_{\pb(\star)}   }{27 \xi_{\sbf(\pb(\star))} \xi_{ \sbf(\star)} \xi_{\star1} \xi_{\star2}} 
\ind_{ \substack{
|\xi_{\sbf(\pb(\star))}| , |\xi_{\sbf(\star)}|, 
|\xi_{\star1}|, |\xi_{\star2}| > 1 \\
1 \ge  |\xi_{\pb(\star)}| 
}}
\prod_{\bul\in \TT^0_\star} K_\bul \ind_{A_\bul^c}
, 
\\
\mul_3^{(3)}(\TT)
&
=
-  \frac{ 1  }{9  \xi_{\sbf(\pb(\star))} \xi_{\star1} \xi_{\star2}} 
\ind_{\substack{
|\xi_{\sbf(\pb(\star))}|, |\xi_{\star1}| , |\xi_{\star2}|, 
|\xi_{\pb(\star)}| > 1  \\
1 \ge |\xi_{\sbf(\star)}| 
}}
\prod_{\bul\in \TT^0_\star} K_\bul \ind_{A_\bul^c}
, 
\label{mul3-new}
\end{align}
where $\TT^{0,\ff}=\{\star\}$, $\TT^0_\star : = \TT^0\setminus\{ \pb(\pb(\star)) , \pb(\star) , \star\}$ and $K_\bul$ and $A_\bul^c$ are as in \eqref{Kmul} and \eqref{defA}. 

\smallskip 

\noi{\rm(ii)}
The following cancellation holds for any $j\ge3$:
\begin{align}\label{cancel}
\sum_{\TT\in \Tr_j} \big( \Ncal_3^{(1)}(\TT;z)  +  \Ncal_4(\TT; z) \big) =0 , 
\end{align}
where $\Ncal_4$ is as in \eqref{Ncal} with Lemma~\ref{LEM:wint}. 
\end{lemma}

From the cancellation observed in \eqref{cancel}, we can rewrite the gauged normal form equation in \eqref{z1} as 
\begin{align}
\dt z + \dx^3 z 
&
= 
\Ncal_1 (\<1>; z) 
+
\sum_{\TT\in\Tr_2} \big( \Ncal_2(\TT;z) + \Ncal_3(\TT;z) \big)
\\
&
\quad 
+
\sum_{j = 3  }^\infty \sum_{\TT\in\Tr_j} \big( \Ncal_2(\TT;z) + \Ncal^{(2)}_3(\TT;z) +  \Ncal^{(3)}_3(\TT;z)  \big) , 
\label{z2}
\end{align}
namely, we no longer see contributions of Type $\IV$ and the Type $\III$ contributions for $j\ge3$ are only of the form $\Ncal^{(k)}_3$ in \eqref{Ncal3} for $k=2,3$. We claim that, having exploited this ultimate cancellation, the multilinear estimates akin to \eqref{eq:bourg_bi_lp} now hold for even lower regularities, namely $s>-1+\frac1{2p}$.
To state this properly,
analogously to \eqref{Ncal}, define the multilinear version of the nonlinear terms $\Ncal_3^{(k)}, \ k=2,3,$ as
\begin{align}
\label{Ncal2}
    \Ft_x\big[ \Ncal_3^{(k)} (\TT; z_1, \ldots, z_{j+1} ) \big](\xi)
    &
    =
    i \int_{\G(\TT)} \mul_3^{(k)}(\TT)(\vec\xi) \prod_{k=1}^{j+1} \ft{z}_k(t, \xi_{\bul_k}) d\vec\xi,\quad k=2,3.
\end{align}

\begin{proposition}[Improved $\FL^p$-estimates]
\label{PROP:Finfty-new}
Let $0< T \le 1$, $j\ge 3$, $\TT\in \Tr_j$ with $|\TT^{0,\ff}| = 1$, and $\Ncal_3^{(k)}(\TT)$ as in \eqref{Ncal2}, $k=2,3$. 
Then, given $2 \le p \le \infty$, $s>-1 + \frac{1}{2p}$, and $b,\wt b$ as in Propositions~\ref{prop:bourgLinfty}-\ref{prop:bourgLp}, we have 
\begin{align}
 & \bigg\| \int_0^t e^{-(t-t')\dx^3}\Ncal_3^{(k)}(\TT; z_1,\dots, z_{j+1})(t')dt' 
 \bigg\|_{X^{s,b}_{\infty, \infty} (0,T)\cap X^{s,\wt b}_{\infty, 1}(0,T)} 
 \\
        & \hspace{2cm} \le T^{0+} C^j\prod_{k=1}^{j+1}\|z_k\|_{X^{s,b}_{\infty, \infty} (0,T) \cap X^{s,\wt{b}}_{\infty,1}(0,T)}, && \text{for } p=\infty, 
        \label{eq:Linftybourg1_improv}
\\
 & \bigg\|\int_0^t e^{-(t-t')\dx^3} \Ncal_3^{(k)}(\TT; z_1,\dots, z_{j+1})(t') dt'  \bigg\|_{X^{s,b}_{p}(0,T)}
 \\
        & \hspace{2cm} \le T^{0+} C^j\prod_{k=1}^{j+1}\|z_k\|_{X^{s,b}_{p}(0,T)}, 
        && \text{for }  2 \le p < \infty, 
        \label{eq:bourg_bi_lp_improv}
\end{align}
for some constant $C>0$ independent of $j$ and $\TT$.

\end{proposition}

The remainder of this section is devoted to the proofs of Lemma~\ref{LEM:cancel} and  Proposition~\ref{PROP:Finfty-new}.

\begin{proof}[Proof of Lemma~\ref{LEM:cancel}]

We first consider the quartic terms, namely  $\TT\in \Tr_3$, $\Ncal_3(\TT; z, z, z, z)$ and $\Ncal_4(\TT;z, z, z, z)$ (defined in \eqref{Ncal}). 
For $\TT= \<31>$, from \eqref{wint0c}, the corresponding multiplier is given by
\begin{align}
\mul_3(\TT) & = -  \frac{ 1   }{9 \xi_2 \xi_{111} \xi_{112}} \ind_{|\xi_1|, |\xi_2|, |\xi_{111}|, |\xi_{112}| > 1}
\\
&
= 
-  \frac{ \xi_{12}   }{9 \xi_2 \xi_{ 12} \xi_{111} \xi_{112}} 
\ind_{|\xi_{111}|, |\xi_{112}|, |\xi_2| > 1 }
\Big( 
\ind_{   |\xi_{12}| > 1 } - \ind_{ |\xi_{12}|> 1 \ge  |\xi_1|}  + \ind_{|\xi_1| > 1 \ge |\xi_{12}|}  
\Big) 
\\
&
=: 
\wt \mul_3^{(1)}(\TT)
+
\wt \mul_3^{(2)}(\TT)
+
\wt \mul_3^{(3)}(\TT)
, 
\end{align}
and we denote by $\Ncal_3^{(k)}(\TT;z)$, $k=1,2,3$, the portions of $\Ncal_3(\TT;z)$ with multiplier $i \wt \mul^{(k)}_3(\TT)$ defined above. 
From \eqref{Ncal}, we note that
\begin{align}
\F_x[\Ncal_3^{(1)}(\TT;z)] (\xi)
&
= 
-\frac{i}{9} 
\int_{\substack{\xi = \xi_2 + \xi_{12} + \xi_{111} + \xi_{112} \\ |\xi_2|, |\xi_{12}|, |\xi_{111}|, |\xi_{112}| >1 }} \frac{ \xi_{12}   }{ \xi_2 \xi_{ 12} \xi_{111} \xi_{112}} \prod_{k=1}^4 \ft z (t, \xi_{\bul_k}) d\vec \xi
\\
&
=
- \frac{ i }{4\times 9} \int_{\substack{\xi = \xi_2 + \xi_{12} + \xi_{111} + \xi_{112} \\ |\xi_2|, |\xi_{12}|, |\xi_{111}|, |\xi_{112}| >1 }} \frac{ \xi_2 + \xi_{12} +  \xi_{111} + \xi_{112}   }{ \xi_2 \xi_{ 12} \xi_{111} \xi_{112}} \prod_{k=1}^4 \ft z (t, \xi_{\bul_k})  d\vec \xi
\\
&
=
- \frac{ i }{4\times 9} \int_{\substack{\xi = \xi_2 + \xi_{12} + \xi_{111} + \xi_{112} \\ |\xi_2|, |\xi_{12}|, |\xi_{111}|, |\xi_{112}| >1 }} \frac{ \xi   }{ \xi_2 \xi_{ 12} \xi_{111} \xi_{112}} \prod_{k=1}^4 \ft z (t, \xi_{\bul_k})  d\vec \xi
\end{align} 
where we argued by symmetry in $\xi_\bul \in \TT^\infty$, 
 and used the notation $\{\bul_1, \cdots, \bul_4\} = \TT^\infty$. 
Proceeding in a similar fashion for $\Ncal_3^{(2)}(\TT;z)$, via symmetrization in $(\xi_{12}, \xi_{111}, \xi_{112})$,  we obtain
\begin{align*}
\F_x[\Ncal_3^{(2)}(\TT;z)] (\xi)
&
= 
\frac{i}{9} 
\int_{\substack{\xi = \xi_2 + \xi_{12} + \xi_{111} + \xi_{112} \\ |\xi_2|, |\xi_{12}|, |\xi_{111}|, |\xi_{112}| >1 \ge |\xi - \xi_2|}} \frac{ \xi_{12}   }{ \xi_2 \xi_{ 12} \xi_{111} \xi_{112}} \prod_{k=1}^4 \ft z (t, \xi_{\bul_k}) d\vec \xi
\\
&
= 
\frac{i}{3 \times 9} 
\int_{\substack{\xi = \xi_2 + \xi_{12} + \xi_{111} + \xi_{112} \\ |\xi_2|, |\xi_{12}|, |\xi_{111}|, |\xi_{112}| >1 \ge |\xi - \xi_2|}} \frac{ \xi_{12}  + \xi_{111} + \xi_{112} }{ \xi_2 \xi_{ 12} \xi_{111} \xi_{112}} \prod_{k=1}^4 \ft z (t, \xi_{\bul_k}) d\vec \xi
\\
&
= 
\frac{i}{3 \times 9} 
\int_{\substack{\xi = \xi_2 + \xi_{12} + \xi_{111} + \xi_{112} \\ |\xi_2|, |\xi_{12}|, |\xi_{111}|, |\xi_{112}| >1 \ge |\xi - \xi_2|}} \frac{ \xi_1 }{ \xi_2 \xi_{ 12} \xi_{111} \xi_{112}} \prod_{k=1}^4 \ft z (t, \xi_{\bul_k}) d\vec \xi .
\end{align*} 
No simplification is needed for $\Ncal^{(3)}_3(\TT;z)$ since its multiplier in \eqref{mul3-new} already has a small numerator. This shows \eqref{Ncal3} with \eqref{mul3-new} holds for $j=3$, which combined with \eqref{wint0c} gives the general result for $j>3$.

Next, we show \eqref{cancel}. 
For $j=3$, from \eqref{wint0c}, the l.h.s. of \eqref{cancel} reduces to 
\begin{align}
\text{l.h.s. } \eqref{cancel}  \vert_{j=3}
= 
 \Ncal_3^{(1)}(\<31>;z)
+
 \Ncal_3^{(1)}(\<32>;z)
+
 \Ncal_3^{(1)}(\<34>;z)
+
 \Ncal_3^{(1)}(\<35>;z)
+
 \Ncal_4(\<33>;z), 
\end{align}
which together with 
\eqref{Ncal}, \eqref{wint0c}, \eqref{Ncal3}, and \eqref{mul3-new}, gives
\begin{align}
&\Ft_x \bigg[ \sum_{\TT\in \Tr_3}  \big( \Ncal_3^{(1)}(\TT;z)  + \Ncal_4(\TT; z) \big)  \bigg](\xi)
 \\
&=
-
\frac{i}{9} \int_{\substack{ \xi=\xi_1 + \xi_2 + \xi_3+ \xi_4 \\ |\xi_1| , |\xi_{2}|, |\xi_{3}| , |\xi_{4}| > 1   }} \frac{\xi}{ \xi_1 \xi_{2} \xi_{3} \xi_{4} } 
\prod_{k=1}^4 \ft{z}(t, \xi_{k}) d\vec\xi
+ 
\frac{i}{9} \int_{\substack{ \xi=\xi_1 + \xi_2 + \xi_3+ \xi_4 \\ |\xi_1| , |\xi_{2}|, |\xi_{3}| , |\xi_{4}| > 1   }} \frac{\xi}{ \xi_1 \xi_{2} \xi_{3} \xi_{4} } 
\prod_{k=1}^4 \ft{z}(t, \xi_{k}) d\vec\xi
\\
& =0 
\label{cancel2}
\end{align}
as intended, 
where we relabelled the frequencies, for clarity. 
For $j\ge 4$, we must find the correct combination of Type $\III$ and Type $\IV$ terms to see the cancellation in \eqref{cancel}. 
To compare the two types of terms, we define the sets of trees $A^j_{\textup{III}}$ and $A^{j}_{\textup{IV}}$ underlying terms of Type $\III$ and $\IV$ with $(j+1)$-factors, see Remark~\ref{rm:tree_type}:
\begin{align}
A^j_{\textup{III}} : \!& =  \big\{  \TT\in \Tr_j : \ \TT^{0,\ff} =\{\star_1, \star_2\} \text{ s.t. } \star_1 = \sbf(\star_2) \big\}, 
\\
A^j_{\textup{IV}}
:\!&
 =  \big\{  \TT\in \Tr_j : \ |\TT^{0,\ff}| = 1 \big\}
 \\
& = \big\{ \TT \in \Tr_j : \ \exists \wt\TT \in A^{j-3}_{\textup{III}} , \exists \star \in \wt\TT^\infty \text{ s.t. } \TT = \wt \TT \overset{\star}{\oplus} \<33>  \}
, 
\end{align}
see Definition~\ref{DEF:tree3} for the notion of concatenation of trees.
Then, we have the following
\begin{align}
&
\text{l.h.s. } \eqref{cancel} 
\\
&
=
\sum_{\TT \in A^j_{\textup{III}} } \Ncal_3^{(1)}(\TT; z) + \sum_{\TT \in A^j_{\textup{IV}} } \Ncal_4(\TT; z)
\\
&
= 
\sum_{\wt \TT \in A^{j-3}_{\textup{III}}  } 
\sum_{\star \in \wt \TT^\infty }
\bigg( 
 \Ncal_3^{(1)}( \wt\TT \overset{\star}{\oplus} \<31>;z)
+
 \Ncal_3^{(1)}( \wt\TT \overset{\star}{\oplus}   \<32>;z)
+
 \Ncal_3^{(1)}( \wt\TT \overset{\star}{\oplus}   \<34>;z)
\\
&
\hspace{3.5cm}
+
 \Ncal_3^{(1)}( \wt\TT \overset{\star}{\oplus}  \<35>;z)
\bigg) 
+ 
\sum_{\wt \TT \in A^{j-3}_{\textup{III}}  } 
\sum_{\star \in \wt\TT^\infty}
 \Ncal_4( \wt\TT \overset{\star}{\oplus} \<33>;z)
\\
&
= 
\sum_{\wt \TT \in A^{j-3}_{\textup{III}}  } 
\sum_{\star \in \wt \TT^\infty }
\bigg( 
 \Ncal_3^{(1)}( \wt\TT \overset{\star}{\oplus} \<31>;z)
+
 \Ncal_3^{(1)}( \wt\TT \overset{\star}{\oplus}   \<32>;z)
+
 \Ncal_3^{(1)}( \wt\TT \overset{\star}{\oplus}   \<34>;z)
\\
&
\hspace{3.5cm}
+
 \Ncal_3^{(1)}( \wt\TT \overset{\star}{\oplus}  \<35>;z)
+
 \Ncal_4( \wt\TT \overset{\star}{\oplus} \<33>;z)
\bigg) , 
\end{align}
from which the cancellation in \eqref{cancel} follows immediately from the observation that, from \eqref{wint0c}, the portion of the contributions on $\wt\TT$ all agree, combined with the cancellation in \eqref{cancel2} for the bottom of the tree. \qedhere

\end{proof}

\begin{proof}[Proof of Proposition~\ref{PROP:Finfty-new}]

We first note that \eqref{eq:bourg_bi_lp_improv}  follows from the interpolation argument in Subsection~\ref{SUBSEC:bourgainLp}, combining \eqref{eq:Linftybourg1_improv}  with the $L^2$-estimates in Proposition~\ref{prop:bourgL2} for the type $\III$, $j\ge 3$ and type $\IV$ terms, which hold in $X^{s, b_2}_{2}(T)$ for $s>-\frac34$.

We now move to the proof of \eqref{eq:Linftybourg1_improv}. 
We proceed as in the proof of Proposition~\ref{prop:bourgLinfty} in Subsection~\ref{SUB:Linfty}, via a reduction to frequency-restricted estimates.
As before, the multilinear version of our operators $\Ncal_3^{(k)}(\TT;z_1, \ldots, z_{j+1})$ are of the form  $\Ncal(z_1, \ldots, z_{j+1})$ in \eqref{eq:multilinear_geral} with multiplier $m=i \mul_3^{(k)}(\TT)$ as in \eqref{mul3-new}.
 Thus, estimate  \eqref{eq:Linftybourg1_improv}  follows from the inhomogeneous estimates \eqref{lin2a}-\eqref{lin2b}, together with Lemma~\ref{lem:reduction_fre_k2} and Lemma~\ref{lem:reduction_fre_kge3}, once we verify the respective hypotheses.
We only consider the case $s<0$. 
Let $\TT\in \Tr_j$, $j\ge3$, of Type $\III$ with $\TT^{0, \ff} = \{\star\}$. Then, from \eqref{mul3-new}, we have
\begin{align}
\jb{\xi}^s \left(| \mul_{3}^{(2)}(\TT) | + | \mul_{3}^{(3)}(\TT) |  \right) \prod_{\bul \in \TT^\infty} \frac{1}{\jb{\xi_\bul}^s} 
&
\le C^j \frac{ 1  }{\jb{\xi_{\sbf(\pb(\star))}}^{s+1} \jb{\xi_{\sbf(\star)}}^{s+1} \jb{\xi_{\star1}}^{s+1} \jb{\xi_{\star2}}^{s+1} } 
\\&\quad \times \prod_{\bul \in \TT^\infty_\star } \frac{1}{\jb{\xi_\bul}^{s+1} \jb{\xi_{\sbf(\bul)}}} 
\label{mul32}
\end{align}
where $\TT^\infty_\star = \TT^\infty \setminus \{   \sbf(\pb(\star)), \sbf(\star), \star1, \star2\}$.

Since \eqref{eq:estMlinha} holds,  it remains to control
\begin{align}
I = & 
\sup_{\xi, \al} \int_{\xi=\xi_1 + \xi_2 + \xi_3+ \xi_4}
\frac{ 1}{ \jb{\xi_1}^{s+1} \jb{\xi_{2}}^{s+1} \jb{\xi_{3}}^{s+1} \jb{\xi_{4}}^{s+1} } 
\ind_{| \Psi - \al|< M }
d\xi_1 d\xi_2 d\xi_3,
\end{align}
where $\Psi = -\xi^3 + \xi_1^3 + \xi_{2}^3 + \xi_{3}^3 + \xi_{4}^3$.
By symmetry, assume that $|\xi_1| \ge \cdots \ge |\xi_4|$.
Also, if $|\xi_1| \le 1$, the quantity above is trivially bounded for any $s$, thus we assume $|\xi_1| >1$ in the following. 
If $|\xi_1| \not\simeq | \xi_4|$, we argue as in \eqref{typeIIIa}, with the change of variables $\xi_1 \mapsto \phi = \Psi$:
\begin{align}
I & \les \sup_{\xi, \al} \int \frac{ 1 }{ \jb{\xi_1}^{s+3-} \jb{\xi_2}^{s+1} \jb{\xi_3}^{s+1} \jb{\xi_4}^{s+1}} \ind_{|\phi - \al|<M} d\phi d\xi_3 d\xi_4 
\\
&
\les M \int \frac{1}{\jb{\xi_3}^{2(s+1)  + 1 -}  \jb{\xi_4}^{2(s+1) + 1 - } } d\xi_3 d\xi_4 \les M 
\end{align}
for $s>-1$. 
If $|\xi_1| \simeq|\xi_4|$, then $|\xi_1| \simeq \cdots \simeq |\xi_4|$ and we argue via Morse's Lemma (Lemma~\ref{LEM:morse}) and Lemma~\ref{lem:quadraticas}; see~\eqref{typeIIIb}:
\begin{align}
I & \les \sup_{\xi, \al} 
\int 
\ind_{|q_1|, |q_2| \les 1, \, |\Phi(1, q_1,q_2) -\al| < \frac{M}{|\xi_4|^3}}
\jb{\xi_{4}}^{-4(s+1)+2  } dq_1 dq_2 d\xi_4 
\\
&
\les M^{1-} \int \jb{\xi_4}^{-4(s+1)-1} d\xi_4 \les M^{1-}, 
\end{align}
for $s>-1$. 
This completes the proof of Proposition~\ref{PROP:Finfty-new}. 
\qedhere

\end{proof}

\section{Proof of the main results}\label{sec:lwp}

The main goal of this section is to present the proofs of Theorems~\ref{thm:lwp_rkdv} and \ref{thm:wp_kdv_intro} and Proposition  \ref{prop:ill_intro}.
We start by showing some properties of the gauge transform $\Gc$ in \eqref{gauge}, see Subsection~\ref{sec:gauge_prop}.
Then, in Subsection~\ref{SUBSEC:lwp}, we prove well-posedness of \eqref{rkdv} (Theorem~\ref{thm:lwp_rkdv}) via a contraction mapping argument on the gauged normal form equation \eqref{z1}, which together with the equivalence between \eqref{z1} and \eqref{kdv} for smooth solutions (Proposition~\ref{PRO:equiv}), gives well-posedness of the original equation \eqref{kdv} by a density argument (Theorem~\ref{thm:wp_kdv_intro}).
Lastly, in Subsection~\ref{SUBSEC:ill}, we establish the ill-posedness results for \eqref{kdv}.

\subsection{Properties of the gauge transform}\label{sec:gauge_prop}
In this section, we study the mapping properties of the gauge transform $\Gc$ (and its inverse).

For the sake of clarity, let us recall the formal definition of the gauge transform \eqref{gauge} via its inverse: given
$f \in \FL^{s,p}(\R)$,
\begin{align}\label{gauge1}
    \Gc^{-1}[f]
    &
    = f + \Ft_x^{-1}\bigg( - \frac{1}{3} \int_{\xi=\xi_1+\xi_2} \ind_{A^c} \frac{1}{\xi_1\xi_2} \ft{f}(\xi_1) \ft{f}(\xi_2) d\xi_1 \bigg)\\&=f + e^{-t\dx^3} \D[e^{t\dx^3} f , e^{t\dx^3}f]
    \\
    &
     =: f + \DD [f, f]
    ,
\end{align}
for $A^c$ as in \eqref{A} and $\D$ as in \eqref{w} (see also \eqref{eq:D_intro}).
Note that $\Gc^{-1}$ and  $\Gc$ are time-independent.

\begin{lemma}\label{LEM:D}

Let $2<p\le \infty$ and consider the map  $\Gc^{-1}$ defined in \eqref{gauge1}.

\smallskip

\noi {\rm(i)}
Let $s' \ge s>-\frac12 - \frac1p$, $0<\dl\ll \eps \ll1$ sufficiently small, and $B_{\dl}^{s',s}$ given by
\begin{equation}
\label{Bss}
B^{s',s}_\delta=\left\{ f\in \FL^{s',p}(\R): \|f\|_{\FL^{s,p}} \le \delta \right\}.
\end{equation}
Then, the map $\Gc^{-1}: B^{s', s}_\dl \to B^{s',s}_{(1+\eps)\dl}$ is continuous with a continuous inverse $\Gc: B^{s',s}_{(1-\eps)\dl} \to B^{s',s}_\dl$.

\smallskip
\noi {\rm(ii)} The map $\Gc^{-1}$ is unbounded on $\FL^{s,p}(\R)$ for $s\le - \frac12 - \frac1p$.

\end{lemma}

\begin{proof}
\underline{\textit{Step 1. Bilinear estimate for $\DD$.}} Consider the bilinear map $\DD$ defined in \eqref{gauge1} and let $z_1,z_2 \in \FL^{s,p}(\R)$. For $p=\infty$, recalling $A^c$ in \eqref{A}, we have
\begin{align}
&
\|\DD[z_1,z_2]\|_{\mathcal{F}L^{s,\infty}}
\\
&
\les
\sup_\xi \jb{\xi}^{s} \int_{A^c}
\frac{1}{\jap{\xi_{1}}\jap{\xi_{2}}} | \ft{z}_1 (\xi_1)  \ft{z}_2(\xi_2)|   d\xi_{1}
\notag
\\
&\les
\bigg( \sup_{\xi} \int_{|\xi_{1}|>1} \frac{1}{\jap{\xi_{1}}^{2s+2 - s\lor 0 }} d\xi_{1}
\bigg)
 \|z_1\|_{\mathcal{F}L^{s,\infty}}\|z_2\|_{\mathcal{F}L^{s,\infty}}
\notag
\\
& \les
 \|z_1\|_{\mathcal{F}L^{s,\infty}}\|z_2\|_{\mathcal{F}L^{s,\infty}},
\label{gaugeLinfty}
\end{align}
for $s>-\frac12$.
Similarly, for $p=2$, by duality and Cauchy-Schwarz inequality, since $s>-1$, we obtain
\begin{align}
&
\|\DD[z_1,z_2] \|_{H^s}
\\
&
\les
\sup_{\|f\|_{L^2} \le 1 }
	 \int_{\xi, (\xi_1,\xi_2)\in A^c} \frac{\jap{\xi}^{s}|\widehat{f}(\xi)|}{\jap{\xi_{1}}\jap{\xi_{2}}} | \ft{z}_1(\xi_1)  \ft{z}_2(\xi_2) |   d\xi_{1}d\xi
     \notag
    \\
    & \les \sup_{\|f\|_{L^2 } \le 1}  \bigg( \sup_{\xi_{\max }}  \int \frac{1}{\jap{\xi_{\min}}^{2s+4  }}d\xi_{\min}\bigg)^{\frac{1}{2}}
 \|z_1\|_{H^s}
\|z_2\|_{H^{s}}\|f\|_{L^2}
\notag
\\
&
\lesssim
 \|z_1\|_{H^s}
\|z_2\|_{H^{s}} ,
\label{gaugeL2}
\end{align}
where $\xi_{\max}$ (resp. $\xi_{\min}$) correspond to the largest (resp. smallest) frequency over $\xi,\xi_1,\xi_2$, for $s>-1$.
Then, from \eqref{gaugeLinfty} and \eqref{gaugeL2} together with interpolation, we have
\begin{equation}\label{eq:est_D}
		\|\DD[z_1,z_2] \|_{ \mathcal{F}L^{s,p}}
        \le C
\|z_1\|_{ \mathcal{F}L^{s,p}}\|z_2\|_{ \mathcal{F}L^{s,p}}
\end{equation}
for $2 \le p \le \infty$ and $s>-\frac12 - \frac1p$, where $C>0$ depends only on $s$.

\smallskip 
\noi 
\underline{\textit{Step 2. Proof of (i) with $s'=s$.}} By Step 1,  given $z_1,z_2\in B_\delta^{s,s}$ and $\delta C\ll \eps$, we have that
\begin{align*}
\|\DD[z_1,z_1] \|_{ \mathcal{F}L^{s,p}}
        & \le \eps\|z_1\|_{ \mathcal{F}L^{s,p}},
        \\
\|\DD[z_1,z_1] - \DD[z_2,z_2]\|_{ \mathcal{F}L^{s,p}}
        &   \le \eps\|z_1-z_2\|_{ \mathcal{F}L^{s,p}}.
\end{align*}
In particular, $\DD$ is Lipschitz from $B^{s,s}_\dl$ to $B^{s,s}_{\eps\dl}$, with constant $\eps$, meaning that $\Gc^{-1}$ is a small Lipschitz perturbation of the identity map.  The continuity and mapping properties of $\Gc^{-1}$ immediately follow, while $\Gc$ is well-defined as the inverse of $\Gc^{-1}$ via a fixed point argument.

\smallskip
\noi
\underline{\textit{Step 3. Proof of (i) for $s'>s$.}} By Step 2, both $\Gc$ and $\Gc^{-1}$ are well-defined for $s'=s$. It remains to check the claimed boundedness and continuity of $\Gc$ and $\Gc^{-1}$ in the $\FL^{s',p}$ topology.
By \eqref{eq:est_D} (replacing $s$ with $s'$), it follows immediately that $\DD$ (and therefore $\Gc^{-1}$) maps  $\FL^{s',p}$ continuously into itself.

Now we move onto $\Gc$. Performing a bilinear von Neumann expansion (which is valid in the $\FL^{s,p}$ topology by Step 2), we have\footnote{Observe that this corresponds to the sum of all the boundary terms in the integral version of \eqref{NFE}, which is in agreement with the fact that the gauge transform absorbs all the boundary terms.}
$$
\Gc[z]=z+\sum_{j\ge 1}\sum_{\TT\in \Tr_j} \mathcal{I}(\TT,-i\mulB(\TT);z).
$$
Fix $j\ge1$ and $\TT \in \Tr_j$. Since $\mathcal{I}(\TT,-i\mulB(\TT);z_1,\dots, z_{j+1})$ is a $j$-fold composition of $\DD$ (a consequence of the von Neumann expansion), Step 1 readily gives
$$
\|\mathcal{I}(\TT,-i\mulB(\TT);z_1,\dots, z_{j+1})\|_{\FL^{s,p}} \lesssim C^j \prod_{k=1}^{j+1}\|z_k\|_{\FL^{s,p}}.
$$
Since $|\xi|\le (j+1)\max_{\bul\in\TT^\infty}|\xi_\bul|$, we can upgrade this estimate to
$$
\|\mathcal{I}(\TT,-i\mulB(\TT);z_1,\dots, z_{j+1})\|_{\FL^{s',p}} \lesssim C^j \max_{\ell=1,\dots, j+1} \bigg\{ \|z_\ell\|_{\FL^{s',p}}\cdot\prod_{\substack{k=1\\ k\neq \ell}}^{j+1}\|z_k\|_{\FL^{s,p}}\bigg\}.
$$
In particular, over $B^{s',s}_\delta$, for $\delta$ sufficiently small, 
\begin{align*}
\|\Gc[z]\|_{\FL^{s',p}}
&
\le \bigg(1+\sum_{j\ge 1}|\Tr_j|(C\delta)^j\bigg)\|z\|_{\FL^{s',p}} 
\le (1+\eps)\|z\|_{\FL^{s',p}} , 
\\
\|\Gc[z_1]-\Gc[z_2]\|_{\FL^{s',p}}
&
\lesssim \bigg(1+\sum_{j\ge 1}|\Tr_j|(C\delta)^{j-1}\max_{k=1,2}\|z_k\|_{\FL^{s',p}}\bigg) \|z_1-z_2\|_{\FL^{s',p}}\\&\lesssim \big(1+\max_{k=1,2}\|z_k\|_{\FL^{s',p}}\big) \|z_1-z_2\|_{\FL^{s',p}}.
\end{align*}
In particular, $\Gc$ maps continuously $B^{s',s}_{\delta}$ into $B^{s',s}_{(1+\eps)\delta}$, as claimed.

\smallskip

\noi
\underline{\textit{Step 4. Proof of (ii).}}
Let $f \in \FL^{s, p}(\R)$ with
\begin{align}
    \ft{f}(\xi) = \jb{\xi}^{-s-\frac1p}(\log \jap{\xi})^{-\frac12}.
\end{align}
If $2 <p <\infty$, from \eqref{gauge1}, we have
\begin{align}
\| \DD[f,f] \|_{\FL^{s,p}}^p
\ges
\int_{|\xi|\ll 1} \bigg| \int_{A^c} \frac{\ft{f}(\xi_1) \ft{f}(\xi_2) }{\xi_1\xi_2} d\xi_1 \bigg|^p d\xi
\sim \bigg( \int_{|\xi_1| \ge 1} \frac{1}{\jb{\xi_1}^{2s+2 + \frac2p  }\log\jap{\xi_1}} d\xi_1\bigg)^p ,
\end{align}
which diverges for $s\le-\frac12 - \frac1p$.
For $p=\infty$, we have
\begin{align}
\| \DD[f,f]\|_{\FL^{s,\infty }}
\ges
\sup_{|\xi| \ll 1}
\bigg| \int_{A^c} \frac{\ft{f}(\xi_1) \ft{f}(\xi_2) }{\xi_1\xi_2} d\xi_1 \bigg|
\sim
\int_{|\xi_1| \ge 1 } \frac{1}{\jb{\xi_1}^{2s+2}\log\jap{\xi_1}} d\xi_1  ,
\end{align}
which also diverges, allowing us to conclude the unboundedness of $\DD$ on $\FL^{s, p}(\R)$, and thus that of $\Gc^{-1}$.
\qedhere

\end{proof}

\subsection{Local well-posedness}
\label{SUBSEC:lwp}

We are finally in a position to prove Theorem \ref{thm:lwp_rkdv} and Theorem \ref{thm:wp_kdv_intro}.
Using the multilinear estimates of Section \ref{sec:multi}, the local well-posedness of \eqref{rkdv} follows from standard arguments (see, for example, \cite{Grunrock_mkdv, HerrGrunrock_FL}).

\begin{proof}[Proof of Theorem~\ref{thm:lwp_rkdv}]
    The proof is standard, we sketch the $p=\infty$ case.  Fix $b=1-$ and $\wt{b}=0+$ given by  Proposition~\ref{prop:bourgLinfty}.
    For $\delta>0$, define
	$$
	B_{2\delta}=\left\{ z\in \mathcal{S}'(\R^2) : \|z\|_{X^{s,b}_{\infty}(0,T)} + \|z\|_{X^{s,\wt{b}}_{\infty,1}(0,T)}\le 2\delta \right\},
	$$
    endowed with the induced metric,
    and the map
    \begin{align}
		\Theta[z](t) := e^{-t\partial_x^3}z_0
        &+ \int_0^t e^{-(t-t')\partial_x^3} \Ncal_1(\<1>; z)(t')dt'
        \\&+
        \sum_{j=2}^\infty \sum_{\TT \in \Tr_j} \sum_{\l=2}^4
        \int_0^t e^{-(t-t')\partial_x^3} \Ncal_\l (\TT;  z)(t')dt'.
        \label{eq:Theta}
\end{align}
Then, using the linear estimates in $X^{s,b}_{\infty,\infty}\cap X^{s,\wt{b}}_{\infty,1}$ (Lemma \ref{lem:linear_estimates}), the multilinear estimates~\eqref{eq:Linftybourg1}, and \eqref{eq:catalan}, we have
\begin{align*}
		\|\Theta[z]\|_{X^{s,b}_{\infty,\infty}(0,T)\cap X^{s,\wt{b}}_{\infty,1}(0,T)} &\le \|z_0\|_{\mathcal{F}L^{s,p}(\R)} + T^{0+} \sum_{j=1}^\infty \sum_{\TT \in \Tr_{j}}  C^j\|z\|_{X^{s,b}_{\infty,\infty}(0,T)\cap X^{s,\wt{b}}_{\infty,1}(0,T)}^{j+1} \\&\le \|z_0\|_{\mathcal{F}L^{s,p}(\R)} + T^{0+} \sum_{j=1}^\infty \frac{(2j)!}{j!(j+1)!} C^j \delta^{j+1}\le 2\delta,
\end{align*}
for $\delta$ small. A similar computation shows that $\Theta$ is a contraction over $B_{2\delta}$ and the result follows.
\end{proof}

\begin{remark}\label{remark:zs'}
	Additionally to the assumptions of Theorem \ref{thm:lwp_rkdv}, if $z_0\in \FL^{s',p}(\R)$ for some $s'\ge s$, a standard persistence-of-regularity argument shows that the \eqref{rkdv} solution $z$ satisfies
	$$
	\|z\|_{X^{s',b}_p(0,T)} \le 2\|z_0\|_{\FL^{s',p}}, 
	$$
where $T=T(\|z_0\|_{\FL^{s,p}})>0$ from the local well-posedness in $\FL^{s,p}(\R)$ in Theorem~\ref{thm:lwp_rkdv}. 

\end{remark}

We obtain the well-posedness of \eqref{kdv} in Theorem~\ref{thm:wp_kdv_intro} by combining the well-posedness of \eqref{rkdv} in Theorem~\ref{thm:lwp_rkdv} with the equivalence of the two equations at high regularity in the following proposition. For the sake of exposition, we present the details for $p<\infty$, as the $p=\infty$ case follows from suitably replacing the Fourier restriction spaces (as in Subsection \ref{sec:fourier_restr})

\begin{proposition}\label{PRO:equiv}
Let $2 \le p < \infty$, $s' \ge s > - \frac12 - \frac1p$ with $s'\gg1$ sufficiently large, $0<\dl \ll 1$ sufficiently small, and $u_0$ satisfying
\begin{align}
\label{equiv-0}
u_0\in \FL^{s',p}(\R)
\quad
\text{and}
\quad
 \|u_0\|_{\FL^{s,p}}<\delta.
\end{align}
Moreover, let $u \in X^{s',b}_p(0,T)$ be the maximal solution to \eqref{kdv} with initial data $u_0$, for some $T>0$, constructed in
\cite{Grunrock_ldv_hier}, and $z\in X^{s', b}_p(0,1)$ be the solution to
\eqref{rkdv} with initial data $z_0 := \Gc[u_0]$ in Theorem~\ref{thm:lwp_rkdv}.
Then, we must have $T \ge 1$ and the following identity holds:
\begin{align} \label{equiv-1}
z(t)=\Gc[u(t)]
\quad
\text{in}
\quad
C([0,1]; \FL^{s',p}(\R)).
\end{align}

\end{proposition}

We first show Theorem~\ref{thm:wp_kdv_intro} assuming Proposition~\ref{PRO:equiv}.

\begin{proof}[Proof of Theorem \ref{thm:wp_kdv_intro}]

Let $2 \le p <\infty$, $s$ be as in \eqref{eq:reg-kdv}, and $R>0$. For $u_0\in B_R \subset \FL^{s,p}(\R)$, we note that the scaled data $S_\ld[u_0]$, with $0< \ld \le 1$, satisfies
\begin{align}\label{scale-sol}
\| S_\ld[u_0] \|_{\FL^{s,p}}
&
= \ld^{1 + \frac1p + s } \| (\ld^{-2} + \xi^2)^{\frac s2} \ft{u}_0(\xi) \|_{L^p_\xi}
\le
C \ld^{1 + \frac1p + s \land 0} R 
<\dl, 
\end{align}
for some constant $C>0$, given that $s> -1 - \frac1p$, by choosing $\ld=\ld(R, \dl)>0$ sufficiently small. 
Then, for smooth initial data $u_0 \in \FL^{s',p}(\R)$ with $s' \gg1$, from \eqref{scale-sol}, we can choose $\ld>0$ so that $S_\ld[u_0]$ satisfies the hypothesis \eqref{equiv-0} of Proposition~\ref{PRO:equiv}, from which we have that 
\begin{align}
\Gc^{-1} \circ \Theta \circ \Gc \circ S_\ld [u_0] 
\end{align}
agrees with the smooth solution to \eqref{kdv} (say that constructed in \cite{Grunrock_ldv_hier}) with initial data $S_\ld [u_0] $ in $C([0,1]; \FL^{s',p}(\R))$. From the invariance of \eqref{kdv} under scaling, we conclude that $\wt\Theta$ in \eqref{wtTa} agrees with the data-to-solution map for \eqref{kdv} in $C([0,\ld^{3}] ; \FL^{s',p}(\R))$. 
Lastly, by Theorem \ref{thm:lwp_rkdv} and the mapping properties of the gauge transform $\Gc$ (Lemma \ref{LEM:D}), $\wt{\Theta}$ in \eqref{wtTa} (and therefore the \eqref{kdv} flow) has a unique continuous extension to $\FL^{s,p}(\R)$.

\end{proof}

It remains to prove that \eqref{kdv} and \eqref{rkdv} are equivalent for sufficiently smooth solutions, namely Proposition~\ref{PRO:equiv}.
To do so, we consider the frequency-truncated KdV dynamics
\begin{equation}\label{kdvn}\tag{KdV$_n$}
	\partial_t u_n+\partial_x^3 u_n = \P_{n}\partial_x(u_n^2),
\end{equation}
where $\P_n$ denotes the Fourier operator with multiplier $\ind_{[-n,n]}(\xi)$. Following the notations in Subsection~\ref{SUB:tree}, for $j\in \N$, $\TT\in\Tr_j$, and $\bul\in \TT^0$, we introduce the truncated multiplier
\begin{align}\label{Kmuln}
K_{\bul, n}:=-\frac{\mathbbm{1}_{[-n,n]}(\xi_\bul)}{3\xi_{\bul 1}\xi_{\bul2}},
\end{align}
as the truncated version of $K_\bul$ in \eqref{Kmul}, and the truncated gauge transform $\Gc_n$ defined via its inverse $\Gc_n^{-1}$:
   \begin{equation}
	\label{gauge_n}
	\Gc_n^{-1}[f]
	:= f + \Ft_x^{-1}\bigg( \int_{\xi=\xi_1+\xi_2} \ind_{A^c}(\xi_1,\xi_2) K_n(\xi_1,\xi_2) \ft{f}(\xi_1) \ft{f}(\xi_2) d\xi_1 \bigg).
\end{equation}
Note that Lemma~\ref{LEM:D} also applies to $\Gc_n$ and $\Gc_n^{-1}$, where the choices of parameters $\dl, \eps$ are independent of $n$.
Lastly, we introduce the truncated gauged KdV equation
\begin{equation}\label{rkdv_n}\tag{$\Gc$-KdV$_n$}
	\dt z_n
	+ \dx^3 z_n
	=
	\Ncal_{1,n}(\<1>;z_n) + \sum_{j=2}^\infty \sum_{\TT \in \Tr_j} \sum_{\l=2}^4  \Ncal_{\l,n}(\TT;z_n),
\end{equation}
where the nonlinear terms $\Ncal_{\l, n}$ are as in \eqref{Ncal}, replacing $\mul_\l(\TT)$ by the truncated multiplier
\begin{align}
\mul_\l (\TT) (\vec\xi) \prod_{\bul\in \TT^0} \ind_{[-n,n]}(\xi_\bul).
\end{align}
In showing the equivalence in Proposition~\ref{PRO:equiv}, we will argue via the truncated dynamics \eqref{kdvn} and \eqref{rkdv_n}.

We first show multilinear estimates in $\FL^{s,p}(\R)$, needed for the convergence of the series in \eqref{rkdv_n}.
\begin{lemma}\label{lem:estimativas_brutas}
Let $2 \le p \le \infty$, $s'\ge s > - \frac12 - \frac1p$, and $s'\ge s +2$. Then, for $j, n\in\N$, $\TT\in\Tr_j$, $1\le \l \le 4$, and $1 \le k \le j+1$, we have
\begin{equation}\label{eq:boundNn}
		\left\|\mathcal{N}_{\ell,n}(\TT)[f_1,\dots,f_{j+1}] \right\|_{\FL^{s,p}} \le  n^4 C^j \|f_k\|_{\FL^{s,p}}
\max_{\substack{1 \le k' \le  j+1 \\ k' \ne k}}
\bigg(
\|f_{k'}\|_{\FL^{s',p}} \prod_{\substack{k''=1 \\ k'' \ne k , k'}}^{j+1} \|f_{k''}\|_{\FL^{s,p}}
\bigg),
	\end{equation}
for some $C>0$ depending only on $s$, where $\Ncal_{\l, n}$ are the truncated nonlinear terms in~\eqref{rkdv_n}.
\end{lemma}

\begin{proof}
We present the proof for $s<0$, as the result for $s\ge 0$ follows from $\jb{\xi}^s \le j \max_{\bul\in \TT^\infty} \jb{\xi_\bul}^s$ and the estimate in negative regularity.
By \eqref{wint0c} in  Lemma \ref{LEM:wint}, we have
\begin{align*}
| \mul_\l (\TT) (\vec\xi) | \prod_{\bul\in \TT^0} \ind_{[-n,n]}(\xi_\bul)
&
\le C^j n \max_{\bul\in \TT^\infty} |\xi_\bul|^2 \ind_{[-n,n]}(\xi) \bigg( \prod_{\bul\in \TT^0} \frac{1}{\jb{\xi_{\bul1}} \jb{\xi_{\bul2}}} \bigg),
\end{align*}
where we used that $|\xi_\bul| \le n$ for $\bul\in\TT^0$ and $\xi = \sum_{\bul \in \TT^\infty} |\xi_\bul|$.
Let $\star_1, \star_2 \in \TT^\infty$ be the distinct final nodes with $|\xi_{\star_1}|$ and $|\xi_{\star_2}|$ being the largest and second largest frequencies in $\{ |\xi_{\bul}|: \bul \in \TT^\infty\}$. Then, we either have
\begin{align}
|\xi_{\star_1}| \les n
\quad
\text{or}
\quad
|\xi_{\star_1}| \sim |\xi_{\star_2}| \gg n,
\end{align}
which implies the following bound
\begin{align}\label{eq:est_mult_n}
| \mul_\l (\TT) (\vec\xi) | \prod_{\bul\in \TT^0} \ind_{[-n,n]}(\xi_\bul)
&
\le C^j n (n^2 + |\xi_{\star_2}|^2)  \ind_{[-n,n]}(\xi) \bigg( \prod_{\bul\in \TT^0} \frac{1}{\jb{\xi_{\bul1}} \jb{\xi_{\bul2}}} \bigg).
\end{align}

To show \eqref{eq:boundNn} for $2\le p \le \infty$, we show the corresponding estimate for $p=\infty$ and $p=2$, obtaining the result for $2<p<\infty$ via interpolation.
For $p=\infty$, from  \eqref{eq:est_mult_n} we~have
\begin{align*}
	\left\|\mathcal{N}_{\ell,n}(\TT)[f_1,\dots,f_{j+1}] \right\|_{\FL^{s,\infty}}
&
\le C^j n^3 \|f_k\|_{\FL^{s,\infty}}
\max_{\substack{1 \le k' \le j+1 \\ k' \ne k }} \bigg(\|f_{k'}\|_{\FL^{s',\infty}} \prod_{\substack{k''=1 \\ k'' \ne k,k'}}^{j+1} \|f_{k''}\|_{\FL^{s,\infty}}\bigg)
\\
&
\quad
\times
 \sup_{|\xi|\le n} \int_{\Gamma(\TT)}
\bigg(  \prod_{\bul \in \TT^0} \frac{1}{\jap{\xi_{\bul1}}\jap{\xi_{\bul2}}} \bigg)
\bigg( \prod_{\bul \in\TT^\infty}\frac{1}{\jap{\xi_\bul}^s} \bigg)
d\vec{\xi},
\end{align*}
for $s'\ge s+2$, where it remains to bound the last factor. For $|\xi| \le n$ fixed, we have
\begin{align*}
\int_{\Gamma(\TT)}
\bigg(  \prod_{\bul \in \TT^0} \frac{1}{\jap{\xi_{\bul1}}\jap{\xi_{\bul2}}} \bigg)
\bigg( \prod_{\bul \in\TT^\infty}\frac{1}{\jap{\xi_\bul}^s} \bigg)
d\vec{\xi}
&
\le n \int_{\Gamma(\TT)} \bigg( \prod_{\bul \in \TT^0} \frac{1}{\jap{\xi_\bul}} \bigg) \bigg( \prod_{\bul \in \TT^\infty} \frac{1}{\jb{\xi_\bul}^{1+s  }} \bigg) d\vec{\xi}
\\
&
 \le C^j n,
\end{align*}
for $s>-\frac12$, by recursively integrating on final subtrees as follows: given $\bul\in \TT^{0,\ff}$, we have
\begin{align*}
\int_{\xi_\bul = \xi_{\bul1} + \xi_{\bul 2}} \frac{1}{\jb{\xi_{\bul}} \jb{\xi_{\bul1}}^{1+s-} \jb{\xi_{\bul2}}^{1+s-}} d\xi_{\bul1}
\les \frac{1}{\jb{\xi_{\bul}}^{2+2s-}} \le \frac{1}{\jb{\xi_{\bul}}^{1+s-}}
\end{align*}
since $-\frac12<s<0$, where we now consider the smaller tree $\TT \overset{\bul}{\ominus} \<1>$ with $\xi_\bul$ as a final node, where the remaining integrals are of the same form as for the original tree $\TT$. We continue in this way until the remaining tree is $\<1>$, where the same estimate holds.

For $p=2$, with $\xi_{\star1}$ as earlier, by Cauchy-Schwarz inequality and \eqref{eq:est_mult_n}, we have
\begin{align*}
	\left\|\mathcal{N}_{\ell,n}(\TT)[f_1,\dots,f_{j+1}] \right\|_{H^s}
&
\le C^j n^3 \|f_k\|_{H^s}
\max_{\substack{1 \le k' \le j+1 \\ k' \ne k }} \bigg(\|f_{k'}\|_{H^{s'}} \prod_{\substack{k''=1 \\ k'' \ne k,k'}}^{j+1} \|f_{k''}\|_{H^s}\bigg)
\\
&
\quad
\times
 \sup_{\xi_{\star1}} \bigg\{  \int_{|\xi| \le n, \Gamma(\TT)}
\bigg(  \prod_{\bul \in \TT^0} \frac{1}{\jap{\xi_{\bul1}}^2 \jap{\xi_{\bul2}}^2 } \bigg)
\bigg( \prod_{\bul \in\TT^\infty}\frac{1}{\jap{\xi_\bul}^{2s}} \bigg)
d\vec{\xi} d\xi \bigg\}^\frac12,
\end{align*}
for $s'\ge s+2$, as before.
To estimate the last factor, we note that for $|\xi| \le n$, we have
\begin{align*}
&  \int_{|\xi|\le n, \Gamma(\TT)}
\bigg(  \prod_{\bul \in \TT^0} \frac{1}{\jap{\xi_{\bul1}}^2 \jap{\xi_{\bul2}}^2 } \bigg)
\bigg( \prod_{\bul \in\TT^\infty}\frac{1}{\jap{\xi_\bul}^{2s}} \bigg)
d\vec{\xi} d\xi
\\
&
\le \int_{|\xi| \le n, \G(\TT)} \bigg( \prod_{\emptyset \ne \bul\in \TT^0 \setminus \TT^{0,\ff}} \frac{1}{\jb{\xi_\bul}^{2}} \bigg) \bigg( \prod_{\bul\in \TT^{0,\ff}} \frac{1}{\jb{\xi_\bul}^2 \jb{\xi_{\bul1}}^{2+2s} \jb{\xi_{\bul2}}^{2+2s}  }\bigg)
d\vec\xi
d\xi
\\
&
\le C^j n,
\end{align*}
where we first integrate on the final subtrees associated with $\bul \in \TT^{0,\ff}$, by integrating in the smallest variable $\min(|\xi_\bul|, |\xi_{\bul1}|, |\xi_{\bul2}|)$, followed by evaluating the $j - |\TT^{0,\ff}| -1$ integrals in $\xi_\bul$ for $\bul \in \TT^{0} \setminus \TT^{0,\ff}$ and $\bul\ne \emptyset$, and lastly integrating in $\xi$, which leads to the loss of $n$. Note that we dropped all the extra factor  $\jb{\xi_\bul}^{-2(s+1)}$ for $\bul\in \TT^\infty$ which do not have a final sibling, which imposes the restriction $s>-1$.

The estimate \eqref{eq:boundNn} follows by interpolation, completing the proof. \qedhere

\end{proof}

Next we show the well-posedness of the truncated gauged equation \eqref{rkdv_n} in high regularity, with improved uniqueness, as compared to Theorem~\ref{thm:lwp_rkdv}. Namely, uniqueness in $C\big([0,1]; \FL^{s',p}(\R)\big)$ with small $\FL^{s,p}$-norm.

\begin{lemma}\label{lem:zn}
Let $2 \le p < \infty$, $s' \ge s > -\frac12 - \frac1p$ with $s' \ge s+2$,  $0<\dl' \le 1$ sufficiently small, and $n\in\N$.
Given initial data $z_{0,n} \in \FL^{s',p}(\R)$ with
\begin{align}\label{z0n}
\| z_{0,n} \|_{\FL^{s, p}} \le  \dl',
\end{align}
there exists a corresponding solution to \eqref{rkdv_n} $z_n \in X^{s', b}_p(0,1)$. Given $0\le T \le 1$, this solution is unique in
$$
\big\{z\in C([0,T];\FL^{s',p}(\R)): \,  \|z\|_{L^\infty_T \FL^{s,p}_x} \le 10\delta '
\big\}.
$$
Moreover, for $t\in [0,1]$,
\begin{equation}\label{eq:bounds_zn}
	\|z_n(t)\|_{\FL^{s,p}}\le 2\dl'
\quad
\text{and}
\quad
\|z_n(t)\|_{\FL^{s',p}}\le 2\|z_{0,n}\|_{\FL^{s',p}}.
\end{equation}
\end{lemma}
\begin{proof}

Since $z_{0,n} \in \FL^{s,p}(\R)$ with small norm as in \eqref{z0n}, we can repeat the proof of Theorem~\ref{thm:lwp_rkdv} for the truncated gauged equation \eqref{rkdv_n}, since all the multilinear estimates in Section~\ref{sec:multi} apply to nonlinear terms $\NN_{\l, n}$ in \eqref{rkdv_n} by dropping the truncations, uniformly in $n$. Thus, for $0 < \dl' \le 1$ sufficiently small independent of $n$ and $z_{0,n}$, there exists a unique solution $z_n \in X^{s, b}_p(0,1)$ to \eqref{rkdv_n} with initial data $z_{0,n}$, satisfying
\begin{align}
\| z_n \|_{X^{s,b}_p(0,1)} \le 2 \dl ',
\end{align}
from which the first inequality in \eqref{eq:bounds_zn} follows. Similarly, from a standard persistence-of-regularity argument (see Remark~\ref{remark:zs'}), we have that
\begin{align*}
\| z_n \|_{X^{s',b}_p (0,1)} \le 2 \| z_{0,n} \|_{\FL^{s',p}},
\end{align*}
which gives the second estimate in  \eqref{eq:bounds_zn}.

\smallskip

It remains to establish the improved uniqueness. Let $0\le T\le 1$ and $z_n,z'_n\in C\big([0,T], \FL^{s',p}(\R)\big)$ be two solutions of \eqref{rkdv_n} satisfying the corresponding Duhamel formulation, with initial data $z_{0,n}$, and
\begin{align*}
\| z_n \|_{L^\infty_T \FL^{s,p}_x} , \| z_n' \|_{L^\infty_T  \FL^{s,p}_x} \le 10 \dl'
\quad
\text{and}
\quad
\| z_n \|_{L^\infty_T \FL^{s',p}_x} , \| z_n' \|_{L^\infty_T  \FL^{s',p}_x} \le 10 \| z_{0,n} \|_{\FL^{s',p}}.
\end{align*}
Let $0\le T'\le T $, to be chosen later. From Lemma~\ref{lem:estimativas_brutas} and the assumption above, we have
	\begin{align*}
&
\|z_n-z'_n\|_{L^\infty_{T'}\FL^{s,p}_x}
\\
&
\le
T' n^{s'-s} \sum_{j=1}^\infty \sum_{\TT \in \Tr_j} \sum_{\l=1}^4 \big\| \Ncal_{\l, n}(\TT; z_n) - \Ncal_{\l, n} (\TT; z_n') \big\|_{L^\infty_{T'} \FL^{s,p}_x}
\\
&
\les T' n^{4 + s'-s} \|z_n - z_n' \|_{L^\infty_{T'} \FL^{s,p}_x}
\big(
\| z_n \|_{L^\infty_T  \FL^{s',p}_x}
+
\| z_n' \|_{L^\infty_T  \FL^{s',p}_x}
\big)
\sum_{j=0}^\infty (10 C \dl')^j
\\&\les T' n^{4+s'-s}\big(
\| z_n \|_{L^\infty_T  \FL^{s',p}_x}
+
\| z_n' \|_{L^\infty_T  \FL^{s',p}_x}
\big)
 \|z_n - z_n' \|_{L^\infty_T \FL^{s,p}_x}
,
	\end{align*}
for $0<\dl'<1 $ potentially smaller.
For $T'  = T'(n, \|z_{n}\|_{L^\infty_T \FL^{s',p}}, \|z_{n}'\|_{L^\infty_T \FL^{s',p}}) $ small, this shows that $z_n=z'_n$ in $C\big([0,T']; \FL^{s',p}(\R) \big)$. Iterating this argument, guarantees uniqueness in $C\big([0,T] ; \FL^{s',p}(\R)\big)$, completing the proof. \qedhere

\end{proof}

Lastly, before proceeding to the proof of Proposition~\ref{PRO:equiv}, we show the equivalence of the truncated equations \eqref{kdvn} and \eqref{rkdv_n}.

\begin{lemma}\label{lem:equiv_n}
Let $n\in\N$, $2 \le p < \infty$, $s'\ge s> - \frac12 - \frac1p$ with $s'\gg1$ sufficiently large, $0< \dl \le 1$ sufficiently small, and $u_0 \in \FL^{s',p}(\R)$ with
\begin{align}
\| u_0 \|_{\FL^{s,p}} \le \dl.
\end{align}
Then, the following hold:

\noi
{\rm(i)}
There exists a unique solution $u_n \in X^{s', b}_p(0,1)$ to \eqref{kdvn} with
\begin{equation}
\label{eq:bounds_unif_u}
		\|u_n\|_{L^\infty([0,1];\FL^{s,p})}
<3 \delta,
\quad 	\|u_n\|_{L^\infty([0,1]; \FL^{s',p})}
<3\|u_0\|_{\FL^{s',p}},
\end{equation}
\begin{equation}
\| u_n\|_{X^{s', b}_p(0,1)}  \le C,
\label{eq:bound_Xsb}
\end{equation}
for some $C=C(\|u_0\|_{\FL^{s',p}})>0$ independent of $n$.

\noi
{\rm(ii)} Given $z_n \in C\big([0,1]; \FL^{s',p}(\R)\big)$ the solution to \eqref{rkdv_n} with initial data $z_{0,n} := \Gc_n[u_0]$ obtained in Lemma~\ref{lem:zn}, it holds that
\begin{align}
\label{eq:equivn}
z_ n = \Gc_n [u_n]
\quad
\text{in}
\quad
C\big([0,1] ; \FL^{s',p}(\R)\big)
.
\end{align}

\end{lemma}

\begin{proof}

Given $u_0 \in \FL^{s',p}(\R)$, for $s'\gg1$, due to the local well-posedness of \eqref{kdv} at high regularity (say in \cite{Grunrock_ldv_hier}), there exists a unique solution $u_n \in X^{s',b}_p(0,T_n)$ to \eqref{kdvn} where $T_n>0$ is the maximal time of existence. We also define $T_{0,n}$ by
\begin{align}
T_{0,n} :=
\max \big\{
0 < T \le \min(1, T_n) : \, \| u_n\|_{L^\infty_T \FL^{s,p}_x} < 3 \dl
\big\}.
\end{align}

\noi
\underline{Step 1.} We first show that, for $\delta\ll 1$, $z'_n := \Gc_n[u_n]$ solves \eqref{rkdv_n} in $C\big([0,T_{0,n}] ; \FL^{s',p}(\R)\big)$ with
\begin{equation}
	z_n'\in C([0,T_{0,n}];\FL^{s',p}(\R))
\quad
\text{and}
\quad
\|z'_n\|_{L^\infty_{T_{0,n}} \FL^{s,p}_x} < (1+\eps) 3 \delta,
\label{eq:z'}
\end{equation}
for some $0<\eps<1$ sufficiently small.
From Lemma~\ref{LEM:D} (adapted to $\Gc_n, \Gc_n^{-1}$) and the definition of $T_{0,n}$, we obtain \eqref{eq:z'}.
It remains to prove that $z'_n$ solves \eqref{rkdv_n}.

Let $v_n(t) : = e^{t\dx^3} u_n(t)$.
Then, $v_n \in C \big( [0,T_{0,n}] ; \FL^{s',p}(\R) \big)$ and it satisfies the truncated \eqref{NFE} (with $K_\bul$ in \eqref{Kmul} replaced with $K_{\bul, n}$ in \eqref{Kmuln}, in the definition of the nonlinear operators in \eqref{NFE1a}).
This is tantamount to the formalization of the infinite normal form reduction (for the truncated \eqref{kdvn}), which relies on four technical assumptions (see also the discussion in \cite[Section 4]{KOY20}): (a) the differentiation-by-parts in time; (b) the dominated convergence theorem to exchange time derivatives and frequency integrals; (c) the summability of the terms $\NN_n(\TT;v_n)$ and $\BB_n(\TT;v_n)$; and (d) the convergence of the remainders $\MM_n(\TT;v_n)$ to zero (in some topology).

Concerning (a)-(b),  this is guaranteed by the frequency decay of both $\widehat{v}_n$ and $\dt \widehat{v}_n$ (as $s\gg 1$) and the fact that, for $\xi$ fixed, $\ft{v_n}(\xi)$ is $C^1$ in the time variable.
For (c)-(d), we argue as in the proof of Lemma \ref{lem:estimativas_brutas}. For $j\in \N$, $\TT\in \Tr_j$, and $0 \le t < T_{0,n}$, we have
\begin{align*}
&
	\left\|\NN_n(\TT;v_n)(t)\right\|_{\FL^{s',p}} + 	\left\|\BB_n(\TT;v_n)(t)\right\|_{\FL^{s',p}} + 	\left\|\MM_n(\TT;v_n)(t)\right\|_{\FL^{s',p}}
\\
&\les_n \|v_n\|_{L^\infty_{T_{0,n}}\FL^{s',p}_x}  C^j  \|v_n\|_{L^\infty_{T_{0,n}}\FL^{s,p}_x}^j
\\
&\les_n \|v_n\|_{L^\infty_{T_{0,n}}\FL^{s',p}_x} \times (3  \delta C )^j,
\end{align*}
where we use $\NN_n, \BB_n,\MM_n$ to denote the analogues of the operators in \eqref{NFE1a} with the additional frequency truncation $|\xi_\bul| \le n$ at each $\bul \in \TT^0$.
Since the number of trees $|\Tr_j|$ grows at most exponentially with $j$ (see \eqref{eq:catalan}), the estimate above ensures the absolute summability of the terms $\NN_n(\TT;v_n)$ and $\BB_n(\TT;v_n)$ and the convergence of the remainders $\MM_n(\TT;v_n)$ to zero in the $L^\infty_t \FL^{s',p}$ topology, for $0< \dl <1$ sufficiently small, depending only on the universal constant $C$. This shows (c)-(d).
We then conclude that $v_n$ satisfies the (truncated) infinite normal form equation \eqref{NFE}.

From the invertibility properties of $\Gc_n$ in $B^{s',s}_{3\dl}$ in Lemma~\ref{LEM:D}, we know that $\Gc_n^{-1} [z'_n] = u_n$ in $C\big([0,T_{0,n}] ; \FL^{s',p}(\R)\big)$, and we can follow the argument in Subsection~\ref{SUBSEC:gauge} to conclude that $z'_n=\Gc_n[u_n]$ is an integral solution to \eqref{rkdv_n}. In particular, note that the cancellations observed in Lemmas \ref{LEM:wbd} and \ref{LEM:wint} are term-by-term cancellations.

\smallskip
\noindent\underline{Step 2.} From Lemma~\ref{LEM:D}, we have
\begin{align}
\| z_{0,n} \|_{\FL^{s,p}} &  = \| \Gc_n [ u_0] \|_{\FL^{s,p}} \le \tfrac{1}{1-\eps} \dl.
\end{align}
Then,
the solution $z_n \in C\big( [0,1]; \FL^{s',p}(\R)\big)$ to \eqref{rkdv_n} in Part (ii), constructed in Lemma~\ref{lem:zn}, satisfies
\begin{align}
\| z_n \|_{L^\infty([0,1]; \FL^{s,p})} \le \tfrac{2}{1-\eps} \dl
.
\end{align}
Thus, together with \eqref{eq:z'}, by Lemma~\ref{lem:zn}, we conclude that
\begin{align}
z_n = z_n'
\quad
\text{in }
C([0,T_{0,n}]; \FL^{s',p}(\R)),
\end{align}
and from \eqref{eq:bounds_zn}, we get
\begin{align}
\| z_n'\|_{L^\infty_{T_{0,n}} \FL^{s,p}_x}
&=
\| z_n \|_{L^\infty_{T_{0,n}} \FL^{s,p}_x} \le \tfrac{2}{1-\eps}\dl.
,
\\
\| z_n'\|_{L^\infty_{T_{0,n}} \FL^{s',p}_x}
&=
\| z_n \|_{L^\infty_{T_{0,n}} \FL^{s',p}_x} \le 2 \|z_{0,n} \|_{\FL^{s',p}}.
\end{align}
Combining the above with the mapping properties of $\Gc_n$ from Lemma~\ref{LEM:D} gives
\begin{align}
\| u_n \|_{L^\infty_{T_{0,n}} \FL^{s,p}_x}
&
\le (1+\eps)
\| z'_n \|_{L^\infty_{T_{0,n}} \FL^{s,p}_x}
\le
\tfrac{2(1+\eps)}{1-\eps} \dl
< 3 \dl,
\label{uns}
\\
\| u_n \|_{L^\infty_{T_{0,n}} \FL^{s',p}_x}
& \le (1+\eps)
\| z'_n \|_{L^\infty_{T_{0,n}} \FL^{s',p}_x}
\le \tfrac{2 (1+\eps) }{1-\eps} \|u_0 \|_{\FL^{s',p}},
\label{unsprime}
\end{align}
for $0< \eps\ll 1 $ sufficiently small.
From the definition of $T_{0,n}$ and a continuity argument, we conclude that $T_{0,n} = \min(1, T_n)$.
Lastly,
if $T_n < 1$ and from the blow-up alternative, we have
\begin{align}
\lim_{t \to T_n} \| u_n(t) \|_{\FL^{s',p}} = +\infty,
\end{align}
which contradicts the uniform in time bound in \eqref{unsprime}. Therefore, we must have $T_n \ge 1$, and the bounds in \eqref{eq:bounds_unif_u} follow from \eqref{uns}-\eqref{unsprime}.

\smallskip
\noindent\underline{Step 3.}
It remains to show \eqref{eq:bound_Xsb}.

From Steps 1-2 we have that $u_n \in X^{s', b}_p(0,1)$ with the bounds in \eqref{eq:bounds_unif_u}.
Consider the uniform in $n$ bound in \eqref{unsprime}, which we can write as
\begin{align}\label{unsprime2}
\| u_n \|_{L^\infty([0,1]; \FL^{s',p})} < 3 \|u_0\|_{\FL^{s',p}} = : M.
\end{align}
Then, from the well-posedness theory for \eqref{kdv} in \cite{Grunrock_ldv_hier} together with the uniform bound on $[0,1]$ in \eqref{unsprime2}, there exists a local time of existence $0<T = T(M)\le 1 $ independent of $n$ where the solution $u_n$ to \eqref{kdvn} satisfies
\begin{align}
\| u_n \|_{X^{s',b}_p(jT,(j+1)T)} \le 2 M , \qquad 0\le j \le \lceil \tfrac1T \rceil -1 .
\end{align}
Then, \eqref{eq:bound_Xsb} follows from a standard gluing argument, completing the proof. \qedhere

\end{proof}

We complete this subsection with proof of equivalence of \eqref{kdv} and \eqref{rkdv} for smooth solutions.

\begin{proof}[Proof of Proposition~\ref{PRO:equiv}]

Fix $u_0 \in \FL^{s',p}(\R)$ as in \eqref{equiv-0} and let $z \in X^{s',b}_p(0,1)$ be the solution to \eqref{rkdv} with data $\Gc[u_0]$ constructed in Theorem~\ref{thm:lwp_rkdv}.

From the well-posedness of \eqref{kdv} in \cite{Grunrock_ldv_hier}, there exist $T>0$ and a solution $u \in X^{s',b}_p(0,T)$ to \eqref{kdv} with initial data $u_0$, where we assume that $T$ is the maximal time of existence.
We will show that $T\ge 1$ and that $u = \Gc^{-1} [z]$ in $C([0,1]; \FL^{s',p}(\R))$ by arguing via approximation through the truncated dynamics \eqref{kdvn} and \eqref{rkdv_n}.

Combining Lemmas~\ref{lem:zn} and \ref{lem:equiv_n}, there exist $z_n, u_n \in X^{s',b}_p(0,1)$ solutions to \eqref{rkdv_n} and \eqref{kdvn}, respectively, with $z_n (0) = \Gc_n[u_0] $ and $u_n(0) = u_0$.
We claim that the following convergence results hold:
\begin{align}
\label{equiv1a}
z_n & \to z \quad \text{in} \quad  C([0,1]; \FL^{s'-1,p}(\R))
,
\\
\label{equiv1c}
\Gc_n [u_n] & \to  \Gc[u]   \quad \text{in} \quad  C([0,T); \FL^{s'-1,p}(\R)),
\end{align}
as $n\to\infty$.
Assuming \eqref{equiv1a}-\eqref{equiv1c} together with the identity \eqref{eq:equivn} in Lemma~\ref{lem:equiv_n}, we note that
\begin{align}
\| \Gc[u] - z\|_{L^\infty_T \FL^{s'-1,p}_x}
&
 \le
\| \Gc[u] - \Gc_n[u_n] \|_{L^\infty_T \FL^{s'-1,p}_x}
+
\| z_n - z \|_{L^\infty_T \FL^{s'-1,p}_x} \to 0,
\end{align}
as $n\to \infty$
which implies that
$u = \Gc^{-1}[z]$ on $C([0,T); \FL^{s',p}(\R))$.
Moreover, note that for $0\le t < T$, we have
\begin{align}
\| u (t) \|_{\FL^{s',p}}  = \|( \Gc \circ \Gc^{-1}) [u(t)] \|_{\FL^{s',p}} \le \tfrac{1}{1-\eps} \| \Gc^{-1}[u](t) \|_{\FL^{s'p}} \le \tfrac{1}{1-\eps} \|z\|_{L^\infty([0,1); \FL^{s',p})}
\end{align}
where we used Lemma~\ref{LEM:D}, and thus the norm of $u$ cannot blow-up as $t$ approaches $T$, and we must have $T \ge 1$, proving \eqref{equiv-1}.
It remains to show \eqref{equiv1a}-\eqref{equiv1c}.

\smallskip

We start by showing \eqref{equiv1a}.
From \eqref{rkdv} and \eqref{rkdv_n}, we have
\begin{align}
&
z(t) - z_n(t)
\\
&
= e^{-t\dx^3} (\Gc - \Gc_n)[u_0]
+ \sum_{j=1}^\infty \sum_{\TT\in\Tr_j} \sum_{\l=1}^4 \int_0^t e^{-(t-t')\dx^3} [\Ncal_{\l}(\TT; z) - \Ncal_{\l, n}(\TT;z)]  dt'
\\
&
\quad
+ \sum_{j=1}^\infty \sum_{\TT\in\Tr_j} \sum_{\l=1}^4
\int_0^t e^{-(t-t')\dx^3} [\Ncal_{\l, n}(\TT; z) - \Ncal_{\l, n}(\TT;z_n)]  dt'.
\end{align}
From the multilinear estimates in Lemma~\ref{lem:linear_estimates_qfinite} (which also apply to $\NN_{\l,n}$ with constants independent of $n$), the first bound in \eqref{eq:bounds_zn}, and the respective bound for $z$ which follows from the proof of Theorem~\ref{thm:lwp_rkdv}, we have that
\begin{align}
&
\sum_{j=1}^\infty \sum_{\TT\in\Tr_j} \sum_{\l=1}^4
\bigg\|
\int_0^t e^{-(t-t')\dx^3} [\Ncal_{\l, n}(\TT; z) - \Ncal_{\l, n}(\TT;z_n)]  dt' \bigg\|_{X^{s'-1,b}_p(0,1)}
\\
&
\le \sum_{j=1}^\infty ( C \tfrac{1}{1-\eps} \dl )^{j} \| z-z_n\|_{X^{s'-1,b}_p(0,1)} \\
&
\le  \tfrac12 \| z-z_n\|_{X^{s'-1,b}_p(0,1)} ,
\end{align}
for $\dl>0$ sufficiently small, where $0<\eps<1$ is as in Lemma~\ref{LEM:D}, namely the Lipschitz constant of $\DD$ in \eqref{eq:est_D}.
Then, we have
\begin{align}
\label{equiv1aa}
\begin{split}
\| z - z_n \|_{X^{s'-1,b}_p(0,1)}
&
\le
2 \|e^{-t\dx^3} (\Gc - \Gc_n)[u_0]\|_{X^{s'-1,b}_p(0,1)}
\\
&
\quad
+
2
\sum_{j=1}^\infty \sum_{\TT\in\Tr_j} \sum_{\l=1}^4
\bigg\| \int_0^t e^{-(t-t')\dx^3} [\Ncal_{\l}(\TT; z) - \Ncal_{\l, n}(\TT;z)]  dt' \bigg\|_{X^{s'-1,b}_p(0,1)}.
\end{split}
\end{align}

For the first term in \eqref{equiv1aa}, we note that by the definition of $\Gc, \Gc_n$, it holds that
\begin{align}
u_0
=
\Gc[u_0] + \DD\big[\Gc[u_0], \Gc[u_0]\big]
=
\Gc_n[u_0] + \P_n \DD\big[\Gc_n[u_0], \Gc_n[u_0]\big],
\end{align}
for $\DD$ as in \eqref{gauge1}. Then, from \eqref{gaugeLinfty}, the mapping properties of $\Gc$, and \eqref{equiv-0}, we have
\begin{align}
&
\| (\Gc - \Gc_n)[u_0] \|_{\FL^{s'-1,p}}
\\
&
\le \big\| \P_n \big( \DD\big[\Gc_n[u_0], \Gc_n[u_0]\big]- \DD\big[\Gc[u_0], \Gc[u_0]\big]  \big)  \big\|_{\FL^{s'-1,p}}
+
\big\| \P_{>n} \DD\big[\Gc[u_0], \Gc[u_0]\big]  \big\|_{\FL^{s'-1,p}}
\\
&
\le 2  \dl \tfrac{\eps}{1-\eps} \| (\Gc - \Gc_n)[u_0] \|_{\FL^{s'-1,p}} +
\big\| \P_{>n} \DD\big[\Gc[u_0], \Gc[u_0]\big]  \big\|_{\FL^{s'-1,p}}
\end{align}
which implies that
\begin{align}
\| (\Gc - \Gc_n)[u_0] \|_{\FL^{s'-1,p}}
\les \big\| \P_{>n} \DD\big[\Gc[u_0], \Gc[u_0]\big]  \big\|_{\FL^{s'-1,p}}\lesssim \frac{1}{n}\big\| \P_{>n} \DD\big[\Gc[u_0], \Gc[u_0]\big]  \big\|_{\FL^{s',p}}
\end{align}
for $\eps,\dl>0$ sufficiently small, where the implicit constant is independent on $n$. Since $\DD\big[\Gc[u_0], \Gc[u_0]\big] \in \FL^{s',p}(\R)$, the r.h.s. converges to 0 as $n\to\infty$, and we conclude the convergence of the first term on the r.h.s. of \eqref{equiv1aa}.

Similarly, the second contribution in \eqref{equiv1aa} converges to 0, by combining the observation that we can rewrite $\NN_{\l}(\TT) - \NN_{\l, n}(\TT)$ as a sum of terms with a similar structure but with $\ind_{|\xi_{\bul}| > n}$ for some $\bul \in \TT^0$, to which the estimates in Lemma~\ref{lem:linear_estimates_qfinite} apply, together with a factor of $\frac1n$ from upgrading the regularity from $s'-1$ to $s'$.
This completes the proof of \eqref{equiv1a}.

\smallskip

Next we show \eqref{equiv1c}.
As a first step, we show that
\begin{align}
\label{equiv1b}
u_n & \to u \quad \text{in} \quad    C([0,T); \FL^{s'-1,p}(\R)).
\end{align}
From \cite{Grunrock_ldv_hier}, there exists $0< T_* \le T$ depending only on $\| u_0\|_{\FL^{s',p}}$ such that
\begin{align}
\| u\|_{X^{s',b}_p(0,T_*)}, \| u_n\|_{X^{s',b}_p(0,T_*)} \le 2 \| u_0\|_{\FL^{s',p}} =: M.
\end{align}
Then, from the relevant bilinear estimate in \cite{Grunrock_ldv_hier}, we have
\begin{align}
\| u_n - u\|_{X^{s'-1,b}_p(0,T_*)}
&
= \bigg\| \int_0^t e^{-(t-t')\dx^3} \dx [ \P_n (u_n^2) - u^2 ] dt' \bigg\|_{X^{s'-1,b}_p(0,T_*)}
\\
&
\le 2M T_*^\ta \| u_n - u \|_{X^{s'-1,b}_p(0,T_*)} + \bigg\| \int_0^t e^{-(t-t')\dx^3} \dx \P_{>n} (u^2)  dt' \bigg\|_{X^{s'-1,b}_p(0,T_*)}
\\
&
\le \tfrac12 \| u_n - u \|_{X^{s'-1,b}_p(0,T_*)} + \bigg\| \int_0^t e^{-(t-t')\dx^3} \dx \P_{>n} (u^2)  dt' \bigg\|_{X^{s'-1,b}_p(0,T_*)}
\end{align}
for some $\ta>0$ and by taking $T_*$ potentially smaller so that the second inequality holds. Then,
\begin{align}
\| u_n - u\|_{X^{s'-1,b}_p(0,T_*)}
\le \frac2n  \bigg\| \int_0^t e^{-(t-t')\dx^3} \dx \P_{>n} (u^2)  dt' \bigg\|_{X^{s',b}_p(0,T_*)} \to 0 ,
\end{align}
as $n\to\infty$, from the bilinear estimate in \cite{Grunrock_ldv_hier}.
We can iterate the above argument on intervals $[ jT_*, (j+1)T_* )$ for $1 \le j \le \lceil \frac{T}{T_*} \rceil -1 $, noting that from the previous iteration we have $\| u_n( jT_*) - u(jT_*) \|_{\FL^{s',p}} \to 0$ as $n\to\infty$. This concludes the proof \eqref{equiv1b}, which combined with \eqref{eq:bounds_unif_u} also implies that
\begin{align}
\label{equiv1bb}
\| u \|_{L^\infty_T \FL^{s,p}_x} < 3 \dl .
\end{align}

We can now show \eqref{equiv1c}. We proceed as in the proof of convergence of the first term in \eqref{equiv1aa}. From the properties of $\Gc, \Gc_n$ in Lemma~\ref{LEM:D} (and its proof), we have that
\begin{align}
u_n &= \Gc_n[u_n] + \P_n \DD\big[ \Gc_n[u_n], \Gc_n[u_n]\big],
\quad
u = \Gc[u] +  \DD\big[ \Gc[u], \Gc[u]\big],
\end{align}
and thus
\begin{align}
&\| \Gc_n[u_n] - \Gc[u] \|_{L^\infty_T \FL^{s'-1,p}_x}
\\
&
\le \| u_n - u \|_{L^\infty_T \FL^{s'-1,p}_x}
+
\big\| \P_n \big(  \DD\big[ \Gc_n[u_n], \Gc_n[u_n]\big] -  \DD\big[ \Gc[u], \Gc[u]\big]  \big) \big\|_{L^\infty_T \FL^{s'-1,p}_x}
\\
&
\quad
+
\big\| \P_{>n} \DD\big[ \Gc[u], \Gc[u]\big] \big\|_{L^\infty_T \FL^{s'-1,p}_x}.
\end{align}
From \eqref{eq:est_D} and the smallness of $u_n,u$ in $L^\infty_T\FL^{s,p}(\R)$ (see \eqref{eq:bounds_unif_u} and
\eqref{equiv1bb}), we get
\begin{align}
&
\big\| \P_n \big(  \DD\big[ \Gc_n[u_n], \Gc_n[u_n]\big] -  \DD\big[ \Gc[u], \Gc[u]\big]  \big) \big\|_{L^\infty_T \FL^{s'-1,p}_x}
\\
&
\le C (\| \Gc_n[u_n] \|_{L^\infty_T \FL^{s,p}_x} + \| \Gc[u] \|_{L^\infty_T \FL^{s,p}_x}) \| \Gc_n[u_n] - \Gc[u] \|_{L^\infty_T \FL^{s'-1,p}_x}
\\
&
\le 2C \tfrac{1}{1-\eps} \dl  \| \Gc_n[u_n] - \Gc[u] \|_{L^\infty_T \FL^{s'-1,p}_x} < \tfrac12  \| \Gc_n[u_n] - \Gc[u] \|_{L^\infty_T \FL^{s'-1,p}_x}
\end{align}
for some $C>0$ and $0<\dl<1$ sufficiently small, giving
\begin{align}
&\| \Gc_n[u_n] - \Gc[u] \|_{L^\infty_T \FL^{s'-1,p}_x}
\\
&
\le \|u_n-u \|_{L^\infty_T \FL^{s'-1,p}_x} +
\frac1n\big\| \P_{>n} \DD\big[ \Gc[u], \Gc[u]\big] \big\|_{L^\infty_T \FL^{s',p}_x}
\to 0,
\end{align}
as $n\to\infty$, from \eqref{equiv1b} for the first term and \eqref{eq:est_D} and boundedness of $\Gc$ in Lemma~\ref{LEM:D} for the second term.
This completes the proof of
\eqref{equiv1c}, as intended. \qedhere

\end{proof}

\subsection{Ill-posedness}
\label{SUBSEC:ill}

We now move to the proof of Theorem \ref{thm:main}-(ii)-(iii). For the sake of exposition, we record here the precise statement.

\begin{proposition}[Ill-posedness for \eqref{kdv} in $\FL^{s,p}$]\label{prop:ill_intro}
		Fix $2\le p \le \infty$.
        
    \noi {\rm(i)} For $2\le p <\infty$ and   $s<-1+\frac1{2p}$, the \eqref{kdv} flow is not uniformly continuous in the $\FL^{s,p}(\R)$ topology.
        
    \noi {\rm(ii)} For $s<s_{H^{-1}}(p)=-\frac12 - \frac1p$, 
    there is no continuous extension of the \eqref{kdv} flow from $\mathcal{F}L^{s,p}(\R)$ to {$\mathcal{D}'((0,T)\times \R)$}, for any $T>0$.
	\end{proposition}

To prove (i), 
we adopt the same strategy of Kenig, Ponce, and Vega \cite{KPV_ill}, which combines the Miura map in \eqref{eq:miura2} with an explicit family of solutions known as \emph{breathers} \cite{Lamb,Wadati}. These solutions are given (up to a constant scaling factor) by
$$
b_{N,\omega}(t, x)=2\sqrt{6} 
\o
\sech(\omega x +\gamma t)\frac{\cos(Nx+\delta t)-(\omega/N)\sin(Nx+\delta t)\tanh(\omega x+\delta t)}{1+(\omega/N)^2\sin^2(Nx+\delta t)\sech(\omega x + \gamma t)}
$$
where
\begin{equation}\label{breather1}
N,\omega>0,\quad \delta=N(N^2-3\omega^2),\quad \gamma=\omega (3N^2-\omega^2).
\end{equation}
For $N\gg 1$ and $\eps>0$, set
\begin{equation}\label{eq:escolhas_breathers}
\omega=N^{-p(s+1)},\quad N_1=N,\quad N_2=N-\eps\omega.
\end{equation}
and
\begin{equation}\label{bnj}
b_j=b_{N_j,\omega},\quad j=1,2.
\end{equation}

The first ingredient for the proof of Proposition~\ref{prop:ill_intro}~(i) is the fact that, for $N$ large, these breather solutions  start $\varepsilon$-close to each other at time $t=0$ and become separated at any fixed  time $T>0$.

\begin{lemma}\label{lem:breather}
Fix $2\le p < \infty$ and $-1 - \frac1p<s<-1+\frac1{2p}$. For $T>0$ and $0<\varepsilon\ll 1$, there exists $c>0$ such that, for $N$ large, we have
\begin{equation}\label{eq:breathers1}
	\|b_1(0)-b_2(0)\|_{\FL^{s+1,p}}\lesssim \eps, \quad \|b_j(0)\|_{\FL^{s+1,p}} \lesssim 1,\ j=1,2,
\end{equation}
for $b_{j}$ as in \eqref{bnj},
and, for some $0<T'\le T$,
\begin{equation}\label{eq:breathers2}
	\|b_1(T')-b_2(T')\|_{\FL^{s+1,p}}\ge c.
\end{equation}
\end{lemma}
\begin{proof}
	Since $\omega\ll N$, 
	$$
	b_j(t,x)\sim b_j^0(t,x):=\omega \sech(\omega x + \gamma_j t)\cos(N_jx+\delta_jt),\quad j=1,2, 
	$$
    where $\dl_j$ and $\g_j$ are as in \eqref{breather1} with $N=N_j$, $j=1,2$.
     For general $s, N, \omega, N_1, N_2$, with $N_1\sim N_2\sim N$, it was shown in \cite{KPV_ill} that
    $$
\|b_j^0(0)\|_{H^{s+1}}\lesssim N^{s+1}\omega^{\frac12},\quad j=1,2,
    $$
and 
$$
\|b_1^0(0)-b_2^0(0)\|_{H^{s+1}}\lesssim N^{s+1}\omega^{\frac12} \frac{|N_1-N_2|}{\omega}.
$$
On the other hand, we have
$$
	\|b_j^0(0)\|_{\FL^{s+1,\infty}}\le N^{s+1}\|b_j^0(0)\|_{L^1}\sim N^{s+1}
	$$
	and
\begin{align*}
		\|b_1^0(0)-b_2^0(0)\|_{\FL^{s+1,\infty}}&\lesssim N^{s+1}\|b_1^0(0)-b_2^0(0)\|_{L^1}\\&
        \lesssim 
        N^{s+1}\int |\omega \sech (\omega x)||\cos(N_1x)-\cos(N_2x)|dx \lesssim  N^{s+1}\frac{|N_1-N_2|}{\omega}.
\end{align*}
It follows by interpolation that
$$
	\|b_j^0(0)\|_{\FL^{s+1,p}}\lesssim N^{s+1}\omega^{\frac1p},\quad \|b_1^0(0)-b_2^0(0)\|_{\FL^{s+1,p}}\lesssim N^{s+1}\omega^{\frac1p}\frac{|N_1-N_2|}{\omega}.
	$$
    Hence, for the specific choices of $N_1, N_2$ and $\omega$ made in \eqref{eq:escolhas_breathers}, we find \eqref{eq:breathers1}.

To prove \eqref{eq:breathers2}, we consider 
$$
\mathbf{b}_j^0 (t,x) = \omega \sech(\omega x + \gamma_jt)e^{iN_jx+i \delta_jt}, \quad j=1,2.
$$
Then
$$
\widehat{\mathbf{b}_j^0}(t, \xi) \sim e^{it\left(\delta_j+\gamma_j\frac{\xi-N_j}{\omega}\right)}\widehat{\sech}\left(\frac{\xi-N_j}{\omega}\right)
$$
and
\begin{align*}
	\left\| {\mathbf{b}_1^0}(T)  -{\mathbf{b}_2^0}(T) \right\|_{\FL^{s+1,p}}&\ge \left\|\jap{\xi}^{s+1} \left(e^{iT\left(\delta_1+\gamma_1\frac{\xi-N_1}{\omega}\right)}- e^{iT\left(\delta_2+\gamma_2\frac{\xi-N_2}{\omega}\right)}\right)\widehat{\sech}\left(\frac{\xi-N_1}{\omega}\right)\right\|_{L^p}\\&\quad - \left\|\jap{\xi}^{s+1}\left( \widehat{\sech}\left(\frac{\xi-N_1}{\omega}\right)-\widehat{\sech}\left(\frac{\xi-N_2}{\omega}\right)\right)\right\|_{L^p}.
\end{align*}

We write
$$
\big| e^{iT(\delta_1+\gamma_1\frac{\xi-N_1}{\omega})}- e^{iT(\delta_2+\gamma_2\frac{\xi-N_2}{\omega})} \big| = \big|e^{iT(\gamma_2-\gamma_1)\cdot\frac{\xi-N_1}{\omega}}-e^{iT\Theta}\big|, \quad \Theta = \delta_1-\delta_2-\gamma_2\tfrac{N_1-N_2}{\omega}.
$$
\\
Let $\eta>0$ be such that $|\widehat{\sech}(\zeta)|^p\ge \frac{1}{2}|\widehat{\sech}(0)|^p$ for $|\zeta|\le \eta$.  Since $T|\gamma_1-\gamma_2|\sim \varepsilon TN^{1-2(s+1)p}\to \infty$, we conclude that, for $N$ large,
$$
\Big|\Big\{ \zeta\in(-\eta,\eta) : \big|e^{iT(\gamma_1-\gamma_2)\cdot\zeta}-e^{iT\Theta}\big|\ge \tfrac{1}{100}\Big\}\Big|\ge\tfrac{\eta}{2}.
$$
In particular,
\begin{align*}
 &\left\|\jap{\xi}^{s+1} \left(e^{iT\left(\delta_1+\gamma_1\frac{\xi-N_1}{\omega}\right)}- e^{iT\left(\delta_2+\gamma_2\frac{\xi-N_2}{\omega}\right)}\right)\widehat{\sech}\left(\frac{\xi-N_1}{\omega}\right)\right\|_{L^p} \\ \ge &\  N^{s+1}\omega^{\frac1p} \left\| \left(e^{iT(\gamma_1-\gamma_2)\cdot\zeta}-e^{iT\Theta}\right)\widehat{\sech}(\zeta)\right\|_{L^p(-\eta,\eta)}  \ge \left(\frac{\eta }{400}\right)^{\frac1p}|\widehat{\sech}(0)|.
\end{align*}
On the other hand, 
\begin{align*}
	\left\|\jap{\xi}^{s+1}\left( \widehat{\sech}\left(\frac{\xi-N_1}{\omega}\right)-\widehat{\sech}\left(\frac{\xi-N_2}{\omega}\right)\right)\right\|_{L^p} \le \eps\|\widehat{\sech}\|_{W^{1,1}}^{\frac1p}.
\end{align*}
In particular, for $\varepsilon$ small,
$$
	\left\| {\mathbf{b}_1^0}(T)  -{\mathbf{b}_2^0}(T) \right\|_{\FL^{s+1,p}}\ge  \left(\frac{\eta }{800}\right)^{\frac1p}|\widehat{\sech}(0)|=:4c.
$$
Hence, splitting into real and imaginary parts, either
$$
	\left\|b_1^0(T)  -b_2^0(T) \right\|_{\FL^{s+1,p}}\ge 2c
$$
and the proof is finished, or
$$
\left\|\omega \sech(\omega x + \gamma_1 T)\sin(N_1x+\delta_1T) - \omega \sech(\omega x + \gamma_2 T)\sin(N_2x+\delta_2T) \right\|_{\FL^{s+1,p}}\ge  2c
$$
and a simple computation then shows that
$$
	\left\|b_1^0(T')  -b_2^0(T') \right\|_{\FL^{s+1,p}}\ge c,\quad  \mbox{for }T'=T-\tfrac{\pi}{2\delta_1}, 
$$
completing the proof.\qedhere

\end{proof}

\begin{remark}
Lemma \ref{lem:breather} shows that the \eqref{mkdv} flow cannot be uniformly continuous in $\FL^{s,p}(\R)$, for $s<\frac{1}{2p}$. In particular, this shows that the well-posedness results of \cite{GrunrockVega} are sharp if one requires local uniform continuity of the flow.
\end{remark}

We now conclude the proof of Proposition \ref{prop:ill_intro}~(i) by applying the Miura map to the family of breathers considered above.

\begin{proof}[Proof of Proposition \ref{prop:ill_intro}~(i)]
	Given $T',\eps>0$, consider the breathers $b_1, b_2$ given by Lemma \ref{lem:breather} and set 
	$$u_j=\mathbf{M}[b_j]=\dx b_j + \tfrac13 b_j^2,\quad j=1,2.$$ Since $u_j\sim \jap{\dx} b_j$ for $\eps$ small, we conclude that
 	$$
	\|u_1(0)-u_2(0)\|_{\FL^{s,p}}\lesssim \eps, \quad \|u_j(0)\|_{\FL^{s,p}} \lesssim 1,\ j=1,2,
	$$
	and
	$$
	\|u_1(T')-u_2(T')\|_{\FL^{s,p}}\ge c.
	$$
	In particular, the \eqref{kdv} flow cannot be uniformly continuous in the $\FL^{s,p}(\R)$ topology, for $2\le p < \infty$ and $s<-1+\frac1{2p}$.\qedhere
\end{proof}

Finally, we prove Proposition \ref{prop:ill_intro}~(ii), which follows from an adaptation of the ill-posedness result in $H^s(\R)$, $s<-1$, shown by Molinet in \cite{Molinet_ill_line}. The argument is based on the fact that the Miura transform $\MMb[u]$ in \eqref{eq:defiMiura}, which
maps solutions to \eqref{mkdv} to solutions of \eqref{kdv},
is not continuous for $u\in H^s(\R)$, for any $s<0$.

\begin{proof}[Proof of Proposition~\ref{prop:ill_intro}~(ii)]

Since the proof is an adaptation of Molinet's proof for $L^2$-based spaces, we briefly recall the strategy in \cite{Molinet_ill_line}.
 First, construct a sequence of (real-valued) initial data $(u_{0,N})_{N\in\N}\subset H^\infty(\R)$ and $\theta\in H^{-\frac34+}(\R)$ such that
\begin{equation}\label{eq:molinet}
	u_{0,N}\rightharpoonup 0 \mbox{ in }L^2(\R)\ \text{for } s<0,\quad \text{with}
    \quad \|u_{0,N}\|_{L^2}\sim 1,
\end{equation}
\begin{equation}\label{eq:molinet2}
 \MMb[u_{0,N}] \to \theta\neq 0  \mbox{ strongly in }H^{s}(\R)
 \ \text{for }
 s<-1.
\end{equation}
From \eqref{eq:molinet}, since $u_{0,N}$ is smooth with bounded $L^2$ norm, one can use Kato smoothing to prove that, as $N\to \infty$, the corresponding solutions to (mKdV) $u_N$ converge to a distributional solution $u$ which vanishes for $t=0$ (see \cite[Proposition 2.1]{Molinet_ill_line}). This implies that the sequence $\MMb[u_N]$, which solves \eqref{kdv} with initial data $\MMb[u_{0,N}]$, converges to $\MMb[u]$ with zero initial data (see \cite[Proposition 2.2]{Molinet_ill_line}).

From \eqref{eq:molinet2}, define $v$ as the solution to \eqref{kdv} emanating from $\theta\neq 0$ (which exists by the standard local well-posedness theory since $\theta\in H^s(\R)$, $s>-\frac34$).
Then, for $s<-1$,
\begin{align}
\MMb[u_{0, N}] \to \theta \mbox{ in }H^s(\R),
\quad
\text{and}
\quad
\MMb[u_N] \not\to v \mbox{ in  }\mathcal{D}'((0,T)\times \R).
\end{align}
This shows that the data-to-solution map is not continuous in the $H^s$-topology.
To extend the failure of continuity to $\FL^{s,p}(\R)$, it suffices to replace \eqref{eq:molinet2} with an analogous convergence in this topology.
In particular, Proposition \ref{prop:ill_intro} follows from constructing a sequence $(u_{0,N})_{N\in \N} \subset H^\infty(\R)$ satisfying \eqref{eq:molinet} and
\begin{equation}\label{eq:molinet2_v2}
 \MMb[u_{0,N}] \to \theta\neq 0 \mbox{ strongly in }\mathcal{F}L^{s,p}(\R) \ \text{for } s+\tfrac1p<-\tfrac12.
\end{equation}

We devote the remaining of the proof to the construction of the sequence $(u_{0,N})_{N\in\N}$.
	For fixed $N$, take
	$$
	\widehat{u}_{0,N}(\xi)= N^{-\frac12} e^{i\xi^2/N}\mathbbm{1}_{[N, 2N]}(\xi) + N^{-\frac12} e^{-i\xi^2/N}\mathbbm{1}_{[-2N, -N]}(\xi).
	$$
	Then
	$$
	\|u_{0,N}\|_{\mathcal{F}L^{s,p}} \sim N^{-\frac12 +s+ \frac1p}
	$$
	which implies \eqref{eq:molinet} and
	$$
	\dx u_{0,N} \to 0 \mbox{ strongly in }H^{s}(\R), \ s+\tfrac1p<-\tfrac12.
	$$
	Let us now compute $u_{0,N}^2$. Since $u_{0,N}$ is real-valued, it suffices to consider $\xi>0$:
	\begin{align*}
		\widehat{(u_{0,N}^2)}(\xi) &= \frac{2}{N}\int_{N+\xi}^{2N} e^{i\frac{\xi_1^2}{N}-i\frac{(\xi-\xi_1)^2}{N}}d\xi_1 + \frac{1}{N}\int_{\max\{N, \xi-2N\}}^{\min\{2N, \xi-N\}}  e^{i\frac{\xi_1^2}{N}+i\frac{(\xi-\xi_1)^2}{N}}d\xi_1\\&=2\int_{1+\frac\xi N}^{2} e^{iN\eta_1^2-i\frac{(\xi-N\eta_1)^2}{N}}d\eta_1 + \int_{\max\{1, \frac{\xi}{N}-2\}}^{\min\{2, \frac{\xi}{N}-1\}}  e^{iN\eta_1^2+i\frac{(\xi-N\eta_1)^2}{N}}d\eta_1
        .
	\end{align*}

    \noi
The second integral converges pointwise to $0$ and therefore in $\mathcal{F}L^{s,p}(\R)$, as $sp<-1$. For the first integral,  for $\xi$ fixed,
\begin{equation}
	\int_{1+\frac\xi N}^{2} e^{iN\eta_1^2-i\frac{(\xi-N\eta_1)^2}{N}}d\eta_1 = e^{-i\frac{\xi^2}{N}}\int_{1 + \frac{\xi}{N}}^{2} e^{2i\eta_1\xi}d\eta_1  = \frac{e^{4i\xi}-e^{2i\left(1+\frac{\xi}{N}\right)\xi}}{2i\xi}
\end{equation}
and thus, for $s+\frac1p<-\frac12$,
\begin{equation}
	\MMb[u_{0,N}]= \dx u_{0,N} + \tfrac13 (u_{0,N})^2 \to \frac{e^{4i\xi}-e^{2i\xi}}{3i\xi} =: \theta \ne 0
    \quad \mbox{in }\mathcal{F}L^{s,p}(\R),
\end{equation}
showing \eqref{eq:molinet2_v2}, as required.
\end{proof}

\subsection{A remark on well-posedness of (dgBO)}\label{subsec:gbo}

In this section, we describe the necessary modifications on the \eqref{kdv} well-posedness theory in order to encompass \eqref{gBO}.
Our goal here is not to provide any (optimal) result, but rather to stress that the methods developed throughout this work can be extended to other nonintegrable dispersive models. 

\medskip

\noi\textbf{Normal form reduction and gauged equation.} 
Due to the quadratic nonlinearity, the infinite normal form reduction for \eqref{kdv} and the gauge transform in Section~\ref{sec:infr} extend to \eqref{gBO}. 
In particular, we have an analogue of the gauged equation \eqref{z1}, where the only change in the operators $\NN_\ell$, $\ell=1,\dots, 4$, comes from the definition of $K_\star$ and $\Psi(\star)$ appearing in \eqref{wint0c}. In the \eqref{gBO} case, these must be replaced with
\begin{align}
\Psi^{(a)} (\star) 
&
:  = - |\xi_\star|^{a}\xi_\star + |\xi_{\star1}|^{a}\xi_{\star1} + |\xi_{\star2}|^{a}\xi_{\star2}
\quad 
\text{and}
\quad
K^{(a)}_\star 
 : = \frac{\xi_\star}{\Psi^{(a)}(\star)}. 
\end{align}
A direct computation shows that for $1<a<2$
\begin{align}
|\Psi^{(a)}(\star) | \ges (\max(|\xi_\star|, |\xi_{\star1}|, |\xi_{\star2}|)^a \min(|\xi_\star|, |\xi_{\star1}|, |\xi_{\star2}|), 
\end{align}
see \cite[Lemma 3.8]{herr_bo} for a proof. In particular, over the set $A_\star^c$,
\begin{equation}\label{eq:Kalpha}
|K_\star^{(a)}|\lesssim \frac{1}{\jap{\xi_{\star1}}^{a-1}\jap{\xi_{\star2}}}, \quad \text{for }|\xi_{\star1}|\ge |\xi_{\star2}|,
\end{equation}
and therefore $K_\star^{(a)}\ind_{A_\star^c}$ remains a bounded multiplier.

\medskip

\noi
\textbf{Multilinear estimates.} The rather cumbersome adaptation of the \eqref{kdv} theory to the context of \eqref{gBO} arises at this stage, where one must check that the frequency-restricted estimates proven in Section \ref{sec:multi} extend naturally for $1<a<2$. To that end, one must incorporate the following facts:
\begin{itemize}
\item (Multiplier estimate) The multiplier $K_\star^{(a)}\ind_{A_\star^c}$ is controlled as in \eqref{eq:Kalpha}.
\item (Integration over $\TT^\infty_\star$) For $p=\infty$, the analogue of \eqref{eq:estMlinha} is then given by
\begin{equation}
    \int \mathcal{M}_2(\xi_{\bul\in \TT^\infty_\star})  d\xi_{\bul\in\TT^\infty_\star} \lesssim \bigg(
\prod_{\bul \in \TT^\infty_\star }
\int \frac{1}{ \jb{\xi_{\bul}}^{s+a-1}\jb{\xi_{\sbf(\bul)}}} d\xi_{\bul} \bigg)  \le C^j, 
\end{equation}
which imposes the condition $s>1-a$. For $p=2$, the correct adaptation of \eqref{eq:int_M2} is
\begin{align}
 \int_{\G(\TT)} \M_2(\vec\xi) 
\ind_{|\Psi(\TT)-\al| < M}
d\xi_{\bul\in\TT^\infty_\star} 
&
\les 
\int
\prod_{\bul\in \TT^\infty_\star} \frac{1}{\jb{\xi_\bul}^{2(s+a-1)} \jb{\xi_{\sbf(\bul)}}^2 } d\xi_{\bul\in\TT^\infty_\star} 
\les C^j,
\end{align}
which holds again for $s>1-a$.

\item (Stationary vs.  nonstationary) The resonance function $\Psi^{(a)}(\TT)$ has the same stationary points as in the KdV case $a=2$.
\item (Nonstationary analysis) Over the nonstationary regions, the change of variables which replaces the nonstationary frequency with $\phi$ yields a jacobian of order $a$.

For instance, \eqref{IIIb-new} is replaced by the following observation: for $|\xi_1| \not\simeq |\xi_2|$ and $|\xi_1 | \ge |\xi_2|$, we have
\begin{align}
|\partial_{\xi_1} \Psi^{(a)}(\xi, \xi_1,\xi-\xi_1) |  \ges_a \big| |\xi_1|^a - |\xi_2|^a \big| 
= |\xi_1|^a \big( 1 - \tfrac{|\xi_2|^a}{|\xi_1|^a}) \big)
\ges |\xi_1|^a .
\end{align}
Note that the lower bound holds since $|\xi_1| - |\xi_2| \ge \frac{1}{C} |\xi_1| $ for some $C>1$, which implies that $1 - \frac{|\xi_1|^a}{|\xi_2|^a} \ge 1 - (1-\frac1C)^a >0$.

\item (Stationary analysis) In the stationary case, the renormalization procedure yields a factor of order $a+1$. For example, the analogue of \eqref{Phi1} is
$$
\Phi(\xi,\xi_{11},\xi_2) = |\xi|^a\xi\Phi(1, p_{11}, p_2).
$$
After the renormalization, the application of Morse's Lemma is identical to the $a=2$ case.
\end{itemize}
Incorporating these facts into the proof of Propositions \ref{prop:bourgLinfty}-\ref{prop:bourgL2}, 
we find that, for some $s_p(a)$, depending continuously on $a$ and satisfying
$$
s_p(2)=\max\left\{-\frac12-\frac1p, -1+\frac1{2p}\right\},
$$
the multilinear estimate \eqref{eq:bourg_bi_lp} holds for any $s>s_p(a)$, $\TT\in \Tr_j$, $j\ge 2$ and $\ell=2,3,4$.

However, the Type $\I$ term imposes a more serious restriction, as can be seen in \eqref{eq:estL21}. Indeed, in the \eqref{gBO} case, a similar computation would give
\begin{align}
\sup_\xi \int_{\xi=\xi_1+\xi_2} \left(\frac{\jb{\xi}^{s} |\xi| }{\jb{\xi_1}^{s}\jb{\xi_2}^{s}}
    \ind_A \ind_{|\Psi(\TT) - \al| < M } \right)^{2} d\xi_1
    \les \sup_\xi \int |\xi|^2 \ind_{|\phi - \al| < M } \frac{d\phi}{|\xi|^a }
\end{align}
    which is unbounded for $a<2$. This obstruction is exactly the reason why in \cite{herr_bo}, Herr
    needed to include the homogeneous space $\dot H^{\frac12 - \frac1a}(\R)$.
    This difficulty can be bypassed by establishing the estimate for the bilinear term of Type $\I$ directly in $X^{s, b}_p$-spaces for $2<p<\infty$, without arguing via interpolation. In that case, we replace \eqref{eq:estL21} with
    \begin{align}
    \sup_\xi \left(\int_{\xi=\xi_1+\xi_2} \left(\frac{\jb{\xi}^{s} |\xi| }{\jb{\xi_1}^{s}\jb{\xi_2}^{s}}
    \ind_A \ind_{|\Psi(\TT) - \al| < M } \right)^{p'} d\xi_1\right)^{\frac{1}{p'}}
    &\les \sup_\xi \left(\int |\xi|^{p'} \ind_{|\phi - \al| < M } \frac{d\phi}{|\xi|^a }\right)^{\frac1{p'}} \\&\lesssim M^{\frac1{p'}},
\end{align}
which now holds for $p'\le a$, that is, $p\ge\frac{a}{a-1}$. The case $p=\infty$ is excluded, as the estimate yields a full power of $M$.

\medskip

\noi\textbf{Well-posedness for \eqref{gBO}.} Having in hand the multilinear estimates, which hold for
\begin{equation}
s>s_p(a), \quad \tfrac{a}{a-1}\le p <\infty,\label{s_bo}
\end{equation}
the local existence for the gauged \eqref{gBO} equation follows exactly the same steps as in the proof of Theorem \ref{thm:lwp_rkdv}. 

All that is left to do is to invert the gauge transform and build the solution to \eqref{gBO} via approximation by smooth solutions.
For the former, using the estimate \eqref{eq:Kalpha}, the conclusions of Lemma \ref{LEM:D} hold for $s>\frac{1-a}{2}-\frac{1}{p}$  (which plays the role of the $H^{-1}$-scale for \eqref{kdv}). 
Finally, for the density argument, we point out that the proof of  Proposition \ref{PRO:equiv} hinges solely on the multilinear estimates proven for the gauged equation and thus it is immediately extendable to the \eqref{gBO} case. We conclude that well-posedness for \eqref{gBO}  holds for
$$
s>\max\left\{s_p(a),\frac{1-a}{2}-\frac{1}{p} \right\},\quad \tfrac{a}{a-1}\le p <\infty, 
$$
where $s_p(a)$ is as in \eqref{s_bo}.

\begin{remark}
The condition  $p\ge \frac{a}{a-1}$, is imposed by the high$\times$low$\to$high interactions also responsible for the obstruction in \cite{herr_bo}, independently of the regularity parameter $s$. This is also linked to the fact that, for $1 \le a<2$, the \eqref{gBO} flow behaves \emph{quasilinearly} and cannot be $C^2$ in $H^s(\R)$, for any $s\in \R$ \cite{KochTzvetkov_bo, MolTzeSaut_ill}. 
The endpoint case $a=1$ corresponding to BO is excluded, providing yet another indicator of the strong quasilinearity of the Benjamin-Ono equation.
For a fixed $1<a<2$, despite the negative results in \cite{KochTzvetkov_bo, MolTzeSaut_ill} at the $H^s$-scale, our argument shows the analytic well-posedness of \eqref{gBO} in $\FL^{s,p}(\R)$ for sufficiently large $s$ and $p$.
Therefore, we conclude that there exists an interplay between the integrability $p$ and the semilinear/quasilinear nature of \eqref{gBO}.
\end{remark}

\appendix

\section{Well-posedness for KdV in $\F L^{s,\infty}(\R)$, $s>-\frac14$}\label{sec:-14}

In this section, we prove the well-posedness of \eqref{kdv}
 in $\FL^{s, \infty}(\R)$ for $s>-\frac14$, complementing the results of \cite{Grunrock_ldv_hier}, without resorting to the gauge transform.
The proof follows from a standard fixed point argument in the (time-localized) space $X^{s, b}_{\infty, \infty} (0,T) \cap X^{s, \wt b }_{\infty, 1} (0,T)$ (see \eqref{Xsb}) for some  suitable $b,\wt b\in \R$,
together with the following bilinear estimate.

\begin{lemma}\label{LEM:bil0}
Let $s>-\frac14$. Then, there exist $b,\wt b >0 $
 such that
\begin{align}
\label{bilinear0}
\left\|\int_0^te^{-(t-t')\partial_x^3}\partial_x(u_1u_2) (t')dt' \right\|_{X^{s,b}_{\infty, \infty}\cap X^{s,\wt{b}}_{\infty, 1} } & \les
\|u_1\|_{X^{s,b}_{\infty, \infty}\cap X^{s,\wt{b}}_{\infty, 1} }
\|u_2\|_{X^{s,b}_{\infty, \infty}\cap X^{s,\wt{b}}_{\infty, 1} }
\end{align}
\end{lemma}

\begin{proof}

For simplicity, we consider only the case of $s<0$.
From Lemma~\ref{lem:reduction_fre_k2} and Lemma~\ref{lem:eta1}, \eqref{bilinear0} follows once we show the frequency-restricted estimate: for $|\xi_1| \ge |\xi_2|$ (by symmetry) and $\Phi$ is as in \eqref{Phi} with $k=2$, then we want to show that
\begin{equation}
\label{est1}
\sup_{\xi} \int_{\xi=\xi_1+\xi_2}
\frac{|\xi|\jap{\xi}^s}{\jap{\xi_1}^{s}\jap{\xi_2}^s}
\fia d\xi_1
\lesssim
\jap{\al}^{\ld } M^{\be}
\end{equation}
for any $M\ge 1$ and $\al\in\R$, and some $\be,\ld \ge 0$ with $\be+\ld < 1$. More precisely, we will show \eqref{est1} for
$$
\beta=\tfrac34,\quad \lambda=\tfrac{1}{6}-\tfrac s3+.
$$

If $|\xi_1| \le 1$, then all frequencies are $\les1$, and the l.h.s. of \eqref{est1} is trivially bounded. Henceforth, assume that $|\xi_1| \ge 1$.
For the remainder of the proof, fix $\xi$. We begin with the case $|\alpha|\gg M$.

\noi
\underline{\textbf{Case 1:} $|\xi_1^2-\xi_2^2| \gtrsim |\xi_1^2|$.}

In this region, we have $|\xi_1|\sim |\xi|$. If $|\xi_2|\ge 1$, by Cauchy-Schwarz, with the change of variables $\xi_2\mapsto \phi=\Phi$, noting that $|\partial_{\xi_2} \Phi|\sim |\xi_1^2 - \xi_2^2| \ges |\xi|^2$, and using $|\xi_2|^{1-2s+}\lesssim |\Phi|^{\frac{1-2s}{3}+}$, we get
\begin{align*}
\text{l.h.s \eqref{est1}}
& \les \left(\int |\xi|^2\jap{\xi_2}^{1-2s+}\fia d\xi_2\right)^\frac12\\&\lesssim \left(\int |\phi|^{\frac{1-2s}{3}+}\mathbbm{1}_{|\phi-\alpha|<M} d\phi\right)^\frac12\lesssim |\alpha|^{\frac16 - \frac s3 +} M^{\frac12}.
\end{align*}

\noi
If $|\xi_2|<1$, we proceed by Cauchy-Schwarz and the change of variables $\xi_2 \mapsto \phi= \Phi$, as above:
\begin{align*}
\text{l.h.s \eqref{est1}}
\les
\left(\int |\xi|^2\fia d\xi_2\right)^\frac12\lesssim \left(\int \mathbbm{1}_{|\phi-\alpha|<M} d\phi\right)^\frac12\lesssim M^{\frac12}.
\end{align*}

\smallskip

\noi
\underline{\textbf{Case 2:} $|\xi| \ll |\xi_1|$}

In this region, we still perform the change of variables $\xi_2\mapsto \phi=\Phi$, since $|\partial_{\xi_2}\Phi|\gtrsim |\xi\xi_1|$. If $|\xi|>1$, since $|\xi|^{1+2s}|\xi_1|^{-4s+}\lesssim |\Phi|^{\frac{1-2s}{3}+} $,  for $s>-\frac12$,  we have
\begin{align*}
\text{l.h.s. \eqref{est1}}
\lesssim \left(\int |\xi|^{2+2s}|\xi_1|^{1-4s+}\fia d\xi_2\right)^\frac12
& \lesssim \left(\int |\phi|^{\frac{1-2s}{3}+}\fia d\phi\right)^\frac12
\\
&
\lesssim |\alpha|^{\frac16 - \frac s3 +} M^{\frac12}.
\end{align*}
For $|\xi|<1$, using $|\xi||\xi_1|^{-4s}\lesssim |\Phi|^{-2s}$,  for $s>-\frac12$, we get
\begin{align*}
    \text{l.h.s. \eqref{est1}}
    \les
    \left(\int |\xi|^{2}|\xi_1|^{1-4s+}\fia d\xi_2\right)^\frac12
    \lesssim \left(\int |\phi|^{-2s}\fia d\phi\right)^\frac12\lesssim |\alpha|^{-s} M^{\frac12}.
\end{align*}

\smallskip

\noi
\underline{\textbf{Case 3:}
$|\xi|\sim|\xi_1| \sim |\xi_2| \gg |\xi_1 - \xi_2|$}

Due to the smallness of $|\xi_1-\xi_2|$,
the earlier change of variables $\xi_2\mapsto \phi = \Phi$ is not useful here.
Instead,
we make a change of variables around the critical point $\xi_1 = \frac\xi2$ of $\Phi$. More precisely, let $\xi_1 \mapsto q = \frac{\xi_1}{\xi} - \frac12$. Then, $\Phi (\xi, \xi_1,\xi-\xi_1)=  -3 \xi^3 (\frac14 - q^2)$ and
\begin{align*}
\text{l.h.s. \eqref{est1}}
\les
|\xi|^{2-s}\int_{|q| \ll 1} \mathbbm{1}_{|-3 \xi^3 (\frac14 - q^2)-\alpha|<M}dq \les |\xi|^{\frac12-s +} M^\frac12  \ind_{|\xi|^3 \les |\al|}
\les |\alpha|^{\frac16-\frac s3+} M^\frac12.
\end{align*}
where we used \eqref{eq:quadratica}.
This concludes the case $|\alpha|\gg M$.

\smallskip

For $|\alpha|\lesssim M$, we may replace the restriction
$$
\fia \quad \mbox{by }\quad \mathbbm{1}_{|\Phi|<M}.
$$

\smallskip
\noi
\underline{\textbf{Case 1:} $|\xi|\gtrsim |\xi_1|$}

If $|\xi_2|>1$, then
$
|\xi||\xi_2|^{-s+1+}\lesssim |\Phi|^{\frac{2-s}{3}+} \lesssim M^{\frac{2-s}{3}+}
$
and thus
\begin{align}
\text{l.h.s. \eqref{est1}}
\les \int |\xi||\xi_2|^{-s+1+}\mathbbm{1}_{|\Phi|<M}\frac{d\xi_2}{|\xi_2|^{1+}}\lesssim M^{\frac{2-s}{3}+}.
\end{align}
If $|\xi_2|<1$, since $|\xi_2|\lesssim |\Phi||\xi|^{-2}\lesssim M|\xi|^{-2}$,
\begin{align}
\text{l.h.s. \eqref{est1}}
\les \int |\xi|\mathbbm{1}_{|\Phi|<M}d\xi_2\lesssim \left(\int |\xi|^2\mathbbm{1}_{|\Phi|<M} d\xi_2\right)^\frac12\lesssim M^{\frac12}.
\end{align}

\smallskip\noi
\underline{\textbf{Case 2:} $|\xi|\ll |\xi_1|\simeq |\xi_2|$}

We bound
$ |\xi|^{1+s}|\xi_1|^{-2s+1+}\lesssim |\Phi|^{\frac{2-s}{3}+}\lesssim M^{\frac{2-s}{3}+}
$,
for $s>-\frac14$,
which implies that
\begin{align}
\text{l.h.s. \eqref{est1}}
&\les \int |\xi|\jap{\xi}^s|\xi_1|^{-2s+1+}\mathbbm{1}_{|\Phi|<M}\frac{d\xi_1}{|\xi_1|^{1+}}\lesssim  M^{\frac{2-s}{3}+}.
\end{align}

\medskip

Gathering the estimates from the various regions, we find that
\begin{align}
\text{l.h.s. \eqref{est1}}
\lesssim |\alpha|^{\frac16-\frac s3+}M^{\frac12} + M^{\frac12} + |\alpha|^{-s}M^{\frac12} +  M^{\frac{2-s}{3}+} \lesssim \jap{\alpha}^{\frac16-\frac s3+}M^{\frac34}
\end{align}
for $s>-\frac14$. The proof is finished.\qedhere

\end{proof}

\begin{remark}
    The restriction $s>-\frac14$ comes from the requirement that the parameters $\beta,\lambda$ in \eqref{est1} must be chosen uniformly for the different frequency domains. When the stationary region $|\xi_1|\simeq |\xi_2|$ is not present (corresponding to the Type I term $\Ncal_1$ in Proposition \ref{prop:bourgLinfty}), one can choose $\lambda=0$ and $\beta=1-$ as long as $s>-1$.
    This also shows that the stationary regions are responsible not only for the above obstruction $s>-\frac14$ in the proof of the bilinear estimate, but also for the more serious well-posedness obstruction $s\ge -\frac12$.

\end{remark}

\bibliographystyle{plain}

\bibliography{biblio}

\begin{thebibliography}{10}

\bibitem{adams}
Joseph Adams.
\newblock Well-posedness for the {NLS} hierarchy.
\newblock {\em J. Evol. Equ.}, 24(4):Paper No. 88, 52, 2024.

\bibitem{AiLiu}
Albert Ai and Grace Liu.
\newblock The dispersion generalized {B}enjamin-{O}no equation.
\newblock {\em Preprint, arXiv:2407.01472}, 2024.

\bibitem{anjolras}
Phillipe Anjolras.
\newblock Scattering of the 2{D} modified {Z}akharov-{K}uznetsov equation.
\newblock {\em Discrete Contin. Dyn. Syst.}, 51:542--565, 2026.

\bibitem{bit}
Anatoli~V. Babin, Alexei~A. Ilyin, and Edriss~S. Titi.
\newblock On the regularization mechanism for the periodic {K}orteweg-de
  {V}ries equation.
\newblock {\em Comm. Pure Appl. Math.}, 64(5):591--648, 2011.

\bibitem{Bambusi_birkhoff}
Dario Bambusi.
\newblock Birkhoff normal form for some nonlinear {PDE}s.
\newblock {\em Comm. Math. Phys.}, 234(2):253--285, 2003.

\bibitem{taobejenaru}
Ioan Bejenaru and Terence Tao.
\newblock Sharp well-posedness and ill-posedness results for a quadratic
  non-linear {S}chr\"{o}dinger equation.
\newblock {\em J. Funct. Anal.}, 233(1):228--259, 2006.

\bibitem{BonaSmith}
Jerry~L. Bona and Ronald Smith.
\newblock The initial-value problem for the {K}orteweg-de {V}ries equation.
\newblock {\em Philos. Trans. Roy. Soc. London Ser. A}, 278(1287):555--601,
  1975.

\bibitem{bourg1}
Jean Bourgain.
\newblock Fourier transform restriction phenomena for certain lattice subsets
  and applications to nonlinear evolution equations. {I}. {S}chr\"{o}dinger
  equations.
\newblock {\em Geom. Funct. Anal.}, 3(2):107--156, 1993.

\bibitem{bourg2}
Jean Bourgain.
\newblock Fourier transform restriction phenomena for certain lattice subsets
  and applications to nonlinear evolution equations. {II}. {T}he
  {K}d{V}-equation.
\newblock {\em Geom. Funct. Anal.}, 3(3):209--262, 1993.

\bibitem{BO94}
Jean Bourgain.
\newblock Periodic nonlinear {S}chr\"odinger equation and invariant measures.
\newblock {\em Comm. Math. Phys.}, 166(1):1--26, 1994.

\bibitem{bou_ill_kdv}
Jean Bourgain.
\newblock Periodic {K}orteweg de {V}ries equation with measures as initial
  data.
\newblock {\em Selecta Math. (N.S.)}, 3(2):115--159, 1997.

\bibitem{bourg_birkhoff_normal}
Jean Bourgain.
\newblock A remark on normal forms and the ``{$I$}-method'' for periodic {NLS}.
\newblock {\em J. Anal. Math.}, 94:125--157, 2004.

\bibitem{Boussinesq}
Joseph Boussinesq.
\newblock Th\'eorie des ondes et des remous qui se propagent le long d'un canal
  rectangulaire horizontal, en communiquant au liquide contenu dans ce canal
  des vitesses sensiblement pareilles de la surface au fond.
\newblock {\em J. Math. Pures Appl. (2)}, 17:55--108, 1872.

\bibitem{Bruned}
Yvain Bruned.
\newblock Derivation of normal forms for dispersive {PDE}s via arborification.
\newblock {\em Preprint, arXiv:2409.03642v2}, 2024.

\bibitem{BuckmasterKoch}
Tristan Buckmaster and Herbert Koch.
\newblock The {K}orteweg--de {V}ries equation at {$H^{-1}$} regularity.
\newblock {\em Ann. Inst. H. Poincar\'e{} C Anal. Non Lin\'eaire},
  32(5):1071--1098, 2015.

\bibitem{CCR_particle}
Mattia Cafasso, Tom Claeys, and Giulio Ruzza.
\newblock Airy kernel determinant solutions to the {K}d{V} equation and
  integro-differential {P}ainlev\'e{} equations.
\newblock {\em Comm. Math. Phys.}, 386(2):1107--1153, 2021.

\bibitem{CamposLinaresSantos}
Luccas Campos, Felipe Linares, and Thyago S.~R. Santos.
\newblock Sharp well-posedness for the $k$-dispersion generalized
  {B}enjamin-{O}no equations: {S}hort and long time results.
\newblock {\em Preprint, arxiv:2410.17217}, 2024.

\bibitem{Chapouto_FL_mkdv}
Andreia Chapouto.
\newblock A remark on the well-posedness of the modified {KDV} equation in the
  {F}ourier-{L}ebesgue spaces.
\newblock {\em Discrete Contin. Dyn. Syst.}, 41(8):3915--3950, 2021.

\bibitem{CCR-NV}
Andreia Chapouto, Sim\~ao Correia, and Jo\~ao~Pedro Ramos.
\newblock Sharp well-posedness for {Novikov-Veselov} equation.
\newblock {\em In preparation}.

\bibitem{CKV}
Andreia Chapouto, Rowan Killip, and Monica Vi\c{s}an.
\newblock Bounded solutions of {K}d{V}: uniqueness and the loss of almost
  periodicity.
\newblock {\em Duke Math. J.}, 173(7):1227--1267, 2024.

\bibitem{Christ_power}
Michael Christ.
\newblock Power series solution of a nonlinear {S}chr\"odinger equation.
\newblock In {\em Mathematical aspects of nonlinear dispersive equations},
  volume 163 of {\em Ann. of Math. Stud.}, pages 131--155. Princeton Univ.
  Press, Princeton, NJ, 2007.

\bibitem{CCT_ill}
Michael Christ, James Colliander, and Terrence Tao.
\newblock Asymptotics, frequency modulation, and low regularity ill-posedness
  for canonical defocusing equations.
\newblock {\em Amer. J. Math.}, 125(6):1235--1293, 2003.

\bibitem{CGKO_normal_form}
Jaywan Chung, Zihua Guo, Soonsik Kwon, and Tadahiro Oh.
\newblock Normal form approach to global well-posedness of the quadratic
  derivative nonlinear {S}chr\"odinger equation on the circle.
\newblock {\em Ann. Inst. H. Poincar\'e{} C Anal. Non Lin\'eaire},
  34(5):1273--1297, 2017.

\bibitem{CKSTT_dnls}
James Colliander, Markus Keel, Gigliola Staffilani, Hideo Takaoka, and Terence
  Tao.
\newblock Global well-posedness for {S}chr\"odinger equations with derivative.
\newblock {\em SIAM J. Math. Anal.}, 33(3):649--669, 2001.

\bibitem{CKSTT_sharp_kdv}
James Colliander, Markus Keel, Gigliola Staffilani, Hideo Takaoka, and Terence
  Tao.
\newblock Sharp global well-posedness for {K}d{V} and modified {K}d{V} on
  {$\Bbb R$} and {$\Bbb T$}.
\newblock {\em J. Amer. Math. Soc.}, 16(3):705--749, 2003.

\bibitem{CollianderKwonOh_birkhoff}
James Colliander, Soonsik Kwon, and Tadahiro Oh.
\newblock A remark on normal forms and the ``upside-down'' {$I$}-method for
  periodic {NLS}: growth of higher {S}obolev norms.
\newblock {\em J. Anal. Math.}, 118(1):55--82, 2012.

\bibitem{COS}
Sim\~ao Correia, Filipe Oliveira, and Jorge~Drumond Silva.
\newblock Sharp local well-posedness and nonlinear smoothing for dispersive
  equations through frequency-restricted estimates.
\newblock {\em SIAM J. Math. Anal.}, 56(4):5604--5633, 2024.

\bibitem{CS20}
Sim\~ao Correia and Jorge~Drumond Silva.
\newblock Nonlinear smoothing for dispersive {PDE}: a unified approach.
\newblock {\em J. Differential Equations}, 269(5):4253--4285, 2020.

\bibitem{CC24}
Simão Correia and Raphaël Côte.
\newblock Sharp blow-up stability for self-similar solutions of the modified
  {K}orteweg-de {V}ries equation.
\newblock {\em Preprint, arXiv:2402.16423}, 2024.

\bibitem{CLS}
Simão Correia, Felipe Linares, and Jorge Silva.
\newblock Sharp local well-posedness for the {S}chrödinger-{K}orteweg-de
  {V}ries system.
\newblock {\em Math. Res. Lett.}, 32:1791--1844, 2025.

\bibitem{Crighton}
David~G. Crighton.
\newblock Applications of {K}d{V}.
\newblock volume~39, pages 39--67. 1995.
\newblock KdV '95 (Amsterdam, 1995).

\bibitem{ErdoganTzirakis_kdv}
Mehmet~B. Erdögan and Nikolaos Tzirakis.
\newblock Global smoothing for the periodic {K}d{V} evolution.
\newblock {\em Int. Math. Res. Not. IMRN}, (20):4589--4614, 2013.

\bibitem{ErdoganTzirakis_nls}
Mehmet~B. Erdögan and Nikolaos Tzirakis.
\newblock Talbot effect for the cubic non-linear {S}chr\"odinger equation on
  the torus.
\newblock {\em Math. Res. Lett.}, 20(6):1081--1090, 2013.

\bibitem{FarahLinaresPastor_bo}
Luiz~G. Farah, Felipe Linares, and Ademir Pastor.
\newblock Global well-posedness for the {$k$}-dispersion generalized
  {B}enjamin-{O}no equation.
\newblock {\em Differential Integral Equations}, 27(7-8):601--612, 2014.

\bibitem{Fordy}
Allan~P. Fordy.
\newblock {\em Soliton theory: a survey of results}.
\newblock Manchester Univ. P., 1990.

\bibitem{Forlano}
Justin Forlano.
\newblock {\em On the deterministic and probabilistic Cauchy problem of
  nonlinear dispersive partial differential equations}.
\newblock Phd thesis, the University of Edinburgh, Edinburgh, UK, December
  2020.

\bibitem{GermainMasmoudi_euler_maxwell}
Pierre Germain and Nader Masmoudi.
\newblock Global existence for the {E}uler-{M}axwell system.
\newblock {\em Ann. Sci. \'Ec. Norm. Sup\'er. (4)}, 47(3):469--503, 2014.

\bibitem{GMS_nls_3d}
Pierre Germain, Nader Masmoudi, and Jalal Shatah.
\newblock Global solutions for 3{D} quadratic {S}chr\"odinger equations.
\newblock {\em Int. Math. Res. Not. IMRN}, (3):414--432, 2009.

\bibitem{GMS_gravity}
Pierre Germain, Nader Masmoudi, and Jalal Shatah.
\newblock Global solutions for the gravity water waves equation in dimension 3.
\newblock {\em Ann. of Math. (2)}, 175(2):691--754, 2012.

\bibitem{Ginibre_interaction_rep}
Jean Ginibre.
\newblock An introduction to nonlinear {S}chr\"odinger equations.
\newblock In {\em Nonlinear waves ({S}apporo, 1995)}, volume~10 of {\em GAKUTO
  Internat. Ser. Math. Sci. Appl.}, pages 85--133. Gakkotosho, Tokyo, 1997.

\bibitem{Grunrock_mkdv}
Axel Gr\"unrock.
\newblock An improved local well-posedness result for the modified {K}d{V}
  equation.
\newblock {\em Int. Math. Res. Not.}, (61):3287--3308, 2004.

\bibitem{Grunrock_ldv_hier}
Axel Gr\"unrock.
\newblock On the hierarchies of higher order m{K}d{V} and {K}d{V} equations.
\newblock {\em Cent. Eur. J. Math.}, 8(3):500--536, 2010.

\bibitem{HerrGrunrock_FL}
Axel Gr\"unrock and Sebastian Herr.
\newblock Low regularity local well-posedness of the derivative nonlinear
  {S}chr\"odinger equation with periodic initial data.
\newblock {\em SIAM J. Math. Anal.}, 39(6):1890--1920, 2008.

\bibitem{GrunrockVega}
Axel Gr\"unrock and Luis Vega.
\newblock Local well-posedness for the modified {K}d{V} equation in almost
  critical {$\widehat{H^r_s}$}-spaces.
\newblock {\em Trans. Amer. Math. Soc.}, 361(11):5681--5694, 2009.

\bibitem{Guo_gwp}
Zihua Guo.
\newblock Global well-posedness of {K}orteweg-de {V}ries equation in
  {$H^{-3/4}(\Bbb R)$}.
\newblock {\em J. Math. Pures Appl. (9)}, 91(6):583--597, 2009.

\bibitem{GuoKwonOh_nls}
Zihua Guo, Soonsik Kwon, and Tadahiro Oh.
\newblock Poincar\'e-{D}ulac normal form reduction for unconditional
  well-posedness of the periodic cubic {NLS}.
\newblock {\em Comm. Math. Phys.}, 322(1):19--48, 2013.

\bibitem{hadac}
Martin Hadac.
\newblock Well-posedness for the {K}adomtsev-{P}etviashvili {II} equation and
  generalisations.
\newblock {\em Trans. Amer. Math. Soc.}, 360(12):6555--6572, 2008.

\bibitem{HHK_kpii}
Martin Hadac, Sebastian Herr, and Herbert Koch.
\newblock Well-posedness and scattering for the {KP}-{II} equation in a
  critical space.
\newblock {\em Ann. Inst. H. Poincar\'e{} C Anal. Non Lin\'eaire},
  26(3):917--941, 2009.

\bibitem{Herr_dnls}
Sebastian Herr.
\newblock On the {C}auchy problem for the derivative nonlinear {S}chr\"odinger
  equation with periodic boundary condition.
\newblock {\em Int. Math. Res. Not.}, pages Art. ID 96763, 33, 2006.

\bibitem{herr_bo}
Sebastian Herr.
\newblock Well-posedness for equations of {B}enjamin-{O}no type.
\newblock {\em Illinois J. Math.}, 51(3):951--976, 2007.

\bibitem{Hormander}
Lars H\"ormander.
\newblock {\em The analysis of linear partial differential operators. {III}}.
\newblock Classics in Mathematics. Springer, Berlin, 2007.
\newblock Pseudo-differential operators, Reprint of the 1994 edition.

\bibitem{IfrimTataru}
Mihaela Ifrim and Daniel Tataru.
\newblock Well-posedness and dispersive decay of small data solutions for the
  {B}enjamin-{O}no equation.
\newblock {\em Ann. Sci. \'Ec. Norm. Sup\'er. (4)}, 52(2):297--335, 2019.

\bibitem{Jockel}
Niklas J\"ockel.
\newblock Well-posedness of the periodic dispersion-generalized
  {B}enjamin-{O}no equation in the weakly dispersive regime.
\newblock {\em Nonlinearity}, 37(8):Paper No. 085002, 37, 2024.

\bibitem{Kappeler_rough}
Thomas Kappeler.
\newblock Solutions to the {K}orteweg-de {V}ries equation with irregular
  initial profile.
\newblock {\em Comm. Partial Differential Equations}, 11(9):927--945, 1986.

\bibitem{KMMT16}
Thomas Kappeler, Alberto Maspero, Jan Molnar, and Peter Topalov.
\newblock On the convexity of the {K}d{V} {H}amiltonian.
\newblock {\em Comm. Math. Phys.}, 346(1):191--236, 2016.

\bibitem{KappelerMolnar18}
Thomas Kappeler and Jan Molnar.
\newblock On the wellposedness of the {K}d{V} equation on the space of
  pseudomeasures.
\newblock {\em Selecta Math. (N.S.)}, 24(2):1479--1526, 2018.

\bibitem{KPST_miura}
Thomas Kappeler, Peter Perry, Mikhail Shubin, and Peter Topalov.
\newblock The {M}iura map on the line.
\newblock {\em Int. Math. Res. Not.}, (50):3091--3133, 2005.

\bibitem{KT06}
Thomas Kappeler and Peter Topalov.
\newblock Global wellposedness of {K}d{V} in {$H^{-1}(\Bbb T,\Bbb R)$}.
\newblock {\em Duke Math. J.}, 135(2):327--360, 2006.

\bibitem{katopusateri}
Jun Kato and Fabio Pusateri.
\newblock A new proof of long-range scattering for critical nonlinear
  {S}chr\"odinger equations.
\newblock {\em Differential Integral Equations}, 24(9-10):923--940, 2011.

\bibitem{Kato_kdv}
Tosio Kato.
\newblock Quasi-linear equations of evolution, with applications to partial
  differential equations.
\newblock In {\em Spectral theory and differential equations ({P}roc.
  {S}ympos., {D}undee, 1974; dedicated to {K}onrad {J}\"orgens)}, volume Vol.
  448 of {\em Lecture Notes in Math.}, pages 25--70. Springer, Berlin-New York,
  1975.

\bibitem{Kato_kdv_disp}
Tosio Kato.
\newblock On the {C}auchy problem for the (generalized) {K}orteweg-de {V}ries
  equation.
\newblock In {\em Studies in applied mathematics}, volume~8 of {\em Adv. Math.
  Suppl. Stud.}, pages 93--128. Academic Press, New York, 1983.

\bibitem{KenigMartelRobbiano}
Carlos~E. Kenig, Yvan Martel, and Luc Robbiano.
\newblock Local well-posedness and blow-up in the energy space for a class of
  {$L^2$} critical dispersion generalized {B}enjamin-{O}no equations.
\newblock {\em Ann. Inst. H. Poincar\'e{} C Anal. Non Lin\'eaire},
  28(6):853--887, 2011.

\bibitem{KPV91}
Carlos~E. Kenig, Gustavo Ponce, and Luis Vega.
\newblock Well-posedness of the initial value problem for the {K}orteweg-de
  {V}ries equation.
\newblock {\em J. Amer. Math. Soc.}, 4(2):323--347, 1991.

\bibitem{KPV_gkdv}
Carlos~E. Kenig, Gustavo Ponce, and Luis Vega.
\newblock Well-posedness and scattering results for the generalized
  {K}orteweg-de {V}ries equation via the contraction principle.
\newblock {\em Comm. Pure Appl. Math.}, 46(4):527--620, 1993.

\bibitem{kpv_bilin}
Carlos~E. Kenig, Gustavo Ponce, and Luis Vega.
\newblock A bilinear estimate with applications to the {K}d{V} equation.
\newblock {\em J. Amer. Math. Soc.}, 9(2):573--603, 1996.

\bibitem{KPV_ill}
Carlos~E. Kenig, Gustavo Ponce, and Luis Vega.
\newblock On the ill-posedness of some canonical dispersive equations.
\newblock {\em Duke Math. J.}, 106(3):617--633, 2001.

\bibitem{KillipVisan_kdv}
Rowan Killip and Monica Vi\c{s}an.
\newblock Kd{V} is well-posed in {$H^{-1}$}.
\newblock {\em Ann. of Math. (2)}, 190(1):249--305, 2019.

\bibitem{Kishimoto_lwp_crit}
Nobu Kishimoto.
\newblock Well-posedness of the {C}auchy problem for the {K}orteweg-de {V}ries
  equation at the critical regularity.
\newblock {\em Differential Integral Equations}, 22(5-6):447--464, 2009.

\bibitem{kishimoto}
Nobu Kishimoto.
\newblock Unconditional uniqueness of solutions for nonlinear dispersive
  equations.
\newblock {\em Preprint, arxiv:1911.04349}, 2019.

\bibitem{Kishimoto_nls}
Nobu Kishimoto.
\newblock Unconditional local well-posedness for periodic {NLS}.
\newblock {\em J. Differential Equations}, 274:766--787, 2021.

\bibitem{KishimotoTsutsumi_kdnls}
Nobu Kishimoto and Yoshio Tsutsumi.
\newblock Gauge transformation for the kinetic derivative nonlinear
  {S}chr\"odinger equation on the torus.
\newblock {\em J. Differential Equations}, 453:Paper No. 113792, 47, 2026.

\bibitem{KlainermanMachedon}
Sergiu Klainerman and Matei Machedon.
\newblock Space-time estimates for null forms and the local existence theorem.
\newblock {\em Comm. Pure Appl. Math.}, 46(9):1221--1268, 1993.

\bibitem{KochTzvetkov_bo}
Herbert Koch and Nikolay Tzvetkov.
\newblock Nonlinear wave interactions for the {B}enjamin-{O}no equation.
\newblock {\em Int. Math. Res. Not.}, (30):1833--1847, 2005.

\bibitem{KortewegDeVries}
Diederik~J. Korteweg and Gustav de~Vries.
\newblock On the change of form of long waves advancing in a rectangular canal,
  and on a new type of long stationary waves.
\newblock {\em Philos. Mag. (5)}, 39(240):422--443, 1895.

\bibitem{Kuznetsov_Nakoryakov_Pokusaev_Shreiber_1978}
Vladimir~V. Kuznetsov, Vladimir~E. Nakoryakov, Boris~G. Pokusaev, and Isaac~R.
  Shreiber.
\newblock Propagation of perturbations in a gas-liquid mixture.
\newblock {\em Journal of Fluid Mechanics}, 85(1):85–96, 1978.

\bibitem{KwonOh_uu_mkdv}
Soonsik Kwon and Tadahiro Oh.
\newblock On unconditional well-posedness of modified {K}d{V}.
\newblock {\em Int. Math. Res. Not. IMRN}, (15):3509--3534, 2012.

\bibitem{KOY20}
Soonsik Kwon, Tadahiro Oh, and Haewon Yoon.
\newblock Normal form approach to unconditional well-posedness of nonlinear
  dispersive {PDE}s on the real line.
\newblock {\em Ann. Fac. Sci. Toulouse Math. (6)}, 29(3):649--720, 2020.

\bibitem{Lamb}
George~L. Lamb, Jr.
\newblock {\em Elements of soliton theory}.
\newblock Pure and Applied Mathematics. John Wiley \& Sons, Inc., New York,
  1980.
\newblock A Wiley-Interscience Publication.

\bibitem{linares_ponce_book}
Felipe Linares and Gustavo Ponce.
\newblock {\em Introduction to nonlinear dispersive equations}.
\newblock Universitext. Springer, New York, {S}econd edition, 2015.

\bibitem{McConnel_nls}
Ryan McConnell.
\newblock Nonlinear smoothing for the periodic generalized nonlinear
  {S}chr\"odinger equation.
\newblock {\em J. Differential Equations}, 341:353--379, 2022.

\bibitem{Miura}
Robert~M. Miura.
\newblock Korteweg-de {V}ries equation and generalizations. {I}. {A} remarkable
  explicit nonlinear transformation.
\newblock {\em J. Mathematical Phys.}, 9:1202--1204, 1968.

\bibitem{Molinet_ill_line}
Luc Molinet.
\newblock A note on ill posedness for the {K}d{V} equation.
\newblock {\em Differential Integral Equations}, 24(7-8):759--765, 2011.

\bibitem{MolinetPilod_bo}
Luc Molinet and Didier Pilod.
\newblock The {C}auchy problem for the {B}enjamin-{O}no equation in {$L^2$}
  revisited.
\newblock {\em Anal. PDE}, 5(2):365--395, 2012.

\bibitem{MolTzeSaut_ill}
Luc Molinet, Jean-Claude Saut, and Nikolay Tzvetkov.
\newblock Ill-posedness issues for the {B}enjamin-{O}no and related equations.
\newblock {\em SIAM J. Math. Anal.}, 33(4):982--988, 2001.

\bibitem{Nezlin}
Mikhail~V. Nezlin.
\newblock Rossby solitary vortices, on giant planets and in the laboratory.
\newblock {\em Chaos: An Interdisciplinary Journal of Nonlinear Science},
  4(2):187--202, 06 1994.

\bibitem{OS12}
Seungly Oh and Atanas Stefanov.
\newblock On quadratic {S}chr\"odinger equations in {${\bf R}^{1+1}$}: a normal
  form approach.
\newblock {\em J. Lond. Math. Soc. (2)}, 86(2):499--519, 2012.

\bibitem{OlsonChristensen}
Peter Olson and Ulrich Christensen.
\newblock Solitary wave propagation in a fluid conduit within a viscous matrix.
\newblock {\em Journal of Geophysical Research: Solid Earth},
  91(B6):6367--6374, 1986.

\bibitem{Russell}
John~Scott Russell.
\newblock The wave of translation in the oceans of water, air, and ether.
\newblock {\em Trübner and Co.}, 1885.

\bibitem{SautTemam}
Jean-Claude Saut and Roger Temam.
\newblock Remarks on the {K}orteweg-de {V}ries equation.
\newblock {\em Israel J. Math.}, 24(1):78--87, 1976.

\bibitem{StewartCorones}
Richard~W. Stewart and James Corones.
\newblock {The breakup of soliton like pulses on a nonlinear, nonuniform
  electrical lattice}.
\newblock {\em Rocky Mountain Journal of Mathematics}, 8(1-2):227--236, 1978.

\bibitem{takaoka_dnls}
Hideo Takaoka.
\newblock Well-posedness for the one-dimensional nonlinear {S}chr\"odinger
  equation with the derivative nonlinearity.
\newblock {\em Adv. Differential Equations}, 4(4):561--580, 1999.

\bibitem{tao_bo}
Terence Tao.
\newblock Global well-posedness of the {B}enjamin-{O}no equation in {$H^1({\bf
  R})$}.
\newblock {\em J. Hyperbolic Differ. Equ.}, 1(1):27--49, 2004.

\bibitem{tao_book}
Terence Tao.
\newblock {\em Nonlinear dispersive equations}, volume 106 of {\em CBMS
  Regional Conference Series in Mathematics}.
\newblock Conference Board of the Mathematical Sciences, Washington, DC; by the
  American Mathematical Society, Providence, RI, 2006.
\newblock Local and global analysis.

\bibitem{TsutsumiMukasa}
Masayoshi Tsutsumi and Toshio Mukasa.
\newblock Parabolic regularizations for the generalized {K}orteweg-de {V}ries
  equation.
\newblock {\em Funkcial. Ekvac.}, 14:89--110, 1971.

\bibitem{Tsutsumi_kdv_measures}
Yoshio Tsutsumi.
\newblock The {C}auchy problem for the {K}orteweg-de {V}ries equation with
  measures as initial data.
\newblock {\em SIAM J. Math. Anal.}, 20(3):582--588, 1989.

\bibitem{Tzvetkov_ill}
Nikolay Tzvetkov.
\newblock Remark on the local ill-posedness for {K}d{V} equation.
\newblock {\em C. R. Acad. Sci. Paris S\'er. I Math.}, 329(12):1043--1047,
  1999.

\bibitem{Wadati}
Miki Wadati.
\newblock The modified {K}orteweg-de {V}ries equation.
\newblock {\em J. Phys. Soc. Japan}, 34:1289--1296, 1973.

\bibitem{WashimiTaniuti}
Haruichi Washimi and Tosiya Taniuti.
\newblock Propagation of ion-acoustic solitary waves of small amplitude.
\newblock {\em Phys. Rev. Lett.}, 17:996--998, Nov 1966.

\bibitem{ZakharovFadeev}
Vladimir~E. Zakharov and Ludvig~D. Faddeev.
\newblock The {K}orteweg-de {V}ries equation is a fully integrable
  {H}amiltonian system.
\newblock {\em Funkcional. Anal. i Prilo\v zen.}, 5(4):18--27, 1971.

\end{thebibliography}

\end{document}